\documentclass[a4paper,11pt]{article}
\usepackage[utf8]{inputenc}

\usepackage{boxedminipage}
\usepackage{amsfonts}
\usepackage{amsmath}
\usepackage{amssymb}
\usepackage{graphicx}
\usepackage{amsthm}
\usepackage{t1enc}
\usepackage{subfig}
\usepackage{cancel}
 \usepackage{textcmds}
 \usepackage{mathabx}
 \usepackage{verbatim}

\newtheorem{theorem}{Theorem}[section]
\newtheorem{lemma}[theorem]{Lemma}

\newtheorem{definition}[theorem]{Definition}

\newtheorem{question}{Question}
\newtheorem*{theorem*}{Theorem}
\newtheorem{claim}{Claim}

\newcommand{\forceP}{\mathbb{P}}
\newcommand{\forceQ}{\mathbb{Q}}

\newcommand{\ZFC}{\mathsf{ZFC}}

\newcommand{\ZFP}{\mathsf{ZF}^-}

\newcommand{\CH}{\mathsf{CH}}

\def\undertilde#1{\mathord{\vtop{\ialign{##\crcr
$\hfil\displaystyle{#1}\hfil$\crcr\noalign{\kern1.5pt\nointerlineskip}
$\hfil\tilde{}\hfil$\crcr\noalign{\kern1.5pt}}}}}

\usepackage{array}
\newif\ifdraftnotes \draftnotesfalse

\usepackage{lmodern,mathtools,microtype}
\usepackage[margin=29mm]{geometry}
\usepackage{enumitem}
\usepackage{float}
\usepackage{tikz}
\usetikzlibrary{arrows.meta}
\usepackage[hidelinks,pdfencoding=auto,psdextra]{hyperref}
\hypersetup{pdftitle={Global Sigma-separation and fixed-level Pi-uniformization in forcing extensions of L},
  pdfauthor={Stefan Hoffelner},pdfdisplaydoctitle=true}
\numberwithin{equation}{section}
\newtheorem{mainthm}{Theorem}

\newtheorem{prop}[theorem]{Proposition}
\newtheorem{lem}[theorem]{Lemma}
\newtheorem{cor}[theorem]{Corollary}
\newtheorem{defn}[theorem]{Definition}
\newcommand{\R}{\mathbb R}
\newcommand{\Pp}{\mathbb P}
\newcommand{\Qq}{\mathbb Q}
\newcommand{\Tt}{\mathbb T}
\newcommand{\Cc}{\mathbb C}
\newcommand{\A}{\mathcal A}
\newcommand{\B}{\mathcal B}
\newcommand{\K}{\mathcal K}
\newcommand{\UP}{\operatorname{Unif}}
\newcommand{\Sep}{\operatorname{Sep}}
\newcommand{\Red}{\operatorname{Red}}
\newcommand{\Pos}{\operatorname{Poss}}
\newcommand{\cert}{\operatorname{Cert}}
\newcommand{\Code}{\operatorname{Code}}
\newcommand{\dom}{\operatorname{dom}}
\newcommand{\htR}{\operatorname{ht}}
\newcommand{\rank}{\operatorname{rk}}
\newcommand{\res}{\mathbin{\upharpoonright}}
\newcommand{\pa}[1]{\langle #1\rangle}
\newcommand{\cat}{\mathbin{{}^\frown}}
\newcommand{\RRrel}[2]{\boldsymbol{#1}^{1}_{#2}}
\newcommand{\resk}[1]{\mathbin{\upharpoonright_{#1}}}
\newcommand{\parm}{\operatorname{par}}
\newcommand{\lc}{\operatorname{lc}}
\newcommand{\da}{\mathfrak a}

\title{Global \(\Sigma\)-separation and fixed-level \(\Pi\)-uniformization\\
in forcing extensions of \(L\)}
\author{Stefan Hoffelner\thanks{Supported by the Austrian Science Fund (FWF), grant
  \href{https://doi.org/10.55776/P37228}{10.55776/P37228}.}\\[2pt]{\normalsize TU Wien}}
\date{25 November 2026}

\begin{document}
\maketitle
\begin{abstract}
We construct a cardinal-preserving set-forcing extension of \(L\) in which,
simultaneously for all \(n\geq3\) and all reals \(a\), any two disjoint
\(\Sigma^1_n(a)\) sets of reals can be separated by a \(\Delta^1_n(a)\) set.
This gives a simultaneous version of the positive \(\Sigma\)-separation
statement in Harrington's 1975 announcement. For each fixed \(n\geq3\) we
also construct an extension of \(L\) in which \(\Pi^1_n\)-uniformization
holds, and hence boldface \(\Pi^1_n\)-uniformization and
\(\Pi^1_n\)-reduction. 
\end{abstract}

\ifdraftnotes
\section*{Notes on this draft}

{\small
This file contains the sections on the coding machinery for the planned paper on $\Sigma^1_n$-separation. The text is taken from ``A Failure of $\Pi^1_3$-Reduction in the Presence of $\Sigma^1_3$-Separation'' \cite{SepRed}, mainly from its Sections 3 and 4. These notes are not part of the paper.

\paragraph{Second version.} The first version was checked against every use of the coding in the $\Sigma^1_4$ draft; the result is the second table below. The following was added: part (c) of Lemma \ref{lem:names_in_M}, part (4) of Lemma \ref{lem:symmetry}, the examples after Definition \ref{def:alpha_allowable_extension}, the last two sentences of Lemma \ref{productallowable} in their present form, and the last paragraph of Section \ref{sec:long}. Nothing was removed, and the copied text is unchanged.

\paragraph{Third version.} The second version was checked once more against the $\Sigma^1_4$ draft, proof by proof, including its hierarchy and the recursion above rank $0$. One statement was too narrow: the normal form gave codes only for forcings over $W$. The remarks on codes after Lemma \ref{lem:names_in_M} and the sentences after Claim \ref{claim:localization} are new. Nothing else was changed.

\paragraph{What is copied.}
\begin{itemize}
\item The preamble, Section 1.1, Section 2, Section 3 and Section 4 up to Lemma \ref{nounwantedcodes}. All numbers in these parts agree with \cite{SepRed}.
\item Section \ref{sec:alpha}, from Section 4.1 of \cite{SepRed}: its first sentences, Definition \ref{def:alpha_allowable_extension}, the paragraph after it, and Lemma \ref{iteration_of_alpha_allowable} with its proof, read for rank $0$. Their numbers agree with \cite{SepRed} as well.
\item Section \ref{sec:products}, from Section 5 of \cite{SepRed}: its first paragraph and Claim \ref{claim:localization} with its proof. The statements (1) to (3) of Lemma \ref{lem:symmetry} are the sentences at the beginning of the proof of Lemma 5.1. Lemma \ref{productallowable} is Lemma 5.2 for rank $0$, with the part of its proof which treats this case.
\item Section \ref{sec:long}, from Section 10 of \cite{SepRed}: its first lines and Claim \ref{claim:reflection}, which is Claim 4 there.
\end{itemize}

\paragraph{What is not copied.} The rank rules of Section 4.1 of \cite{SepRed}, that is the list of formulas, the restrictions $E$, Cases 1 and 2, the limit stages and the remarks, and Lemmas 5.1 to 5.3 for positive rank. They belong to the proof of $\Sigma^1_3$-separation. The hierarchy for $\Sigma^1_4$-separation takes their place and can be built on Section \ref{sec:alpha}. The paragraphs which are commented out in the source of Section 4 are left out.

\paragraph{What is new.}
\begin{itemize}
\item Lemma \ref{lem:fresh_real}: a coding forcing adds a new real, also below a condition.
\item Lemma \ref{lem:names_in_M}: the bookkeeping function in $L$ can be replaced by names in $M$, so allowability over $M$ is definable in $M$. Over a suitable pair $(M,J)$ it only depends on $G^0 \times G^1$ and $J$. This is the base case of the recursion for the hierarchy.
\item The remarks after Lemma \ref{lem:names_in_M} and after Claim \ref{claim:localization}: the forcings which are allowable over a suitable pair $(M,J)$ have hereditarily countable codes, and these codes belong to $N_J$.
\item Lemma \ref{lem:limits}: countable limits of allowable forcings over a suitable pair.
\item The proof of Lemma \ref{lem:symmetry}, the second half of its statement (1), and its statement (4), which moves the cover at a later stage.
\item Lemma \ref{lem:long_exactness}: $\Phi_0$ and $\Phi_1$ are exact for iterations of length $\omega_1$. In \cite{SepRed} this argument is used for single reals, in Section 9.2 for the statement (P) and in the proof of Claim 5.
\item Single sentences: all cardinals of $L$ are preserved in $W$ (Section \ref{sec:coding}); the examples after Definition \ref{def:alpha_allowable_extension}; countable sequences of ordinals in Claim \ref{claim:localization}; the last sentences of Lemma \ref{productallowable}; the dense subset of size $\aleph_1$ and the last paragraph in Section \ref{sec:long}.
\end{itemize}

\paragraph{Changes in the copied text.}
\begin{enumerate}
\item Section 2: the sentence on the role of the trees no longer speaks of the failure of $\Pi^1_3$-reduction.
\item Section 2: $\vec{S}^0$ and $\vec{S}^1$ are sequences and $S^l_{\beta}$ is defined. In \cite{SepRed} they are introduced as sets, and $S^l_{\beta}$ is used later.
\item Section 2: ``$\mathsf{ZF}^-$-models $M$ with,'' is completed to ``with $\omega_1 \in M$''. This is a guess at the missing condition.
\item Section 3: the coding area $\eta$ is added as an argument of $\operatorname{Code}$ in three places.
\item Section 3: $\rho$ is defined on all countable sets, as it is applied to every $g \cap \alpha$, $\alpha < \omega_1$.
\item Section 3: $\mathbb{P}^0$ in place of $\mathbb{P}_0$; ``Theorem \ref{DefinitionIndependentSequence}'' in place of ``lemma \ref{DefinitionIndependentSequence}''; ``Definition \ref{definitionadcoding}'' in place of ``Definition 2.5''.
\item Section 3: $(\ast)_1$ and $\sigma_1$ carry the index $1$ throughout. In the display of $(\ast)_1$ a comma and a second full stop are removed.
\item Section 3: in $\Psi_1$ the condition $\omega_1^M \in \bar{N}$ is added. The statement $\sigma_1$ has it, and $\bar{N}$ needs it to compute $\vec{S}$ correctly by Theorem \ref{DefinitionIndependentSequence}.
\item Section 3: the reference to Section 10 before Lemma \ref{lem:local_decoding} is removed.
\item Section 4: in the footnote to the proof of Lemma \ref{FirstPropertiesOfAllowableForcings} the lengths are called $\delta_0$ and $\delta_1$, as in the main text.
\item Section 4: $\Phi_0$ gets its argument in the display before Lemma \ref{nounwantedcodes}.
\item Sections \ref{sec:alpha} and \ref{sec:products}: ``allowable'' in place of ``$\alpha$-allowable'', and ``with respect to $E$'' and the restrictions $E$ are removed. In the proof of Lemma \ref{iteration_of_alpha_allowable} the two sentences on the rules are removed.
\item Section \ref{sec:products}: the first paragraph says that an allowable forcing is not tied to the coding areas it uses. In \cite{SepRed} it says this of the existence of an $\alpha$-allowable forcing which forces some $\Sigma^1_3$-statement.
\item Lemma \ref{lem:symmetry}: the sentences of \cite{SepRed} are arranged as the parts (1) to (3) of a lemma. The words which introduce $(M,J)$, $H'$ and $\pi$ are new, and so is part (4).
\item Lemma \ref{productallowable}: the two forcings are called $\mathbb{Q}^1$ and $\mathbb{Q}^2$, as $\mathbb{P}^1$ is the product of the Cohen forcings.
\item Section \ref{sec:long}: in Claim \ref{claim:reflection} the words on the intermediate models are removed, and after its proof the words ``but we will not need this''.
\end{enumerate}

\paragraph{Left as in the source.} The statement $\sigma_1$ speaks of models of $\mathsf{ZF}^-$, the formula $\Psi_1$ of models of $\mathsf{ZFC}^-$, and Lemma \ref{lem:local_decoding} of models of $\mathsf{T}$. The first paragraph of Section 3 writes $\operatorname{Code}(x)$. The branch $b_{\beta} \subset S^1_{\beta}$ in Section 3 is indexed in $\vec{S}^1$, the branch $b_{\xi}$ in Lemma \ref{DefinabilityInProperInnerModels} in $\vec{S}$.

\paragraph{Dictionary for Section 2 of the $\Sigma^1_4$ draft.}\

\medskip
\renewcommand{\arraystretch}{1.15}%
\noindent
\begin{tabular}{@{}>{\raggedright\arraybackslash}p{0.43\textwidth}>{\raggedright\arraybackslash}p{0.53\textwidth}@{}}
\emph{$\Sigma^1_4$ draft} & \emph{here} \\[2pt]
$W = L[\vec{b}][\vec{g}]$, cardinals, $\mathsf{CH}$ & Section \ref{sec:coding}, first paragraphs; $\vec{b}$ and $\vec{g}$ are given by $G^0$ and $G^1$ \\
word $w$, area $\eta$, $\operatorname{Code}(w,\eta)$ & the real $(w,0,0)$ on the side $1$, coding area $\eta$, $\operatorname{Code}(w,0,0,1,\eta)$. So every real is a word, see the examples after Definition \ref{def:alpha_allowable_extension}. The side $0$ is not needed \\
an active factor adds a new real, also below a condition & Lemma \ref{lem:fresh_real} \\
inactive factor & the trivial forcing in Definition \ref{def:allowable} \\
coding iteration; names for word, area, activity and restricting condition & allowable forcing, Definitions \ref{def:allowable} and \ref{def:allowable_extension}; the condition $\dot{c}$ in the remarks after Definition \ref{def:allowable} \\
ccc, dense set of size $\aleph_1$, $\mathsf{CH}$, up to length $\omega_1$ & Lemma \ref{FirstPropertiesOfAllowableForcings}, the normal form after it, Section \ref{sec:long} \\
every real occurs at a bounded stage & Section \ref{sec:long}, Claim \ref{claim:reflection} \\
$\Phi$ and equation (1) & $\Phi(w) := \Phi_1((w,0,0))$; Lemma \ref{nounwantedcodes} for countable length, Lemma \ref{lem:long_exactness} for length $\omega_1$ \\
$U$, $N = L[\vec{b}, \vec{g} \upharpoonright U, h]$, $M = N[H]$ & the cover $J$, $N_J = L[G^0][G^1 \upharpoonright J][G]$, Claim \ref{claim:localization}; $G$ and the sequence of its generic reals are definable from each other over $L[G^0][G^1 \upharpoonright J]$ \\
prefixes and cones & Lemma \ref{iteration_of_alpha_allowable} (2) and (3) \\
products in either order & Lemma \ref{FirstPropertiesOfAllowableForcings} (4), Lemma \ref{productallowable} \\
recursion for named stages, rank $0$ & Lemmas \ref{lem:names_in_M} and \ref{lem:limits} \\
symmetry of the reservoir & Lemma \ref{lem:symmetry} \\
\end{tabular}

\paragraph{Check against Sections 3 to 9 of the $\Sigma^1_4$ draft.} Every place where these sections use the coding is listed, with the statement here which covers it.

\medskip
\renewcommand{\arraystretch}{1.15}%
\noindent
\begin{tabular}{@{}>{\raggedright\arraybackslash}p{0.55\textwidth}>{\raggedright\arraybackslash}p{0.41\textwidth}@{}}
\emph{use in the $\Sigma^1_4$ draft} & \emph{here} \\[2pt]
stage $s = (M,U,B)$, without $B$ & suitable pair, Definition \ref{def:alpha_allowable_extension} (1) \\
block over $s$, sets $U_{\gamma}$, extra unavailable areas & Definition \ref{def:alpha_allowable_extension} (2), Lemma \ref{lem:names_in_M} (c) \\
one step block which writes a word; padding by inactive steps & the examples after Definition \ref{def:alpha_allowable_extension} \\
tags: the generic real of a step is new (Lemma 3.3, Proposition 5.1, Section 9) & Lemma \ref{lem:fresh_real} \\
$\operatorname{Safe}$ is a first order test in the stage (Definition 3.1); rank $0$ of the set recursion and its parameters (Lemma 4.2) & Lemma \ref{lem:names_in_M} \\
nice names and real codes for blocks (proof of Lemma 4.2) & the normal form before Lemma \ref{DefinabilityInProperInnerModels}; the remarks on codes after Lemma \ref{lem:names_in_M} \\
cones, prefixes, quotients, two-step compositions (Lemmas 3.2 and 4.2) & Lemma \ref{iteration_of_alpha_allowable} \\
countable concatenations, covers at limits (Lemma 4.2, Sections 6, 7 and 9) & Lemma \ref{lem:limits} \\
$M = N[H]$, all parameters and the codes of blocks lie in $N$ (Lemma 4.3) & Claim \ref{claim:localization} and the sentences after it \\
rewriting of the areas by a permutation, also at an internal prefix; the automorphism argument (Lemma 4.3) & Lemma \ref{lem:symmetry} \\
products in both orders, mixed prefixes, same generic reals (Proposition 5.1) & Lemma \ref{productallowable} \\
reals and $\Sigma^1_3$-statements at bounded stages; tails of the long iteration (Lemmas 3.2 and 7.1, Section 8) & Section \ref{sec:long}, Claim \ref{claim:reflection} \\
(7): a word is protected if and only if $\lnot\Phi$ holds in the end & Lemmas \ref{nounwantedcodes} and \ref{lem:long_exactness} \\
cardinals, $\mathsf{CH}$, a dense set of size $\aleph_1$ in $L$ (Section 7) & Sections \ref{sec:coding} and \ref{sec:long}, Lemma \ref{FirstPropertiesOfAllowableForcings} \\
\end{tabular}

\medskip
\noindent
Not in this file: the protected words, the declarations and the hierarchy $A_{\rho}$, $\operatorname{Poss}_{\rho}$, $\operatorname{side}_{\nu}$ of the $\Sigma^1_4$ draft, and the proof of its Lemma 4.2 above rank $0$. The fact that two mutually generic extensions have only the reals of the ground model in common, which is used in Proposition 5.1, is not special to the coding and is not repeated here.
}

\bigskip
\hrule
\bigskip
\fi

\section{Introduction}

Separation, reduction and uniformization are the three classical
structural properties of the projective pointclasses. Write
\(\R=2^\omega\), with its usual recursive coding of finite tuples. For a
projective pointclass \(\Gamma\), let \(\check\Gamma\) be its dual. The
\emph{separation property} \(\Sep(\Gamma)\) says that for any two disjoint
sets \(A_0,A_1\in\Gamma\) there is a set \(C\in\Gamma\cap\check\Gamma\)
with \(A_0\subseteq C\subseteq\R\setminus A_1\). The \emph{reduction
property} \(\Red(\Gamma)\) says that for all \(A_0,A_1\in\Gamma\) there
are disjoint sets \(B_0,B_1\in\Gamma\) with \(B_i\subseteq A_i\) for
\(i<2\) and \(B_0\cup B_1=A_0\cup A_1\). The \emph{uniformization
property} \(\UP(\Gamma)\) says that every \(A\subseteq\R^2\) in \(\Gamma\)
contains the graph of a function in \(\Gamma\) whose domain is
\(\{x:\exists y\, A(x,y)\}\). We write \(\Sigma^1_n(a)\) and
\(\Pi^1_n(a)\) for the effective pointclasses relative to a real \(a\).
Omitting \(a\) denotes the lightface pointclass, and boldface symbols
allow arbitrary real parameters.

The three properties are linked. Uniformization for \(\Gamma\) implies
reduction for \(\Gamma\), reduction for \(\Gamma\) implies separation for
\(\check\Gamma\), and reduction for \(\Gamma\) excludes separation for
\(\Gamma\) itself (Proposition~\ref{prop:duality}). So at each level at
most one of \(\Sigma^1_n\) and \(\Pi^1_n\) has the reduction property, and
if one of them has it, separation holds on the other side and fails on
this one.

In \(\ZFC\) the first two levels are settled. Lusin's separation theorem
\cite{lusin} gives \(\Sep(\Sigma^1_1)\), and Kond\^o's theorem \cite{kondo}
gives \(\UP(\Pi^1_1)\) and hence \(\UP(\Sigma^1_2)\) \cite[4E]{mosch}. So reduction holds for
\(\Pi^1_1\) and for \(\Sigma^1_2\), and \(\Sigma^1_2\)-separation and
\(\Pi^1_2\)-uniformization fail. Above the second level \(\ZFC\) decides
none of these properties, and the two canonical extensions of \(\ZFC\)
give different answers. If \(V=L\), the good \(\Sigma^1_2\) wellorder of
the reals yields \(\UP(\Sigma^1_n)\) for every \(n\ge2\)
\cite{addison}, \cite[5A]{mosch}. Reduction then holds on the
\(\Sigma\)-side at every level, and \(\Sigma^1_n\)-separation fails for
every \(n\ge2\). Under projective determinacy, Moschovakis' second
periodicity theorem yields \(\UP(\Pi^1_{2k+1})\) and
\(\UP(\Sigma^1_{2k+2})\) for every \(k\) \cite{mosch71}, \cite[6C]{mosch},
\cite[\S39]{kechris}.
Reduction then alternates between the two sides with the parity of the
level, and \(\Sigma^1_n\)-separation fails for every even \(n\). The
reduction property for these pointclasses was obtained first, by Martin
\cite{martin68} and by Addison and Moschovakis \cite{addisonmosch}.
Projective determinacy follows from the existence of infinitely many
Woodin cardinals \cite{martinsteel}; see \cite{msw} for the correspondence
between its levels and inner models with finitely many Woodin cardinals.
In both
pictures a single global structure, the definable wellorder of the reals
in one case and the scales given by determinacy in the other, decides
the pattern at all levels at once. Under the axiom of determinacy the
pattern is rigid: for every nonselfdual pointclass \(\Gamma\) which is
closed under continuous preimages, exactly one of \(\Gamma\) and
\(\check\Gamma\) has the separation property \cite{vanwesep,steel81}.

This raises the question which other patterns are consistent. Mathias
asked whether \(\Sigma^1_3\)-separation can hold in a forcing extension
\cite[Problem~3029]{mathias}. A more general result announced by Harrington
in 1975 is recorded by Sami \cite[Theorem~1.21]{sami}; see also
\cite[Chapter~V, \S3, p.~230]{hinman}. As reproduced in
\cite[\S1, footnote~1]{klh}, it includes, for each fixed
\(n\geq3\), an extension of \(L\) which satisfies
\(\Sep(\RRrel{\Sigma}{n})\) and \(\neg\Sep(\RRrel{\Pi}{n})\). The May 1974
addendum to \cite{harr} treats the failure of separation for both
pointclasses.

At the third level, forcing over \(L\) has answered several of these
questions. \(\Sigma^1_3\)-separation holds in a forcing extension of \(L\)
\cite{sep3}; \(\Pi^1_3\)-reduction can hold while
\(\Pi^1_3\)-uniformization fails \cite{Ho3}; \(\Sigma^1_3\)-separation can
hold while \(\Pi^1_3\)-reduction fails \cite{SepRed}; and
\(\Pi^1_3\)-uniformization can be forced over \(L\) \cite{up}. These
constructions use a feature of the third level. A \(\Sigma^1_3\)
statement consists of an existential real quantifier followed by a
\(\Pi^1_2\) matrix, and the matrix is absolute between the models of the
iteration by Shoenfield's theorem \cite{shoenfield}, so a witness found at
some stage
remains a witness. At a level \(n\ge4\) the matrix is \(\Pi^1_{n-1}\), and
a witness found at an intermediate stage can be destroyed later. So far
the constructions were lifted to a level \(n>3\) only by working over a
canonical inner model with \(n-3\) Woodin cardinals \cite{steel95} instead
of \(L\) \cite[\S4]{up}, \cite{SepRed}. For the different setting of considering only projective subsets of \(\omega\) instead of projective sets of the reals, Harrington
announced in a 1975 addendum to \cite{harr}, with a sketch of the proof, an
extension of \(L\) in which lightface separation holds for both
\(\Sigma^1_3\) and \(\Pi^1_3\); Kanovei and Lyubetsky gave a full proof and
extended the result to any fixed level \cite[Theorem~1.1]{kl}. 

\subsection{Main results}
We show that this barrier can be overcome over \(L\) itself, without
large cardinals, and for separation at all levels within a single
forcing extension.

\begin{mainthm}\label{main:sep}
There is a cardinal-preserving set-forcing extension \(V_{\mathrm{sep}}\)
of \(L\) satisfying \(\CH\) and
\[
  \forall n\geq3\ \forall a\in\R\quad \Sep\bigl(\Sigma^1_n(a)\bigr).
\]
\end{mainthm}

\begin{mainthm}\label{main:up}
For each integer \(n\geq3\) there is a cardinal-preserving set-forcing
extension \(V_n\) of \(L\) satisfying \(\CH\) and \(\UP(\Pi^1_n)\).
\end{mainthm}

As is well known, lightface \(\Pi^1_n\)-uniformization implies
\(\Pi^1_n(a)\)-uniformization for every real \(a\), and hence boldface
\(\RRrel{\Pi}{n}\)-uniformization. Indeed, fix a lightface universal
relation \(A_n(\pa{m,a,x},y)\) for \(\Pi^1_n(a)\) relations, and regard
\(\pa{m,a,x}\) as the input. A lightface \(\Pi^1_n\) graph which
uniformizes \(A_n\) yields, for all \(m\) and \(a\), a \(\Pi^1_n(a)\)
graph which uniformizes the relation with code \(m\) and parameter \(a\).
So \(V_n\) satisfies \(\UP(\Pi^1_n(a))\) for every real \(a\). For
separation this argument is not available, as the two sets of a
universal pair need not be disjoint for every parameter; this is why
Theorem~\ref{main:sep} quantifies over the parameter. Taking
\(a=0^\omega\) there gives the lightface version, and as any two boldface
\(\Sigma^1_n\) sets are \(\Sigma^1_n(a)\) for a single real \(a\), it
also gives the boldface version. The bound \(n\geq3\) is optimal, as the first two levels are
settled in \(\ZFC\).

Theorem~\ref{main:sep} realizes a pattern which neither of the two
classical pictures allows. In \(V_{\mathrm{sep}}\),
\(\Sigma^1_n\)-separation holds at every level \(n\ge3\), whereas it
fails at every level \(n\ge2\) if \(V=L\), and at every even level under
projective determinacy. Moreover, it holds without the uniformization
property which produces separation in both pictures: in
\(V_{\mathrm{sep}}\), \(\Sigma^1_k\)- and \(\Pi^1_k\)-uniformization fail
for every \(k\ge3\) (Corollary~\ref{cor:models}). Theorem~\ref{main:sep}
is a simultaneous version of the positive \(\Sigma\)-part of the result
announced by Harrington.

By Proposition~\ref{prop:duality}, the model \(V_n\) of
Theorem~\ref{main:up} also satisfies \(\Red(\Pi^1_n)\),
\(\Sep(\Sigma^1_n)\) and \(\UP(\Sigma^1_{n+1})\), and
\(\Pi^1_n\)-separation fails there. So \(V_n\) realizes the separation
pattern announced by Harrington at the level \(n\), together with
\(\Pi^1_n\)-uniformization. At the level \(n\), the model \(V_n\) differs
from \(L\) for every \(n\). For odd \(n\), the levels \(n\) and \(n+1\) of
\(V_n\) look as under projective determinacy, although \(V_n\) is an
extension of \(L\) without large cardinals; for even \(n\), they are the
opposite of the pattern under projective determinacy. By
\cite{woodin82,hauserschindler}, large cardinal strength is needed if every
projective relation is to be uniformized by a projective function while
all projective sets are Lebesgue measurable and have the Baire property.
Two adjacent
\(\Pi\)-levels can never both have the uniformization property
(Section~\ref{sec:consequences}), so a global version of
Theorem~\ref{main:up} could concern at most every other level.

\subsection{The method}
Many of the methods which were developed at level three play a crucial role for the proofs of both theorems. Key techniques for this article were introduced in \cite{sep3}, \cite{Ho3}, \cite{up}, \cite{bpfa} and \cite{mixed}. 
Both constructions of this paper start from the same ground model \(W\), an extension
of \(L\) by generic branches through an independent sequence of
constructible Suslin trees and by Cohen subsets of \(\omega_1\), and code
reals into the trees by almost disjoint coding \cite{jensensolovay}. Thus the ground model is basically the same as in \cite{Ho3} and \cite{up}. This
yields a lightface
\(\Sigma^1_3\) predicate which is exact: in the final model it holds of
exactly the reals which were coded. From it we define coding predicates
\(D_j\) at all levels by the recursion
\(D_{j+1}(d)\leftrightarrow\exists t\,\neg D_j(e_j(d,t))\)
(Section~\ref{sec:common}). To make \(D_{j+1}(d)\) true, we protect one
real \(e_j(d,t)\) from being coded at level \(j\); to make it false, we
code all of them at level \(j\), also those for reals \(t\) which are
added later.

The substantial difficulty is to decide, during an iteration of length \(\omega_1\),
which statements to make true. For \(n\ge4\), the truth of a
\(\Sigma^1_n\) statement at an intermediate stage is a poor guide. A
witness has a \(\Pi^1_{n-1}\) matrix, so it can be refuted later and new
witnesses can appear, and the forcing which acts on the current answer
can itself reverse it. Section~\ref{sec:four} begins with an example in
which this happens infinitely often.

We therefore replace truth at a stage by a forcing-theoretic notion. A
witness is \emph{certified} at a stage if no continuation of the
iteration which respects the coding restrictions in effect can refute its
matrix, and a statement is \emph{possible} if some such continuation
produces a certified witness. Both notions, and the protection of a
real \(e_j(d,t)\) with a new real \(t\), have a prototype at level four
in \cite[\S7.5]{mixed}: there a value of a \(\Sigma^1_4\)-uniformizing
function is selected once a forcing of the current fixed-point class adds
a witness whose \(\Pi^1_3\) matrix no further forcing of this class can
refute, and a new real is added whose extensions are never coded. For
separation, we do not take these notions relative to the fixed point of a
thinning-out process. The restrictions are recorded in assignments, and
certification, possibility and the assignments are defined by one
recursion through all finite levels. The heart of the paper is that possibility does not
depend on the stage: it can neither arise nor disappear under further
allowable forcing (Proposition~\ref{s:product}). The proof moves
witnessing forcings to disjoint coding areas by a symmetry argument and
compares the two presentations of their product. This is a considerably more complicated version of several similar lemmas which were central for the separation, reduction or uniformization proofs on the third level. It needs complete
restrictions at all lower levels along the product, which are
constructed by induction on the level. Along an iteration which carries
out every possible witness, that is which follows the greedy paradigm already present in \cite{sep3}, possibility coincides with truth in the final
model (Lemma~\ref{s:truth}). 

For separation, each pair of \(\Sigma^1_n(a)\) sets is first tested for
a possible common point. If none is possible, the restrictions under
which this was established are stored with the pair, and every real is
placed on the side whose membership is possible under the stored
restrictions, on a fixed default side if neither is. By stage
independence this decision is never revised, and the coding predicate at
level \(n\) turns the placement into the separating set. All levels are
treated within one iteration, which proves Theorem~\ref{main:sep}.

For uniformization at a fixed level \(n\), certificates of the same kind
show that a counterexample to membership in a section, once found,
remains a counterexample at all later stages. This is what the rank
construction of \cite[\S2.6]{up}, given there for the third level,
needs. With it we select an element of each nonempty \(\Pi^1_n\) section
over \(L\), without the inner models with Woodin cardinals used in
\cite[\S4]{up}.

\subsection{How to read this paper}\label{sec:howtoread}
The main constructions of this paper are in
Sections~\ref{sec:four}--\ref{sec:up}. The coding machinery on which they
rest is to a large extent known from earlier work, and the paper is
arranged so that these sections can be read without it.

\emph{The coding machinery.}
Sections~\ref{sec:trees}--\ref{sec:allowable} contain the coding machinery
of \cite{sep3,Ho3,up,mixed}: independent Suslin trees and almost disjoint
coding, the coding forcings, and allowable forcings over suitable pairs,
their products and iterations of length \(\omega_1\). They follow
\cite[\S\S2--4]{SepRed} closely. Section~\ref{sec:common} introduces the
higher coding predicates and serves as the interface to the rest of the
paper: Sections~\ref{sec:four}--\ref{sec:up} use the coding only through
Definition~\ref{def:alpha_allowable_extension}, which defines suitable
pairs and allowability over them, the format of the iterations of length
\(\omega_1\) at the end of Section~\ref{sec:long}, the forcings
\(\Code(w,\eta)\), and Lemma~\ref{lem:interface}, which lists the
properties of the coding used later. A reader who takes
Lemma~\ref{lem:interface} for granted can therefore skip
Sections~\ref{sec:trees}--\ref{sec:allowable}, apart from
Definition~\ref{def:alpha_allowable_extension} and the last paragraph of
Section~\ref{sec:long}, and start with Section~\ref{sec:common}.

\emph{Separation.} Section~\ref{sec:four} carries out the method at level
four alone. There are no requirements at the lower levels, so every step
can be checked directly, and the section is meant as an introduction to
the general construction; we recommend reading it before
Section~\ref{sec:sep}. Section~\ref{f:toy} explains the ideas informally,
starting from an example in which the truth of a membership at the
intermediate stages changes infinitely often, and
Section~\ref{f:examples} collects examples which show why the definitions
take their form. Section~\ref{sec:sep} proves Theorem~\ref{main:sep}. It
begins with a list of the four changes compared with
Section~\ref{sec:four}; apart from them, the definitions and the arguments
follow Section~\ref{sec:four}. The requirements on the forcings of
Section~\ref{sec:sep} are collected in
\ref{rule:base}--\ref{rule:finite}, Table~\ref{tab:rules} records where
each of them is verified, and Figure~\ref{fig:deps} shows how the parts of
the recursion and of the proof of Proposition~\ref{s:product} depend on
each other. On a first reading, one may skip the definability of the
recursion, (D1)--(D4), and take Proposition~\ref{s:product} and
Lemma~\ref{s:extend} for granted.

\emph{Uniformization.} Section~\ref{sec:up} proves Theorem~\ref{main:up}
for \(n\ge4\); the case \(n=3\) is taken from \cite{up}. The section
begins with an example at level four, and Section~\ref{sec:uptransfer}
explains why the coding iterations of \cite{up} can be replaced by
allowable forcings over \(W\). The rank and product arguments follow
\cite[\S\S2.6--3]{up}. They are proved in full, and
Section~\ref{sec:ranks} contains a dictionary to the notation of
\cite{up}. Section~\ref{sec:up} does not use the constructions of
Sections~\ref{sec:four} and~\ref{sec:sep}. From Section~\ref{sec:sep} it
only takes the notation \ref{N:models}--\ref{N:par} and a few auxiliary
facts, each of which is cited where it is used.

\subsection{Notation}
Our notation is mostly standard. A forcing iteration of length $\gamma$ is written $\forceP=(\forceP_{\alpha} : \alpha < \gamma)$, with initial segments $\forceP_{\alpha}$, and we also write $\forceP_{\gamma}$ for $\forceP$. Where we have to refer to the iterand at stage $\alpha$, we write $\forceP(\alpha)$; it is a $\forceP_{\alpha}$-name for a partial order, but we omit the dot. For $\alpha \le \beta \le \gamma$, $\forceP_{\alpha,\beta}$ denotes the part of the iteration on the interval $[\alpha,\beta)$, so that $\forceP_{\beta} \cong \forceP_{\alpha} \ast \forceP_{\alpha,\beta}$; once a $\forceP_{\alpha}$-generic filter is fixed, we regard $\forceP_{\alpha,\beta}$ as a forcing over the corresponding extension. A two-step iteration is written $\forceP \ast \dot{\forceQ}$, where $\dot{\forceQ}$ is a $\forceP$-name for a partial order.

For a set $X$, $\Sigma_n(X)$ denotes the set of $\Sigma_n$-formulas with parameters from $X$.

We write $\mathbf 1 \Vdash \varphi$ if every condition forces $\varphi$. For $p \in \forceP$ let $\forceP \upharpoonright p := \{ q \in \forceP : q \le p \}$ be the restriction of $\forceP$ below $p$. If $p \Vdash \varphi$, we often work with $\forceP \upharpoonright p$ in place of $\forceP$, which reduces the notational load. For an allowable forcing, $\forceP \upharpoonright p$ is again presented as an allowable iteration, see the remarks after Definition~\ref{def:allowable}.

\subsection{Background and declaration of AI usage}

The investigation of a possible proof of both theorems of this paper was the driving force of a majority of my research since 2019, in a direct (obviously \cite{sep3}, \cite{Ho3} and \cite{up}) or indirect way. As an example, the paper \cite{bpfa}, even though it pulls in the opposite direction from the results presented here, was born out of failed attempts to force $\Pi^1_4$-uniformization over $L$; similarly, the paper \cite{mixed} grew out of investigating what a possible forcing for $\Pi^1_4$-uniformization could even look like. Note here that such a forcing must not look like the symmetry-based coding forcings from \cite{Ho3}, it must not look like the final construction from \cite{up}, as $\Pi^1_3$-uniformization implies $\Sigma^1_4$-uniformization, and it certainly must not resemble anything close to the constructions from \cite{bpfa}. Some of the ideas of \cite{mixed} were key in developing the forcing construction for $\Sigma^1_4$-separation; in particular, a first form of what I call \emph{possible} membership (Definition~\ref{f:poss}, motivated in Section~\ref{f:toy}) appears in \cite[\S7.5]{mixed}. Even the very recently finished \cite{pi2n} is an offspring of a wild idea to construct a model of $\Pi^1_n$-uniformization by forming an inverse iteration, which failed in the projective case but goes through in the second-order case. Building on these, a successful attack on forcing $\Sigma^1_4$-separation and $\Pi^1_4$-uniformization over just $L$ was possible without any AI assistance, together with an idea for tackling $\Sep(\Sigma^1_4)$ and $\Sep(\Sigma^1_5)$ jointly. However, for the higher projective levels and their combination I profited immensely from working closely with ChatGPT Astra 6.2 and Claude Fable 5.2; without them, these levels would have been out of reach for me, at least for the next year or two (this is hard to estimate), due to the huge technical obstacles this proof posed. I used both models extensively to help with the writing, for checking proofs, and in developing the proofs at the higher levels. I have checked all arguments in this paper myself and take full responsibility for its content.

\section{Independent Suslin trees in $L$, almost disjoint coding}\label{sec:trees}
We use the coding method of \cite{sep3,Ho3,up,mixed}: reals are written into an independent sequence of Suslin trees by almost disjoint coding, and the trees can be destroyed generically, independently of each other. Almost disjoint coding is due to Jensen and Solovay \cite{jensensolovay}; for its use in obtaining projectively definable sets of reals in forcing extensions of $L$ see also \cite{harrington77}. The use of constructible Suslin trees for this purpose goes back to Jensen and Johnsbr\aa ten \cite{jensenjohnsbraten}, and Abraham and Shelah \cite{abrahamshelah} code by products of Aronszajn trees, some of which are made special while others remain Suslin. Coding areas are also used in \cite{friedmanschrittesser}. Sections~\ref{sec:trees}--\ref{sec:allowable} follow \cite[\S\S2--4]{SepRed} closely. New are Lemma~\ref{lem:names_in_M} with the codes of allowable forcings, Lemma~\ref{lem:limits}, and part~(4) of Lemma~\ref{lem:symmetry} and its proof.
\begin{definition}
 Let $\vec{T} = (T_{\alpha} \, : \, \alpha < \kappa)$ be a sequence of Suslin trees. We say that the sequence is an
 independent family of Suslin trees if for every finite set of pairwise distinct indices $e= \{e_0, e_1,...,e_n\} \subset \kappa$ the product $T_{e_0} \times T_{e_1} \times \cdot \cdot \cdot \times T_{e_n}$
 is a Suslin tree again.
\end{definition}
If $\vec{T} = (T_{\alpha} \mid \alpha < \kappa)$ is an independent sequence of Suslin trees and $A \subset \kappa$, then forcing with the finite support product $\prod_{i \in A} T_i$, where each tree is ordered by its tree order, keeps every $T_{\alpha}$ with $\alpha \notin A$ Suslin, as $(\prod_{i \in A} T_i) \times T_{\alpha}$ has the ccc; for a more general argument of this type see \cite[Lemma~51]{Ho}.

One can easily force the existence of independent sequences of Suslin trees with products of Jech's or Tennenbaum's forcing, or with just products of ordinary Cohen forcing. On the other hand independent sequences of length $\omega_1$ already exist in $L$.

\begin{theorem}[{\cite[Theorem~2.2]{SepRed}}]\label{DefinitionIndependentSequence}
There are a sequence $\vec{T} = ( T_{\eta} \mid \eta < \omega_1) \in L$, which in $L$ is an independent sequence of Suslin trees whose nodes are countable ordinals, and a $\Sigma_1$-formula $\phi(v_0,v_1)$ with the following property. If $\aleph_1= \aleph_1^L$ and $(M,\in)$ is a transitive model of $\ZFC^{-}$ with $\omega_1 \in M$, then
\[ \{ (\eta,s,t) \mid \eta < \omega_1, \ s <_{T_{\eta}} t \} = \{ u \in M \mid M \models \phi(u,\omega_1) \}. \]
So $\vec{T}$ is uniformly $\Sigma_1(\{\omega_1\})$-definable over $M$; in particular $\vec{T}$ can be defined correctly in any transitive, uncountable model $M$ of $\ZFC^-$ with $\omega_1 \in M$, using $\phi$ and $\omega_1$ as a parameter of $\phi$.
\end{theorem}

The sequence of \cite[Theorem~2.2]{SepRed} is obtained from the canonical $\diamondsuit$-sequence of $L$ \cite{jensen72} by the standard construction of a Suslin tree \cite{devlinj}, carried out for all trees simultaneously, so that the maximal antichains of all finite products of the trees are sealed. Its trees are normal trees of height $\omega_1$: every node has infinitely many immediate successors and successors on every higher level, and at limit levels distinct nodes have distinct sets of predecessors. In particular every level is countably infinite. Replacing the nodes of height $\gamma$ by the ordinals in $[\omega\gamma,\omega\gamma+\omega)$, via the $<_L$-least bijection, we obtain isomorphic trees. The new sequence is again independent, and as the renaming is computed in $L_{\omega_1}$, which every $M$ as in Theorem \ref{DefinitionIndependentSequence} computes correctly, it is uniformly $\Sigma_1(\{\omega_1\})$-definable in the same sense. We assume from now on that the nodes of height $\gamma$ are exactly the ordinals in $[\omega\gamma,\omega\gamma+\omega)$. So if $\omega\alpha = \alpha$, then the nodes below $\alpha$ are exactly the nodes of height below $\alpha$, and we write $T_{\eta} \upharpoonright \alpha$ for the set of these nodes.

From now on we write $\vec{S} = (S_{\eta} \mid \eta < \omega_1)$ for the sequence $\vec{T}$ of Theorem \ref{DefinitionIndependentSequence}.

The trees from $\vec{S}$ will later be used to define two $\Sigma^1_3$-predicates on the reals, which will be manipulated with the coding forcings of the next section. For this we split $\vec{S}$ into two independent sequences
\[ \vec{S}^0 := ( S_{\alpha} \mid \alpha \text{ is even}) \]
and
\[\vec{S}^1:= ( S_{\alpha} \mid \alpha \text{ is odd}). \]
We write $S^l_{\beta}$ for the $\beta$-th tree of $\vec{S}^l$, so $\vec{S}^l = (S^l_{\beta} \mid \beta < \omega_1)$ for $l \in \{0,1\}$.

We briefly introduce the almost disjoint coding forcing due to R. Jensen and R. Solovay \cite{jensensolovay}. We will identify subsets of $\omega$ with their characteristic function and will use the word reals for elements of $2^{\omega}$ and subsets of $\omega$ respectively.
Let $D=\{d_{\alpha} \, : \, \alpha < \aleph_1 \}$ be a family of almost disjoint subsets of $\omega$, i.e. a family such that if $r, s \in D$ then
$r \cap s$ is finite. Let $X\subseteq \omega_1$ be a set of ordinals. Then there
is a ccc forcing, the almost disjoint coding $\mathbb{A}_D(X)$ which adds
a new real $x$ which codes $X$ relative to the family $D$ in the following way
$$\alpha \in X \text{ if and only if } x \cap d_{\alpha} \text{ is finite.}$$
\begin{definition}\label{definitionadcoding}
 The almost disjoint coding $\mathbb{A}_D(X)$ relative to an almost disjoint family $D$ consists of
 conditions $(r, R) \in [\omega]^{<\omega} \times D^{<\omega}$ and
 $(s,S) < (r,R)$ holds if and only if
 \begin{enumerate}
  \item $r \subset s$ and $R \subset S$.
  \item If $\alpha \in X$ and $d_{\alpha} \in R$ then $r \cap d_{\alpha} = s \cap d_{\alpha}$.
 \end{enumerate}
\end{definition}
We shall briefly discuss the $L$-definable, $\aleph_1^L$-sized almost disjoint family of reals $D$  we will use throughout this article. The family $D$ is the canonical almost disjoint family one obtains when recursively adding the $<_L$-least real $x_{\beta}$ not yet chosen and replace it with $d_{\beta} \subset \omega$ where this $d_{\beta}$  is the real which codes the initial segments of $x_{\beta}$ using some recursive bijections between $\omega$ and $\omega^{<\omega}$. The definition of $D$ is uniform over any uncountable, transitive $\ZFP$-models $M$ with $\omega_1 \in M$, as we can correctly compute $L$ up to $\aleph_1^L$ inside $M$ and then apply the above definition inside $L$'s version of $M$. Even more is true, if $M$ is a countable, transitive model of $\ZFP+``$ $\aleph_1$ exists and $\aleph_1=\aleph_1^L"$, then $M$ will compute $D \upharpoonright \omega_1^M$ in a correct way. The reason is again, that $M$ can define an initial segment of $L$ correctly which suffices to calculate $D \upharpoonright \omega_1^M$.

Finally, as $\mathbb{A}_D(X)$ has the Knaster property, it preserves Suslin trees.

\section{Coding machinery}\label{sec:coding}
We continue with the construction of the appropriate notions of forcing which we want to use in our proof. The goal is to first define a coding forcing $\operatorname{Code} (x)$ for reals $x$, which will force for $x$ that a certain $\Sigma^1_3$-formula $\Phi(x)$ becomes true in the resulting generic extension.

In a first step we force over $L$ to destroy all members of $\vec{S}= \vec{S}^0 \cup \vec{S}^1$ via generically adding an $\omega_1$-branch, that is we  form \[\forceP^0:=\prod_{\alpha \in \omega_1} S_{\alpha}\] with finite support. Note that this is an $\aleph_1$-sized, ccc forcing over $L$.

In a second step
we add $\omega_1$-many $\omega_1$-Cohen subsets with a countably supported product, for technical reasons we want the second forcing to be defined as in $L$ and not as in $L^{\forceP^0}$. We let $$\forceP^1:= (\prod_{\alpha< \omega_1} \mathbb{C} (\omega_1))^L.$$
This forcing is $\sigma$-closed only in $L$, and as it is defined in $L$, the product $\forceP^0 \times \forceP^1$ sits densely in the two step iteration $\forceP^0 \ast \forceP^1$. By Easton's Lemma (see Lemma 15.19 in \cite{Jech}), $\forceP^1$ adds no countable sets of ordinals and every forcing which is ccc in $L$ remains ccc in $L^{\forceP^1}$. So $\vec{S}$ is still an independent sequence of Suslin trees in $L^{\forceP^1}$, and $\forceP^0 \times \forceP^1$, which has size $\aleph_1$ in $L$, preserves all cardinals of $L$ and $\CH$.

 We use $W$ to denote this generic extension of $L$, we fix a filter $G^0 \times G^1$ which is generic for $\forceP^0 \times \forceP^1$  and let \[W=L[G^0][G^1]. \]

 Let $x \in W$ be a real,  and let $m,k \in \omega$ and let $\eta <\omega_1$. We simply write $(x,m,k)$ for a real $w$ which codes the triple $(x,m,k)$ in a recursive way. The forcing $\operatorname{Code}(x,m,k,1,\eta)$ \footnote{The other coding forcing,  $\operatorname{Code}(x,m,k,0,\eta)$, which uses the trees from the $\vec{S}^0$-sequence instead is defined in the analogous way}  which codes the triple $(x,m,k)$ into $\vec{S^1}$ is defined as the almost disjoint coding forcing of a specific set $Y \subset \omega_1$, that is
 \[ \operatorname{Code}(x,m,k,1,\eta) := {\mathbb{A}}({Y}).\]  We will define the crucial set $Y \subset \omega_1$ now.

 To ease notation we let $g \subset \omega_1$ be $g_{\eta}$ for $\eta < \omega_1$, where $g_{\eta}$ is the $\eta$-th coordinate of the $\prod_{\alpha < \omega_1} \mathbb{C}(\omega_1)$-generic filter over $L^{\forceP^0}$. We let $\rho: ([\omega_1]^{\le\omega})^{L} \rightarrow \omega_1$ be some canonically definable, constructible bijection between
these two sets. We use $\rho$ and $g$ to define the set $h \subset \omega_1$, which eventually shall be the set of indices of $\omega$-blocks of $\vec{S}$, where we ``code up the characteristic function of the real $(x,m,k)"$, the latter slogan will be made precise in a moment. Let \[h:= \{\rho( g \cap \alpha) \,: \, \alpha < \omega_1 \}\] and let
\begin{align*}
A:= &\{ \omega \gamma +2n \mid \gamma \in h, n \notin (x,m,k) \} \cup \\& \{\omega \gamma + 2n+1 \mid \gamma \in h, n \in (x,m,k) \}.
\end{align*}

Let $X \subset \omega_1$ be chosen such that it canonically codes the following objects, with $A$ as the set of its even entries:
\begin{itemize}
\item The set $A \subset \omega_1$.
\item The set $\{ b_{\beta} \subset S^1_{\beta} \mid \beta \in A \}$ of generic $\omega_1$-branches we added with $\forceP^0$.

\end{itemize}

For this coding over $W$, every $\gamma\in h$ gives the following branch pattern in $L[X]$.
\begin{align*}
   (\ast)_1({\gamma},(x,m,k)) = {} & \parbox[t]{0.75\textwidth}{$n \in (x,m,k)$ if and only if $S^1_{\omega \cdot \gamma +2n+1}$ has an $\omega_1$-branch, and $n \notin (x,m,k)$ if and only if $S^1_{\omega \cdot \gamma +2n}$ has an $\omega_1$-branch.}
\end{align*}

The branch prescribed by $A$ is recovered from $X$. For the other tree in the pair, let $\xi$ be its index in the full sequence $\vec S$ and put
\[
 N_\xi=L[G^1][G^0\upharpoonright(\omega_1\setminus\{\xi\})].
\]
By independence in $L[G^1]$, $S_\xi$ is Suslin in $N_\xi$. Since $W=N_\xi[b_\xi]$ and forcing with a normal Suslin tree adds no reals, the input real belongs to $N_\xi$. The set $g_\eta$ and the sequence of all the prescribed branches also belong to $N_\xi$, so the canonical code $X$ belongs to $N_\xi$. Thus $L[X]$ has no cofinal branch through $S_\xi$.

We note that we can apply an argument resembling David's trick \footnote{see \cite{David} for the original argument, where the strings in Jensen's coding machinery \cite{bjw} are altered such that certain unwanted universes are destroyed. This destruction is emulated in our context as seen below.} in this situation. We rewrite the information of $X \subset \omega_1$ as a subset $Y \subset \omega_1$ using the following line of reasoning. All choices in this construction (of the model $N$, its code $X_N$, the ordinal $\beta$, the club $C$ and the models $M_{\alpha}$ below) are made canonically inside $L[X]$, using its canonical well-order; in particular, $Y \in L[X]$ is determined by $X$.
In $L[X]$ choose $N$ as the transitive collapse of an elementary submodel of $H_{\omega_2}^{L[X]}$ of size $\aleph_1$ containing $\omega_1\cup\{\omega_1,X\}$. The collapse fixes $X$ and $\omega_1$, so $N$ recovers $h$, the coded real and the prescribed branches from $X$, and computes the trees correctly by Theorem~\ref{DefinitionIndependentSequence}. For the coding over $W$ just considered, $N$ therefore satisfies $(\ast)_1(\gamma,(x,m,k))$ for every $\gamma\in h$.
Consequently, if we code the structure $(N,\in,X)$, with $X$ as a distinguished constant, as a set $X_N \subset \omega_1$ (say $X_N := \{ \Gamma(\xi,\zeta) \mid f(\xi) \in f(\zeta) \}$, where $\Gamma$ is G\"odel's pairing function and $f\colon \omega_1 \to N$ is the $<_{L[X]}$-least bijection with $f(0) = X$), then for any $\beta>\omega_1$ such that $L_{\beta}[X_N] \models \ZFP$:
\begin{align*}
    L_{\beta}[X_N] \models & \text{\ldq The model decoded out of }X_N \text{ satisfies $(\ast)_{1}(\gamma,(x,m,k))$} \\& \text{for every $\gamma \in h$\rdq.}
\end{align*}

The strict inequality $\beta>\omega_1$ ensures that $X_N \in L_{\beta}[X_N]$, so that $L_{\beta}[X_N]$ can decode the structure coded by $X_N$; this may fail for $\beta=\omega_1$.

For the latter assertion to be meaningful, we need that the model decoded out of $X_N$ has $h$ as an element but this is clear as $h$ can be read off from $A$ which in turn is coded into $X$ which belongs to $N$ as an element.
Taking an elementary submodel of size $\aleph_1$ of a sufficiently large level of $L[X_N]$ containing $\omega_1\cup\{\omega_1,X_N\}$ and applying condensation gives a $\beta$ with $\omega_1<\beta<\omega_2$ and $L_\beta[X_N]\models\ZFP$. We take the least such $\beta$. This choice also defines $\beta$ when $X$ is formed at a later coding stage. Fix a club $C\subset\omega_1$ and a sequence $(M_{\alpha} \, : \, \alpha \in C)$ of countable elementary submodels of $(L_{\beta} [X_N],\in,X_N)$ such that
\[\forall \alpha \in C (M_{\alpha} \prec L_{\beta}[X_N] \land M_{\alpha} \cap \omega_1 = \alpha)\]
Now let the set $Y\subset \omega_1$ code the pair $(C, X_N)$ such that the odd entries of $Y$ should code $X_N$ and if $E(Y)$ denotes  the set of even entries of $Y$ and $\{c_{\alpha} \, : \, \alpha < \omega_1\}$ is the enumeration of $C$ then
\begin{enumerate}
\item $E(Y) \cap \omega$ codes a well-ordering of type $c_0$.
\item $E(Y) \cap [\omega, c_0) = \emptyset$.
\item For all $\beta$, $E(Y) \cap [c_{\beta}, c_{\beta} + \omega)$ codes a well-ordering of type $c_{\beta+1}$.
\item For all $\beta$, $E(Y) \cap [c_{\beta}+\omega, c_{\beta+1})= \emptyset$.
\end{enumerate}
We claim that the following assertion is true, which can be seen as a local version of $(\ast)_{1}(\gamma, (x,m,k))$.

\begin{align*}
\sigma_1 (x,m,k) := {} & \parbox[t]{0.75\textwidth}{For any countable transitive model $M$ of ``$\ZFP$ and $\aleph_1$ exists'' such that $\omega_1^M=(\omega_1^L)^M$ and $Y \cap \omega_1^M \in M$, $M$ can construct its version of the universe $L[Y \cap \omega_1^M]$, and the latter will see that there is an $\aleph_1^M$-sized, transitive model $\bar{N} \in L[Y \cap \omega_1^M]$, $\bar{N} \models \ZFC^- \land \omega_1^M \in \bar{N}$ which models $(\ast)_{1} ({\gamma,(x,m,k)})$ for $\aleph_1^M$-many $\gamma$.}
\end{align*}
Put $\delta=\omega_1^M$. The even entries of $Y\cap\delta$ imply $\delta\in C$, and its odd entries code $(N\cap M_\delta,\in,X)$ on $\delta$. The collapse belongs to $L[Y\cap\delta]^M$. Since $N\cap M_\delta\prec N$, the collapse sends $\omega_1$ to $\delta$ and $h$ to $h\cap\delta$, which is unbounded in $\delta$. For this coding over $W$, it therefore satisfies $(\ast)_1(\gamma,(x,m,k))$ for every $\gamma\in h\cap\delta$. The decoding is verified in Lemma~\ref{lem:local_decoding}.

We have finally defined the desired set $Y$ and now we use
\[\operatorname{Code} (x,m,k,1,\eta):= \mathbb{A} (Y) \]
 relative to our previously defined, almost disjoint family of reals $D \in  L $ (see the paragraph after Definition \ref{definitionadcoding})  to code the set $Y$ into one real $r$. This forcing only depends on the subset of $\omega_1$ we code, thus $\mathbb{A}_D(Y)$ will be independent of the surrounding universe in which we define it, as long as it has the right $\omega_1$ and contains the set $Y$.

The effect of the coding forcing $\operatorname{Code} (x,m,k,1,\eta)$ is that it generically adds a real $r$ such that
this formula becomes true:

\begin{align*}
\Psi_1 (r, (x,m,k)) := {} & \parbox[t]{0.75\textwidth}{For any countable, transitive model $M$ of ``$\ZFC^-$ and $\aleph_1$ exists'', such that $\omega_1^M=(\omega_1^L)^M$ and $r \in M$, $M$ can construct its version of $L[r]$, denoted by $L[r]^M$, which in turn thinks that there is a transitive $\ZFC^-$-model $\bar{N}$ of size $\aleph_1^M$ with $\omega_1^M \in \bar{N}$ such that $\bar{N}$ believes $(\ast)_{1}({\gamma,(x,m,k)})$ for an $\aleph_1^M$-sized set of ordinals $\gamma$.}
\end{align*}

Indeed, for $\delta=\omega_1^M$, both $M$ and $L[r]^M$ compute $D\upharpoonright\delta$ correctly, as discussed after Definition~\ref{definitionadcoding}. Thus $Y\cap\delta\in L[r]^M$, and the preceding argument gives the required $\bar N$ in $L[r]^M$.

The assertion $\Psi_1(r,(x,m,k))$ is $\Pi^1_2$. A failure is witnessed by a real coding a countable model and its satisfaction relation. Well-foundedness of the membership relation is $\Pi^1_1$, and the requirements on the theory, the parameters and the failure of the displayed conclusion are arithmetical in these codes. We say in the above situation that the real $(x,m,k)$ \emph{ is written into $\vec{S}^1$}, or that $(x,m,k)$ \emph{is coded into} $\vec{S^1}$. To summarize our discussion, given an arbitrary real of the form $(x,m,k)$, then our forcing $\operatorname{Code} (x,m,k,1,\eta)$, when applied over $W$, will add a real $r$ which will turn the $\Pi^1_2$-formula $ \Psi_1 (r, (x,m,k))$ into a true statement in $W^{\operatorname{Code} (x,m,k,1,\eta)}$.

The coding forcing which codes a given real $(x,m,k)$ into the $\vec{S}^0$, denoted by $\operatorname{Code} {(x,m,k,0,\eta)}$ is defined in the same way and so are the according formulas $\sigma_0$ and $\Psi_0$.

We make the local decoding precise. Let $\mathsf{T}$ be the theory $\ZFP$ + ``$\aleph_1$ exists and $\aleph_1 = \aleph_1^L$''. Let $M$ be a transitive model of $\mathsf{T}$, let $r \in M$ be a real, and let $\delta := \omega_1^M$. In $L[r]^M$ let $Y^M_r := \{ \xi < \delta \mid r \cap d_{\xi} \text{ is finite} \}$. If the odd entries of $Y^M_r$ code, via $\Gamma$ as above, a well-founded and extensional relation, let $(\bar{N}^M_r,\in,\bar{X}^M_r)$ be its transitive collapse together with the image of $0$, let $\bar{A}^M_r$ be the set of even entries of $\bar{X}^M_r$, and let
\[ a^M(r) := \{ \gamma < \delta \mid \omega \gamma + j \in \bar{A}^M_r \text{ for some } j \in \omega \} \]
be the \emph{coding area decoded from $r$ in $M$}. If one of these decodings fails, then $a^M(r)$ is undefined.

\begin{lemma}\label{lem:local_decoding}
Let $r$ be the real added by $\operatorname{Code}(x,m,k,1,\eta)$, either over $W$ or at a stage of an allowable iteration as in Definition~\ref{def:allowable}. Let $h$, $A$, $X$, $N$, $X_N$, $f$, $C$ and $(M_{\alpha} \mid \alpha \in C)$ be as in its definition. For a countable transitive model $M$ of $\mathsf{T}$ with $r \in M$, put $\delta := \omega_1^M$.
\begin{enumerate}
\item $Y^M_r = Y \cap \delta$, and $\delta \in C$.
\item $\bar{N}^M_r$ is defined and is the transitive collapse of $N \cap M_{\delta}$. It is a model of $\ZFC^-$ with $\delta = \omega_1^{\bar{N}^M_r} \in \bar{N}^M_r$, we have $\bar{X}^M_r = X \cap \delta$, and $\bar{N}^M_r$ computes the trees $S_{\xi}$, $\xi < \delta$, as $S_{\xi} \upharpoonright \delta$.
\item $a^M(r)=h\cap\delta$ is unbounded in $\delta$. If the coding is applied directly over $W$, then $\bar N^M_r$ satisfies $(\ast)_1(\gamma,(x,m,k))$ for every $\gamma\in a^M(r)$.
\end{enumerate}
The same holds for $\operatorname{Code}(x,m,k,0,\eta)$, with $(\ast)_0$ in place of $(\ast)_1$.
\end{lemma}
\begin{proof}
(1) As discussed after Definition \ref{definitionadcoding}, $M$ computes $D \upharpoonright \delta$ correctly. As $r$ is generic for $\mathbb{A}_D(Y)$, $r \cap d_{\xi}$ is finite if and only if $\xi \in Y$, so $Y^M_r = Y \cap \delta$. Suppose $\delta \notin C$. As $C$ is closed, $C \cap \delta$ is empty or has a largest element $c_{\nu}$. Then $E(Y) \cap \omega$, respectively $E(Y) \cap [c_{\nu},c_{\nu}+\omega)$, codes a well-ordering of type $c_0 > \delta$, respectively $c_{\nu+1} > \delta$. This well-ordering belongs to $M$, as $Y \cap \delta \in M$, so $\delta$ would be countable in $M$, a contradiction.

(2) As $\delta \in C$, $M_{\delta} \prec (L_{\beta}[X_N],\in,X_N)$ and $M_{\delta} \cap \omega_1 = \delta$. The bijection $f$ is the inverse of the transitive collapse of the relation coded by $X_N$, so it is definable from $X_N$ in $L_{\beta}[X_N]$ and $f \in M_{\delta}$. Hence $f[\delta] = N \cap M_{\delta}$, $\delta$ is closed under $\Gamma$, and the odd entries of $Y \cap \delta$ code $X_N \cap \delta = \{ \Gamma(\xi,\zeta) \mid \xi,\zeta < \delta, \ f(\xi) \in f(\zeta) \}$, i.e.\ a copy of $(N \cap M_{\delta},\in,X)$. As $N$ is definable from $X_N$ in $L_{\beta}[X_N]$, we have $N \cap M_{\delta} \prec N$. Let $\pi\colon N \cap M_{\delta} \to \bar{N}$ be the transitive collapse. Transitive collapses of well-founded relations are computed correctly in transitive models of $\ZFP$, so $\bar{N}^M_r = \bar{N}$ and $\bar{X}^M_r = \pi(X)$. We have $\pi(\xi) = \xi$ for $\xi < \delta$, $\pi(\omega_1) = \delta$, and $\pi(Z) = Z \cap \delta$ for $Z \in N \cap M_{\delta}$ with $Z \subseteq \omega_1$; in particular $\pi(X) = X \cap \delta$. By elementarity, $\bar{N} \models \ZFC^-$ and $\delta = \omega_1^{\bar{N}}$. As $\omega_1 \in N$, $N$ computes $\vec{S}$ correctly by Theorem \ref{DefinitionIndependentSequence}, so $\bar{N}$ computes $\pi(S_{\xi}) = S_{\xi} \upharpoonright \delta$ for $\xi < \delta$; recall that the nodes of $S_{\xi}$ below $\delta$ are exactly its nodes of height below $\delta$.

(3) The sets $A$ and $h$ are definable from $X$ in $N$, so $\pi(A) = A \cap \delta$ is the set of even entries of $\pi(X)$, and $\pi(h) = h \cap \delta$. As $\delta$ is closed under $\gamma \mapsto \omega\gamma + \omega$, we obtain $a^M(r) = h \cap \delta$. As $h$ is unbounded in $\omega_1$ and $h \in M_{\delta}$, the set $h \cap \delta$ is unbounded in $\delta$. For coding directly over $W$, the preceding argument in $L[X]$ gives $N\models(\ast)_1(\gamma,(x,m,k))$ for every $\gamma\in h$. By elementarity the same holds in $\bar N$ for every $\gamma\in h\cap\delta$. Lemma~\ref{nounwantedcodes} treats the branch pattern at later stages.
\end{proof}

We record a consequence of $\Psi_l(r,w)$ for an arbitrary witness $r$ and either side $l$. Suppose $\omega_1=\omega_1^L$ and $\Psi_l(r,w)$ holds. Then for unboundedly many $\gamma<\omega_1$,
\[
 \text{for every }n<\omega,\qquad
 L[r]\models\text{\ldq $S^l_{\omega\gamma+2n+w(n)}$ has a cofinal branch\rdq,}
\]
where $w(n)$ is the characteristic function of $w$ evaluated at $n$. To see this, let $M=H_\theta^{L[r,w]}$ for a regular $\theta>\omega_2$ of $L[r,w]$. Its $\omega_1$ is the true $\omega_1$. If the conclusion required of $M$ in $\Psi_l(r,w)$ failed, the transitive collapse of a countable elementary submodel of $M$ containing $r,w,\omega_1$ would witness the same failure. Hence $L[r]^M$ contains a transitive $\bar N\models\ZFC^-$ with $\omega_1\in\bar N$ which satisfies $(\ast)_l(\gamma,w)$ for unboundedly many $\gamma$. By Theorem~\ref{DefinitionIndependentSequence}, $\bar N$ computes the trees correctly. Each prescribed cofinal branch in $\bar N$ therefore belongs to $L[r]$.

Finally, a coding forcing always adds a new real, also below a condition: the set $Y$ is infinite, as its odd entries code $X_N$, so a routine density argument shows that the real added by $\operatorname{Code}(x,m,k,l,\eta)\upharpoonright c$ is not in the model over which we force. The same holds for the coding forcings over generic extensions of $W$ with the same $\omega_1$, which are defined in the next section.

\section{Allowable forcings}\label{sec:allowable}

Next we define the set of forcings which we will use in our proof.
We aim to iterate the coding forcings we defined in the last section.

\begin{definition}\label{def:allowable}
Let $W=L[G^0 \times G^1]$ be our ground model. Let $\alpha < \omega_1$ and let $F\in L$, $F: \alpha \rightarrow H(\omega_2)^L$ be a bookkeeping function.
A finite support iteration $\forceP=(\forceP_{\beta}\,:\, {\beta< \alpha})$ is called allowable (relative to the bookkeeping function $F$) if the function $F$ determines $\forceP$ inductively as follows:
 \begin{itemize}
 \item[] We assume that $\beta \ge 0$ and $\forceP_{\beta}$ is defined.
 We let $G_{\beta}$ be a $\forceP_{\beta}$-generic filter over $W$ and assume that $F(\beta)=(\dot{x},\dot{m},\dot{k}, \dot{l}, \dot{\eta})$ is a tuple of $(\forceP^0 \times\forceP^1) \ast \forceP_{\beta}$-names. We assume that $\dot{x}^{(G^0 \times G^1)\ast G_{\beta}}=:x$ is a real, $\dot{m}^{(G^0 \times G^1)\ast G_{\beta}}=:m $ and $\dot{k}^{(G^0 \times G^1)\ast G_{\beta}}=:k$ are natural numbers, $\dot{l}^{(G^0 \times G^1)\ast G_{\beta}}=:l \in \{0,1\}$ and $\dot{\eta}^{(G^0 \times G^1)\ast G_{\beta}}=:\eta$ is an ordinal $< \omega_1$. If $F(\beta)$ is not a tuple of names as above, or if its entries are not interpreted as indicated, then $\forceP(\beta)^{G_{\beta}}$ is the trivial forcing; otherwise we apply the following condition.

 \paragraph{The Freshness Condition:} We say that the \emph{Freshness Condition} fails at stage $\beta$ if there is a prior stage $\gamma < \beta$ and a $\forceP_{\gamma}$-name of a tuple $(\dot{a}', \dot{m}', \dot{k}', \dot{l}' ,\dot{\eta}') $ such that $\dot{a}'^{(G^0 \times G^1)\ast G_{\gamma}}= a \in \omega^{\omega}$, $\dot{m}'^{(G^0 \times G^1)\ast G_{\gamma}}= m' \in \omega$, $\dot{k}'^{(G^0 \times G^1)\ast G_{\gamma}}=k' \in \omega$, $\dot{l}'^{(G^0 \times G^1)\ast G_{\gamma}}=l' \in \{0,1\}$, and $\dot{\eta}'^{(G^0 \times G^1)\ast G_{\gamma}} = \eta$, and the forcing at stage $\gamma$ was defined as $\forceP(\gamma)^{G_{\gamma}} = \operatorname{Code}(a, m', k', l', \eta)$. In other words, the condition fails if the coding area $\eta$ has already been used for coding at an earlier stage of the iteration.

 Then we split into two cases:
 \begin{itemize}
 \item If the Freshness Condition fails at stage $\beta$, then we force with the trivial forcing.
 \item If the Freshness Condition holds at stage $\beta$, then let $\forceP(\beta)^{G_{\beta}}:= \operatorname{Code} (x,m,k,l,\eta)$.
 \end{itemize}
 \end{itemize}
\end{definition}

If $\forceP \in W$ is a forcing such that there is an $\alpha < \omega_1$ and an $F \in L$, $F: \alpha \rightarrow H(\omega_2)^L$ such that $\forceP$ is allowable with respect to $F$, then we often just drop the $F$ and simply say that $\forceP \in W$ is allowable.

Here $H(\omega_2)^L$ is computed in $L$. This is no restriction: $\forceP^0 \times \forceP^1$ has size $\aleph_1$ in $L$, and every allowable forcing of countable length has a dense subset of size $\aleph_1$ consisting of hereditarily countable conditions (see the normal form before Lemma \ref{DefinabilityInProperInnerModels}); so all names which we need can be chosen hereditarily of size at most $\aleph_1$.

We also allow that $F(\beta)$ contains, as an additional entry, a name $\dot{c}$ for a condition of $\operatorname{Code}(x,m,k,l,\eta)$. In this case $\forceP(\beta)^{G_{\beta}}$ is the restriction $\operatorname{Code}(x,m,k,l,\eta)\upharpoonright c$ of this coding forcing below $c$. Such a stage counts as a coding forcing which codes the same real $(x,m,k)$ into $\vec{S}^l$, using the same coding area $\eta$, and everything we say about coding forcings applies to it as well. Consequently, if $\forceP$ is allowable and $p \in \forceP$, then $\forceP\upharpoonright p$ can be presented as an allowable forcing: restrict the iterands at the finitely many stages in the support of $p$.

Informally speaking, the bookkeeping $F$ hands us at every step a real of the form $(x,m,k)$, decides whether to use $\vec{S}^0$ or $\vec{S}^1$ and uses a $\mathbb{C}(\omega_1)$-set (this set we will call a \emph{coding area}) which gives rise to a subset of $h \subset \omega_1$ which corresponds to the places where we code up the relevant branches through $\vec{S}^l$ to compute $(x,m,k)$ using the coding mechanism described in the previous section. The definition ensures that the sets $h$ determined by the coding areas are an almost disjoint family of $\omega_1$, i.e. two distinct such sets have countable intersection. The definition of allowable demands that each such coding area is used at most once in an allowable forcing, which will help to ensure that we will not accidentally code an unwanted $(x,m,k)$ into $\vec{S}$.

\begin{definition}
Let $\forceP= ((\forceP_{\alpha}, \forceP({\alpha})) \mid \alpha < \delta)$ be an allowable forcing. Let $G \subset \forceP$ be a generic filter over $W$. Then
\begin{align*}
C^G:= \{ \eta < \omega_1\mid \exists \beta < \delta \exists \dot{x},&\dot{m},\dot{k}, \dot{l}, \dot{\eta} \in W^{\forceP_{\beta}}  \\&( {\forceP}({\beta}))^{G_{\beta}} = \operatorname{Code} (\dot{x}^{G_{\beta}},\dot{m}^{G_{\beta}},\dot{k}^{G_{\beta}},\dot{l}^{G_{\beta}},\dot{\eta}^{G_{\beta}}=\eta )\}
\end{align*}
is the set of coding areas of $\forceP$ relative to $G$.
We also let
\[ C^{\forceP} := \{\eta < \omega_1 \mid \exists p \in \forceP ( p \Vdash \eta  \in C^{\dot{G}}) \}. \]
\end{definition}
The set $C^G$ is countable, as $\forceP$ has countable length. Also $C^{\forceP}$ is countable, by the ccc: $\forceP$ is a finite support iteration of almost disjoint coding forcings and trivial forcings, hence ccc, so for every $\beta<\delta$ a maximal antichain of $\forceP_{\beta}$ which decides whether a coding forcing is used at stage $\beta$, and if so, its coding area, is countable. Thus $C^{\forceP}$ is contained in a countable union of countable sets. Next we derive some properties of allowable forcings.
\begin{definition}\label{def:allowable_extension}
Let $\mathbb{P} = (\mathbb{P}_{\xi} : \xi < \delta_1)$ be an allowable forcing over $W$ with bookkeeping function $F_{\mathbb{P}} : \delta_1 \to H(\omega_2)^L$, and let $G_{\mathbb{P}}$ be a $\mathbb{P}$-generic filter over $W$. A finite support iteration $\mathbb{Q} = (\mathbb{Q}_{\beta} : \beta < \delta_2)$ in $W[G_{\mathbb{P}}]$ is called \emph{allowable over the extension $W[G_{\mathbb{P}}]$} if there exists a bookkeeping function $F_{\mathbb{Q}} : \delta_2 \to H(\omega_2)^L$ in $L$ such that at every stage $\beta < \delta_2$, if we let $G_{\mathbb{Q}_{\beta}}$ denote the $\forceQ_{\beta}$-generic filter over $W[G_{\forceP}]$, then the iterand $\mathbb{Q}(\beta)^{G_{\mathbb{Q}_\beta}}$ is determined exactly as in Definition \ref{def:allowable} using $F_{\mathbb{Q}}(\beta)$, but subject to the \emph{Global Freshness Condition}:
\begin{quote}
We say the \emph{Global Freshness Condition} fails at stage $\beta$ if the proposed coding area $\eta$ has already been used at a prior stage in $\mathbb{Q}$ (i.e., some $\alpha < \beta$) or it was already used at any stage in the forcing $\mathbb{P}$ (i.e., some $\xi < \delta_1$).
\end{quote}
If the Global Freshness Condition fails, $\mathbb{Q}(\beta)^{G_{\mathbb{Q}_\beta}}$ is defined to be the trivial forcing. Otherwise, we force with the corresponding coding forcing determined by $F_{\forceQ} (\beta)$ and the generic filters $G_{\forceP}$ and $G_{\forceQ_{\beta}}$.
\end{definition}

Here the names listed by $F_{\forceQ}$ are read as $((\forceP^0 \times \forceP^1)\ast \forceP)\ast \forceQ_{\beta}$-names, i.e.\ they are evaluated with the whole generic filter $(G^0 \times G^1)\ast G_{\forceP} \ast G_{\forceQ_{\beta}}$. We call $W[G_{\forceP}]$ an \emph{allowable extension} of $W$. Whenever we speak of an allowable forcing over an allowable extension $W[G]$, we mean an allowable forcing over $W[G]$ in the sense of Definition \ref{def:allowable_extension}; in particular, the coding areas used in $W[G]$ cannot be used again. By Lemma \ref{iteration_of_allowable} below, an allowable extension of an allowable extension of $W$ is again an allowable extension of $W$.

\begin{lemma}\label{iteration_of_allowable}
Let $\mathbb{P}$ be an allowable forcing over $W$ of length $\delta_1$, and suppose $\mathbb{P} \Vdash ``\dot{\mathbb{Q}} \text{ is an allowable forcing over } W[\dot{G}_{\mathbb{P}}] \text{ of length } \dot{\delta}_2<\omega_1"$. Then there is a countable ordinal $\delta_2$ such that the two-step iteration $\mathbb{P} \ast \dot{\mathbb{Q}}$, where $\dot{\mathbb{Q}}$ is extended by trivial iterands to length $\delta_2$, is an allowable forcing over $W$ of length $\delta_1 + \delta_2$. If the length of $\dot{\mathbb{Q}}$ is a fixed ordinal, it can be taken as $\delta_2$.
\end{lemma}

\begin{proof}
As $\mathbb{P}$ is ccc, there is a countable ordinal $\delta_2$ with $\mathbb{P} \Vdash \dot{\delta}_2 \le \check{\delta}_2$. Extending $\dot{\mathbb{Q}}$ by trivial iterands to length $\delta_2$ does not change the generic extensions; these iterands are permitted by Definition \ref{def:allowable}, using entries which are not of the required form. Let $F_1$ be the bookkeeping function for $\mathbb{P}$, and let $\dot{F}_2$ be a $\mathbb{P}$-name for a bookkeeping function witnessing that $\dot{\mathbb{Q}}$ is allowable over $W[\dot{G}_{\mathbb{P}}]$. As $W=L[G^0 \times G^1]$, $\mathbb{P}$-names in $W$ are given by $(\forceP^0 \times \forceP^1)\ast\mathbb{P}$-names in $L$. Hence we can define the combined bookkeeping function $F \in L$ with domain $\delta_1 + \delta_2$ by $F \restriction \delta_1 := F_1$, and by letting $F(\delta_1 + \beta)$ be the $(\forceP^0 \times \forceP^1)\ast\mathbb{P}\ast\mathbb{Q}_{\beta}$-name canonically obtained from the $\mathbb{P}$-name $\dot{F}_2(\beta)$ for a $((\forceP^0 \times \forceP^1)\ast\mathbb{P})\ast\mathbb{Q}_{\beta}$-name. In particular, $F$ is chosen before fixing a generic filter for $\mathbb{P}$.

Let $\mathbb{R} = (\mathbb{R}_{\xi} : \xi < \delta_1 + \delta_2)$ be the finite support iteration generated by $F$ over $W$. We show $\mathbb{R} \cong \mathbb{P} \ast \dot{\mathbb{Q}}$ by induction on $\xi$. For $\xi \le \delta_1$, $F \restriction \delta_1 = F_1$ yields $\mathbb{R}_{\xi} = \mathbb{P}_{\xi}$. For $\xi = \delta_1 + \beta$, the entry $F(\xi)$ has the same interpretation as $\dot{F}_2(\beta)$, and by Definition \ref{def:allowable_extension}, the Global Freshness Condition for $\dot{\mathbb{Q}}$ at stage $\beta$ holds if and only if the proposed coding area $\eta$ satisfies $\eta \notin C^{G_{\mathbb{P}}} \cup C^{G_{\mathbb{Q}_\beta}}$, which is the Freshness Condition for $\mathbb{R}$ at stage $\delta_1 + \beta$. Consequently, the iterands coincide, and the limit steps use finite support. We conclude $\mathbb{R} \cong \mathbb{P} \ast \dot{\mathbb{Q}}$, proving that the two-step iteration is allowable over $W$.
\end{proof}

 \begin{lemma} \label{FirstPropertiesOfAllowableForcings}

\begin{enumerate}
\item If $\forceP=(\forceP(\beta) \, : \, \beta < \delta) \in W$ is allowable then for every $\beta < \delta$, $\forceP_{\beta} \Vdash| \forceP(\beta)|\le \aleph_1$, thus every factor of $\forceP$ is forced to have size at most $\aleph_1$.
\item Every allowable forcing over $W$ is ccc and thus preserves cardinals.
\item Every allowable forcing over $W$ preserves $\CH$. Furthermore, if $\forceP= (\forceP(\alpha) \, : \, \alpha < \omega_1) \in W$ is an $\omega_1$-length iteration such that each initial segment of the iteration is allowable over $W$, then $W^{\forceP} \models \CH$.
\end{enumerate}
\end{lemma}
\begin{proof}
The assertions follow immediately from the definition modulo some well-known results.
\end{proof}
Every allowable forcing is determined by its bookkeeping function $F$. The function $F$ was not required to be a list of names for reals, but it can be replaced by a bookkeeping function in $L$ which yields the same forcing and whose entries are tuples of names for a real, two natural numbers, an element of $\{0,1\}$, a countable ordinal and a condition: replace each entry by names for its interpreted values, and where the entry is not as required in Definition \ref{def:allowable}, use a fixed tuple which prescribes the trivial forcing. We always assume that $F$ is of this form.

Moreover, we fix the following presentation. Let $\Lambda := [\omega]^{<\omega} \times [\omega_1]^{<\omega}$. Identifying $d_{\xi} \in D$ with $\xi$, every almost disjoint coding forcing $\mathbb{A}_D(Y)$ is an ordering of $\Lambda$; its restriction below a condition $c$ is presented as a preordering of $\Lambda$ by identifying every element of $\Lambda$ which is not below $c$ with $c$; and the trivial forcing is presented as the preordering of $\Lambda$ in which all elements are equivalent. So all iterands of an allowable forcing $\forceP$ of length $\delta$ are presented on $\Lambda$, and $\forceP$ is equivalent to the set $\operatorname{Fn}(\delta,\Lambda)$ of finite partial functions from $\delta$ to $\Lambda$, ordered by $q \le p$ if and only if $\operatorname{dom}(p) \subseteq \operatorname{dom}(q)$ and $q \upharpoonright \beta \Vdash_{\forceP_{\beta}} q(\beta) \le p(\beta)$ for every $\beta \in \operatorname{dom}(p)$. As $\forceP$ is ccc, every entry of $F$ is, in $W$, given by nice $\forceP_{\beta}$-names with respect to this presentation, i.e.\ by countable sets of hereditarily countable objects. Fixing a countable ordinal $\kappa$ containing all ordinals involved in these names and a bijection between $\omega$ and $\kappa$ in $L$, the whole list of these names is coded by a single real $z_{\forceP} \in W$. The real $z_{\forceP}$ does not code the sets $Y$ used by the coding forcings.

We shall utilize this normal form in showing that every allowable forcing can in fact be defined already in a proper inner model of $W$.

\begin{lemma}\label{DefinabilityInProperInnerModels}
Let $\forceP \in W$ be an allowable forcing of length $\delta<\omega_1$. Then there is an uncountable, co-uncountable subset $I \subset \omega_1$ and a countable subset $J \subset \omega_1$ such that $\forceP$ can successfully be defined already in an inner model of $W$ of the form $L[(G^0 \upharpoonright I) \times (G^1 \upharpoonright J)]$, where $G^0 \upharpoonright I$ is just the restriction of the generic $G^0$ to coordinates which are in $I$ and $G^1 \upharpoonright J$ being defined similarly. More precisely, if $J \in L$ is countable with $C^{\forceP} \subseteq J$, if $K \in L$ is countable with $z_{\forceP} \in L[G^0 \upharpoonright K]$, and if $H_J$ denotes the set of indices of all trees from $\vec{S}$ in the $\omega$-blocks indexed by $\bigcup_{\eta \in J} h_{\eta}$, where $h_{\eta} = \{\rho(g_{\eta} \cap \alpha) : \alpha < \omega_1\}$ as in Section \ref{sec:coding}, then $I := H_J \cup K$ works.

Moreover, let $G \subset \forceP$ be generic over $W$ and let $\xi \notin K$ be such that no coding forcing of $G$ uses the branch $b_{\xi}$ through $S_{\xi}$, i.e.\ $b_{\xi}$ is not among the branches coded into the set $X$ of a coding forcing used by $G$ (see Section \ref{sec:coding}). Put $N_{\xi} := L[G^1][G^0 \upharpoonright (\omega_1 \setminus \{\xi\})]$. Then $W[G] = N_{\xi}[G][b_{\xi}]$, the tree $S_{\xi}$ is Suslin in $N_{\xi}[G]$, and $W[G]$ and $N_{\xi}[G]$ have the same reals. In particular, there is no real $r \in W[G]$ such that $L[r]$ contains a cofinal branch through $S_{\xi}$.
\end{lemma}
\begin{proof}
We first show:
\begin{itemize}
\item[$(\ast)$] Let $J' \in L$, $I' \in L[G^1 \upharpoonright J']$ be subsets of $\omega_1$, let $N := L[G^1 \upharpoonright J'][G^0 \upharpoonright I']$ and assume that $z_{\forceP} \in N$. Let $p \in \forceP$ force that every coding forcing of $\forceP$ has its coding area in $J'$ and uses only branches $b_{\zeta}$ with $\zeta \in I'$. Then there is a finite support iteration $\forceQ \in N$ of length $\delta$, with iterands presented on $\Lambda$, such that $\forceP \restriction p = \forceQ \restriction p$ as preorders on $\operatorname{Fn}(\delta,\Lambda)$.
\end{itemize}
To see $(\ast)$, define $\forceQ$ in $N$ by recursion on $\beta \le \delta$. The entry of the bookkeeping at stage $\beta$ is read off from $z_{\forceP}$ as a nice $\forceQ_{\beta}$-name. If its value is a valid entry whose coding area $\eta$ satisfies the Freshness Condition and lies in $J'$, and whose set $X$ (see Section \ref{sec:coding}) only involves branches $b_{\zeta}$ with $\zeta \in I'$, then $\forceQ(\beta)$ is the indicated (possibly restricted) coding forcing, computed in $N$; note that then $X \in N$ and that $Y$ is computed canonically in $L[X]$. Otherwise $\forceQ(\beta)$ is trivial. We show by induction on $\beta \le \delta$ that for all $r, q \in \operatorname{Fn}(\beta,\Lambda)$:
\begin{align*}
& r \le_{\forceP_{\beta}} p \upharpoonright \beta \ \Longleftrightarrow\ r \le_{\forceQ_{\beta}} p \upharpoonright \beta, \tag{$\dagger_1$}\\
& \text{if } r \le_{\forceP_{\beta}} p \upharpoonright \beta, \text{ then } \bigl( r \le_{\forceP_{\beta}} q \ \Longleftrightarrow\ r \le_{\forceQ_{\beta}} q \bigr). \tag{$\dagger_2$}
\end{align*}
For $\beta=\delta$, $(\dagger_1)$ says that the cones below $p$ in $\forceP$ and in $\forceQ$ consist of the same conditions, and $(\dagger_2)$ that their orders agree; this is $(\ast)$. Both assertions only involve finitely many coordinates, so the limit steps follow from the definition of the order on $\operatorname{Fn}(\beta,\Lambda)$. Assume $(\dagger_1)$ and $(\dagger_2)$ at $\beta$. Then the cones below $p \upharpoonright \beta$ in $\forceP_{\beta}$ and in $\forceQ_{\beta}$ coincide as preorders, and for a filter $H$ which is generic for this cone, the upward closures of $H$ in $\forceP_{\beta}$ and in $\forceQ_{\beta}$ coincide by $(\dagger_2)$; so every nice name coded in $z_{\forceP}$ has the same value, whether it is read as a $\forceP_{\beta}$-name or as a $\forceQ_{\beta}$-name. Moreover, $p \upharpoonright \beta$ forces in $\forceP_{\beta}$ the assumption on $p$ for the stage $\beta$, as an extension of $p \upharpoonright \beta$ forcing the opposite could be combined with $p \upharpoonright [\beta,\delta)$. So below $p \upharpoonright \beta$, the entry at stage $\beta$, its freshness and its set $X$, and hence the iterand at stage $\beta$, are the same for both iterations. Consequently, if $r \upharpoonright \beta \le p \upharpoonright \beta$ and $u,v \in \Lambda$, then $r \upharpoonright \beta$ forces $u \le v$ in the iterand at stage $\beta$ of $\forceP$ if and only if it does so for $\forceQ$. As $r \le q$ in $\operatorname{Fn}(\beta+1,\Lambda)$ means $r \upharpoonright \beta \le q \upharpoonright \beta$ and, if $\beta \in \operatorname{dom}(q)$, that $\beta \in \operatorname{dom}(r)$ and $r \upharpoonright \beta$ forces $r(\beta) \le q(\beta)$, we obtain $(\dagger_1)$ and $(\dagger_2)$ at $\beta+1$.

For the first assertion, note that $K$ exists: by Easton's Lemma (see Lemma 15.19 from \cite{Jech}) the $\sigma$-closed forcing $\forceP^1$ adds no reals over $L[G^0]$, so $z_{\forceP} \in L[G^0]$, and a nice name for $z_{\forceP}$ with respect to the ccc forcing $\forceP^0$ involves only countably many coordinates. Now apply $(\ast)$ with $p=1$, $J' := J$ and $I' := H_J \cup K$: every coding forcing of $\forceP$ has its coding area in $C^{\forceP} \subseteq J$ and only uses branches with indices in $H_J$. So $\forceP = \forceQ$ is defined in $L[G^1 \upharpoonright J][G^0 \upharpoonright I]$. In other words, over $U := L[G^1 \upharpoonright J]$ the forcing producing $W[G]$, for $G \subset \forceP$, factors as
\[ \bigl[(\forceP^0 \upharpoonright I) \ast \forceP\bigr] \times \bigl(\forceP^0 \upharpoonright (\omega_1 \setminus I)\bigr) \times \bigl(\forceP^1 \upharpoonright (\omega_1 \setminus J)\bigr)^L, \]
where $\forceP^0 \upharpoonright (\omega_1 \setminus I)$ is a finite support product. We may assume $J \ne \emptyset$, so $I$ is uncountable. To see that $\omega_1 \setminus I$ is uncountable, choose $\eta^* \notin J$. As $J$ is countable and the sets $h_{\eta}$ are pairwise almost disjoint, $h_{\eta^*} \cap \bigcup_{\eta \in J} h_{\eta}$ is countable, so uncountably many $\omega$-blocks indexed by $h_{\eta^*}$ are disjoint from $H_J$, and removing the countable set $K$ leaves uncountably many indices outside of $I$.

For the second assertion, choose $p \in G$ which forces that no coding forcing uses $b_{\xi}$, and apply $(\ast)$ with $N := N_{\xi}$, i.e.\ $J' = \omega_1$ and $I' = \omega_1 \setminus \{\xi\}$; note that $z_{\forceP} \in N_{\xi}$ as $\xi \notin K$. So $\forceP \restriction p = \forceQ \restriction p$ for some $\forceQ \in N_{\xi}$. As $\forceP^1$ is $\sigma$-closed in $L$, $\vec{S}$ is still an independent sequence of Suslin trees in $L[G^1]$, so $S_{\xi}$ is Suslin in $N_{\xi}$, which is a generic extension of $L[G^1]$ by the finite support product of the other trees. The forcing $\forceQ$ is a finite support iteration of $\sigma$-centered forcings and hence has the Knaster property, so $S_{\xi} \times (\forceQ \restriction p)$ is ccc in $N_{\xi}$. As $W = N_{\xi}[b_{\xi}]$ and $G$ is $(\forceQ \restriction p)$-generic over $W$, the pair $(b_{\xi}, G)$ is $S_{\xi} \times (\forceQ \restriction p)$-generic over $N_{\xi}$. Hence $S_{\xi}$ is Suslin in $N_{\xi}[G]$, and $W[G] = N_{\xi}[G][b_{\xi}]$ is obtained by forcing with the Suslin tree $S_{\xi}$, which adds no reals. Finally, if $L[r]$ contained a cofinal branch through $S_{\xi}$ for a real $r \in W[G]$, this branch would belong to $N_{\xi}[G]$, which is impossible.
\end{proof}

Let $\forceP= (\forceP(\beta) \, : \, \beta < \delta)$ be an allowable forcing with respect to some $F \in L$.
The set of  (names of) reals which are enumerated by $F$, and where the second case in the definition of allowable applies, i.e. we actually used a coding forcing, is dubbed \emph{the set of reals which are coded by $\forceP$}. That is, for every $\beta$ at which a coding forcing is used, if we let $\dot{x}_{\beta}$ be the (name) of a real listed by $F(\beta)$ and if we let $G \subset \forceP$ be a generic filter over $W$ and finally if we let
$ \dot{x}_{\beta}^G =:x_{\beta}$,  then we say that
$\{ x_{\beta} \, : \, \beta < \delta \}$ is the set of reals coded by $\forceP$ and $G$ (though we will suppress the $G$), and that $x_{\beta}$ is coded into $\vec{S}^{l_{\beta}}$ at stage $\beta$, where $l_{\beta}$ is the side listed by $F(\beta)$.

Next we show, that allowable forcings will not add accidentally new elements to the set of reals defined by the $\Sigma^1_3$-formula \[\Phi_1((x,m,k)):= \exists r \Psi_1 (r,(x,m,k)), \]
where $\Psi_1$ is defined in the previous section and by the same argument will also not add non-intended members to the set defined by
\[ \Phi_0((x,m,k)) := \exists r \Psi_0( r,(x,m,k)).\] We let
\[ \Phi(x) \equiv \Phi_0 (x) \lor \Phi_1 (x). \]
To utilize already established jargon we will say that $``x$ is coded into $\vec{S^i}$'' whenever $\Phi_i (x)$ is true. Likewise, if $\lnot \Phi_i (x)$ holds true, then we say that ``$x$ is not coded into $\vec{S^i}$''.

\begin{lemma}\label{nounwantedcodes}
If $\forceP \in W$ is allowable, $\forceP=(\forceP_{\beta} \, : \, \beta < \delta)$, $G \subset \forceP$ is generic over $W$ and $\{ x_{\beta} \, : \, \beta < \delta\}$ is the set of reals which are coded by $\forceP$. Let $\Phi(v_0)$ be the distinguished formula from above. Then in $W[G]$, for $i \in \{0,1\}$, a real satisfies $\Phi_i(v_0)$ if and only if it is coded into $\vec{S}^i$ at some stage of $\forceP$. In particular, the set of reals which satisfy $\Phi(v_0)$ is exactly $\{ x_{\beta} \, : \, \beta < \delta\}$.
\end{lemma}
\begin{proof}
Let $B\subseteq\delta$ be the set of coding stages. For $\beta\in B$ write $l_\beta$ for its side, $\eta_\beta$ for its coding area and $h_\beta=\{\rho(g_{\eta_\beta}\cap\alpha):\alpha<\omega_1\}$. By the Freshness Condition the $\eta_\beta$ are distinct, so the $h_\beta$ are pairwise almost disjoint.

Fix $\beta\in B$ and work in $W[G_\beta]$, where $G_\beta=G\cap\forceP_\beta$. Put $w=x_\beta$ and $l=l_\beta$. Apply Lemma~\ref{DefinabilityInProperInnerModels} to $\forceP_\beta$ and take the countable set $K\in L$ given there. Let $I^l_\gamma$ be the set of indices in the full sequence $\vec S$ of the trees $S^l_{\omega\gamma+j}$, $j<\omega$. The set
\[
 E_\beta=\left(h_\beta\cap\bigcup_{\alpha\in B\cap\beta}h_\alpha\right)
          \cup\{\gamma\in h_\beta:I^l_\gamma\cap K\ne\varnothing\}
\]
is countable. Indeed $\beta$ is countable, each intersection $h_\beta\cap h_\alpha$ is countable for $\alpha\in B\cap\beta$, and the blocks $I^l_\gamma$ are disjoint.

Let $\gamma\in h_\beta\setminus E_\beta$ and $n<\omega$. Write $\xi$ for the full index of $S^l_{\omega\gamma+2n+(1-w(n))}$. No earlier coding uses $b_\xi$, and $\xi\notin K$. By Lemma~\ref{DefinabilityInProperInnerModels}, $S_\xi$ is Suslin in $N_\xi[G_\beta]$ and this model has the same reals as $W[G_\beta]$. It contains $w$, $g_{\eta_\beta}$ and the sequence of branches with coordinate $\xi$ omitted. The code $X$ for this step is formed there from these objects, since the branch at $\xi$ is not prescribed. Thus $L[X]$ has no cofinal branch through $S_\xi$. The prescribed branch through $S^l_{\omega\gamma+2n+w(n)}$ is recovered from $X$. This proves $(\ast)_l(\gamma,w)$ in $L[X]$ for every $\gamma\in h_\beta\setminus E_\beta$.

The model $N\in L[X]$ used to define $Y$ recovers $h_\beta,w$ and all prescribed branches from $X$, and computes the trees correctly. Thus the set of $\gamma\in h_\beta$ at which $N$ satisfies $(\ast)_l(\gamma,w)$ contains $h_\beta\setminus E_\beta$. In particular
\[
 N\models\text{\ldq $\{\gamma\in h_\beta:(\ast)_l(\gamma,w)\}$
 is unbounded in $\omega_1$\rdq.}
\]
This statement uses only $h_\beta,w$ as parameters. Now let $r_\beta$ be the generic real at stage $\beta$, also allowing a restriction below a condition. For every countable transitive $M$ as in $\Psi_l$ containing $r_\beta$, put $\delta_M=\omega_1^M$. By Lemma~\ref{lem:local_decoding}(1)--(3), $L[r_\beta]^M$ contains the collapse $\bar N^M_{r_\beta}$ of the elementary submodel $N\cap M_{\delta_M}$ of $N$. The collapse sends $\omega_1$ to $\delta_M$, $h_\beta$ to $h_\beta\cap\delta_M$, and fixes $w$. By elementarity, the set of $\gamma$ satisfying the pattern in the collapse is unbounded in $\delta_M$. This set belongs to the collapse and has size $\aleph_1^M$ in $L[r_\beta]^M$. Hence $\Psi_l(r_\beta,w)$ holds in $W[G_{\beta+1}]$ and remains true in $W[G]$ by Shoenfield absoluteness.

Conversely, suppose that $\Psi_l(r,w)$ holds in $W[G]$, but that $w \ne x_{\beta}$ for every $\beta \in B$ with $l_{\beta} = l$. By the discussion at the end of Section \ref{sec:coding}, there is an unbounded set $h\subseteq\omega_1$ such that $L[r]$ contains a cofinal branch through $S^l_{\omega\gamma+2n+w(n)}$ for every $\gamma \in h$ and $n \in \omega$, where $w(n)=1$ if $n \in w$ and $w(n)=0$ otherwise. The set $Z := \bigcup_{\beta \ne \beta' \in B} (h_{\beta} \cap h_{\beta'})$ is countable. Let $\gamma \in h \setminus Z$. Then $\gamma \in h_{\beta}$ for at most one $\beta \in B$. If there is such a $\beta$ with $l_{\beta} = l$, choose $n$ with $w(n) \ne x_{\beta}(n)$; otherwise let $n:=0$. In both cases the branch through $S^l_{\omega\gamma+2n+w(n)}$ is not used by any coding forcing of $G$. As distinct $\gamma$ give distinct trees, we obtain uncountably many such trees, and we can pick one of them, $S_{\xi}$, with $\xi \notin K$, where $K$ is as in Lemma \ref{DefinabilityInProperInnerModels}. By the second part of Lemma \ref{DefinabilityInProperInnerModels}, $L[r]$ cannot contain a cofinal branch through $S_{\xi}$. This is the desired contradiction.
\end{proof}

\subsection{Allowable forcings over suitable pairs}\label{sec:alpha}

The forcings used in Sections~\ref{sec:four}--\ref{sec:up} are allowable over suitable pairs, as defined in this section, and satisfy further restrictions which are introduced there.
Since the restrictions will be introduced along our iteration, we define allowability over generic extensions of $W$. To obtain the closure under products in Section \ref{sec:products}, we also keep track of the coding areas which must not be used any more, by a countable set in $L$ which covers them.

\begin{definition}\label{def:alpha_allowable_extension}
\begin{enumerate}
\item A \emph{suitable pair} is a pair $(M,J)$, where $M = W[G]$ for a filter $G$ which is generic over $W$ for an allowable forcing, and $J \in L$ is a countable set of coding areas with $C^G \subseteq J$. We let $N_J := L[G^0][G^1 \upharpoonright J][G]$.
\item Let $(M,J)$ be a suitable pair. A forcing $\forceQ = (\forceQ_{\gamma} : \gamma < \delta)$ is \emph{allowable over $(M,J)$} if it is allowable over $M$ (Definition \ref{def:allowable_extension}), and there are a countable set $K \in L$ with $K \cap J = \emptyset$ and $\forceQ_{\gamma}$-names $J_{\gamma}$, $\gamma \le \delta$, for countable sets in $L$ such that in every $\forceQ$-generic extension of $M$:
\begin{itemize}
\item[(i)] $J_0 = J$ and $J_{\gamma} \subseteq J_{\gamma'} \subseteq J \cup K$ for $\gamma \le \gamma' \le \delta$;
\item[(ii)] if $\forceQ(\gamma)$ is a coding forcing with coding area $\eta$, then $\eta \in K \setminus J_{\gamma}$ and $\eta \in J_{\gamma+1}$.
\end{itemize}
We say that $\forceQ$ \emph{uses coding areas from $K$}, and we call $J_{\gamma}$ the \emph{cover} at stage $\gamma$ and $J_{\delta}$ the \emph{final cover}. We regard $K$ and $(J_{\gamma})_{\gamma \le \delta}$ as part of $\forceQ$.
\end{enumerate}
\end{definition}

By (i) and (ii), the cover $J_{\gamma}$ contains $J$ and all coding areas used by $\forceQ$ before stage $\gamma$. So $(M[G_{\gamma}], J_{\gamma})$ is again a suitable pair, and the Global Freshness Condition holds automatically. A coding area in $J_{\gamma}$ can never be used at a later stage, even if it was not used before.

For example, let $(M,J)$ be a suitable pair, let $x \in M$ be a real, $m,k \in \omega$, $l \in \{0,1\}$ and $\eta \in \omega_1 \setminus J$. Then the one step forcing $\operatorname{Code}(x,m,k,l,\eta)$, and also its restriction below a condition, is allowable over $(M,J)$, with $K = \{\eta\}$ and the covers $J_0 = J$ and $J_1 = J \cup \{\eta\}$. So over a suitable pair every real of the form $(x,m,k)$ can be coded, on every coding area outside of the cover. The iteration of trivial forcings of any countable length is allowable over $(M,J)$ as well, with $K = \emptyset$ and all covers equal to $J$.

The two step iteration of allowable forcings over suitable pairs results in an allowable forcing again, and allowability is inherited by initial segments, by the remaining parts of the iteration and by restrictions below conditions.

\begin{lemma}\label{iteration_of_alpha_allowable}
Let $(M,J)$ be a suitable pair.
\begin{enumerate}
\item Let $\forceP$ be an allowable forcing over $(M,J)$ of length $\delta_1$ which uses coding areas from $K_1$ and has final cover $B$, and suppose $\forceP \Vdash ``\dot{\forceQ}$ is allowable over $(M[\dot{G}_{\forceP}],B)$ of length $\dot{\delta}_2 < \omega_1$''. Then there is a countable ordinal $\delta_2$ such that the two step iteration $\forceP \ast \dot{\forceQ}$, where $\dot{\forceQ}$ is extended by trivial iterands to length $\delta_2$, is allowable over $(M,J)$ of length $\delta_1 + \delta_2$, with the covers of $\forceP$ followed by those of $\dot{\forceQ}$.
\item If $\forceQ$ is allowable over $(M,J)$, with covers $(J_{\gamma})_{\gamma \le \delta}$ and coding areas from $K$, and $\gamma_0 \le \delta$, then the initial segment $\forceQ_{\gamma_0}$ is allowable over $(M,J)$, and for every $\forceQ_{\gamma_0}$-generic filter $G_{\gamma_0}$ over $M$ the remaining forcing $\forceQ_{\gamma_0,\delta}$, with the covers $(J_{\gamma})_{\gamma_0 \le \gamma \le \delta}$ and coding areas from $K \setminus J_{\gamma_0}$, is allowable over $(M[G_{\gamma_0}],J_{\gamma_0})$.
\item If $\forceQ$ is allowable over $(M,J)$ and $q \in \forceQ$, then $\forceQ \upharpoonright q$, presented as in the remarks after Definition \ref{def:allowable}, is allowable over $(M,J)$ with the same covers.
\end{enumerate}
\end{lemma}
\begin{proof}
For (1), by the ccc there are countable ordinals $\delta_2$ and $\theta$ such that $\forceP$ forces $\dot{\delta}_2 \le \delta_2$ and that $\dot{\forceQ}$ uses coding areas from a subset of $\theta$; let $K := K_1 \cup (\theta \setminus J)$. As in the proof of Lemma \ref{iteration_of_allowable} we obtain a bookkeeping function for $\forceP \ast \dot{\forceQ}$, and conditions (i) and (ii) of Definition \ref{def:alpha_allowable_extension} hold for $K$ and the indicated covers. The proof of (2) is the same, and (3) holds as restricting iterands below conditions changes neither the intermediate models, nor the covers, nor the coded reals, sides and coding areas.
\end{proof}

The bookkeeping function of an allowable forcing over an allowable extension $M$ of $W$ is an element of $L$. The next lemma shows that it can be replaced by a sequence of names in $M$. So allowability over $M$ is definable in $M$, and we can form limits of allowable forcings. For a suitable pair $(M,J)$ the Global Freshness Condition can be replaced by a condition on the covers, so that allowability over $(M,J)$ only depends on $G^0 \times G^1$ and $J$.

\begin{lemma}\label{lem:names_in_M}
Let $M = W[G]$ be an allowable extension of $W$ and let $\forceQ = (\forceQ_{\beta} : \beta < \delta) \in M$ be a finite support iteration of countable length. Then the following are equivalent:
\begin{itemize}
\item[(a)] $\forceQ$ is allowable over $M$.
\item[(b)] There is a sequence $(\tau_{\beta} \mid \beta < \delta) \in M$ such that $\tau_{\beta}$ is a $\forceQ_{\beta}$-name for every $\beta < \delta$, and at every stage $\beta < \delta$ the iterand $\forceQ(\beta)^{G_{\forceQ_{\beta}}}$ is determined as in Definition \ref{def:allowable_extension}, using the value $\tau_{\beta}^{G_{\forceQ_{\beta}}}$ in place of the value of $F_{\forceQ}(\beta)$.
\end{itemize}
Now let $(M,J)$ be a suitable pair, let $K \in L$ be countable with $K \cap J = \emptyset$, and let $J_{\gamma}$, $\gamma \le \delta$, be $\forceQ_{\gamma}$-names for countable sets in $L$ which satisfy (i) of Definition \ref{def:alpha_allowable_extension}. Then $\forceQ$ is allowable over $(M,J)$ with this set $K$ and these covers if and only if the following holds:
\begin{itemize}
\item[(c)] There is a sequence $(\tau_{\beta} \mid \beta < \delta) \in M$ such that $\tau_{\beta}$ is a $\forceQ_{\beta}$-name for every $\beta < \delta$, and at every stage $\beta < \delta$: if the value of $\tau_{\beta}$ is a tuple as required in Definition \ref{def:allowable} whose coding area $\eta$ belongs to $K \setminus J_{\beta}$, then $\forceQ(\beta)^{G_{\forceQ_{\beta}}}$ is the (possibly restricted) coding forcing which is given by this tuple, and $\eta \in J_{\beta+1}$; otherwise $\forceQ(\beta)^{G_{\forceQ_{\beta}}}$ is the trivial forcing.
\end{itemize}
Consequently, the forcings which are allowable over $M$ form a class which is definable in $M$ from the parameters $G^0 \times G^1$ and $C^G$, by a formula which does not depend on $M$. The forcings which are allowable over $(M,J)$, together with their sets $K$ and their covers, form a class which is definable in $M$ from the parameters $G^0 \times G^1$ and $J$, again by a formula which does not depend on $M$.
\end{lemma}
\begin{proof}
Let $\forceP$ be the allowable forcing with $G \subset \forceP$. If $F_{\forceQ} \in L$ witnesses (a), let $\tau_{\beta} \in M$ be the $\forceQ_{\beta}$-name which is obtained from $F_{\forceQ}(\beta)$ by evaluating it partially with $(G^0 \times G^1) \ast G$. Then the sequence of these names is as in (b). Conversely, let $(\tau_{\beta} \mid \beta < \delta)$ be as in (b). As $M = L[G^0 \times G^1][G]$, this sequence has a $(\forceP^0 \times \forceP^1) \ast \forceP$-name in $L$. As in the proof of Lemma \ref{iteration_of_allowable}, we let $F_{\forceQ}(\beta)$ be the $((\forceP^0 \times \forceP^1) \ast \forceP) \ast \forceQ_{\beta}$-name which is canonically obtained from this name. Then $F_{\forceQ} \in L$, and the value of $F_{\forceQ}(\beta)$ under $(G^0 \times G^1) \ast G \ast G_{\forceQ_{\beta}}$ is $\tau_{\beta}^{G_{\forceQ_{\beta}}}$. By the remarks after Definition \ref{def:allowable}, the entries of $F_{\forceQ}$ can be chosen in $H(\omega_2)^L$. So $F_{\forceQ}$ witnesses (a).

Now let $(M,J)$, $K$ and the names $J_{\gamma}$ be as in the second part. Assume first that $\forceQ$ is allowable over $(M,J)$ with $K$ and these covers, and let $(\tau_{\beta} \mid \beta < \delta)$ be as in (b). Let $\tau'_{\beta}$ be a $\forceQ_{\beta}$-name for the value of $\tau_{\beta}$, if a coding forcing is used at stage $\beta$, and for a fixed tuple which is not as required in Definition \ref{def:allowable} otherwise. If a coding forcing with coding area $\eta$ is used at stage $\beta$, then $\eta \in K \setminus J_{\beta}$ and $\eta \in J_{\beta+1}$ by (ii) of Definition \ref{def:alpha_allowable_extension}. So the sequence $(\tau'_{\beta} \mid \beta < \delta)$ is as in (c). Conversely, let $(\tau_{\beta} \mid \beta < \delta)$ be as in (c). By (i) and by induction on $\beta$, the cover $J_{\beta}$ contains $J \supseteq C^G$ and all coding areas which are used before stage $\beta$. So the Global Freshness Condition holds at every stage $\beta$ at which the value of $\tau_{\beta}$ is a tuple as required whose coding area belongs to $K \setminus J_{\beta}$. Let $\tau''_{\beta}$ be a $\forceQ_{\beta}$-name for the value of $\tau_{\beta}$ in this case, and for a fixed tuple which is not as required otherwise. Then the sequence $(\tau''_{\beta} \mid \beta < \delta)$ is as in (b), so $\forceQ$ is allowable over $M$, and (ii) holds by (c).

For the last assertions, note that (b) and (c) are statements about $\forceQ$ in $M$. The coding forcings are defined from the coded real, the side, the coding area and the generic branches, and the latter two are given by $G^0 \times G^1$. The Global Freshness Condition at a stage only refers to $C^G$ and to the coding areas which are used at the earlier stages of $\forceQ$, and (c) does not refer to it at all. Finally, the set $K$ and the covers in Definition \ref{def:alpha_allowable_extension} are elements of $L$ and names for elements of $L$, and $L$ is definable in $M$.
\end{proof}

The normal form before Lemma \ref{DefinabilityInProperInnerModels} applies to the forcings which are allowable over a suitable pair as well. Let $\forceQ$ be allowable over $(M,J)$, of length $\delta$, with the set $K$ and the covers $(J_{\gamma})_{\gamma \le \delta}$. In $M$ the iterands of $\forceQ$ are presented on $\Lambda$, and $\forceQ$ is equivalent to $\operatorname{Fn}(\delta,\Lambda)$ with the ordering described there. As $\forceQ$ is ccc, the names $\tau_{\beta}$ in (c) can be chosen as tuples of nice $\forceQ_{\beta}$-names with respect to this presentation, and each cover $J_{\gamma}$ is given by a countable maximal antichain of $\forceQ_{\gamma}$ together with a countable set in $L$ for each of its elements. We call the triple which consists of $K$, of the sequence $(\tau_{\beta} \mid \beta < \delta)$ and of the sequence of these names for the covers a \emph{code} of $\forceQ$. A code is a hereditarily countable set, and it determines $\forceQ$, by recursion on $\beta$ as in (c). So the forcings which are allowable over $(M,J)$ are given by codes in $H(\omega_1)^M$, and the set of these codes is definable in $M$ from the parameters $G^0 \times G^1$ and $J$, by a formula which does not depend on $M$. Codes for the forcings in Lemma \ref{iteration_of_alpha_allowable} are obtained in the usual way: for an initial segment we restrict the sequences, for the remaining forcing we evaluate the names partially with the generic filter of the initial segment, and for a two step iteration we compose the names as in the proof of Lemma \ref{iteration_of_allowable} and replace them by nice names.

\begin{lemma}\label{lem:limits}
Let $(M,J)$ be a suitable pair, let $\lambda < \omega_1$ be a limit ordinal and let $\forceQ = (\forceQ_{\gamma} : \gamma < \lambda) \in M$ be a finite support iteration. Assume that there are $\forceQ_{\gamma}$-names $J_{\gamma}$, $\gamma < \lambda$, such that for every $\gamma_0 < \lambda$ the initial segment $\forceQ_{\gamma_0}$ is allowable over $(M,J)$ with the covers $(J_{\gamma})_{\gamma \le \gamma_0}$. Then there is a countable ordinal $\theta$ such that $\forceQ$ forces $J_{\gamma} \subseteq \theta$ for every $\gamma < \lambda$. For every such $\theta$, the forcing $\forceQ$ is allowable over $(M,J)$, with the covers $(J_{\gamma})_{\gamma < \lambda}$ and the final cover $J_{\lambda} := \theta$, and it uses coding areas from $\theta \setminus J$.
\end{lemma}
\begin{proof}
The forcing $\forceQ$ is ccc, as it is a finite support iteration all of whose proper initial segments are ccc. In every $\forceQ$-generic extension of $M$ the union of the sets $J_{\gamma}$, $\gamma < \lambda$, is a countable set of countable ordinals. So by the ccc there is a countable ordinal $\theta$ as desired. Let $K := \theta \setminus J$. Then $K \in L$ is countable, $K \cap J = \emptyset$ and $\theta = J \cup K$, as $J = J_0 \subseteq \theta$.

Each stage of $\forceQ$ is a stage of one of the forcings $\forceQ_{\gamma_0}$, $\gamma_0 < \lambda$, with the same intermediate model and the same cover. So $J_0 = J$ and $J_{\gamma} \subseteq J_{\gamma'} \subseteq \theta = J \cup K$ for $\gamma \le \gamma' \le \lambda$, which is (i) of Definition \ref{def:alpha_allowable_extension}. If $\forceQ(\gamma)$ is a coding forcing with coding area $\eta$, then $\eta \notin J_{\gamma}$ and $\eta \in J_{\gamma+1}$ by (ii) for $\forceQ_{\gamma+1}$. Hence $\eta \in \theta \setminus J_{\gamma} \subseteq K \setminus J_{\gamma}$, which is (ii).

It remains to show that $\forceQ$ is allowable over $M$. For $\beta < \lambda$ we choose, in $M$, a $\forceQ_{\beta}$-name $\tau_{\beta}$ which is as in Lemma \ref{lem:names_in_M} (b) for the stage $\beta$ of the allowable forcing $\forceQ_{\beta+1}$. The Global Freshness Condition at stage $\beta$ only refers to the stages before $\beta$. So the sequence $(\tau_{\beta} \mid \beta < \lambda)$ is as in Lemma \ref{lem:names_in_M} (b) for $\forceQ$.
\end{proof}

\subsection{Products and the symmetry of the coding areas}\label{sec:products}

We show that allowable forcings over a suitable pair with disjoint sets of coding areas can be combined into a product. The main point is that an allowable forcing is not tied to the coding areas it uses, as they are generically added sets. So we can shift the coding areas of such a forcing, using automorphisms of $(\prod_{\omega_1} \mathbb{C} (\omega_1))^L$. We first show that a suitable pair $(M,J)$ is a generic extension of $N_J$ by the Cohen subsets of $\omega_1$ with indices outside of $J$.

\begin{claim}\label{claim:localization}
Let $(M,J)$ be a suitable pair, $M = W[G]$. Then $G^1 \upharpoonright (\omega_1 \setminus J)$ is $(\forceP^1 \upharpoonright (\omega_1 \setminus J))^L$-generic over $N_J$, $M = N_J[G^1 \upharpoonright (\omega_1 \setminus J)]$, and $M$ and $N_J$ have the same reals and the same countable sequences of ordinals.
\end{claim}
\begin{proof}
Let $\forceP$ be the allowable forcing with $G \subset \forceP$, and let $p \in G$ force that all coding areas used by the generic filter belong to $J$. By the proof of Lemma \ref{DefinabilityInProperInnerModels}, $z_{\forceP} \in L[G^0]$, so we can apply $(\ast)$ from that proof with $N := L[G^0][G^1 \upharpoonright J]$, i.e.\ $J' := J$ and $I' := \omega_1$, and obtain $\forceQ \in N$ with $\forceP \restriction p = \forceQ \restriction p$. Now $W = N[G^1 \upharpoonright (\omega_1 \setminus J)]$, where $G^1 \upharpoonright (\omega_1 \setminus J)$ is generic over $N$, and $G \cap (\forceQ \restriction p)$ is generic over $W$ for $\forceQ \restriction p \in N$. By the product lemma, $G^1 \upharpoonright (\omega_1 \setminus J)$ is generic over $N[G] = N_J$, and $M = N_J[G^1 \upharpoonright (\omega_1 \setminus J)]$. Finally, $N_J$ is an extension of $L[G^1 \upharpoonright J]$ by a ccc forcing, and $(\forceP^1 \upharpoonright (\omega_1 \setminus J))^L$ is $\sigma$-closed in $L[G^1 \upharpoonright J]$. So by Easton's Lemma it adds no reals and no countable sequences of ordinals over $N_J$.
\end{proof}

Every hereditarily countable set is coded by a real, so $M$ and $N_J$ also have the same hereditarily countable sets. In particular the codes of the forcings which are allowable over $(M,J)$ (see the remarks after Lemma \ref{lem:names_in_M}) belong to $N_J$; the forcings which they determine depend on $G^1 \upharpoonright (\omega_1 \setminus J)$ as well.

\begin{lemma}\label{lem:symmetry}
Let $(M,J)$ be a suitable pair, let $\mathbb{C} := (\forceP^1 \upharpoonright (\omega_1 \setminus J))^L$ and $H := G^1 \upharpoonright (\omega_1 \setminus J)$, so $M = N_J[H]$ by Claim \ref{claim:localization}.
\begin{enumerate}
\item If $H'$ is $\mathbb{C}$-generic over $N_J$, then $N_J[H']$ is of the same form as $M$, with the coding areas with indices outside of $J$ computed from $H'$. More precisely, $G^0 \times ((G^1 \upharpoonright J) \times H')$ is $\forceP^0 \times \forceP^1$-generic over $L$, and if $W'$ denotes the model which it generates over $L$, then $N_J[H'] = W'[G]$, where $G$ is generic over $W'$ for a forcing which is allowable over $W'$ and whose coding areas belong to $J$. So $(N_J[H'],J)$ is a suitable pair with respect to $W'$ in place of $W$, and its model $N_J$ is the same.
\item Every permutation $\pi \in L$ of $\omega_1$ with $\pi \upharpoonright J = \operatorname{id}$ induces an automorphism $\pi^*$ of $\mathbb{C}$, given by $\pi^*(c)(\pi(\eta)) = c(\eta)$. We have $N_J[\pi^*[H']] = N_J[H']$, and the coding area with index $\pi(\eta)$ computed from $\pi^*[H']$ is the coding area with index $\eta$ computed from $H'$.
\item For $H'$ and $\pi$ as in (1) and (2), a forcing is allowable over $(N_J[H'],J)$ and uses coding areas from $K$, computed from $H'$, if and only if it is allowable over $(N_J[\pi^*[H']],J)$ and uses coding areas from $\pi[K]$, computed from $\pi^*[H']$ (with the covers transformed by $\pi$).
\item Let $H'$ and $\pi$ be as in (1) and (2), let $\forceQ'$ be allowable over $(N_J[H'],J)$, computed from $H'$, with final cover $V$, and let $G'$ be $\forceQ'$-generic over $N_J[H']$. Then $(N_J[H'][G'],\pi[V])$ is a suitable pair with respect to $\pi^*[H']$. A forcing is allowable over $(N_J[H'][G'],V)$ and uses coding areas from $K$, computed from $H'$, if and only if it is allowable over $(N_J[H'][G'],\pi[V])$ and uses coding areas from $\pi[K]$, computed from $\pi^*[H']$ (with the covers transformed by $\pi$).
\end{enumerate}
\end{lemma}

\begin{proof}
(1) Let $\forceP$ be the allowable forcing with $G \subset \forceP$ and let $N := L[G^0][G^1 \upharpoonright J]$. Let $p \in G$ and $\forceQ \in N$ be as in the proof of Claim \ref{claim:localization}, so $\forceP \restriction p = \forceQ \restriction p$. By $(\dagger_2)$ in the proof of Lemma \ref{DefinabilityInProperInnerModels}, $G$ is the upward closure of $G \cap (\forceQ \restriction p)$ in $\forceQ$. So $G$ is $\forceQ$-generic over $W$, and $N_J = N[G]$. As $H'$ is $\mathbb{C}$-generic over $N[G]$, the product lemma shows that $H'$ is $\mathbb{C}$-generic over $N$ and that $G$ is $\forceQ$-generic over $N[H']$. As $J \in L$, the forcing $\forceP^1$ is isomorphic in $L$ to $(\forceP^1 \upharpoonright J)^L \times \mathbb{C}$. So $G^0 \times ((G^1 \upharpoonright J) \times H')$ is $\forceP^0 \times \forceP^1$-generic over $L$, $W' = N[H']$ and $N_J[H'] = W'[G]$. By its definition in the proof of Lemma \ref{DefinabilityInProperInnerModels}, $\forceQ$ is the finite support iteration which is generated by the names coded in $z_{\forceP}$, where a coding forcing is only used if its coding area belongs to $J$. Its coding forcings are computed from the coding areas with indices in $J$ and from the branches given by $G^0$. These objects are the same for $W$ and for $W'$, and almost disjoint coding forcings are absolute. So over $W'$, $\forceQ$ is the iteration which is generated by these names, where the entries with a coding area outside of $J$ are replaced by entries which prescribe the trivial forcing. By Lemma \ref{lem:names_in_M}, applied to $W'$ and the trivial extension of $W'$, this forcing is allowable over $W'$, and all its coding areas belong to $J$.

(2) As $\pi \in L$ maps $\omega_1 \setminus J$ onto itself, $\pi^*$ is an automorphism of $\mathbb{C}$ which belongs to $L$. So $\pi^*[H']$ is $\mathbb{C}$-generic over $N_J$ and $N_J[\pi^*[H']] = N_J[H']$. The coordinate of $\pi^*[H']$ with index $\pi(\eta)$ is $\{ \pi^*(c)(\pi(\eta)) \mid c \in H' \} = \{ c(\eta) \mid c \in H' \}$.

(3) By (2), the coding forcing $\operatorname{Code}(x,m,k,l,\pi(\eta))$ computed from $\pi^*[H']$ is the coding forcing $\operatorname{Code}(x,m,k,l,\eta)$ computed from $H'$; recall that $\pi(\eta) = \eta$ for $\eta \in J$. Let $\forceQ = (\forceQ_{\beta} : \beta < \delta)$ be allowable over $(N_J[H'],J)$, using coding areas from $K$, computed from $H'$, and let $(\tau_{\beta} \mid \beta < \delta)$ be a sequence of names as in Lemma \ref{lem:names_in_M} (b). Let $\tau'_{\beta}$ be a $\forceQ_{\beta}$-name for the tuple which is obtained from the value of $\tau_{\beta}$ by replacing its coding area $\eta$ by $\pi(\eta)$, if this value is a tuple as required in Definition \ref{def:allowable}, and for the value of $\tau_{\beta}$ otherwise. As $\pi$ is injective and the identity on $J \supseteq C^G$, the Global Freshness Condition holds for $\tau_{\beta}$ with respect to $H'$ if and only if it holds for $\tau'_{\beta}$ with respect to $\pi^*[H']$. So by induction on $\beta$, the sequence $(\tau'_{\beta} \mid \beta < \delta)$ generates the same iteration $\forceQ$, now with respect to $\pi^*[H']$, and $\forceQ$ is allowable over $N_J[\pi^*[H']]$ by Lemma \ref{lem:names_in_M}. Conditions (i) and (ii) of Definition \ref{def:alpha_allowable_extension} hold for $K$ and the covers $J_{\gamma}$ with respect to $H'$ if and only if they hold for $\pi[K]$ and the covers $\pi[J_{\gamma}]$ with respect to $\pi^*[H']$. The other direction follows in the same way, with $\pi^{-1}$ in place of $\pi$.

(4) By (3), $\forceQ'$ is allowable over $(N_J[\pi^*[H']],J)$, computed from $\pi^*[H']$, with final cover $\pi[V]$. So by (1), applied to $\pi^*[H']$, and by the remark after Definition \ref{def:alpha_allowable_extension}, the pair $(N_J[H'][G'],\pi[V])$ is suitable with respect to $\pi^*[H']$; recall that $N_J[\pi^*[H']] = N_J[H']$. Now let $\forceQ = (\forceQ_{\beta} : \beta < \delta)$ be allowable over $(N_J[H'][G'],V)$, using coding areas from $K$ and with the covers $(V_{\gamma})_{\gamma \le \delta}$, computed from $H'$. Let $(\tau_{\beta} \mid \beta < \delta)$ be a sequence of names as in Lemma \ref{lem:names_in_M} (c), and define $\tau'_{\beta}$ as in the proof of (3). A coding area $\eta$ belongs to $K \setminus V_{\beta}$ if and only if $\pi(\eta)$ belongs to $\pi[K] \setminus \pi[V_{\beta}]$, and $\eta \in V_{\beta+1}$ if and only if $\pi(\eta) \in \pi[V_{\beta+1}]$. So by the first sentence of the proof of (3), the sequence $(\tau'_{\beta} \mid \beta < \delta)$ is as in Lemma \ref{lem:names_in_M} (c) for the same iteration $\forceQ$, for $\pi[K]$ and the covers $\pi[V_{\gamma}]$, computed from $\pi^*[H']$. Condition (i) of Definition \ref{def:alpha_allowable_extension} is preserved as well. The other direction follows with $\pi^{-1}$ in place of $\pi$.
\end{proof}

\begin{lemma}\label{productallowable}
Let $(M,J)$ be a suitable pair. Let $\forceQ^1$ and $\forceQ^2$ be allowable over $(M,J)$, of lengths $\delta_1$ and $\delta_2$, using coding areas from $K_1$ and $K_2$, with covers $(B^1_{\xi})_{\xi \le \delta_1}$ and $(B^2_{\xi})_{\xi \le \delta_2}$, and assume that $K_1 \cap K_2 = \emptyset$. If $G_1$ is $\forceQ^1$-generic over $M$, then, in $M[G_1]$, $\forceQ^2$ is allowable over $(M[G_1],B^1_{\delta_1})$, using coding areas from $K_2$ and with covers $(B^1_{\delta_1} \cup B^2_{\xi})_{\xi \le \delta_2}$. Consequently, $\forceQ^1 \times \forceQ^2$ densely embeds into the forcing $\forceQ^1 \ast \check{\forceQ}^2$, which is allowable over $(M,J)$.
The same holds with the roles of $\forceQ^1$ and $\forceQ^2$ interchanged. So the product can be presented as an allowable forcing over $(M,J)$ in either order, and both presentations lead to the same model, with the same reals added at the stages of $\forceQ^1$ and of $\forceQ^2$ and with the same final cover $B^1_{\delta_1} \cup B^2_{\delta_2}$. By Lemma \ref{iteration_of_alpha_allowable} (2), all this applies to the initial segments $\forceQ^1_{\xi_1}$ and $\forceQ^2_{\xi_2}$ as well, with the final cover $B^1_{\xi_1} \cup B^2_{\xi_2}$.
\end{lemma}
\begin{proof}
The iterands of $\forceQ^2$ are computed in the same way in $M[G_2 \upharpoonright \xi]$ and in $M[G_1][G_2 \upharpoonright \xi]$, as almost disjoint coding forcings are absolute. So in $M[G_1]$, $\forceQ^2$ is the iteration generated by the same bookkeeping function. As $B^1_{\delta_1} \subseteq J \cup K_1$ is disjoint from $K_2$, conditions (i) and (ii) of Definition \ref{def:alpha_allowable_extension} hold for $K_2$ and the covers $B^1_{\delta_1} \cup B^2_{\xi}$; in particular, the Freshness Condition holds. Finally, $(p,q) \mapsto (p,\check{q})$ is a dense embedding of $\forceQ^1 \times \forceQ^2$ into $\forceQ^1 \ast \check{\forceQ}^2$, and the latter is allowable over $(M,J)$ by Lemma \ref{iteration_of_alpha_allowable} (1). The last assertions follow by symmetry, as both presentations are equivalent to the product, and from Lemma \ref{iteration_of_alpha_allowable} (2).
\end{proof}

\subsection{Iterations of length $\omega_1$}\label{sec:long}

The forcing which we will eventually use is an iteration of length $\omega_1$ whose initial segments of countable length are allowable. So let $\forceP_{\omega_1} = (\forceP_{\beta} : \beta < \omega_1) \in W$ be a finite support iteration such that $\forceP_{\beta}$ is allowable over $W$ for every $\beta < \omega_1$. As for allowable forcings, we say that a real is coded into $\vec{S}^l$ at stage $\beta$ of $\forceP_{\omega_1}$ if this holds for the allowable forcing $\forceP_{\beta+1}$, and we speak of the reals which are coded by $\forceP_{\omega_1}$.

We let $G_{\omega_1}$ be a $\forceP_{\omega_1}$-generic filter over $W$, and we let $G_{\beta} := G_{\omega_1} \cap \forceP_{\beta}$ for $\beta < \omega_1$.

We collect some standard observations. The normal form before Lemma \ref{DefinabilityInProperInnerModels} applies to $\forceP_{\omega_1}$ as well: $\forceP_{\omega_1}$ is equivalent to $\operatorname{Fn}(\omega_1,\Lambda)$ with the ordering described there, so it is a ccc forcing with a dense subset of size $\aleph_1$. Hence $W[G_{\omega_1}]$ has the same cardinals as $L$, and it satisfies $\CH$ by Lemma \ref{FirstPropertiesOfAllowableForcings}. As $\forceP_{\omega_1}$ is a finite support iteration of ccc forcings, every real of $W[G_{\omega_1}]$ belongs to $W[G_{\beta}]$ for some $\beta < \omega_1$. So $\Sigma^1_3$-statements reflect to countable stages: if a $\Sigma^1_3$-statement with real parameters from $W[G_{\beta}]$, $\beta < \omega_1$, holds in $W[G_{\omega_1}]$, then there is $\gamma$ with $\beta \le \gamma < \omega_1$ such that it holds in every model $M$ with $W[G_{\gamma}] \subseteq M \subseteq W[G_{\omega_1}]$, as a witness belongs to some $W[G_{\gamma}]$ and its $\Pi^1_2$ matrix is absolute by Shoenfield absoluteness. Together with Lemma \ref{nounwantedcodes} for the allowable forcings $\forceP_{\gamma}$, $\gamma < \omega_1$, and Shoenfield absoluteness for the $\Pi^1_2$ formulas $\Psi_0$ and $\Psi_1$, this shows that $\Phi_0$ and $\Phi_1$ remain exact: in $W[G_{\omega_1}]$, for $i \in \{0,1\}$, a real satisfies $\Phi_i(v_0)$ if and only if it is coded into $\vec{S}^i$ at some stage of $\forceP_{\omega_1}$, so the set of reals which satisfy $\Phi(v_0)$ is exactly the set of reals which are coded by $\forceP_{\omega_1}$.

In our application the initial segments of $\forceP_{\omega_1}$ are moreover allowable over a suitable pair, with covers which do not depend on the initial segment. We regard $W$ as the extension of $W$ by the trivial forcing, so $(W,\emptyset)$ is a suitable pair. Assume that there are $\forceP_{\beta}$-names $J_{\beta}$, $\beta < \omega_1$, such that for every $\gamma < \omega_1$ the forcing $\forceP_{\gamma}$ is allowable over $(W,\emptyset)$ with the covers $(J_{\beta})_{\beta \le \gamma}$. Then for $\beta < \gamma < \omega_1$ the pair $(W[G_{\beta}],J_{\beta})$ is suitable and $\forceP_{\beta,\gamma}$ is allowable over it, with the covers $(J_{\xi})_{\beta \le \xi \le \gamma}$, by Lemma \ref{iteration_of_alpha_allowable} (2). Such an iteration is built by recursion: Lemma \ref{iteration_of_alpha_allowable} (1) is used at the successor steps, where the next allowable forcing may be given by a name, and Lemma \ref{lem:limits} at the limit steps, where the cover is a sufficiently large countable ordinal.

\section{The coding predicates}\label{sec:common}

In what follows we only code into \(\vec S^1\), and only reals of the form
\((w,0,0)\). For a real \(w\) and a coding area \(\eta<\omega_1\) we write
\[
  \Code(w,\eta):=\Code(w,0,0,1,\eta),
\]
and we say that \(w\) \emph{is coded} at stage \(\beta\) if \((w,0,0)\) is
coded into \(\vec S^1\) at stage \(\beta\). By Lemma~\ref{nounwantedcodes}
and its extension to iterations of length \(\omega_1\) in
Section~\ref{sec:long}, in the final model of an allowable forcing,
and of an iteration of length \(\omega_1\) whose initial segments are
allowable,
\begin{equation}\label{eq:exact}
  \Phi_1((w,0,0))\ \longleftrightarrow\ w\text{ is coded at some stage.}
\end{equation}
As noted at the end of Section~\ref{sec:coding}, the real which is added by a coding forcing
does not belong to the model over which we force, also below a condition.

For separation, and for uniformization at the levels \(n\ge4\), we work
over \(W=L[G^0][G^1]\),
with the allowable forcings over suitable pairs of Section~\ref{sec:alpha}
and their products from Section~\ref{sec:products}. The final iterations
have length \(\omega_1\), and their initial segments are allowable over
\((W,\varnothing)\), as at the end of Section~\ref{sec:long}; so
\eqref{eq:exact} holds in their final models. Followed by a name for such
an iteration, \(\Pp^0\times\Pp^1\) is a set forcing in \(L\). For
uniformization at the levels \(n\ge4\), Section~\ref{sec:uptransfer}
explains why this presentation can replace the one of \cite[\S2.4]{up},
where the Cohen subsets of \(\omega_1\) are added along the iteration. The
case \(n=3\) of Theorem~\ref{main:up} is taken directly from
\cite[Theorem~3.1 and Corollary~3.2]{up}.

\paragraph{The interface.}
Sections~\ref{sec:four}--\ref{sec:up} use the coding machinery of
Sections~\ref{sec:trees}--\ref{sec:allowable} only through
Definition~\ref{def:alpha_allowable_extension}, the format of the
iterations of length \(\omega_1\) at the end of Section~\ref{sec:long}, the
forcings \(\Code(w,\eta)\) and the following lemma. So they can be checked
on the basis of Lemma~\ref{lem:interface}, without looking into the
coding.

\begin{lem}\label{lem:interface}
Let \((M,J)\) be a suitable pair, and let \(\Qq\) and \(\mathbb R\) be
allowable over \((M,J)\), of lengths \(\delta\) and \(\delta'\), using
coding areas from \(K\) and \(K'\), presented as in the normal form
before Lemma~\ref{DefinabilityInProperInnerModels}. Let \(G\) be
\(\Qq\)-generic over \(M\), let \(\Cc_J:=(\Pp^1\res(\omega_1\setminus J))^L\)
and \(\bar G:=G^1\res(\omega_1\setminus J)\).
\begin{enumerate}[leftmargin=*,topsep=.4em,itemsep=.2em,label=\textup{(\arabic*)}]
\item\label{T:loc} \emph{Localization.} \(\bar G\) is \(\Cc_J\)-generic
over \(N_J\), \(M=N_J[\bar G]\), and \(M\) and \(N_J\) have the same reals,
the same countable sequences of ordinals and the same hereditarily
countable sets.
\item\label{T:seg} \emph{Initial segments.} For \(\gamma\le\delta\),
\(\Qq_\gamma\) is allowable over \((M,J)\), and \(\Qq_{\gamma,\delta}\) is
allowable over \((M[G\res\gamma],J_\gamma)\), with the covers
\((J_\xi)_{\gamma\le\xi\le\delta}\) and coding areas from
\(K\setminus J_\gamma\). If \(p\in\Qq\) and \(q\le p\res\gamma\) in
\(\Qq_\gamma\), then \(q\cup p\res[\gamma,\delta)\le p,q\).
\item\label{T:lift} \emph{Restrictions and lifts.} For \(q\in\Qq\),
\(\Qq\res q\) is allowable over \((M,J)\) with the same covers. If
\(c\in\Qq\) forces that \(\dot\Tt\) is allowable over
\((M[\dot G],\dot J_\delta)\), of length \(\check\delta''\) and using
coding areas from \(\check K''\), then there is a \(\Qq\)-name
\(\dot\Tt^\uparrow\) with
\[
 c\Vdash\dot\Tt^\uparrow=\dot\Tt,\qquad
 c'\Vdash\dot\Tt^\uparrow\text{ is trivial of length }\delta''
 \text{ with constant covers}\quad(c'\perp c),
\]
and \(\Qq*\dot\Tt^\uparrow\) is allowable over \((M,J)\).
\item\label{T:closure} \emph{Iterations, single codings and limits.} If
\(\Qq\) forces that \(\dot\Tt\) is allowable over
\((M[\dot G],\dot J_\delta)\) and has countable length, then there is a
countable ordinal \(\delta''\) such that \(\Qq*\dot\Tt\), where \(\dot\Tt\)
is extended by trivial iterands to the length \(\delta''\), is allowable
over \((M,J)\), with the covers of \(\Qq\) followed by those of
\(\dot\Tt\). For a real \(w\in M\) and \(\eta\in\omega_1\setminus J\), the
forcing \(\Code(w,\eta)\) of length one, and its restriction below a
condition, is allowable over \((M,J)\), with coding areas from
\(\{\eta\}\) and the covers \(J\) and \(J\cup\{\eta\}\). An iteration of
trivial forcings of countable length is allowable over \((M,J)\) with
constant covers. Finally, let \(\lambda<\omega_1\) be a limit ordinal and
let \(\Tt=(\Tt_\gamma:\gamma<\lambda)\in M\) be a finite support iteration
with \(\Tt_\gamma\)-names \(J_\gamma\), \(\gamma<\lambda\), such that every
proper initial segment \(\Tt_{\gamma_0}\) is allowable over \((M,J)\) with
the covers \((J_\gamma)_{\gamma\le\gamma_0}\). Then there is a countable
ordinal \(\theta\) such that \(\Tt\) forces \(J_\gamma\subseteq\theta\) for
every \(\gamma<\lambda\), and for every such \(\theta\), \(\Tt\) is
allowable over \((M,J)\) with the final cover \(\theta\).
\item\label{T:prod} \emph{Products.} If \(K\cap K'=\varnothing\), then
\(\Qq\times\mathbb R\) can be presented as an allowable forcing over
\((M,J)\) in either order. Both presentations lead to the same model, with
the same reals added at the stages of \(\Qq\) and of \(\mathbb R\), and
with the union of the two final covers as final cover; the same holds for
the initial segments \(\Qq_\gamma\times\mathbb R_{\gamma'}\). If moreover
\(\gamma\le\delta\) and \(\dot\Tt\) is a \(\Qq_\gamma\)-name for a forcing
which is allowable over \((M[G\res\gamma],J_\gamma)\) and uses coding
areas from \(K''\in L\) with \(K''\cap K'=\varnothing\), then
\[
 (\mathbb R\times\Qq_\gamma)*\dot\Tt\ \simeq\
 \mathbb R\times(\Qq_\gamma*\dot\Tt),
\]
and both sides are allowable over \((M,J)\).
\item\label{T:trans} \emph{Translations of the coding areas.} Let
\(\pi\in L\) be a permutation of \(\omega_1\) with \(\pi\res J=\mathrm{id}\),
let \(\pi^*\) be the automorphism of \(\Cc_J\) given by
\(\pi^*(c)(\pi(\eta))=c(\eta)\), and let \(H'\) be \(\Cc_J\)-generic over
\(N_J\). Then \(G^0\times((G^1\res J)\times H')\) is
\(\Pp^0\times\Pp^1\)-generic over \(L\), and \((N_J[H'],J)\) is a suitable
pair with respect to the model \(W'\) which it generates in place of
\(W\), with the same model \(N_J\). Moreover \(N_J[\pi^*[H']]=N_J[H']\),
the coding area with index \(\pi(\eta)\) computed from \(\pi^*[H']\) is the
coding area with index \(\eta\) computed from \(H'\), and for a later
suitable pair with cover \(J'\), the model \(N_{\pi[J']}\) computed from
\(\pi^*[H']\) is the model \(N_{J'}\) computed from \(H'\). A forcing is
allowable over \((N_J[H'],J)\) and uses coding areas from \(K''\), computed
from \(H'\), if and only if it is allowable over \((N_J[H'],J)\) and uses
coding areas from \(\pi[K'']\), computed from \(\pi^*[H']\), with the covers
transformed by \(\pi\). The same holds over every suitable pair
\((N_J[H'][G'],V)\) which arises from \((N_J[H'],J)\) by an allowable
forcing with final cover \(V\), computed from \(H'\), with \(\pi[V]\) in
place of \(V\) on the right side.
\item\label{T:disj} \emph{Disjoint coding areas.} Let
\(\varphi(\dot\Qq,\dot q)\) be a formula whose parameters belong to \(N_J\)
and which is invariant under the translations of \ref{T:trans}, and let
\(Z\in L\) be countable. If \(M\models\varphi(\Qq,q)\) for some
\(q\in\Qq\), then \(M\models\varphi(\Qq',q')\) for some \(\Qq'\) which is
allowable over \((M,J)\) with coding areas disjoint from \(Z\) and some
\(q'\in\Qq'\). Here \(\varphi\) may assert that \(q\) forces a statement
over \(M\) and that \(\Qq\) belongs to a class which is defined by an
invariant formula.
\item\label{T:codes} \emph{Codes.} \(\Qq\) is given by a code in the sense
of the remarks after Lemma~\ref{lem:names_in_M}, a hereditarily countable
set in \(M\) which determines \(\Qq\) together with its coding areas and
its covers. The set of the codes of all forcings which are allowable over
\((M,J)\) is definable in \(M\) from \(G^0\times G^1\) and \(J\), by a
formula which does not depend on \(M\). Codes for initial segments,
remaining forcings and two step iterations are obtained from the given
codes in a canonical way.
\item\label{T:size} \emph{Size and reals.} \(\Qq\) is ccc and has a dense
subset of size \(\aleph_1\). The models \(M\) and \(M[G]\) have the same
cardinals as \(L\) and satisfy \(\CH\). The same holds for the final model
of an iteration \(\Pp_{\omega_1}\) of length \(\omega_1\) whose initial
segments are allowable over \((W,\varnothing)\), as at the end of
Section~\ref{sec:long}, and every real of this final model belongs to the
model of an initial segment of countable length.
\item\label{T:fresh} \emph{Fresh reals.} If \(\Qq(\gamma)\) is a coding
forcing, possibly restricted below a condition, then the real which it adds
does not belong to \(M[G\res\gamma]\).
\item\label{T:exact} \emph{Exactness.} \eqref{eq:exact} holds in \(M[G]\)
and in the final model of an iteration \(\Pp_{\omega_1}\) as in
\ref{T:size}.
\item\label{T:stage} \emph{Stage models.} \(\R^{M}=\R^{N_J}\). If
\(B\subseteq\R^M\) is defined over \(M\) by a formula as in
\ref{T:disj}, then \(B\in N_J\). In the situation of \ref{T:prod}, for
incomparable pairs \(a,b\) of stages of \(\Qq\times\mathbb R\), the models
\(N_{J_a}\) and \(N_{J_b}\) are mutually generic extensions of
\(N_{J_{a\wedge b}}\), and hence
\[
 \mathcal P(\R)^{N_{J_a}}\cap\mathcal P(\R)^{N_{J_b}}
 =\mathcal P(\R)^{N_{J_{a\wedge b}}}.
\]
\end{enumerate}
\end{lem}
\begin{proof}
\ref{T:loc} is Claim~\ref{claim:localization} together with the remark
after it. \ref{T:seg} is Lemma~\ref{iteration_of_alpha_allowable}~(2); for
the last assertion note that \(\Qq\) is a finite support iteration on
\(\operatorname{Fn}(\delta,\Lambda)\). The first part of \ref{T:lift} is
Lemma~\ref{iteration_of_alpha_allowable}~(3). For the lift, let
\((\tau_\beta)_{\beta<\delta''}\) and \((\dot J_\beta)_{\beta\le\delta''}\)
be names as in Lemma~\ref{lem:names_in_M}~(c) for \(\dot\Tt\) below \(c\).
Let \(\tau^\uparrow_\beta\) be \(\tau_\beta\) below \(c\) and a fixed tuple
which is not as required in Definition~\ref{def:allowable} below every
\(c'\perp c\), and let \(\dot J^\uparrow_\beta\) be \(\dot J_\beta\) below
\(c\) and \(\dot J_\delta\) below every \(c'\perp c\). Condition (c) of
Lemma~\ref{lem:names_in_M} only refers to the single stages, so it holds
for these names, and \(\Qq*\dot\Tt^\uparrow\) is allowable over \((M,J)\)
by Lemma~\ref{iteration_of_alpha_allowable}~(1).

\ref{T:closure} consists of Lemma~\ref{iteration_of_alpha_allowable}~(1),
the examples after Definition~\ref{def:alpha_allowable_extension} and
Lemma~\ref{lem:limits}.

The first part of \ref{T:prod} is Lemma~\ref{productallowable}. For the
second part, \(\Qq_\gamma*\dot\Tt\) is allowable over \((M,J)\) by
Lemma~\ref{iteration_of_alpha_allowable}~(1) and uses coding areas from
\(K\cup K''\), which is disjoint from \(K'\). So
Lemma~\ref{productallowable} applies to \(\mathbb R\) and
\(\Qq_\gamma*\dot\Tt\). The iterands of \(\dot\Tt\) are computed in the
same way in \(M[G\res\gamma]\) and in \(M[G_{\mathbb R}][G\res\gamma]\), as
almost disjoint coding forcings are absolute, which gives the displayed
equivalence.

\ref{T:trans} is Lemma~\ref{lem:symmetry}~(1)--(4). For the models \(N_{\pi[J']}\) note that \(\pi\in L\),
and that the coordinates of \(\pi^*[H']\) with indices in
\(\pi[J'\setminus J]\) are those of \(H'\) with indices in
\(J'\setminus J\), by Lemma~\ref{lem:symmetry}~(2).

For \ref{T:disj} we use the symmetry of the coding areas; the same
argument appears in the proof of Lemma~\ref{s:symmetry}. By the ccc, see
\ref{T:size}, \(\Qq\) uses coding areas from \(J\cup Z_0\) for a countable
\(Z_0\in L\). Some condition \(c\) of \(\Cc_J\) forces over \(N_J\) that
\(\varphi(\Qq,q)\) holds for some code of a forcing \(\Qq\) allowable over
\((M,J)\), with coding areas from \(Z_0\), and some \(q\in\Qq\). Given any
condition of \(\Cc_J\), a permutation \(\pi\in L\) of \(\omega_1\) with
\(\pi\res J=\mathrm{id}\) moves \(Z_0\) and the support of \(c\) above \(Z\)
and above the support of the given condition. By \ref{T:trans} and the
invariance of \(\varphi\), \(\pi^*(c)\) forces the same with coding areas
from \(\pi[Z_0]\), and it is compatible with the given condition. Density
gives \ref{T:disj}.

\ref{T:codes} is Lemma~\ref{lem:names_in_M} together with the remarks
after it. \ref{T:size} is Lemma~\ref{FirstPropertiesOfAllowableForcings}
together with Section~\ref{sec:coding}, the normal form and
Section~\ref{sec:long}. \ref{T:fresh} is the remark at the end of Section~\ref{sec:coding}, and
\ref{T:exact} is Lemma~\ref{nounwantedcodes} together with its extension
to iterations of length \(\omega_1\) in Section~\ref{sec:long}.

For \ref{T:stage}, \(\R^M=\R^{N_J}\) by \ref{T:loc}.
Let \(B=\{y\in\R^M:M\models\psi(y,p)\}\) with \(p\in N_J\) and \(\psi\) as
in \ref{T:disj}. For \(y\in\R^{N_J}\) and conditions \(c,c'\) of \(\Cc_J\),
a permutation \(\pi\) as above makes \(\pi^*(c)\) compatible with \(c'\),
and \(\pi^*(c)\) forces \(\psi(y,p)\) if \(c\) does. So the trivial
condition of \(\Cc_J\) decides \(\psi(y,p)\), and
\(B=\{y\in\R^{N_J}:\mathbf 1\Vdash_{\Cc_J}\psi(\check y,\check p)\}\in N_J\).
Now let \(a=(\beta',\gamma')\) and \(b=(\beta,\gamma)\) with
\(\beta<\beta'\) and \(\gamma'<\gamma\), so \(a\wedge b=(\beta,\gamma')\),
and let \(C_1:=J_a\setminus J_{a\wedge b}\subseteq K\) and
\(C_2:=J_b\setminus J_{a\wedge b}\subseteq K'\), which are disjoint. The
covers are given by names, so we take the factorization below conditions
of \(\Qq_{\beta'}\times\mathbb R_\gamma\) which decide \(J_a\), \(J_b\) and
\(J_{a\wedge b}\); on such a cone \(C_1\) and \(C_2\) are fixed disjoint
sets in \(L\), and the two localized tails use disjoint sets of new Cohen
coordinates. As in
the proof of Claim~\ref{claim:localization}, applied to the presentation
of Lemma~\ref{productallowable}, \(N_{J_a\cup J_b}\) is an extension of
\(N_{J_{a\wedge b}}\) by the product of
\((\Pp^1\res C_1)^L*\Qq_{\beta,\beta'}\) and
\((\Pp^1\res C_2)^L*\mathbb R_{\gamma',\gamma}\), and \(N_{J_a}\),
\(N_{J_b}\) are the extensions by the two factors. By a theorem of Solovay
\cite[Part~I, Lemma~2.5 and the remark after its proof]{solovay}, mutually
generic extensions have only the sets of ordinals of the ground model in
common.
A set of reals in both models is a subset of \(\R^{N_{J_{a\wedge b}}}\),
and it is coded by a set of ordinals via a bijection between an ordinal
and \(\R^{N_{J_{a\wedge b}}}\) in \(N_{J_{a\wedge b}}\).
\end{proof}

Fix recursive injections \(e_j\colon\R\times\R\to\R\), \(j\geq3\), whose
ranges are pairwise disjoint and whose inverses on these ranges are
recursive. Further recursive injections, with ranges disjoint from these
and from each other, are used below for the inputs and for an auxiliary
real \(w_*\). The higher coding predicates are defined by the recursion
\begin{equation}\label{eq:common-decoder}
  D_3(d)\ \longleftrightarrow\ \Phi_1((d,0,0)),\qquad
  D_{j+1}(d)\ \longleftrightarrow\
       \exists t\,\neg D_j(e_j(d,t)).
\end{equation}
Each \(D_j\) is lightface \(\Sigma^1_j\). This is the successor step used
in \cite[\S3]{bpfa}, where the argument of \(D_3\) is the whole real which
is coded. Informally, we say that \(d\) is coded at level \(j\) if
\(D_j(d)\) holds. To make \(D_{j+1}(d)\) true, the forcing can reserve a
new real \(t\) and make sure that \(e_j(d,t)\) is never coded at level
\(j\). To make \(D_{j+1}(d)\) false, it has to code all reals \(e_j(d,t)\)
at level \(j\), including those for reals \(t\) which are added later. The
bookkeeping in both constructions takes care of these new reals.

\section{Separation at level four}\label{sec:four}

This section proves separation at the first level at which the matrix of a
\(\Sigma^1_n\) statement is not absolute. Its result is not new: by
Theorem~\ref{main:up} and Corollary~\ref{cor:models}, the model \(V_4\)
satisfies \(\Sep(\Sigma^1_4(a))\) for every real \(a\). The point is the
method. The construction is a simplified form of that of
Section~\ref{sec:sep}: there are no requirements at the lower levels and no
nodes of higher level, and the index of a pair is a constructible real
instead of a natural number. Each step can be checked directly.
Section~\ref{sec:sep} begins with a list of the changes.

\begin{prop}\label{f:main}
There is a cardinal-preserving set-forcing extension \(V\) of \(L\)
satisfying \(\CH\) with the following property. For every real \(y\) there
is a constructible real \(d_y\) such that any two disjoint
\(\Sigma^1_4(y)\) sets of reals are separated by a \(\Delta^1_4(y,d_y)\) set,
and \(d_{0^\omega}=0^\omega\). In particular,
\(V\models\Sep(\Sigma^1_4)\land\Sep(\RRrel{\Sigma}{4})\).
\end{prop}

\subsection{A toy problem}\label{f:toy}

Let \(A_0,A_1\) be two \(\Sigma^1_4\) sets. If they are disjoint in the
final model \(V\), we want two complementary \(\Sigma^1_4\) sets which
contain \(A_0\) and \(A_1\), respectively. Every real \(x\) then has to be
put on exactly one side, also if it belongs to neither set. But a single
real has no correct side before we know whether the pair is disjoint, so
the construction starts from the pair and not from \(x\). Write
\(x\in A_i\) as \(\exists b\,Q_i(x,b)\), with \(Q_i\) a \(\Pi^1_3\)
matrix.

The truth of \(x\in A_i\) at an intermediate stage is a poor guide, as it
can change infinitely often. A witness \(b\) can be refuted by a later
forcing, and new witnesses can appear. A refutation of \(Q_i(x,b)\) is a
\(\Sigma^1_3\) fact and persists by Shoenfield absoluteness, so a fixed
witness is refuted at most once, and \(x\in A_i\) holds in \(V\) if and
only if some witness is never refuted. Moreover, the forcing which we use
to place \(x\) can itself produce the next witness.

Here is a pair for which we can see this happen. Write
\(\langle i,z,t,j\rangle\) for a fixed recursive coding of tuples by reals
whose range avoids \(w_*\) and the ranges of the injections \(e_j\) and of
the tuple function \(c\) of Section~\ref{f:leaves}. For \(i<2\) let
\[
 P_i(z)\ \longleftrightarrow\
 \exists t\,\Bigl[
   \bigl(\forall n<\omega\;D_3(\langle i,z,t,n+1\rangle)\bigr)
   \ \land\ \neg D_3(\langle i,z,t,0\rangle)\Bigr].
\]
Think of \(t\) as an attempt to put \(z\) on side \(i\). We build it by
coding the reals \(\langle i,z,t,n+1\rangle\), \(n<\omega\), and we cancel
it by coding the real \(\langle i,z,t,0\rangle\). The predicate \(P_i\)
is \(\Sigma^1_4\), since the witnesses for
the countably many positive \(D_3\) assertions can be coded by one real.
Take \(A_0=\{z:P_1(z)\}\) and \(A_1=\{z:P_0(z)\}\), so that each set reads
our attempts to put a real on the other side.

Suppose we place \(x\) according to the truth at the current stage, with
\(P_0,P_1\) as provisional sides. Start over \(W\) on side \(0\) and code
the reals \(\langle 0,x,t_0,n+1\rangle\) one at a time on fresh coding
areas. After these \(\omega\) steps, \(P_0(x)\) holds, so we have put
\(x\) into \(A_1\). Switch to side \(1\) by coding
\(\langle 0,x,t_0,0\rangle\) and starting again with \(1\) and a fresh
\(t_1\). Now the same procedure puts \(x\) into \(A_0\). Repeat this for
\(\omega\) rounds, coding only the indicated reals. At the limit
\(\omega^2\) every attempt we began has been cancelled. For every other
\((i,t)\), including those for reals first appearing at the limit, the
real \(\langle i,x,t,1\rangle\) was never coded. By \eqref{eq:exact}
no \(t\) witnesses \(P_i(x)\). Thus \(x\notin A_0\cup A_1\), and the two
sets are disjoint since all our moves concerned \(x\). Either side would
now be correct, yet neither provisional side contains \(x\), and the truth
value of \(x\in A_0\) has changed infinitely often.

\paragraph{Certified witnesses.}
We therefore decide by whether a witness can be certified, not by what
holds at the current stage. The notions below are relative to rules which
restrict the codings. They are recorded in a set of \emph{protected}
leaves, whose reals may not be coded, and in a register of entries for
pairs (Definitions~\ref{f:forcings} and~\ref{f:poss}), and a continuation
is \emph{allowable} if it respects them. A witness \(b\) for \(x\in A_i\)
is \emph{certified} at a stage if no continuation which respects the
protected leaves can refute \(Q_i(x,b)\) (Definition~\ref{f:cert}). The
membership \(x\in A_i\) is \emph{possible} at a stage if some allowable
continuation produces a certified witness (Definition~\ref{f:poss}). Both
notions go back to \cite[\S7.5]{mixed}, where they select, without these
names, the values of a \(\Sigma^1_4\)-uniformizing function; there they
are taken relative to the fixed point of a thinning-out process instead
of a register. The rules are essential here. A \(\Pi^1_3\) fact \(\neg D_3(w)\) is certified
only if the rules forbid coding \(w\) (Lemma~\ref{f:persist}), and in the
example no real \(\langle i,x,t,0\rangle\) is ever forbidden. So no witness
for \(P_i(x)\) is ever certified.

Three facts drive the construction.

\emph{Possibility does not depend on the stage}, as long as the register is
kept fixed. It cannot arise later: if an allowable continuation \(\Tt\)
leads to a stage at which \(\Qq\) produces a certified witness, then
\(\Tt\) followed by \(\Qq\) produces one already before \(\Tt\). It cannot
be lost: let \(\Qq\) produce a certified witness, and let \(\Tt\) be any
allowable continuation. Move \(\Qq\) to coding areas which are disjoint
from those of \(\Tt\) (Lemma~\ref{f:basic}), and form the product. Read
as \(\Qq\) followed by \(\Tt\), the product keeps the certificate, since
certificates persist through allowable forcing (Lemma~\ref{f:persist}).
Read as \(\Tt\) followed by \(\Qq\), it shows that \(\Qq\) is still
allowable after \(\Tt\): the leaves which \(\Tt\) protects have new last
coordinates, the extensions by \(\Qq\) and by \(\Tt\) have only the reals
of the initial stage in common, and the side on which \(\Qq\) puts a real
is decided by the same statement before and after \(\Tt\), for a shorter
register. This is Proposition~\ref{f:product}, and it is the heart of the
proof. At level four the persistence of certificates is immediate. In
Section~\ref{sec:sep} it is the same proposition one level down, so
Proposition~\ref{s:product} is proved by induction on the level, and its
products have to carry complete assignments at all lower levels.

\emph{The iteration refutes whatever can be refuted}, by task~(2) in
Section~\ref{f:iteration}. So a \(\Pi^1_3\) statement holds in \(V\) if and
only if it is certified at a stage whose model contains its parameters
(Lemma~\ref{f:truth}).

\emph{Rules are only added.} A register is only extended, and the rule for
an entry refers to the register before the entry, which stays fixed
(Lemma~\ref{f:parts}). A membership which holds in \(V\) is possible with
respect to every register, since its witness is certified at a later stage
and the trivial continuation is allowable.

\paragraph{Placing reals.}
In general the two sets are \(\Sigma^1_4(y)\) for a real \(y\), the
\emph{parameter} of the pair, and for each real \(y\) there are countably
many pairs with parameter \(y\), one for each code of two \(\Sigma^1_4\)
formulas. The iteration \emph{treats} every real \(y\) exactly once, at a
stage \(\sigma_y\) at which a segment of the iteration begins and whose
model contains \(y\) (Section~\ref{f:iteration}). The treatment deals with
all pairs with parameter \(y\) at once. First it chooses an index \(d\)
for their separators which is \emph{fresh}: before \(\sigma_y\), no real
has been put on a side under the index \(d\) for a pair with parameter
\(y\). Then it tests these pairs by a recursion on \(\nu\), which extends
the register built before \(\sigma_y\), call it \(H_0\), to registers
\(H_\nu\). Given \(H_\nu\), a pair with parameter \(y\) which has no entry
in \(H_\nu\) \emph{fails the test} if the \(\Sigma^1_4\) statement
\(\exists x\,(x\in A_0\wedge x\in A_1)\) about its two sets is not
possible at \(\sigma_y\) with respect to \(H_\nu\). If some pair fails the
test, let \(p\) be the one with the least code, and let \(H_{\nu+1}\) be
\(H_\nu\) followed by the entry \(r=(p,d,H_\nu)\). So the register
\(H_r:=H_\nu\) which is stored in \(r\) is the one with respect to which
\(p\) failed the test. At limits we take unions, and if no pair fails the
test, the recursion stops. Thus a pair without an entry is tested again
at every later step of the recursion, with respect to the larger
register. This is necessary: an entry restricts the continuations which
are allowable, so a common point which is possible with respect to
\(H_\nu\) need not be possible with respect to \(H_{\nu+1}\)
(Section~\ref{f:examples}). Each step appends an entry for a new pair, so
the recursion stops at a countable ordinal \(\nu^*\), and we let
\(H^{\sigma_y}:=H_{\nu^*}\). At this point every pair with parameter
\(y\) either has an entry \(r\), and a common point is not possible at
\(\sigma_y\) with respect to \(H_r\), or it has no entry, and a common
point is possible at \(\sigma_y\) with respect to \(H^{\sigma_y}\).

The segment which begins at \(\sigma_y\) takes care of the pairs without
an entry. For each of them, one after the other, it carries out a
continuation which is allowable with respect to \(H^{\sigma_y}\) and at
whose end there are \(x,b_0,b_1\) such that both \(Q_0(x,b_0)\) and
\(Q_1(x,b_1)\) are certified; by stage independence such a continuation
still exists after the earlier ones. By Lemma~\ref{f:persist}\,(iii),
\(x\in A_0\cap A_1\) in \(V\), so these pairs need no separation. For a
pair with an entry \(r\), every real \(x\), whether it exists at
\(\sigma_y\) or appears later, is put on side \(1\) if \(x\in A_1\) is
possible with respect to \(H_r\) at the stage at which \(x\) is placed,
and on side \(0\) otherwise. This is the rule \eqref{f:index}.

The placement is correct. Fix a pair \(p\) with an entry \(r\). From
\(\sigma_y\) on, the iteration is allowable with respect to \(H_r\), as
\(H_r\) is an initial segment of every later register
(Lemma~\ref{f:parts}), so stage independence applies to possibility with
respect to \(H_r\). First, at no stage after \(\sigma_y\) are
\(x\in A_0\) and \(x\in A_1\) both possible with respect to \(H_r\)
(Lemma~\ref{f:incompat}). Otherwise, move the two witnessing forcings to
disjoint coding areas and take their product. Both certificates survive
in the product, so \(x\) and the two witnesses show that a common point
is possible with respect to \(H_r\) at this stage. As \(x\) is part of
the witness, the only real parameter of this statement is \(y\), so by
stage independence a common point was possible with respect to \(H_r\)
at \(\sigma_y\), even if \(x\) appeared later. But this is the test which
\(p\) failed when \(r\) was appended. Now if \(x\in A_1\) holds in \(V\),
then it is possible with respect to \(H_r\), and \(x\) goes to side
\(1\). If \(x\in A_0\) holds in \(V\), then \(x\in A_0\) is possible, so
\(x\in A_1\) is not, and \(x\) goes to side \(0\). The register \(H_r\)
never changes, so by stage independence the answer for \(x\) is the same
at all stages after \(\sigma_y\) whose model contains \(x\), and the side
of \(x\) is never revised. In the example no membership is ever possible.
So a common point is not possible, the pair receives an entry, and every
real goes to side \(0\); in \(V\) both sets are empty, see
Section~\ref{f:examples}.

Three further points are separate from the oscillation. The index \(d\)
is fresh because steps for the pair may have occurred before it was
tested, for instance inside a witnessing forcing in which its parameter
first appears, and such steps obey no rule for the pair. The rule for
\(x\) refers to the stored register \(H_r\), which does not contain the
entry \(r\), and not to the register in effect, which does; otherwise
the rule would refer to its own outcome. Finally, a pair may be defined
from the candidate separators themselves and so read the side which we
are about to choose. Such a pair is dealt with by its test, before any
side is chosen; see Section~\ref{f:examples}.

\paragraph{Carrying out the choice.}
Let \(c_i\) code the pair, its index \(d\), the real \(x\) and the side
\(i\). On the chosen side we keep \(e_3(c_i,t)\) uncoded for a fresh real
\(t\). This use of a fresh real is taken from \cite[\S7.5]{mixed}, where a
new real is added whose extensions are never coded. The real
\(e_3(c_i,t)\) is the real of a leaf which is protected, so the rules forbid
coding it, and \(\neg D_3(e_3(c_i,t))\) is certified. On the opposite side
every real \(e_3(c_{1-i},t)\) is eventually coded, including those for
reals \(t\) which are added later. By \eqref{eq:common-decoder}, exactly
the chosen side satisfies \(D_4\) in the final model. The rest of this
section carries this out.

\subsection{Leaves and certificates}\label{f:leaves}

We work over \(W=L[G^0][G^1]\) and construct a finite support iteration
\(\Pp_{\omega_1}=(\Pp_\beta:\beta<\omega_1)\) as in the last paragraph of
Section~\ref{sec:long}: there are \(\Pp_\beta\)-names \(J_\beta\),
\(\beta<\omega_1\), such that for every \(\gamma<\omega_1\) the forcing
\(\Pp_\gamma\) is allowable over the suitable pair \((W,\emptyset)\) with
the covers \((J_\beta)_{\beta\le\gamma}\). We use suitable pairs \((M,J)\),
the models \(N_J\) and allowability over \((M,J)\) as in
Definition~\ref{def:alpha_allowable_extension}, and the forcing \(\Cc_J\)
of Lemma~\ref{lem:interface}. Apart from these definitions and the format
of \(\Pp_{\omega_1}\), this section and Section~\ref{sec:sep} use the
coding only through Lemma~\ref{lem:interface}. By
Lemma~\ref{lem:interface}\,\ref{T:loc},
\(M=N_J[G^1\res(\omega_1\setminus J)]\), and \(M\) and \(N_J\) have the
same reals, the same countable sequences of ordinals and the same
hereditarily countable sets.

As the letter \(K\) will denote the first component of an assignment in
Section~\ref{sec:sep}, we write \(Z\) for the set from which an allowable
forcing over \((M,J)\) uses its coding areas, i.e.\ for the set \(K\) of
Definition~\ref{def:alpha_allowable_extension}(2) and of
Lemma~\ref{lem:interface}. A coding area in a cover
\(J_\gamma\) can never be used at a later stage, even if it was not used
before. The models \(N_J\) increase along the iteration as the covers
increase, and if two independent forcings are carried out in either order,
then the model \(N_J\), computed from both generic filters and the union of
both covers, is the same. By
Lemma~\ref{lem:interface}\,\ref{T:seg}--\ref{T:closure}, allowability over
suitable pairs is preserved by two step iterations, and it is inherited by
initial segments, by the remaining forcings and by restrictions below
conditions; trivial iterands are allowed, and countable limits are
allowable, with a countable ordinal as final cover.

\paragraph{Pairs and leaves.}
A \emph{pair} is \(p=(e,y)\), where \(e\in\omega\) codes two \(\Sigma^1_4\)
formulas and \(y\) is a real, the \emph{parameter} \(\parm(p):=y\). Let
\(A^4_{p,0},A^4_{p,1}\) be the two \(\Sigma^1_4(y)\) sets of reals which
these formulas define from \(y\). We present them in normal form,
\[
 x\in A^4_{p,i}\ \longleftrightarrow\ \exists b\,P_{p,i}(x,b)\qquad(i<2),
\]
with \(\Pi^1_3\) matrices \(P_{p,i}\), which have the parameter \(y\). In
general, a \(\Sigma^1_4\) assertion in normal form is \(\exists b\,P(b)\),
where \(P\) is a finite conjunction of \(\Pi^1_3\) assertions with real
parameters. The assertion \(I^4_p\) that \(A^4_{p,0}\) and \(A^4_{p,1}\)
have a common element is presented with the witness \((x,b_0,b_1)\) and the
matrix \(P_{p,0}(x,b_0)\land P_{p,1}(x,b_1)\).

Fix a recursive coding \(c(j,p,d,x,i)\) by reals of tuples of natural
numbers and reals, with a recursive inverse on its range, whose range is
disjoint from the ranges fixed in Section~\ref{sec:common}. In this section
\(j=4\), and the index \(d\) is a real. A \emph{leaf} is a tuple
\[
 v=(4,p,d,x,i,t),\qquad p\text{ a pair},\ i<2,\ d,x,t\in\R,
\]
and its real is \(w(v):=e_3(c(4,p,d,x,i),t)\). We write \(\lc(v):=t\) for
its last coordinate. Applying \(e_3^{-1}\) and the inverse of \(c\) to
\(w(v)\) recovers all coordinates of \(v\), so \(v\le_T w(v)\), and
\(w(v)=w(u)\) only if \(v=u\). The reals of leaves are different from
\(w_*\). For a pair \(p\), reals \(d,x\) and \(i<2\) let
\begin{equation}\label{f:sep}
 S^4_i(p,d,x)\ :\longleftrightarrow\
 \exists t\ \neg D_3\bigl(w((4,p,d,x,i,t))\bigr)
 \ \longleftrightarrow\ D_4\bigl(c(4,p,d,x,i)\bigr),
\end{equation}
by \eqref{eq:common-decoder}. For fixed \(e\) and \(d\), both sets
\(S^4_i(p,d,\cdot)\), \(i<2\), are \(\Sigma^1_4(y,d)\). For a pair which is
disjoint in the final model we shall find a real \(d\), which only depends
on \(y\), such that these two sets are complements of each other and
contain \(A^4_{p,0}\) and \(A^4_{p,1}\), respectively.

\begin{defn}\label{f:forcings}
A \emph{stage} \(s\) consists of a suitable pair \((M_s,J_s)\) and a
countable set \(B_s\in M_s\) of leaves, the \emph{protected} leaves. The
\emph{initial stage} is \((W,\emptyset,\emptyset)\). A \emph{level-four
forcing} over a stage \(s\) is a forcing \(\Qq=(\Qq_\gamma:\gamma<\delta)\)
of countable length which is allowable over \((M_s,J_s)\), together with
\(\Qq_\gamma\)-names \(B_\gamma\), \(\gamma\le\delta\), for countable sets of
leaves, such that \(B_0=B_s\),
\(B_\gamma=\bigcup_{\gamma'<\gamma}B_{\gamma'}\) for limit \(\gamma\), and
the following hold for every \(\gamma<\delta\), where \(G\) is
\(\Qq\)-generic over \(M_s\).
\begin{enumerate}[leftmargin=*,topsep=.4em,itemsep=.2em,label=\textup{(\roman*)}]
\item Either \(B_{\gamma+1}=B_\gamma\), or \(\Qq(\gamma)\) is a coding
forcing for \(w_*\) and \(B_{\gamma+1}=B_\gamma\cup\{v\}\) for a leaf
\(v\) whose last coordinate is the real added at \(\gamma\) and whose other
coordinates belong to \(M_s[G\res\gamma]\). In the second case \(v\) is
\emph{protected} at \(\gamma\), and \(\Qq\) \emph{makes a step for
\((p,d)\) at \(\gamma\), concerning \(x\) and \(i\)}, where
\(v=(4,p,d,x,i,\lc(v))\).
\item If \(\Qq(\gamma)\) is a coding forcing for the real \(w(v)\) of a leaf
\(v\), then \(v\notin B_\gamma\).
\end{enumerate}
Here a coding forcing for a real \(w\) is \(\Code(w,\eta)\) or its
restriction below a condition. The stage \(s[G\res\gamma]\) consists of
\((M_s[G\res\gamma],J_\gamma)\) and \(B_\gamma\), and \(s[G]:=s[G\res\delta]\)
is the stage \emph{reached} by \(\Qq\). A level-four forcing without steps is
\emph{ordinary}.
\end{defn}

If the steps of a level-four forcing \(\Qq\) over \(s\) are
\emph{forgotten}, i.e.\ if every \(B_\gamma\) is replaced by \(B_s\), we
obtain an ordinary forcing over \(s\) with the same iterands, generic
extensions and covers; (ii) still holds, as \(B_s\subseteq B_\gamma\). The
trivial forcing is ordinary, and if \(w\) is a real of \(M_s\) which is not
the real of a leaf in \(B_s\) and \(\eta\notin J_s\), then the forcing
\(\Code(w,\eta)\) of length one is ordinary over \(s\), by
Lemma~\ref{lem:interface}\,\ref{T:closure}.

\begin{defn}\label{f:cert}
Let \(s\) be a stage and \(\theta\) a \(\Pi^1_3\) assertion with parameters
from \(M_s\). Then \(\cert_3(s,\theta)\) holds if no condition of an
ordinary forcing over \(s\) forces \(\neg\theta\). For a finite conjunction
of \(\Pi^1_3\) assertions, \(\cert_3\) means \(\cert_3\) for each conjunct.
\end{defn}

As the trivial forcing is ordinary, a certified assertion holds in \(M_s\).
If \(\theta\) is \(\Pi^1_2\), then \(\cert_3(s,\theta)\) holds if and only
if \(\theta\) holds in \(M_s\), by Shoenfield absoluteness. In
Section~\ref{sec:sep}, \(\cert_3\) is defined from possibility at level three
with respect to the complete level-three assignment. Here there are no
requirements at level three, and the forcings which are allowable at level
three are the ordinary ones.

\begin{lem}\label{f:persist}
Let \(s\) be a stage and \(\theta\) a \(\Pi^1_3\) assertion with parameters
from \(M_s\).
\begin{enumerate}[leftmargin=*,topsep=.4em,itemsep=.2em,label=\textup{(\roman*)}]
\item If \(\cert_3(s,\theta)\) holds and \(s'\) is reached from \(s\) by a
level-four forcing, then \(\cert_3(s',\theta)\).
\item Every leaf which is protected by a level-four forcing over \(s\) has
its last coordinate outside \(M_s\). If \(s\) is reached from the initial
stage by a level-four forcing, then for every real \(w\in M_s\),
\[
 \cert_3(s,\neg D_3(w))\ \longleftrightarrow\ w=w(v)\ \text{for some }v\in B_s .
\]
\item Let \(\Pp_{\omega_1}\) be an iteration as above, with names
\(B_\gamma\) for countable sets of leaves, such that for all stages
\(\alpha\le\beta<\omega_1\) the part \(\Pp_{\alpha,\beta}\), with the sets
\(B_\gamma\), \(\alpha\le\gamma\le\beta\), is a level-four forcing over the
stage \(\alpha\). If \(\cert_3\) holds for \(\theta\) at a stage \(\alpha\),
then \(\theta\) holds in the final model.
\end{enumerate}
\end{lem}
\begin{proof}
(i) Let \(\Tt\) be the forcing from \(s\) to \(s'=s[G]\), and suppose that a
condition \(q\) of an ordinary forcing \(\Qq'\) over \(s'\) forces
\(\neg\theta\). Choose \(t\in G\) and \(\Tt\)-names \(\dot\Qq',\dot q\) such
that \(t\) forces that \(\dot\Qq'\) is ordinary over \(s[\dot G]\) and that
its condition \(\dot q\) forces \(\neg\theta\). By
Lemma~\ref{lem:interface}\,\ref{T:lift} and~\ref{T:closure},
\((\Tt\res t)*\dot\Qq'\) is allowable over \((M_s,J_s)\). We forget the steps
of \(\Tt\). The result is an ordinary forcing over \(s\), since a coding of
the real of a leaf \(v\) occurs at a stage at which \(v\) is not protected,
so that \(v\notin B_s\). Its condition \((t,\dot q)\) forces \(\neg\theta\),
as forgetting the steps does not change the forcing. So
\(\cert_3(s,\theta)\) fails.

(ii) The real added at a stage \(\gamma\) does not belong to
\(M_s[G\res\gamma]\supseteq M_s\), by
Lemma~\ref{lem:interface}\,\ref{T:fresh}. This proves the first assertion.
If \(w\) is not the real of a leaf in \(B_s\), take \(\eta\notin J_s\). The
forcing \(\Code(w,\eta)\) is ordinary over \(s\), and it forces \(D_3(w)\) by
Lemma~\ref{lem:interface}\,\ref{T:exact}. Now let \(v\in B_s\). The forcing
which reaches \(s\) from the initial stage protects \(v\) at some stage
\(\gamma_v\). As \(w(v)\) computes \(\lc(v)\), which is the real added at
\(\gamma_v\), \(w(v)\) does not belong to the model at \(\gamma_v\). The real
which is coded at a stage belongs to the model at that stage, and the real
coded at \(\gamma_v\) is \(w_*\). So \(w(v)\) is not coded at a stage
\(\le\gamma_v\). Later \(v\) is protected, so \(w(v)\) is coded neither
before \(s\) nor by an ordinary forcing over \(s\). By
Lemma~\ref{lem:interface}\,\ref{T:exact}, \(\neg D_3(w(v))\) holds in every
extension of \(M_s\) by an ordinary forcing, and so
\(\cert_3(s,\neg D_3(w(v)))\).

(iii) Write \(\neg\theta\) as \(\exists z\,R(z)\) with \(R\in\Pi^1_2\). If
\(\neg\theta\) held in the final model, a witness \(z\) would belong to the
model at some stage \(\beta\ge\alpha\), by
Lemma~\ref{lem:interface}\,\ref{T:size}, and \(R(z)\) would hold there, by
Shoenfield absoluteness. But \(\cert_3\) holds for \(\theta\) at \(\beta\), by
(i), and the trivial forcing is ordinary.
\end{proof}

\subsection{Registers and possibility}\label{f:registers}

\begin{defn}\label{f:poss}
A \emph{register} is a sequence \(H=(r_\nu:\nu<\lambda)\), \(\lambda<\omega_1\),
of \emph{entries} \(r_\nu=(p_\nu,d_\nu,H\res\nu)\), where the \(p_\nu\) are
pairwise different pairs, the \(d_\nu\) are reals, and
\(H\res\nu:=(r_\mu:\mu<\nu)\). For an entry \(r=r_\nu\) we call
\(H_r:=H\res\nu\) its \emph{stored register}. By recursion on the length of
\(H\), simultaneously for all stages \(s\), let:
\begin{enumerate}[leftmargin=*,topsep=.4em,itemsep=.2em,label=\textup{(\roman*)}]
\item \(\B_4(s;H)\) be the class of the level-four forcings \(\Qq\) over
\(s\) such that every step of \(\Qq\) for \((p,d)\) at a stage \(\gamma\),
concerning \(x\) and \(i\), where \((p,d,H_r)\) is an entry of \(H\),
satisfies
\begin{equation}\label{f:index}
 i=1\ \Longleftrightarrow\ \Pos_4\bigl(s[G\res\gamma];H_r,\,x\in A^4_{p,1}\bigr);
\end{equation}
\item for a \(\Sigma^1_4\) assertion \(A=\exists b\,P(b)\) in normal form with
parameters from \(M_s\),
\begin{equation}\label{f:possdef}
 \Pos_4(s;H,A)\ :\Longleftrightarrow\ \exists\Qq\in\B_4(s;H)\ \
 \Qq\Vdash\exists b\ \cert_3\bigl(s[\dot G],P(b)\bigr).
\end{equation}
\end{enumerate}
\end{defn}

Clause (i) for \(H\) only refers to \(\Pos_4\) for the stored registers of the
entries of \(H\), which are shorter than \(H\), and (ii) refers to
\(\B_4(s;H)\) and to \(\cert_3\). A register is a hereditarily countable
set. A step for \((p,d)\) for which \(H\) has no entry is unrestricted, also
if \(H\) has an entry \((p,d',H_r)\) with \(d'\neq d\). The coding of the real
of a leaf which is not protected is never restricted, also if its pair has an
entry. Every ordinary forcing over \(s\) belongs to every class
\(\B_4(s;H)\), and forgetting the steps of a forcing in \(\B_4(s;H)\) gives an
ordinary forcing over \(s\). If \(H'\) is an initial segment of \(H\), then
\(\B_4(s;H)\subseteq\B_4(s;H')\), as an entry of \(H'\) has the same stored
register in \(H\); hence \(\Pos_4(s;H,A)\) implies \(\Pos_4(s;H',A)\). By the
following lemma, \(\Pos_4(s;H,A)\) holds as soon as a condition of a forcing
in \(\B_4(s;H)\) forces the statement in \eqref{f:possdef}, in particular if
\(A\) has a witness \(b\in M_s\) with \(\cert_3(s,P(b))\).

\begin{lem}\label{f:basic}
Let \(s\) be a stage and \(H\) a register.
\begin{enumerate}[leftmargin=*,topsep=.4em,itemsep=.2em,label=\textup{(\roman*)}]
\item \(\B_4(s;H)\) is closed under restrictions below conditions, under
initial segments and remaining parts, with the inherited stages, and under
two step iterations: if \(\Qq\in\B_4(s;H)\) forces that
\(\dot\Tt\in\B_4(s[\dot G];H)\), then \(\Qq*\dot\Tt\in\B_4(s;H)\). The same
holds for a countable sequence of forcings, each over the stage reached by
the preceding ones, with the unions of the sets of protected leaves and a
countable ordinal as cover at limits.
\item The class \(\B_4(s;H)\) and the predicates \(\Pos_4(s;H,\cdot)\) and
\(\cert_3(s,\cdot)\) are defined over \(M_s\), uniformly in \(s\), by
formulas of set theory with the parameters \(G^0\times G^1\), \(J_s\),
\(B_s\) and \(H\). The same holds for their versions given by names over
initial segments of allowable forcings.
\item Let \(Z'\in L\) be countable. A witness to \(\Pos_4(s;H,A)\), and an
ordinary forcing over \(s\) with a condition which forces a given
projective assertion with parameters from \(M_s\), can be chosen with
coding areas disjoint from \(Z'\).
\end{enumerate}
\end{lem}
\begin{proof}
(i) The requirements of Definitions~\ref{f:forcings} and~\ref{f:poss}
concern the single stages of a forcing and the data at these stages. Under
the indicated operations every stage keeps its iterand, its model and its
data, and allowability is preserved by
Lemma~\ref{lem:interface}\,\ref{T:seg}--\ref{T:closure}. So all
requirements carry over. Unlike in Section~\ref{sec:sep}, no requirement
concerns the final stage, so nothing has to be checked at limits.

(ii) A forcing which is allowable over \((M_s,J_s)\) is given by a code, and
the set of these codes is definable in \(M_s\) from \(G^0\times G^1\) and
\(J_s\), by Lemma~\ref{lem:interface}\,\ref{T:codes}. The requirements at a
stage \(\gamma\) are read through the forcing relation of the initial
segment \(\Qq_\gamma\), and \eqref{f:index} at \(\gamma\) refers to
\(\Pos_4\) for a shorter register, in the extension by \(\Qq_\gamma\). By
recursion on the length of \(H\) this gives the definitions. That the
versions given by names agree with the definitions evaluated in the
extensions is shown as in (D1)--(D4) in Section~\ref{sec:sep}, with the
length of the register in place of the index \((k,\eta)\) used there.

(iii) Let \(\pi\in L\) be a permutation of \(\omega_1\) with
\(\pi\res J_s=\mathrm{id}\). By Lemma~\ref{lem:interface}\,\ref{T:trans},
the translation by \(\pi\) keeps the iterands, the generic reals, the
sets of protected leaves and the registers, and it preserves allowability,
with translated coding areas and covers. By induction on the length of
\(H\), as for \eqref{s:inv} in the proof of Lemma~\ref{s:symmetry}, it
preserves membership in \(\B_4(s;H)\) and the truth of \(\Pos_4\) and
\(\cert_3\). The data of the stage and the register are hereditarily
countable, so they belong to \(N_{J_s}\) by
Lemma~\ref{lem:interface}\,\ref{T:loc}. So
Lemma~\ref{lem:interface}\,\ref{T:disj} applies.
\end{proof}

\subsection{Products}\label{f:products}

Two level-four forcings \(\Qq^1,\Qq^2\) over a stage \(s\) are
\emph{disjoint} if they use coding areas from sets \(Z_1,Z_2\in L\) with
\(Z_1\cap Z_2=\emptyset\). Let \(\delta_i\) be the length of \(\Qq^i\) and
\(B^i_\xi\) its sets of protected leaves, let \(G_1\times G_2\) be
\(\Qq^1\times\Qq^2\)-generic over \(M:=M_s\), and let
\(M_{\beta,\gamma}:=M[G_1\res\beta][G_2\res\gamma]\). By
Lemma~\ref{lem:interface}\,\ref{T:prod}, for \(\beta\le\delta_1\) and
\(\gamma\le\delta_2\) the product \(\Qq^1_\beta\times\Qq^2_\gamma\) can be
presented as an allowable forcing over \((M,J_s)\) in either order, with
the same reals added at the stages of \(\Qq^1\) and of \(\Qq^2\) and with the
union of the covers. We write \(\Qq^1_\beta*\Qq^2_\gamma\) for the
presentation which runs through \(\Qq^1\) up to \(\beta\) and then, over
\(M_{\beta,0}\), through \(\Qq^2\) up to \(\gamma\), and we give its stage
\((\beta',\gamma')\) the protected leaves \(B^1_{\beta'}\cup B^2_{\gamma'}\).
The presentation \(\Qq^2_\gamma*\Qq^1_\beta\) is defined symmetrically. Both
reach the model \(M_{\beta,\gamma}\) with the protected leaves
\(B^1_\beta\cup B^2_\gamma\). The extensions \(M[G_1]\) and \(M[G_2]\) are
mutually generic over \(M\), so they have only the reals of \(M\) in common
\cite{solovay}.

\begin{prop}\label{f:product}
For every register \(H\) and every stage \(s\) the following hold.
\begin{enumerate}[leftmargin=*,topsep=.4em,itemsep=.2em,label=\textup{(\roman*)}]
\item If \(\Qq^1,\Qq^2\in\B_4(s;H)\) are disjoint, then
\(\Qq^1_\beta*\Qq^2_\gamma\) and \(\Qq^2_\gamma*\Qq^1_\beta\) belong to
\(\B_4(s;H)\), for all \(\beta\le\delta_1\) and \(\gamma\le\delta_2\).
\item If \(\Tt\in\B_4(s;H)\) reaches the stage \(s'\), then for every
\(\Sigma^1_4\) assertion \(A\) in normal form with parameters from \(M_s\),
\begin{equation}\label{f:inv}
 \Pos_4(s;H,A)\ \longleftrightarrow\ \Pos_4(s';H,A).
\end{equation}
\end{enumerate}
\end{prop}
\begin{proof}
By induction on the length of \(H\), simultaneously for all stages. For
\(H\) we first prove (i), using (i) and (ii) for the stored registers of
the entries of \(H\), and then (ii), using (i) for \(H\).

(i) We consider \(\Qq^1_\beta*\Qq^2_\gamma\); the other presentation is
symmetric. At a stage \((\beta',0)\) the iterand, the step and the data are
those of \(\Qq^1\) at \(\beta'\). Consider a stage \((\beta,\gamma')\),
\(\gamma'<\gamma\), at which \(\Qq^2\) acts. The real which \(\Qq^2\) codes
at \(\gamma'\), and the coordinates other than the last one of a leaf which
it protects there, belong to \(M[G_2\res\gamma']\subseteq M[G_2]\).

\emph{Protected leaves.} Let \(v\in B^1_\beta\setminus B_s\). Then
\(\lc(v)\in M[G_1]\setminus M\) by Lemma~\ref{f:persist}\,(ii), so neither
\(v\) nor \(w(v)\), which computes \(v\), belongs to \(M[G_2]\). So if
\(\Qq^2\) codes the real \(w(u)\) of a leaf \(u\) at \(\gamma'\), then
\(u\notin B^1_\beta\setminus B_s\), and \(u\notin B^2_{\gamma'}\) by
Definition~\ref{f:forcings}\,(ii) for \(\Qq^2\); so \(u\) is not protected
at \((\beta,\gamma')\). If \(\Qq^2\) makes a step at \(\gamma'\), the last
coordinate of the new leaf is the real added at \(\gamma'\), which is the
same in the product, and its other coordinates belong to
\(M[G_2\res\gamma']\subseteq M_{\beta,\gamma'}\). So
Definition~\ref{f:forcings} holds at \((\beta,\gamma')\).

\emph{The rule \eqref{f:index}.} Let \(\Qq^2\) make a step at \(\gamma'\)
for \((p,d)\), concerning \(x\) and \(i\), where \(r=(p,d,H_r)\) is an
entry of \(H\). Let \(s_2:=s[G_2\res\gamma']\) be the stage of \(\Qq^2\)
alone, and \(s_{12}\) the stage \((\beta,\gamma')\) of the product. By
\eqref{f:index} for \(\Qq^2\),
\[
 i=1\ \Longleftrightarrow\ \Pos_4(s_2;H_r,x\in A^4_{p,1}).
\]
As \(H_r\) is an initial segment of \(H\), \(\Qq^1\) and \(\Qq^2\) belong to
\(\B_4(s;H_r)\), and \(H_r\) is shorter than \(H\). By (i) for \(H_r\),
\(\Qq^2_{\gamma'}*\Qq^1_\beta\in\B_4(s;H_r)\). Its remaining part after
\(\gamma'\) is \(\Qq^1_\beta\), carried out over \(s_2\); it belongs to
\(\B_4(s_2;H_r)\) by Lemma~\ref{f:basic}\,(i), and it reaches \(s_{12}\). As
\(x\in M[G_2\res\gamma']\), (ii) for \(H_r\) gives
\[
 \Pos_4(s_2;H_r,x\in A^4_{p,1})\ \longleftrightarrow\
 \Pos_4(s_{12};H_r,x\in A^4_{p,1}).
\]
So the step satisfies \eqref{f:index} in the product. If \(H\) is empty,
there is no such step.

(ii) Write \(A=\exists b\,P(b)\), and let \(G_\Tt\) be the generic filter
with \(s'=s[G_\Tt]\).

(\(\rightarrow\)) Let \(\Qq\in\B_4(s;H)\) witness \(\Pos_4(s;H,A)\). By
Lemma~\ref{f:basic}\,(iii) we may assume that \(\Qq\) and \(\Tt\) are
disjoint. By (i), \(\Tt*\Qq\in\B_4(s;H)\), so its remaining part after
\(\Tt\), which is \(\Qq\) carried out over \(s'\), belongs to \(\B_4(s';H)\).
Let \(G_\Qq\) be \(\Qq\)-generic over \(M[G_\Tt]\). Then \(G_\Tt\times G_\Qq\)
is \(\Tt\times\Qq\)-generic over \(M\), so there is \(b\) with
\(\cert_3(s[G_\Qq],P(b))\). By (i), \(\Qq*\Tt\) belongs to \(\B_4(s;H)\) and
reaches the same stage as \(\Tt*\Qq\). Its remaining part after \(\Qq\) is
a level-four forcing over \(s[G_\Qq]\), so \(\cert_3\) holds for \(P(b)\) at
the common final stage, by Lemma~\ref{f:persist}\,(i). Hence \(\Qq\),
carried out over \(s'\), witnesses \(\Pos_4(s';H,A)\).

(\(\leftarrow\)) Choose \(t\in G_\Tt\) and a \(\Tt\)-name \(\dot\Qq\) such
that \(t\) forces that \(\dot\Qq\in\B_4(s[\dot G];H)\) witnesses
\(\Pos_4(s[\dot G];H,A)\). Then \((\Tt\res t)*\dot\Qq\in\B_4(s;H)\) by
Lemma~\ref{f:basic}\,(i), and it witnesses \(\Pos_4(s;H,A)\).
\end{proof}

\begin{cor}\label{f:certinv}
Let \(s'\) be reached from a stage \(s\) by a level-four forcing, and let
\(\theta\) be a \(\Pi^1_3\) assertion with parameters from \(M_s\). Then
\(\cert_3(s,\theta)\leftrightarrow\cert_3(s',\theta)\).
\end{cor}
\begin{proof}
Write \(\neg\theta\) as \(\exists z\,R(z)\) with \(R\in\Pi^1_2\), a
\(\Sigma^1_4\) assertion in normal form. We claim that at every stage
\(s''\) whose model contains the parameters of \(\theta\),
\[
 \neg\cert_3(s'',\theta)\ \longleftrightarrow\
 \Pos_4(s'';\emptyset,\exists z\,R(z)).
\]
Every level-four forcing over \(s\) belongs to \(\B_4(s;\emptyset)\), so the
corollary follows from the claim and Proposition~\ref{f:product}\,(ii) for
the empty register. If a condition \(q\) of an ordinary forcing \(\Qq\) over
\(s''\) forces \(\exists z\,R(z)\), then \(\Qq\res q\) witnesses
\(\Pos_4(s'';\emptyset,\exists z\,R(z))\), as \(\cert_3\) for a true
\(\Pi^1_2\) assertion holds by Shoenfield absoluteness. Conversely, let
\(\Qq\in\B_4(s'';\emptyset)\) witness
\(\Pos_4(s'';\emptyset,\exists z\,R(z))\). A witness \(z\) with \(\cert_3\)
for \(R(z)\) satisfies \(R(z)\). So \(\Qq\), with its steps forgotten, is an
ordinary forcing over \(s''\) which forces \(\exists z\,R(z)\), and
\(\cert_3(s'',\theta)\) fails.
\end{proof}

\begin{lem}\label{f:incompat}
Let \(r=(p,d,H_r)\) be an entry of a register, let \(\sigma\) be a stage
with \(\parm(p)\in M_\sigma\) at which \(\Pos_4(\sigma;H_r,I^4_p)\) fails,
and let \(s\) be reached from \(\sigma\) by a forcing in
\(\B_4(\sigma;H_r)\). Then for no \(x\in M_s\) do both
\(\Pos_4(s;H_r,x\in A^4_{p,0})\) and \(\Pos_4(s;H_r,x\in A^4_{p,1})\) hold.
\end{lem}
\begin{proof}
Suppose that both hold, witnessed by \(\Qq^0\) and \(\Qq^1\), which are
disjoint by Lemma~\ref{f:basic}\,(iii). By Proposition~\ref{f:product}\,(i),
\(\Qq^0*\Qq^1\in\B_4(s;H_r)\), and it reaches the same stage as
\(\Qq^1*\Qq^0\). In the extension by \(\Qq^0\) alone there is \(b_0\) with
\(\cert_3\) for \(P_{p,0}(x,b_0)\), and this persists to the common final
stage along the remaining part of \(\Qq^0*\Qq^1\), by
Lemma~\ref{f:persist}\,(i). Likewise there is \(b_1\) with \(\cert_3\) for
\(P_{p,1}(x,b_1)\) at the common final stage, by the other presentation. So
\(\Qq^0*\Qq^1\) witnesses \(\Pos_4(s;H_r,I^4_p)\), with the witness
\((x,b_0,b_1)\). The real \(x\) is part of the witness, so the only real
parameter of \(I^4_p\) is \(\parm(p)\in M_\sigma\). By
Proposition~\ref{f:product}\,(ii), \(\Pos_4(\sigma;H_r,I^4_p)\) holds, a
contradiction.
\end{proof}

\subsection{The iteration}\label{f:iteration}

The iteration \(\Pp_{\omega_1}\) is built in \(\omega_1\) segments of
countable length. The stages at which the segments begin, and their
limits, are the \emph{boundary stages}. For a stage \(\alpha\) of
\(\Pp_{\omega_1}\) we write \(s_\alpha\) for the stage which consists of
the model and the cover at \(\alpha\) and of the set \(B_\alpha\) of
protected leaves, which we keep along the iteration, with \(B_0=\emptyset\).
At each boundary stage \(\sigma\) we have a register \(H^\sigma\), the
register \emph{in effect} at the stages of the segment which begins at
\(\sigma\), and this segment belongs to \(\B_4(s_\sigma;H^\sigma)\).
Registers are extended only at boundary stages, and only by appending
entries. At a limit of boundary stages we take the unions of the sets of
protected leaves and of the registers, and a countable ordinal as cover,
by Lemma~\ref{lem:interface}\,\ref{T:closure}. Some reals \(y\) are
\emph{treated}, each at a single boundary stage \(\sigma_y\); there \(y\)
receives a constructible real \(d_y\), and the entries for the pairs with
parameter \(y\) are appended.

\paragraph{Treating a parameter.}
Let \(\sigma\) be a boundary stage, let \(H\) be the union of the registers
at the earlier boundary stages, and let \(y\in M_\sigma\) be a real which
has not been treated. Let \(d_y\) be the \(<_L\)-least constructible real
such that no leaf \((4,p,d_y,x,i,t)\) with \(\parm(p)=y\) belongs to
\(B_\sigma\). It exists, as \(B_\sigma\) is countable and
\(\omega_1^L=\omega_1\). We append entries to \(H\). Let \(H_0:=H\). If there
is \(e\in\omega\) such that the pair \(p=(e,y)\) has no entry in \(H_\nu\)
and
\begin{equation}\label{f:test}
 \neg\Pos_4\bigl(s_\sigma;H_\nu,I^4_p\bigr),
\end{equation}
we take the least such \(e\) and put \(H_{\nu+1}:=H_\nu\cat(p,d_y,H_\nu)\);
at limits we take unions. Each step uses a new \(e\), so the recursion stops
at a countable ordinal \(\nu^*\), and we let \(H^\sigma:=H_{\nu^*}\). By
Lemma~\ref{f:basic}\,(ii) the recursion is carried out in \(M_\sigma\). Every
entry \(r=(p,d_y,H_r)\) appended at \(\sigma\) satisfies
\(\neg\Pos_4(s_\sigma;H_r,I^4_p)\), and every pair \(p=(e,y)\) without an
entry in \(H^\sigma\) satisfies \(\Pos_4(s_\sigma;H^\sigma,I^4_p)\).

The segment which begins at \(\sigma\) produces common elements for the
pairs without entry. Let \((p_n:n<N)\), \(N\le\omega\), enumerate these
pairs, and let \(s^0:=s_\sigma\). Given the stage \(s^n\), which is
reached from \(s_\sigma\) by a forcing in \(\B_4(s_\sigma;H^\sigma)\),
\(\Pos_4(s^n;H^\sigma,I^4_{p_n})\) holds by
Proposition~\ref{f:product}\,(ii). Let \(\Qq_n\in\B_4(s^n;H^\sigma)\)
witness it, and let \(s^{n+1}\) be the stage reached by \(\Qq_n\). The
segment is the concatenation of the forcings \(\Qq_n\), \(n<N\), which
belongs to \(\B_4(s_\sigma;H^\sigma)\) by Lemma~\ref{f:basic}\,(i). At the
final stage of \(\Qq_n\) there are \(x,b_0,b_1\) such that \(\cert_3\)
holds for \(P_{p_n,0}(x,b_0)\) and for \(P_{p_n,1}(x,b_1)\).

\paragraph{Tasks.}
Every segment carries out one of the following tasks, at its first stage
\(\sigma\). In (2) and (3), \(H^\sigma\) is the union of the registers at
the earlier boundary stages.
\begin{enumerate}[leftmargin=*,topsep=.4em,itemsep=.2em,label=\textup{(\arabic*)}]
\item For a real \(y\in M_\sigma\) which has not been treated, treat \(y\) at
\(\sigma\) as above. At \(\sigma=0\) we treat \(y=0^\omega\); there
\(B_0=\emptyset\), so \(d_{0^\omega}=0^\omega\).
\item For a \(\Pi^1_3\) assertion \(\theta\) with parameters from
\(M_\sigma\): if \(\cert_3(s_\sigma,\theta)\) fails, the segment is an
ordinary forcing over \(s_\sigma\), restricted below a condition which
forces \(\neg\theta\).
\item For an entry \(r=(p,d,H_r)\) of \(H^\sigma\) and a real
\(x\in M_\sigma\), the segment makes one step for \((p,d)\), concerning \(x\)
and the index \(i\) with \(i=1\leftrightarrow\Pos_4(s_\sigma;H_r,x\in
A^4_{p,1})\), on a coding area outside the cover.
\end{enumerate}
If there is nothing to do, the segment is trivial. The segments of (2) are
ordinary, and the segment of (3) belongs to \(\B_4(s_\sigma;H^\sigma)\), as
its step satisfies \eqref{f:index} by the choice of \(i\). So every segment
belongs to \(\B_4(s_\sigma;H^\sigma)\).

At each boundary stage we choose a name for an enumeration of its reals in
order type \(\omega_1\). As in Section~\ref{sec:iteration}, every kind of
task, with every tuple of indices into the enumerations at earlier
boundary stages and every code of a formula, is scheduled cofinally often.
Names for the segments are chosen by the forcing theorem, their lengths are
bounded by countable ordinals by the ccc, and every segment has at least
one stage. So \(\Pp_{\omega_1}\) has length \(\omega_1\), and every real of
the final model belongs to the model at a boundary stage. Every real is
treated, task (2) considers every \(\Pi^1_3\) assertion at arbitrarily late
boundary stages, and task (3) makes a step for every entry and every real.
The iteration preserves all cardinals and \(\CH\), by
Lemma~\ref{lem:interface}\,\ref{T:size}, and, presented as a name over
\(\Pp^0\times\Pp^1\), it yields a set forcing over \(L\). Let
\(V:=W[G_{\omega_1}]\) be the final model.

\begin{lem}\label{f:parts}
Let \(\alpha\le\beta<\omega_1\) be stages of \(\Pp_{\omega_1}\), and let
\(H'\) be an initial segment of the register in effect at \(\alpha\). Then
the part \(\Pp_{\alpha,\beta}\), with the sets \(B_\gamma\),
\(\alpha\le\gamma\le\beta\), belongs to \(\B_4(s_\alpha;H')\).
\end{lem}
\begin{proof}
Let a step for \((p,d)\) be made at a stage \(\gamma\in[\alpha,\beta)\),
where \((p,d,H_r)\) is an entry of \(H'\). It is made in a segment whose
register extends \(H'\), and the entry, with its stored register, is the
same in this register, so the step satisfies \eqref{f:index} at \(\gamma\).
The other requirements hold at the single stages, as every segment is a
level-four forcing. At limits of boundary stages we argue as in
Lemma~\ref{f:basic}\,(i).
\end{proof}

\begin{lem}\label{f:truth}
For every stage \(\alpha<\omega_1\) and every \(\Pi^1_3\) assertion
\(\theta\) with parameters from \(M_\alpha\),
\(V\models\theta\leftrightarrow\cert_3(s_\alpha,\theta)\).
\end{lem}
\begin{proof}
By Lemma~\ref{f:parts} for the empty register, every part of
\(\Pp_{\omega_1}\) between two stages is a level-four forcing. So the
implication from right to left is Lemma~\ref{f:persist}\,(iii). Conversely,
suppose that \(\cert_3(s_\alpha,\theta)\) fails. By
Corollary~\ref{f:certinv}, it fails at every later stage, in particular at a
later boundary stage at which task (2) considers \(\theta\). There the
segment forces the \(\Sigma^1_3\) assertion \(\neg\theta\), which holds in
\(V\) by Shoenfield absoluteness.
\end{proof}

\begin{lem}\label{f:exact}
Let \(B^\infty\) be the union of the sets \(B_\gamma\), \(\gamma<\omega_1\).
For every leaf \(v\) of \(V\),
\(V\models\neg D_3(w(v))\leftrightarrow v\in B^\infty\).
\end{lem}
\begin{proof}
Let \(\alpha\) be a stage with \(v\in M_\alpha\). A leaf of \(M_\alpha\)
which is not protected at \(\alpha\) is never protected later, by the first
part of Lemma~\ref{f:persist}\,(ii); so \(v\in B^\infty\) if and only if
\(v\in B_\alpha\). By Lemma~\ref{f:persist}\,(ii) and as \(w(v)\) determines
\(v\), \(v\in B_\alpha\) if and only if \(\cert_3(s_\alpha,\neg D_3(w(v)))\).
Lemma~\ref{f:truth} completes the proof.
\end{proof}

\subsection{The separators}\label{f:separators}

\begin{proof}[Proof of Proposition~\ref{f:main}]
Let \(y\) be a real of \(V\), treated at the boundary stage
\(\sigma:=\sigma_y\) with the constructible real \(d:=d_y\), and let
\(p=(e,y)\) be a pair such that \(A^4_{p,0}\cap A^4_{p,1}=\emptyset\) in
\(V\). If \(p\) had no entry in \(H^\sigma\), then the segment which begins
at \(\sigma\) would produce \(x,b_0,b_1\) such that \(\cert_3\) holds for
\(P_{p,0}(x,b_0)\) and for \(P_{p,1}(x,b_1)\) at a stage, and \(x\) would be
a common element of the two sets in \(V\), by Lemma~\ref{f:truth}. So \(p\)
has an entry \(r=(p,d,H_r)\), appended at \(\sigma\), and
\(\Pos_4(s_\sigma;H_r,I^4_p)\) fails.

For a stage \(\alpha\ge\sigma\) and a real \(x\in M_\alpha\) let
\(j_\alpha(x):=1\) if \(\Pos_4(s_\alpha;H_r,x\in A^4_{p,1})\), and
\(j_\alpha(x):=0\) otherwise. If \(\sigma\le\alpha\le\beta\), then
\(\Pp_{\alpha,\beta}\in\B_4(s_\alpha;H_r)\) by Lemma~\ref{f:parts}, as \(H_r\)
is an initial segment of the register in effect at \(\alpha\). So
\(j_\alpha(x)=j_\beta(x)\) by Proposition~\ref{f:product}\,(ii), and we write
\(j(x)\).

\emph{The steps for \((p,d)\).} By the choice of \(d\), no leaf
\((4,p,d,x,i,t)\) belongs to \(B_\sigma\). Every later step for \((p,d)\),
inside a segment or by task (3), is made at a stage \(\gamma\ge\sigma\) in
a segment whose register contains \(r\). So it concerns some \(x\) and the
index \(j(x)\), by \eqref{f:index}. Task (3) makes such a step for every
\(x\). By Lemma~\ref{f:exact}, for every real \(x\) of \(V\) there is \(t\)
with \(V\models\neg D_3(w((4,p,d,x,j(x),t)))\), and
\(V\models D_3(w((4,p,d,x,1-j(x),t)))\) for every \(t\). By \eqref{f:sep},
\[
 V\models S^4_{j(x)}(p,d,x)\land\neg S^4_{1-j(x)}(p,d,x).
\]
So the sets \(S^4_0(p,d,\cdot)\) and \(S^4_1(p,d,\cdot)\) are complements of
each other in \(V\). Both are \(\Sigma^1_4(y,d)\), so both are
\(\Delta^1_4(y,d)\).

\emph{Correctness.} Let \(x\in A^4_{p,i}\) in \(V\), witnessed by \(b\), and
let \(\alpha\ge\sigma\) be a stage with \(x,b\in M_\alpha\). By
Lemma~\ref{f:truth}, \(\cert_3(s_\alpha,P_{p,i}(x,b))\), so the trivial
forcing witnesses \(\Pos_4(s_\alpha;H_r,x\in A^4_{p,i})\). If \(i=1\), then
\(j(x)=1\). If \(i=0\), then \(\Pos_4(s_\alpha;H_r,x\in A^4_{p,1})\) fails
by Lemma~\ref{f:incompat}, applied to
\(\Pp_{\sigma,\alpha}\in\B_4(s_\sigma;H_r)\), so \(j(x)=0\). Hence
\(A^4_{p,i}\subseteq S^4_i(p,d,\cdot)\) for \(i<2\), and the
\(\Delta^1_4(y,d)\) set \(S^4_0(p,d,\cdot)\) separates \(A^4_{p,0}\) from
\(A^4_{p,1}\).

Every two disjoint \(\RRrel{\Sigma}{4}\) sets are given by a pair with some
parameter \(y\). For \(y=0^\omega\) we have \(d_y=0^\omega\), so two disjoint
lightface \(\Sigma^1_4\) sets are separated by a lightface \(\Delta^1_4\)
set.
\end{proof}

\subsection{Examples}\label{f:examples}

\paragraph{The toy problem revisited.}
Let \(p\) be the lightface pair of the toy problem, with
\(A^4_{p,0}=\{z:P_1(z)\}\) and \(A^4_{p,1}=\{z:P_0(z)\}\). The reals
\(\langle i,z,t,j\rangle\) are not reals of leaves, so by
Lemma~\ref{f:persist}\,(ii) \(\cert_3\) fails at every stage for
\(\neg D_3(\langle i,z,t,j\rangle)\), and \(V\models
D_3(\langle i,z,t,j\rangle)\) by Lemma~\ref{f:truth}. So
\(V\models\neg P_i(z)\) for all \(i\) and \(z\), and
\(A^4_{p,0}=A^4_{p,1}=\emptyset\). The matrix of a witness for \(P_i(z)\)
contains the conjunct \(\neg D_3(\langle i,z,t,0\rangle)\), which is never
certified. So no membership is ever possible, \(p\) receives an entry when
\(0^\omega\) is treated, and every real is put on the side \(0\).

\paragraph{A pair which reads its own separator.}
By the recursion theorem \cite[Chapter~7]{mosch} there is \(e\) such that
the pair
\(p=(e,0^\omega)\) satisfies
\[
 A^4_{p,0}=\{x:S^4_1(p,0^\omega,x)\},\qquad
 A^4_{p,1}=\{x:S^4_0(p,0^\omega,x)\},
\]
with the witnesses \(t\) of \eqref{f:sep} and the matrices
\(\neg D_3(w(v))\). So \(x\) belongs to \(A^4_{p,0}\) as soon as it is put on
the side \(1\) of the candidate with the index \(d_{0^\omega}=0^\omega\), which
\(p\) would receive. Consider the recursion \eqref{f:test} when
\(0^\omega\) is treated at the initial stage, and let \(H_\nu\) be one of its
registers without an entry for \(p\). Let \(\Qq\) consist of two steps for
\((p,0^\omega)\), concerning the same \(x\) and the indices \(0\) and \(1\),
which protect leaves \(v_0\) and \(v_1\). Then
\(\Qq\in\B_4(s_0;H_\nu)\), as no entry restricts these steps. At its final
stage \(\cert_3\) holds for \(\neg D_3(w(v_0))\) and for
\(\neg D_3(w(v_1))\), by Lemma~\ref{f:persist}\,(ii), and these are the
matrices of \(x\in A^4_{p,1}\) and of \(x\in A^4_{p,0}\). So \eqref{f:test}
fails for \(p\) at every \(\nu\), \(p\) receives no entry, and the segment
produces a common element of the two sets. A pair which can read the side
we are about to choose is thus dealt with before any side is chosen for it.

\paragraph{A particular witness.}
In \eqref{f:possdef}, \(\cert_3\) is applied to the matrix of a particular
witness, not to the assertion \(A\) itself. To see that this matters, let
\(p\) be the lightface pair with \(A^4_{p,0}=\R\) and
\(A^4_{p,1}=\emptyset\), presented with a true and a false matrix. It
receives an entry with the index \(0^\omega\) at the initial stage, and
\(j(x)=0\) for every \(x\), so \(V\models\neg S^4_1(p,0^\omega,0^\omega)\).
But at every stage, and after every ordinary forcing, only countably many
reals have been coded, so \(S^4_1(p,0^\omega,0^\omega)\) holds there. So no
condition of an ordinary forcing forces the negation of this
\(\Sigma^1_4\) assertion, at any stage, although it fails in \(V\).
Equivalently, the \(\Pi^1_4\) sentence \(\neg S^4_1(p,0^\omega,0^\omega)\) is
true in \(V\) and false at every stage.

\paragraph{Why the test is repeated.}
Let \(p\) be the pair of the previous paragraph, and let \(p'=(e',0^\omega)\) be
a pair whose two sets both consist of the reals \(x\) with
\(S^4_1(p,0^\omega,0^\omega)\). In \(V\) both sets are empty. As long as \(p\)
has no entry, a single step for \((p,0^\omega)\), concerning \(0^\omega\) and the
index \(1\), protects a leaf \(v\), and \(\cert_3\) holds for
\(\neg D_3(w(v))\); so a common element of the two sets of \(p'\) is possible.
Once the entry for \(p\) has been appended, every step for \((p,0^\omega)\) has
the index \(0\), and no leaf \((4,p,0^\omega,0^\omega,1,t)\) is protected at the
initial stage. By Lemma~\ref{f:persist}\,(ii), \(S^4_1(p,0^\omega,0^\omega)\) is
then no longer possible, and neither is a common element for \(p'\). So the
outcome of the test for \(p'\) depends on whether \(p\) already has an entry.
This is why \eqref{f:test} is evaluated with the current register \(H_\nu\),
and why the recursion only stops when every remaining pair passes the test
with respect to \(H^\sigma\). Here both pairs receive entries.

\paragraph{Separators which are not \(\Delta^1_3\).}
Let \(C\) be a lightface \(\Pi^1_3\) set which is not \(\Sigma^1_3\). Then
\(C\) and \(\R\setminus C\) are disjoint lightface \(\Sigma^1_4\) sets whose
union is \(\R\), so \(C\) is their only separator, and \(\R\setminus C\) is the
only separator of the pair in the other order. So the separators cannot be
\(\Delta^1_3\) in general, and they cannot be given by the coding at level
three on both sides, as at the third level in \cite{sep3}. This is why
\eqref{f:sep} uses \(D_4\) on both sides.

\section{Global separation}\label{sec:sep}

We now separate the pairs of all levels within one iteration, in the
setting of Section~\ref{f:leaves}. Compared with Section~\ref{sec:four},
the construction changes in four ways.
\begin{enumerate}[leftmargin=*,topsep=.4em,itemsep=.2em,label=\textup{(\alph*)}]
\item \emph{Indices.} The separator of a pair with parameter \(y\) should
only use \(y\). So the index \(d\) is a natural number, which occurs in the
definition of the separator as a numeral. To keep a fresh index available,
every pair receives a finite set \(K(p)\) of admissible indices when its
level is first completed after its parameter appears; these sets form the
first component of an assignment (Definition~\ref{def:assignment}). All
pairs are classified at once, and a possible common point is produced
later, by task~(1) in Section~\ref{sec:iteration}. So possibility has to be
preserved when the assignment grows by data for new parameters. This is
the change of assignment in Proposition~\ref{s:product}\,(ii).
\item \emph{Lower levels.} At a level \(k\ge4\), the certificates are
those of level \(k-1\), which are given by possibility at level \(k-1\) with
respect to the complete assignment at level \(k-1\). A \(k\)-allowable forcing may change the assignments at the
lower levels, at the stages of the sets \(C_\ell\)
(Definition~\ref{def:mallowable}). The persistence of certificates, which
is Lemma~\ref{f:persist}\,(i) at level four, becomes the change of
assignment one level down. So Proposition~\ref{s:product} is proved by
induction on \(k\), and its products have to carry complete assignments at
all lower levels (part~C of its proof).
\item \emph{Nodes.} For \(j\ge5\), the separator at level \(j\) unfolds
through nodes of the levels \(j,j-1,\ldots,4\), and a step puts a node into
one of the sets \(B_k\). The leaves of Section~\ref{sec:four} are the nodes
of level four.
\item \emph{Truth.} Possibility agrees with truth in the final model at
every level (Lemma~\ref{s:truth}), since the certificates at the next level
depend on it.
\end{enumerate}
The registers of Section~\ref{sec:four} correspond to the second component
\(R\) of an assignment. Apart from these changes, the definitions and the
arguments follow Section~\ref{sec:four}.

\subsection{The candidate separators}\label{sec:separators}

For \(j\geq3\) let \(p=(e,y)\), where \(e\in\omega\) codes two
\(\Sigma^1_j\) formulas and \(y\) is a real, and let
\(A^j_{p,0},A^j_{p,1}\) be the two \(\Sigma^1_j(y)\) sets which these
formulas define from \(y\). A natural number \(d\) will distinguish the
reals reserved for this pair from those which were used before it was
known whether the two sets intersect. For a real \(x\) and \(i<2\) let
\(c=c(j,p,d,x,i)\), where \(c\) is the tuple function of
Section~\ref{f:leaves}. For \(j=3\)
the real which is coded for \((p,d,x,i)\) is \(w_3(p,d,x,i):=c\). For
\(j\geq4\) we use finite sequences of reals, which we call \emph{nodes}.
The \emph{roots} are the sequences \(v=(j,p,d,x,i,t)\), \(t\in\R\), and
they have level \(j\). A node \(v\) of level \(k>4\) has the children
\(v\cat u\), \(u\in\R\), which have level \(k-1\). The nodes of level four
are the \emph{leaves}. To each node \(v\) we assign a real \(a(v)\): for a
root let \(a(v):=e_{j-1}(c,t)\), and let
\[
                a(v\cat u):=e_{k-2}(a(v),u)
                     \qquad(v\text{ of level }k>4).
\]
The real which is coded for a leaf \(v\) is \(w(v):=a(v)\). The reals
\(w(v)\), the reals \(w_3(p,d,x,i)\) and the auxiliary real
\(w_{\rm aux}:=w_*\) lie in pairwise disjoint ranges. If a leaf \(v\)
lies below a root of level \(j\), then applying
\(e_3^{-1},e_4^{-1},\ldots,e_{j-1}^{-1}\) successively to \(w(v)\), each time
to the first component and retaining the second, recovers all coordinates
of \(v\) and finally \(c(j,p,d,x,i)\). So \(v\le_T w(v)\) for every
leaf~\(v\).

Besides the coded reals, the iteration carries sets \(B_k\), \(k\geq4\),
of nodes of level \(k\). These sets increase along the iteration, and at
limit stages we take unions. A node is put into \(B_k\) only in the
following way: at a single stage we code the auxiliary real
\(w_{\rm aux}\) and put one node \(v\) of level \(k\) into \(B_k\), where
all coordinates of \(v\) except the last one are chosen before this stage,
and the last coordinate of \(v\) is the generic real added by this coding
forcing. A node may be put into \(B_k\) only if its parent does not belong
to \(B_{k+1}\), and the real \(w(v)\) of a leaf \(v\) may be coded only if
\(v\notin B_4\). When a node is put into \(B_k\), none of its children
belongs to \(B_{k-1}\), and if it is a leaf, then its real has not been
coded before. Indeed, every such earlier step would compute the new last
coordinate.

For nodes \(v\) and for \(p,d,x\) let
\begin{align}
 T_4(v)&\ \longleftrightarrow\ \neg D_3(w(v)),\notag\\
 T_{k+1}(v)&\ \longleftrightarrow\ \forall u\,\neg T_k(v\cat u),
                                                    \label{s:decoder}\\
 S^3_i(p,d,x)&\ \longleftrightarrow\ D_3(w_3(p,d,x,i)),\notag\\
 S^j_i(p,d,x)&\ \longleftrightarrow\
      \exists t\,T_j((j,p,d,x,i,t))\quad(j\ge4).     \label{s:separators}
\end{align}
By induction, using \eqref{eq:common-decoder}, we have
\(T_k(v)\longleftrightarrow\neg D_{k-1}(a(v))\), and hence
\(S^j_i(p,d,x)\longleftrightarrow D_j(c(j,p,d,x,i))\) for every
\(j\geq3\). For fixed \(e\) and \(d\), both sets \(S^j_i(p,d,\cdot)\),
\(i<2\), are \(\Sigma^1_j(y)\). For a pair which is disjoint in the final
model we shall find some \(d\) for which these two sets are complements of
each other and contain \(A^j_{p,0}\) and \(A^j_{p,1}\), respectively.

\subsection{Assignments and allowable forcings}

Fix \(j\), \(p\) and \(d\). A \emph{step for \((p,d)\) at level three} is
the coding of a real of the form \(w_3(p,d,x,i)\). A \emph{step for
\((p,d)\) at level \(k\geq4\)} is a stage at which a root
\((k,p,d,x,i,t)\) is put into \(B_k\). Putting a node which is not a root
into \(B_k\) is also called a step of level \(k\), but such steps are not
subject to the restrictions on pairs \((p,d)\) below. Let \(U_k(p)\) be
the set of all \(d\) such that a step for \((p,d)\) at level \(k\) has
occurred so far; for \(k=3\), every coding of a real of the form
\(w_3(p,d,x,i)\) counts. The countable set of reals coded so far, the
sets \(B_k\), the model \(M\) and the cover \(J\) form the \emph{coding
data}, on which the definitions below depend. If two independent forcings
are carried out in either order, then the coding data at the common final
model are the same. Next we define the assignments which regulate how the
coding data may change.

\begin{defn}\label{def:assignment}
Let \((M,J)\) be a suitable pair and \(k\geq3\). A \emph{level-\(k\)
assignment} is a triple \(H=(K,R,F)\in H(\omega_2)^{N_J}\) with the
following properties.
\begin{enumerate}[leftmargin=*,topsep=.4em,itemsep=.2em,label=\textup{(\roman*)}]
\item \(K\) is a partial function on pairs \(p\) whose values are finite
subsets of \(\omega\).
\item \(R\) is a set of \emph{entries} \(r=(p,d,H_r)\), at most one for
each \(p\), where \(p\in\dom K\), \(d\in K(p)\) and
\(H_r=(K_r,R_r,F_r)\) is a level-\(k\) assignment with \(H_r\preceq H\).
\item \(F\) is a countable set of nodes of level \(k+1\). For \(k=3\)
these are leaves whose reals may not be coded. For \(k\geq4\), no child of
an element of \(F\) may be put into \(B_k\).
\end{enumerate}
Here \(H'\preceq H\) means that \(H'\) is contained in \(H\)
componentwise, where \(K\) is regarded as a set of pairs and entries are
compared together with all their data. We require that the relation
\(r'\in R_r\) on entries is well-founded, of height below \(\omega_2\),
and that every set \(R_r\) is closed under predecessors with respect to
this relation.
\end{defn}

A step for \((p,d)\) at level \(k\) is allowed only if \(d\in K(p)\),
provided that \(p\in\dom K\). If moreover \(r=(p,d,H_r)\in R\), then a
step for \((p,d)\) at level \(k\) which concerns the input \(x\) and the
index \(i\), i.e.\ the coding of \(w_3(p,d,x,i)\) or putting a root
\((k,p,d,x,i,t)\) into \(B_k\), is allowed only if
\begin{equation}\label{s:index}
 i=1\ \Longleftrightarrow\ \Pos_k(a;H_r,\,x\in A^k_{p,1}),
\end{equation}
so that \(x\) is put on side \(0\) if this instance of \(\Pos_k\) fails.
Here, at a stage \(\gamma\), \(a=\da_k(\gamma)\) denotes the data
with respect to which \(\Pos_k\) is evaluated, see \ref{N:data}
below: the current coding data through level \(k\)
(the coded reals and the sets \(B_\ell\), \(\ell\leq k\)), together with
the current complete assignments at the levels below \(k\). The assignment
\(H_r\), including \(F_r\), stays fixed. In particular, \eqref{s:index}
does not take into account restrictions at levels above \(k\) which are
added later.

\paragraph{The predicates \(\Pos_k\) and \(\cert_k\).}
We present a \(\Sigma^1_k\) assertion as an existential real quantifier
followed by a \(\Pi^1_{k-1}\) matrix, and a \(\Pi^1_k\) assertion as the
negation of such an assertion, and we call these presentations
\emph{normal forms}. We allow the matrix to be a finite conjunction of
\(\Pi\)-assertions of the same level in normal form. This includes the
assertion \(I^k_p\) which says that \(A^k_{p,0}\) and \(A^k_{p,1}\) have a
common element: its witness is a triple \((x,b_0,b_1)\), and its two
matrices are those of \(x\in A^k_{p,0}\) and of \(x\in A^k_{p,1}\). For a
conjunction, \(\cert\) means \(\cert\) for each conjunct. The negation of
an existential assertion in normal form can in turn serve as a
\(\Pi\)-assertion in normal form. Ordinary prenex formulas are included.
Let \(\cert_2(a,P)\) mean that the \(\Pi^1_2\) matrix \(P\) is true.
Recursively, let
\begin{align}
 \Pos_k(a;H,\exists b\,P(b))
 &\ \longleftrightarrow\
 \exists\Qq\in\B_k(a;H)\ \
     \Qq\Vdash\exists b\,\cert_{k-1}(a[\dot G],P(b)),\label{s:poss}\\
 \cert_k(s,\neg A)
 &\ \longleftrightarrow\
       \neg\Pos_k(a;(K_k,R_k,B_{k+1}),A).\label{s:cert}
\end{align}
In \eqref{s:cert}, level \(k\) is complete at the stage \(s\), and
\(\Pos_k\) is evaluated with the whole current complete assignment and the
current set \(B_{k+1}\). In \eqref{s:poss}, the final stage of \(\Qq\) is
complete through level \(k-1\); for \(k=3\) only the final model is needed
there. The trivial forcing belongs to \(\B_k(a;H)\), and \(\B_k(a;H)\) is
closed under restrictions below conditions.

\paragraph{Complete assignments.}
A stage is \emph{complete through level \(k\)} if complete level-\(\ell\)
assignments, as defined next, are given there for all \(3\leq\ell\leq k\).
No assignments are needed below level three, so every stage is complete
through level two.

\begin{defn}\label{def:complete}
Let \(s\) be a stage with suitable pair \((M,J)\) which is complete
through level \(k-1\). A \emph{complete level-\(k\) assignment} at \(s\)
is a level-\(k\) assignment \(H_k=(K_k,R_k,B_{k+1})\), where \(B_{k+1}\)
is the current set, with \(\dom K_k=\omega\times\R^M\) and
\(U_k(p)\subseteq K_k(p)\) for all \(p\), such that:
\begin{enumerate}[leftmargin=*,topsep=.4em,itemsep=.2em,label=\textup{(\roman*)}]
\item for every \(r=(p,d,H_r)\in R_k\), \(\Pos_k(a;H_r,I^k_p)\) fails;
\item for every \(p\) without an entry in \(R_k\),
\(\Pos_k(a;H_k,I^k_p)\) holds.
\end{enumerate}
\end{defn}

At later stages at which level \(k\) is completed, the set \(K_k(p)\) and
the entry for \(p\), or the absence of an entry, are kept for every old
\(p\). At the first stage at which an entry for \(p\) appears, its index
\(d\) must not belong to \(U_k(p)\) at that stage. For \(m\)-allowable
forcings, both rules are part of Definition~\ref{def:mallowable}(i).
Adding an entry does
not involve a coding forcing. In the definition of \(\Pos_k\), arbitrary
assignments with these properties are allowed. This freedom allows us to
use the same assignments for both orders of a product.

\paragraph{Allowable forcings at level \(m\).}

\begin{defn}\label{def:mallowable}
Let \(m\geq3\), let \(a\) be the data at a stage which is complete through
level \(m-1\), and let \(H=(K,R,F)\) be a level-\(m\) assignment at this
stage. A forcing \(\Qq\) of countable length \(\delta\) is
\emph{\(m\)-allowable} over \(a\) with respect to \(H\), and we write
\(\Qq\in\B_m(a;H)\), if it is allowable over the suitable pair of this
stage, puts nodes only into sets \(B_k\) with \(k\leq m\), keeps the
level-\(m\) assignment \(H\) fixed and changes assignments only at levels
below \(m\), and if it comes with closed sets of stages
\[
 [0,\delta]=C_2\supseteq C_3\supseteq\cdots\supseteq C_{m-1},
 \qquad 0,\delta\in C_k,
\]
(for \(m=3\) only \(C_2\) is used) such that the following hold. Whether a
stage belongs to \(C_k\), and the assignments at that stage, are given by
names in the corresponding initial segment of \(\Qq\), and \(\Qq\) forces
that the sets \(C_k\) are closed and nested as displayed.
\begin{enumerate}[leftmargin=*,topsep=.4em,itemsep=.2em,label=\textup{(\roman*)}]
\item At the stages \(\gamma\in C_k\), \(3\le k<m\), level \(k\) is
completed: the level-\(k\) assignment at \(\gamma\) is complete at
\(\gamma\) (Definition~\ref{def:complete}); for every earlier stage
\(\gamma'\in C_k\) and every pair \(p\) whose parameter belongs to the
model at \(\gamma'\), the set \(K_k(p)\) and the entry for \(p\), or the
absence of an entry, are the same at \(\gamma\) as at \(\gamma'\); and if
\(\gamma>0\) and the assignment at \(\gamma\) has an entry
\((p,d,H_r)\), but \(p\) has no entry at any earlier stage in \(C_k\),
then \(d\) does not belong to \(U_k(p)\) at \(\gamma\). The components
\(K_k,R_k\) of the resulting assignment remain in effect until the next
stage in \(C_k\). The initial assignments at the levels below \(m\) are
part of \(a\).
\item A step of level \(k\), \(4\le k\le m\), begins and ends at stages in
\(C_{k-1}\), so the levels below \(k\) are complete when its restrictions
are checked. The coding of a real \(w_3(p,d,x,i)\) has this property with
respect to \(C_2\). All steps satisfy the conditions on the last coordinate
and on parents from Section~\ref{sec:separators}. At level \(m\), the set
\(F\) is used as the set of forbidden parents, or, for \(m=3\), of
forbidden leaves. At lower levels the current sets \(B_{k+1}\) are used.
Steps for \((p,d)\) respect \(K(p)\) and \eqref{s:index} for the assignment
in effect. Codings of other reals are unrestricted, except for the
condition on leaves in \(B_4\) whenever it applies.
\item Every \(U_k(p)\), \(3\le k\le m\), is finite at the end of \(\Qq\),
and hence at every stage of \(\Qq\). At limit stages the coding data are
the unions of the earlier ones. The stages in the sets \(C_k\), the
iterands and the choice of the data may depend on the generic filter up
to the stage in question.
\end{enumerate}
\end{defn}

The predicates \(\cert\) at lower levels use only the coding data through
their level, the assignments at lower levels, and the relevant fixed or
current assignment. For \(m=3\), the fixed set \(F\) is exactly the
restriction on coding reals of leaves. For \(m>3\), a coding has to
respect \(B_4\), which is part of the data of \(\Qq\) itself. Restriction
below a condition preserves all requirements which are forced. Two
\(m\)-allowable forcings with the same fixed assignment can be
concatenated. The requirements of Definition~\ref{def:mallowable}(i) on
old pairs and on new entries carry over to the concatenation, as an old
pair keeps its data along both forcings. By the ccc, a name for the length of an \(m\)-allowable
forcing is bounded by a countable ordinal. We extend the forcing beyond
its original end by trivial iterands and repeat its final assignments at
all added stages of the sets \(C_k\). This also works for names in initial
segments for the sets \(C_k\). The concatenation of two \(m\)-allowable
forcings preserves the finiteness of the sets \(U_k(p)\). For countable
concatenations, finiteness at the end and completeness at the new limit
stages have to be proved.

\paragraph{Notation.}
Let \(\Qq=(\Qq_\gamma:\gamma<\delta)\in\B_m(a;H)\) be an \(m\)-allowable
forcing over a stage with suitable pair \((M,J)\), and let \(G\) be
\(\Qq\)-generic over \(M\). For \(\alpha\le\beta\le\delta\) and
\(\gamma\le\delta\) we use the following notation. The same notation is
used for the stages of the main iteration \(\Pp_{\omega_1}\) of
Section~\ref{sec:iteration}.
\begin{enumerate}[leftmargin=*,topsep=.4em,itemsep=.2em,label=\textup{(N\arabic*)}]
\item\label{N:models}
\(G\res\gamma:=G\cap\Qq_\gamma\) and \(M_\gamma:=M[G\res\gamma]\). The
remaining forcing between \(\alpha\) and \(\beta\) is
\(\Qq_{\alpha,\beta}\), so that \(\Qq_\beta\cong\Qq_\alpha*\Qq_{\alpha,\beta}\);
it is a forcing over \(M_\alpha\). If \(\Qq\) uses coding areas from \(Z\),
then by Lemma~\ref{lem:interface}\,\ref{T:seg}, applied to
\(\Qq_\beta\), the forcing \(\Qq_{\alpha,\beta}\) is allowable over
\((M_\alpha,J_\alpha)\), with the covers \((J_\gamma)_{\alpha\le\gamma\le\beta}\)
and coding areas from \(Z\setminus J_\alpha\).

\item\label{N:coding}
\(J_\gamma\) is the cover at \(\gamma\), \(X_\gamma\) is the set of all
reals coded before \(\gamma\), including those coded before the initial
stage of \(\Qq\), and \(B^\gamma_\ell\) and \(U^\gamma_\ell(p)\) are the
sets \(B_\ell\) and \(U_\ell(p)\) at \(\gamma\). We read
\(U^\gamma_3(p)\) from \(X_\gamma\) and, for \(\ell\ge4\),
\(U^\gamma_\ell(p)\) from the roots in \(B^\gamma_\ell\).
\item\label{N:assign}
For \(3\le\ell<m\) let \(\gamma^*_\ell:=\max(C_\ell\cap[0,\gamma])\), which
exists since \(C_\ell\) is closed and \(0\in C_\ell\). The level-\(\ell\)
assignment in effect at \(\gamma\) is
\[
 H^\gamma_\ell=(K^\gamma_\ell,R^\gamma_\ell,B^\gamma_{\ell+1}),
\]
where \((K^\gamma_\ell,R^\gamma_\ell)\) are the first two components of
the complete level-\(\ell\) assignment at \(\gamma^*_\ell\).
In fact \(B^\gamma_{\ell+1}=B^{\gamma^*_\ell}_{\ell+1}\). Each addition
to this set occurs at the end of a level-\((\ell+1)\) step, hence at a
stage in \(C_\ell\) by Definition~\ref{def:mallowable}(ii), and at limits
we take unions. At the fixed level put \(H^\gamma_m:=H\).
\item\label{N:data}
The \emph{data at \(\gamma\) through level \(k\le m\)} are
\[
 \da_k(\gamma):=\bigl(M_\gamma,\,J_\gamma,\,X_\gamma,\,
   (B^\gamma_\ell)_{4\le\ell\le k},\,(H^\gamma_\ell)_{3\le\ell<k}\bigr).
\]
For \(k\ge3\), at a stage complete through level \(k-1\), \(\Pos_k\)
uses \(\da_k(\gamma)\) and the specified level-\(k\) assignment.
At a stage complete through level \(k\), \(\cert_k\) uses the whole
current complete assignment in \eqref{s:cert}. We have \(a=\da_m(0)\).
\item\label{N:par}
For a pair \(p=(e,y)\) let \(\parm(p):=y\), and for an entry
\(r=(p,d,H_r)\) let \(\parm(r):=\parm(p)\). For a node \(v\), \(\lc(v)\)
is its last coordinate.
\item\label{N:newpar}
For a model \(N\) and level-\(k\) assignments
\(H=(K,R,F)\preceq H'=(K',R',F')\) we write \(H\preceq_N H'\) if
\(H'\) extends \(H\) only on parameters outside \(N\):
\[
 \parm(p)\notin N\ \ (p\in\dom K'\setminus\dom K),\quad
 \parm(r)\notin N\ \ (r\in R'\setminus R),\quad
 v\notin N\ \ (v\in F'\setminus F).
\]
\end{enumerate}

\paragraph{Locality of completions.}\phantomsection\label{obs:completion-locality}
For \(3\le\ell<m\) and \(\gamma\in C_\ell\), Definition~\ref{def:complete}
tests a proposed
\((K,R,B^\gamma_{\ell+1})\) using \(\da_\ell(\gamma)\) and
\(U^\gamma_\ell\). Along \(C_\ell\), Definition~\ref{def:mallowable}(i)
requires agreement on old pairs, including the absence of an entry.
Freshness \(d\notin U^\gamma_\ell(p)\) is checked for entries first
added at \(\gamma\).
Thus these checks use only data through level \(\ell+1\).

\paragraph{Summary of the requirements.}
For reference we collect what Section~\ref{sec:separators},
Definitions~\ref{def:assignment}--\ref{def:mallowable} and the remarks
after Definition~\ref{def:complete} require of a forcing
\(\Qq\in\B_m(a;H)\) of length \(\delta\), where \(H=(K,R,F)\). We use the
notation \ref{N:models}--\ref{N:newpar}. Recall that for \(k<m\) the
level-\(k\) assignment in effect at a stage \(\gamma\) is
\(H^\gamma_k=(K^\gamma_k,R^\gamma_k,B^\gamma_{k+1})\), and that
\(H^\gamma_m=H\). The \emph{forbidden set} of level \(k\) at \(\gamma\) is
the third component of \(H^\gamma_k\), that is, \(B^\gamma_{k+1}\) for
\(k<m\) and \(F\) for \(k=m\).
\begin{enumerate}[leftmargin=*,topsep=.4em,itemsep=.2em,label=\textup{(\Roman*)}]
\item\label{rule:base} \emph{Base.} \(\Qq\) is allowable over the suitable
pair of its initial stage (Definition~\ref{def:alpha_allowable_extension}).
In particular, every coding forcing uses a coding area outside the current
cover. At limit stages the coding data are the unions of the earlier ones.
\item\label{rule:nodes} \emph{Nodes.} A node is put into a set \(B_k\) only
if \(4\le k\le m\), and only at a stage \(\gamma\) at which
\(w_{\rm aux}\) is coded. At most one node is put into these sets at
\(\gamma\). Its last coordinate is the generic real of the stage
\(\gamma\), and its other coordinates belong to \(M_\gamma\).
\item\label{rule:parents} \emph{Parents.} If a node with a parent is put
into \(B_k\) at \(\gamma\), then its parent does not belong to the
forbidden set of level \(k\) at \(\gamma\).
\item\label{rule:leaves} \emph{Leaves.} If the real \(w(v)\) of a leaf
\(v\) is coded at \(\gamma\), then \(v\) does not belong to the forbidden
set of level three at \(\gamma\).
\item\label{rule:indices} \emph{Indices.} If a step for \((p,d)\) at level
\(k\le m\) is made at \(\gamma\) and \(p\in\dom K^\gamma_k\), then
\(d\in K^\gamma_k(p)\).
\item\label{rule:sides} \emph{Sides.} If moreover
\((p,d,H_r)\in R^\gamma_k\) and the step concerns the input \(x\) and the
index \(i\), then \eqref{s:index} holds with \(a=\da_k(\gamma)\).
\item\label{rule:timing} \emph{Timing.} The sets
\(C_2\supseteq\cdots\supseteq C_{m-1}\) are closed, contain \(0\) and
\(\delta\), and are given by names in the initial segments of \(\Qq\). If
a step of level \(k\ge4\) is made at \(\gamma\), then
\(\gamma,\gamma+1\in C_{k-1}\).
\item\label{rule:assign} \emph{Assignments.} For \(3\le\ell<m\) and
\(\gamma\in C_\ell\), \(H^\gamma_\ell\) is a level-\(\ell\) assignment
over \((M_\gamma,J_\gamma)\) in the sense of
Definition~\ref{def:assignment}. In particular, it belongs to
\(H(\omega_2)^{N_{J_\gamma}}\).
\item\label{rule:complete} \emph{Completeness.} For such \(\ell\) and
\(\gamma\), \(H^\gamma_\ell\) is complete at \(\gamma\)
(Definition~\ref{def:complete}):
\(\dom K^\gamma_\ell=\omega\times\R^{M_\gamma}\),
\(U^\gamma_\ell(p)\subseteq K^\gamma_\ell(p)\) for all \(p\),
\(\Pos_\ell(\da_\ell(\gamma);H_r,I^\ell_p)\) fails for every entry
\((p,d,H_r)\in R^\gamma_\ell\), and
\(\Pos_\ell(\da_\ell(\gamma);H^\gamma_\ell,I^\ell_p)\) holds for every
\(p\) without an entry in \(R^\gamma_\ell\).
\item\label{rule:coherent} \emph{Coherence.} For \(\gamma'<\gamma\) in
\(C_\ell\), \(3\le\ell<m\),
\[
 K^\gamma_\ell\res\dom K^{\gamma'}_\ell=K^{\gamma'}_\ell
 \qquad\text{and}\qquad
 \{r\in R^\gamma_\ell:\parm(r)\in M_{\gamma'}\}=R^{\gamma'}_\ell .
\]
So for a pair \(p\) with \(\parm(p)\in M_{\gamma'}\), the set
\(K_\ell(p)\) and the entry for \(p\), or the absence of an entry, are
kept.
\item\label{rule:fresh} \emph{Fresh indices.} If
\(\gamma\in C_\ell\setminus\{0\}\), \(3\le\ell<m\),
\((p,d,H_r)\in R^\gamma_\ell\) and \(p\) has no entry in
\(R^{\gamma'}_\ell\) for any \(\gamma'\in C_\ell\cap\gamma\), then
\(d\notin U^\gamma_\ell(p)\).
\item\label{rule:finite} \emph{Finiteness.} \(U^\delta_k(p)\) is finite
for all \(p\) and \(3\le k\le m\).
\end{enumerate}
Here \ref{rule:base}, \ref{rule:timing} and \ref{rule:finite} are part of
Definition~\ref{def:mallowable}; \ref{rule:nodes}--\ref{rule:leaves} come
from Section~\ref{sec:separators} and Definition~\ref{def:mallowable}(ii);
\ref{rule:indices} and \ref{rule:sides} from the rules stated with
\eqref{s:index} and from Definition~\ref{def:mallowable}(ii);
\ref{rule:assign} from Definition~\ref{def:assignment};
and \ref{rule:complete}--\ref{rule:fresh} from
Definition~\ref{def:mallowable}(i), together with
Definition~\ref{def:complete}.

All these requirements are statements which \(\Qq\) forces about its
stages. So they are preserved by restrictions below conditions, and by (a)
and (b) in the proof of Lemma~\ref{s:symmetry} they are invariant under
the translations of the coding areas. Table~\ref{tab:rules} lists where
each of them is verified for the forcings which are constructed below:
the restrictions of Lemma~\ref{s:restriction}, the two presentations of a
product in Proposition~\ref{s:product}, and the main iteration
\(\Pp_{\omega_1}\) of Section~\ref{sec:iteration}, which satisfies the
requirements for every \(m\) between any two stages in \(C_{m-1}\) by
Lemma~\ref{s:extend}. The last column names the place where the
requirement is used.

\begin{table}[H]
\centering\footnotesize
\setlength{\tabcolsep}{3pt}\renewcommand{\arraystretch}{1.3}
\begin{tabular}{@{}>{\raggedright\arraybackslash}p{0.19\textwidth}
  >{\raggedright\arraybackslash}p{0.13\textwidth}
  >{\raggedright\arraybackslash}p{0.2\textwidth}
  >{\raggedright\arraybackslash}p{0.165\textwidth}
  >{\raggedright\arraybackslash}p{0.11\textwidth}
  >{\raggedright\arraybackslash}p{0.135\textwidth}@{}}
\hline
 & restrictions & products
 & \multicolumn{2}{l}{main iteration (Lemma~\ref{s:extend})} & used in\\
 & (Lemma~\ref{s:restriction}) & (Proposition~\ref{s:product})
 & segments, single steps & limits & \\
\hline
\ref{rule:base} Base & \ref{N:models}
 & Lemma~\ref{lem:interface}\,\ref{T:prod}; Lemma~\ref{s:symmetry}
 & Lemma \ref{lem:interface}\,\ref{T:closure}
 & Lemma \ref{lem:interface}\,\ref{T:closure} & \eqref{eq:exact}\\
\ref{rule:nodes} Nodes & proof, (ii); Definition \ref{def:restriction}\,(i)
 & (A1)
 & (E1); tasks (2), (3) & --- & Lemma~\ref{s:nodes}, \eqref{s:B3}\\
\ref{rule:parents} Parents & proof, (ii); \eqref{s:restr1}, \eqref{s:restr2}
 & \eqref{s:A1c}, \eqref{s:A1b}
 & (E1), (E2); task (2) & --- & Lemma~\ref{s:nodes}\\
\ref{rule:leaves} Leaves & proof, end of (ii)
 & \eqref{s:A1c}, \eqref{s:A1b}
 & (E1), (E2); task (2) & --- & Lemma~\ref{s:nodes}\\
\ref{rule:indices} Indices & proof, (ii); \eqref{s:restr1}, \eqref{s:restr2}
 & level \(k\): end of (A1); below: (C2), \ref{V:steps}
 & (E1), (E2); task (3) & --- & \eqref{s:E3}\\
\ref{rule:sides} Sides & proof, (ii); \eqref{s:restr1}
 & level \(k\): proof of \eqref{s:A3}; below: (C1), (C2), \ref{V:steps}
 & (E1), (E2); task (3) & --- & proof of Theorem~\ref{main:sep}\\
\ref{rule:timing} Timing & \eqref{s:restr1}; proof, (i), (ii)
 & paths \(\Qq^i[\zeta]\) before part~A; (C6)
 & \eqref{s:E1}; single steps & (E3) & \ref{N:assign}\\
\ref{rule:assign} Assignments & Definition \ref{def:restriction}\,(iii)
 & (B2), (C3)
 & (E1); single steps & (E3) & Lemma~\ref{s:symmetry}, (B2)\\
\ref{rule:complete} Completeness & proof, (i)
 & (B2)--(B4), (C4), \ref{V:complete}
 & (E1); single steps & (E3) & Lemma~\ref{s:truth}, Theorem~\ref{main:sep}\\
\ref{rule:coherent} Coherence & proof, (i)
 & (B3), (C3), (C4), \ref{V:complete}; before part~A
 & (E1); single steps & (E3) & Lemma~\ref{s:truth}, Theorem~\ref{main:sep}\\
\ref{rule:fresh} Fresh indices & proof, (i)
 & (B4), (C5), \ref{V:fresh}; before part~A
 & (E1); single steps & (E3) & proof of Theorem~\ref{main:sep}\\
\ref{rule:finite} Finiteness & proof, (iii)
 & \eqref{s:A2}, \ref{V:finite}
 & (E1) & \eqref{s:E3} & (B1)\\
\hline
\end{tabular}
\caption{Where the requirements \ref{rule:base}--\ref{rule:finite} are
verified. Parts (A1)--(C6) belong to the proof of
Proposition~\ref{s:product}, parts (E1)--(E3) to the proof of
Lemma~\ref{s:extend}, and the tasks and single steps are those of
Section~\ref{sec:iteration}. In the main iteration, (E1) treats the
segments chosen in task~(1) and, together with the paragraph on single
steps, the steps of tasks~(2) and~(3); (E2) transfers the requirements
to any two stages in \(C_{m-1}\). At the initial stage of the main
iteration, \ref{rule:assign} and \ref{rule:complete} hold by the first
paragraph of Section~\ref{sec:iteration}, and \ref{rule:coherent} and
\ref{rule:fresh} are void. A dash means that there is nothing to verify,
as no step is made at a limit stage.}\label{tab:rules}
\end{table}

\subsection{Restriction to lower levels and symmetry}

\begin{defn}\label{def:restriction}
Let \(\Qq\in\B_m(a;H)\), let \(3\le k<m\), let \(\alpha\le\beta\) be stages
in \(C_{k-1}\), and let \(H'\) be a level-\(k\) assignment over
\((M_\alpha,J_\alpha)\). The \emph{restriction}
\(\Qq_{\alpha,\beta}\resk{k}H'\) is the forcing
\(\Qq_{\alpha,\beta}\) over \(M_\alpha\) with the following data:
\begin{enumerate}[leftmargin=*,topsep=.4em,itemsep=.2em,label=\textup{(\roman*)}]
\item the iterands \(\Qq(\gamma)\), \(\alpha\le\gamma<\beta\), given by
their original names after fixing \(G\res\alpha\), where a step of level
\(\ell>k\) counts only as a coding of \(w_{\rm aux}\), i.e.\ its node is
not put into \(B_\ell\);
\item the sets \(B^\gamma_\ell\), \(4\le\ell\le k\), and the sets
\(C_\ell\cap[\alpha,\beta]\), \(2\le\ell\le k-1\);
\item the complete level-\(\ell\) assignments of \(\Qq\) at the stages in
\(C_\ell\cap[\alpha,\beta]\), \(3\le\ell\le k-1\);
\item the fixed level-\(k\) assignment \(H'\).
\end{enumerate}
The endpoint hypotheses are understood in the forcing sense. For fixed
ordinals \(\alpha,\beta\), membership in \(C_{k-1}\) is read by
\(\Qq_\beta\), and we work throughout below a condition there forcing
both memberships when necessary. For endpoints given by ordinal names,
first decide their values below a condition of \(\Qq\) and apply this
convention to the resulting ordinals.
\end{defn}

\begin{lem}\label{s:restriction}
Let \(\Qq\in\B_m(a;H)\), let \(3\le k<m\), and let \(\alpha\le\beta\) be
stages in \(C_{k-1}\). Then for every level-\(k\) assignment \(H'\)
over \((M_\alpha,J_\alpha)\),
\begin{equation}\label{s:restr}
 H'\preceq H^\alpha_k\ \Longrightarrow\
 \Qq_{\alpha,\beta}\resk{k}H'\in\B_k\bigl(\da_k(\alpha);H'\bigr).
\end{equation}
Moreover, for every level-\(m\) assignment \(H_0\preceq H\) at the
initial stage, \(\Qq\in\B_m(a;H_0)\).
\end{lem}
\begin{proof}
Let \(H'=(K',R',F')\preceq H^\alpha_k\) and
\(\Qq':=\Qq_{\alpha,\beta}\resk{k}H'\), and mark the data of \(\Qq'\) by a
prime. Let \(\alpha\le\gamma\le\beta\). For \(3\le\ell<k\) we have
\(\alpha\in C_{k-1}\subseteq C_\ell\), hence
\(\max(C_\ell\cap[\alpha,\gamma])=\gamma^*_\ell\), and therefore
\begin{equation}\label{s:restr1}
 \begin{gathered}
 \Qq'(\gamma)=\Qq(\gamma)\ \ (\gamma<\beta),\qquad
 C'_\ell=C_\ell\cap[\alpha,\beta]\ \ (2\le\ell<k),\\
 \da'_k(\gamma)=\da_k(\gamma),\qquad
 U'^\gamma_\ell(p)=U^\gamma_\ell(p)\ \ (3\le\ell\le k).
 \end{gathered}
\end{equation}
The two forcings use the same generic filters and covers, hence the same
models \(M_\gamma\) and \(N_{J_\gamma}\). The coded reals and the nodes
through level \(k\) also agree. In particular the sets \(U_\ell(p)\) include
the same indices already used before \(\alpha\), by~\ref{N:coding}.
The iterands are the original names, also if data of
level \(>k\) were used to choose them; this is allowed by
Definition~\ref{def:mallowable}(iii). By~\ref{N:models}, \(\Qq'\) is
allowable over \((M_\alpha,J_\alpha)\) with the inherited covers. We verify
(i)--(iii) of Definition~\ref{def:mallowable} for \(\Qq'\), with \(k\) in
place of \(m\).

(i) We have \(\alpha,\beta\in C_{k-1}\subseteq C_\ell\) for \(2\le\ell<k\),
so \(\alpha,\beta\in C'_\ell\), and the sets
\([\alpha,\beta]=C'_2\supseteq\cdots\supseteq C'_{k-1}\) are closed.
The complete assignments at \(\alpha\) are inherited, so their entries
are old at the start of \(\Qq'\). At each later
\(\gamma\in C'_\ell\), \(3\le\ell<k\), the candidate assignment, the data
through level \(\ell+1\le k\), and the earlier assignments on
\([\alpha,\gamma)\) agree by \eqref{s:restr1}. \hyperref[obs:completion-locality]{Locality of completions} gives (i).

(ii) Let \(\gamma\in[\alpha,\beta)\) be the stage of a step of \(\Qq'\) of
level \(\ell\le k\). It is a step of \(\Qq\) of level \(\ell\), so
\(\gamma,\gamma+1\in C_{\ell-1}\cap[\alpha,\beta]=C'_{\ell-1}\), and if it
puts a node \(v\) into \(B_\ell\), then \(\lc(v)\) is the generic real added
at \(\gamma\).
If \(\ell<k\), its restrictions in \(\Qq'\) are given by
\(B^\gamma_{\ell+1}\), \(H^\gamma_\ell\) and \(\da_\ell(\gamma)\), which are
the same as in \(\Qq\) by \eqref{s:restr1}.
Let \(\ell=k\). As \(\alpha\le\gamma\), the sets \(B_{k+1}\) increase, and
the values of \(K\) and the entries of old pairs are kept at later
completions of level \(k\), we have
\begin{equation}\label{s:restr2}
 F'\subseteq B^\alpha_{k+1}\subseteq B^\gamma_{k+1},\qquad
 K'\subseteq K^\alpha_k\subseteq K^\gamma_k,\qquad
 R'\subseteq R^\alpha_k\subseteq R^\gamma_k .
\end{equation}
In \(\Qq\) the step respects \(H^\gamma_k\), and by \eqref{s:restr2} it
respects \(H'\):
\begin{alignat*}{2}
 v\notin B^\gamma_{k+1}&\ \Longrightarrow\ v\notin F'
   &\qquad&\text{(\(v\) the parent of the new node, \(k\ge4\)),}\\
 d\in K^\gamma_k(p)&\ \Longrightarrow\ d\in K'(p)=K^\gamma_k(p)
   &&\text{(a step for \((p,d)\), \(p\in\dom K'\)).}
\end{alignat*}
If moreover \(r=(p,d,H_r)\in R'\), then \eqref{s:index} for \(r\) reads
\(i=1\Longleftrightarrow\Pos_k(\da_k(\gamma);H_r,\,x\in A^k_{p,1})\) in \(\Qq\) as
well as in \(\Qq'\), since \(r\in R^\gamma_k\) has the same stored assignment
\(H_r\) and \(\da'_k(\gamma)=\da_k(\gamma)\).
A step of \(\Qq\) of level \(>k\) is a coding of \(w_{\rm aux}\) in
\(\Qq'\), which is unrestricted. The coding of the real \(w(v)\) of a leaf
\(v\) requires \(v\notin B^\gamma_4=B'^\gamma_4\) if \(k\ge4\); if \(k=3\),
it requires \(v\notin F'\), which follows from
\(v\notin B^\gamma_4\supseteq B^\alpha_4\supseteq F'\).

(iii) For \(3\le\ell\le k\) and every \(p\),
\(U'^\beta_\ell(p)=U^\beta_\ell(p)\subseteq U^\delta_\ell(p)\), which is
finite. The inherited coding data have the same unions at limits.

This proves \eqref{s:restr}. For the second assertion let
\(H_0=(K_0,R_0,F_0)\preceq H=(K,R,F)\). The fixed assignment enters only
into (ii) at level \(m\), and there \(F_0\subseteq F\), \(K_0\subseteq K\)
and \(R_0\subseteq R\), with the same stored assignments. So every restriction
imposed by \(H_0\) is imposed by \(H\), and \(\Qq\in\B_m(a;H_0)\).
\end{proof}

The proof of Lemma~\ref{s:restriction} with \(k=m\) and \(H'=H\), where
\(H^\gamma_m=H\), also shows that remaining parts between stages in
\(C_{m-1}\) stay \(m\)-allowable:
\begin{equation}\label{s:remaining}
 \alpha\le\beta\ \text{stages in }C_{m-1}\ \Longrightarrow\
 \Qq_{\alpha,\beta}\in\B_m\bigl(\da_m(\alpha);H\bigr).
\end{equation}

\paragraph{Definability of the recursion.}

The predicates \(\Pos_k\) and \(\cert_k\), the classes \(\B_k(a;H)\) and
completeness are defined by one recursion also when the
data are given by names. For a level-\(k\) assignment \(H=(K,R,F)\) let
\[
 \rank(r):=\sup_{s\in R_r}\bigl(\rank(s)+1\bigr)\ \ (r\in R),\qquad
 \htR(R):=\sup_{r\in R}\bigl(\rank(r)+1\bigr).
\]
The \emph{index} of \(\Pos_k(a;H,A)\) and of \(\Qq\in\B_k(a;H)\) is
\((k,\htR(R))\in\omega\times\omega_2\), and the indices are ordered
lexicographically.

\emph{(D1) Dependencies.} At the index \((k,\eta)\) we first determine
membership in \(\B_k(a;H)\), and then we evaluate the existential
statement in \eqref{s:poss}. As \(\htR(R_r)=\rank(r)<\htR(R)\) for
\(r\in R\), these clauses only refer to the base or to smaller indices:
\begin{equation}\label{s:D1}
\begin{array}{@{}l@{\qquad}l@{}}
 \text{allowability over a suitable pair:} & \text{the base, by
   Lemma~\ref{lem:interface}\,\ref{T:codes};}\\[2pt]
 \text{completion and \eqref{s:index} at a level }\ell<k: &
   \text{indices }(\ell,\zeta);\\[2pt]
 \text{\eqref{s:index} at level }k\text{ for }r\in R: &
   \text{the index }(k,\rank(r)),\ \rank(r)<\eta;\\[2pt]
 \cert_{k-1}\text{ at the final stage:} & \text{indices }(k-1,\zeta).
\end{array}
\end{equation}
At the bottom, \(\cert_2(a,P)\) is the satisfaction of the code of \(P\) in
the structure \((\omega,\R^M)\) of second order arithmetic, where \(M\) is
the model of \(a\). Figure~\ref{fig:deps}\,(a) shows these dependencies.

\emph{(D2) Auxiliary forcings.} By Lemma~\ref{lem:interface}\,\ref{T:codes},
the forcings which are allowable over a suitable pair
\((M,J)\) are given by their codes in \(H(\omega_1)^M\), presented on
\(\operatorname{Fn}(\delta,\Lambda)\), and the set of these codes is
definable in \(M\) from \(G^0\times G^1\) and \(J\). The data of a stage,
i.e.\ the assignments, the presentations of iterations and the sets
\(C_k\) and \(B_k\), belong to \(H(\omega_2)^{N_J}\). An element \(x\) of
\(H(\omega_2)\) is given by a nice name for a relation
\(E\subseteq D\times D\), \(D\subseteq\omega_1\), with
\((D,E)\cong(\operatorname{trcl}(\{x\}),\in)\); by the ccc, such a name has
size at most \(\aleph_1\), and so has a countable sequence of such names.
Let \(\mathbb F\) be the set of pairs \((\Pp,\nu)\), where \(\Pp\) is a ccc
forcing over \(W\) of density at most \(\aleph_1\), presented as in the
normal form, and \(\nu\) is a countable sequence of such names for data,
conditions and candidate forcings. Then \(\mathbb F\) is a set, and it is
closed under restrictions below conditions and under the two-step
iterations \(\Pp_\sigma*\dot\Qq\) of an initial segment of the iteration
with a name for a candidate forcing. A name for an ordinal which is used
as a height is decided below a stronger condition.

\emph{(D3) The tables.} For an index \((k,\eta)\) and
\((\Pp,\nu)\in\mathbb F\), with names \(\dot a,\dot H,\dot A\) in \(\nu\),
\(\htR(\dot R)=\check\eta\), and \(\dot A=\exists b\,\dot P(b)\), we
define the set \(T_{k,\eta}(\Pp,\nu)\subseteq\Pp\) of conditions which
force \(\Pos_k(\dot a;\dot H,\dot A)\), and likewise the conditions which
force that a candidate \(\dot\Qq\) belongs to \(\B_k(\dot a;\dot H)\);
requirements on initial segments use \(\Pp\) followed by the respective
initial segment. The clauses are
\begin{equation}\label{s:D3}
 \begin{aligned}
 q\in T_{k,\eta}(\Pp,\nu)\ &\longleftrightarrow\ \forall q'\le q\
   \exists q''\le q'\ \exists\dot\Qq,\dot b\ \bigl[\,q''\Vdash_\Pp
   \dot\Qq\in\B_k(\dot a;\dot H)\\
 &\hphantom{{}\longleftrightarrow{}}\qquad\wedge\
   (q'',\mathbf 1_{\dot\Qq})\Vdash_{\Pp*\dot\Qq}\cert_{k-1}
   \bigl(\dot a[\dot G],\dot P(\dot b)\bigr)\bigr],\\
 q\Vdash_\Pp\neg\Pos_k(\dot a;\dot H,\dot A)\ &\longleftrightarrow\
   \neg\exists q'\le q\ \ q'\in T_{k,\eta}(\Pp,\nu).
 \end{aligned}
\end{equation}
Here the forcing statements on \(\cert_{k-1}\), on completions and
on \eqref{s:index} are read in the tables which have already been
constructed for smaller indices, by \eqref{s:D1}, using
\eqref{s:cert} for \(\cert_{k-1}\), and the requirements on the coding
forcings themselves use the ordinary forcing relation. An entry of \(R\)
which is given by a name is decided below a stronger condition, or its
rank \(\zeta<\eta\) is decided and its name is used at the index
\((k,\zeta)\).

\emph{(D4) Interpretation and correctness.} Let \(\varphi\) be one of the
statements
\[
 \dot\Qq\in\B_k(\dot a;\dot H),\qquad \Pos_k(\dot a;\dot H,\dot A),\qquad
 \dot H\ \text{is a complete level-\(k\) assignment at }\dot a ,
\]
with names from \(\nu\), where \((\Pp,\nu)\in\mathbb F\). We write
\(q\Vdash^*_\Pp\varphi\) if \(q\) belongs to the table of \(\varphi\) from
(D3), where the table of the third statement is defined from
Definition~\ref{def:complete} in the same way, and
\(q\Vdash^*_\Pp\neg\varphi\) if no \(q'\le q\) satisfies
\(q'\Vdash^*_\Pp\varphi\). The statement
\(\cert_{k-1}(\dot a[\dot G],\dot P(\dot b))\) in \eqref{s:D3} stands for the
statements \(\neg\Pos_{k-1}\) which \eqref{s:cert} gives for the conjuncts
of \(\dot P(\dot b)\), and for \(k=3\) for the ordinary forcing relation.
For a \(\Pp\)-generic filter \(G\) over \(W\) we say that
\begin{equation}\label{s:D4}
 \varphi\ \text{\emph{holds in} }W[G]\quad:\Longleftrightarrow\quad
 \exists q\in G\ \ q\Vdash^*_\Pp\varphi .
\end{equation}
The following holds for every index, every \((\Pp,\nu)\in\mathbb F\) and
every \(\Pp\)-generic filter \(G\) over \(W\).
\begin{enumerate}[leftmargin=*,topsep=.4em,itemsep=.2em,label=\textup{(\alph*)}]
\item The set of all \(q\) with \(q\Vdash^*_\Pp\varphi\) is open, and it
contains \(q\) if it is dense below \(q\). So exactly one of
\(q\Vdash^*_\Pp\varphi\) and \(q\Vdash^*_\Pp\neg\varphi\) occurs for some
\(q\in G\), and
\[
 q\Vdash^*_\Pp\varphi\ \Longleftrightarrow\ \varphi\text{ holds in }W[G']
 \text{ for every \(\Pp\)-generic }G'\ni q .
\]
\item In \(M:=W[G]\) the clauses \eqref{s:index}, \eqref{s:poss},
\eqref{s:cert} and those of Definitions~\ref{def:complete}
and~\ref{def:mallowable} are true, where a statement about the extension of
\(M\) by a candidate forcing \(\Qq=\dot\Qq^G\), or by one of its initial
segments \(\Qq_\gamma\), is read in the sense of \eqref{s:D4} for
\(\Pp*\dot\Qq\), respectively for \(\Pp*\dot\Qq_\gamma\). In particular,
\(\Pos_k(a;H,A)\) holds in \(M\) if and only if there are \(\Qq\) and
\(q'\in\Qq\) such that \(\Qq\in\B_k(a;H)\) holds in \(M\) and, for every
\(\Qq\)-generic filter \(G'\ni q'\) over \(M\), there is \(b\) such that
\(\cert_{k-1}(a[G'],P(b))\) holds in \(W[G*G']\).
\item Whether \(\varphi\) holds in \(W[G]\) depends only on the model
\(W[G]\) and on the values of the names of \(\nu\) under \(G\); it does not
depend on \(\Pp\), \(G\) and \(\nu\) themselves.
\end{enumerate}
We prove (a)--(c) simultaneously by induction on the index. For (a),
\eqref{s:D3} describes the conditions below which a certain set \(D\) of
conditions is dense. The set \(D\) is open, by (a) for \(\B_k\) at the same
index and for \(\Pos_{k-1}\) at smaller indices, and the tables of \(\B_k\)
and of completeness are given by requirements on all stronger conditions.
This gives the first sentence of (a); the rest follows by genericity.

For (b) we treat \(\Pos_k\); the other clauses are verified in the same
way. Let \(q\in G\) with \(q\Vdash^*_\Pp\Pos_k(\dot a;\dot H,\dot A)\). By
\eqref{s:D3} and genericity there are \(q''\in G\), \(\dot\Qq\) and
\(\dot b\) with
\[
 q''\Vdash^*_\Pp\dot\Qq\in\B_k(\dot a;\dot H),\qquad
 (q'',\mathbf 1)\Vdash^*_{\Pp*\dot\Qq}
    \cert_{k-1}(\dot a[\dot G],\dot P(\dot b)).
\]
So \(\Qq:=\dot\Qq^G\), \(q':=\mathbf 1\) and \(b:=\dot b^{G*G'}\) are as
required. Conversely, let \(\Qq=\dot\Qq^G\) and \(q'=\dot q^G\) be as in
(b). Let \(E\) be the set of all \(\dot r^G\le q'\) in \(\Qq\) for which
there are \(r\in G\) and \(\dot b\) with
\((r,\dot r)\Vdash^*_{\Pp*\dot\Qq}\cert_{k-1}(\dot a[\dot G],\dot P(\dot b))\).
By (a) and (c) at smaller indices, the hypothesis on the filters \(G'\)
says that \(E\) is dense below \(q'\). This is a statement about \(G\) in
\(W[G]\), so it is forced by some \(q_0\in G\), and we may assume
\(q_0\Vdash^*_\Pp\dot\Qq\in\B_k(\dot a;\dot H)\). Let \(q_1\le q_0\). Take a
\(\Pp\)-generic filter through \(q_1\), an element \(\dot r\) of the dense
set, with \(r\) and \(\dot b\), and \(q''\le q_1,r\) in this filter. The
table of \(\B_k\) is closed under restrictions below conditions, so
\[
 q''\Vdash^*_\Pp\dot\Qq\res\dot r\in\B_k(\dot a;\dot H),\qquad
 (q'',\mathbf 1)\Vdash^*_{\Pp*(\dot\Qq\res\dot r)}
   \cert_{k-1}(\dot a[\dot G],\dot P(\dot b)),
\]
the second by (a) and (c) at smaller indices, as the generic extensions by
\(\Pp*(\dot\Qq\res\dot r)\) through \((q'',\mathbf 1)\) are those by
\(\Pp*\dot\Qq\) through \((q'',\dot r)\). So \(q_0\) belongs to the table
\eqref{s:D3}.

For (c), by (b) the truth of \(\varphi\) in \(W[G]\) is characterized by a
recursion which only mentions \(W[G]\), the values of the names, and
statements of smaller index in extensions of \(W[G]\), for which (c) is the
induction hypothesis.

Consequently \(\B_k\), \(\Pos_k\), \(\cert_k\) and completeness are defined
in all models under consideration by formulas of set theory with the
parameter \(G^0\times G^1\), in which \(k\) is a variable. By Replacement,
the sequence \((H^\sigma_k)_{k\ge3}\) of the completions at one stage
\(\sigma\) exists. Proposition~\ref{s:product} and Lemma~\ref{s:extend}
below are statements about these formulas, and their proofs by induction
on \(k\) are inductions inside \(\ZFC\). The truth of projective
assertions in the final model enters only in Lemmas~\ref{s:truth}
and~\ref{s:nodes} and in the proof of Theorem~\ref{main:sep}.

\begin{lem}\label{s:symmetry}
Let \((M,J)\) be the suitable pair at a stage \(\beta<\omega_1\). A
witness to \(\Pos_k(a;H,A)\) can be chosen such that it uses coding areas
from a countable set in \(L\) which is disjoint from \(J\) and from any
prescribed countable set. If the given data belong to \(N_J\), then the
trivial condition of \(\Cc_J\) over \(N_J\) decides whether
\(\Pos_k(a;H,A)\) holds.
\end{lem}
\begin{proof}
Let \(\bar G:=G^1\res(\omega_1\setminus J)\), so \(M=N_J[\bar G]\) by
Lemma~\ref{lem:interface}\,\ref{T:loc}. Let \(\pi\in L\) be a permutation of
\(\omega_1\) with \(\pi\res J=\mathrm{id}\), and let \(\pi^*\) be the
automorphism of \(\Cc_J\) which is induced by \(\pi\). By
Lemma~\ref{lem:interface}\,\ref{T:trans}, the coding area with index \(\pi(\eta)\)
computed from \(\pi^*[\bar G]\) is the coding area with index \(\eta\)
computed from \(\bar G\). So a computation from \(\bar G\) is translated
into a computation from \(\pi^*[\bar G]\) by replacing each coding area
\(\eta\) by \(\pi(\eta)\) and each cover \(J'\) by \(\pi[J']\). The coding
forcings are the same, and for every cover \(J'\) the model
\(N_{\pi[J']}\) computed from \(\pi^*[\bar G]\) is the model \(N_{J'}\)
computed from \(\bar G\). By the same item,
allowability over \((M,J)\) and over the later suitable pairs, together with
the sets of coding areas and the covers, is invariant under this
translation. Under this translation the coded
reals, the sets \(C_k\), the sets \(B_k\) and the assignments are
unchanged, and all fixed old parameters have check names over \(N_J\).

We make the induction on the indices explicit. For a \(\Cc_J\)-generic
filter \(\bar G'\) over \(N_J\), a superscript \(\bar G'\) marks a notion
which is computed in \(N_J[\bar G']\) from
\(G^0\times((G^1\res J)\times\bar G')\), see Lemma~\ref{lem:interface}\,\ref{T:trans}.
For the code \(c\) of a forcing, together with its data as an
\(m\)-allowable forcing, let \(\pi\cdot c\) be obtained by replacing every
coding area \(\eta\) by \(\pi(\eta)\), the set \(Z\) of coding areas by
\(\pi[Z]\) and every cover \(J'\) by \(\pi[J']\); all other data are kept.
For data \(a\) at a stage with cover \(J'\), let \(\pi\cdot a\) be \(a\)
with \(\pi[J']\) in place of \(J'\). By the first part of the proof, \(c\)
computed from \(\bar G\) and \(\pi\cdot c\) computed from \(\pi^*[\bar G]\)
are the same iteration, with the same intermediate models and the same
models \(N_{J'}\). By induction on the index \((k,\eta)\) of (D1) we show
\phantomsection
\begin{equation*}\tag{\(\ast\)}\label{s:inv}
 \begin{aligned}
 c\in\B_k(a;H)^{\bar G}\ &\Longleftrightarrow\
   \pi\cdot c\in\B_k(\pi\cdot a;H)^{\pi^*[\bar G]},\\
 \Pos_k(a;H,A)^{\bar G}\ &\Longleftrightarrow\
   \Pos_k(\pi\cdot a;H,A)^{\pi^*[\bar G]},
 \end{aligned}
\end{equation*}
and likewise for the completeness of \(H\) at \(a\).
\begin{enumerate}[leftmargin=*,topsep=.4em,itemsep=.2em,label=\textup{(\alph*)}]
\item \emph{Clauses without \(\Pos\).} Allowability over the suitable pair
is invariant by Lemma~\ref{lem:interface}\,\ref{T:trans}. The requirements of
Definition~\ref{def:mallowable} on the sets \(C_\ell\), on last
coordinates, parents and leaves, on \(K(p)\) and on the sets \(U_\ell(p)\)
only mention the iterands, the generic reals of the stages, and reals,
nodes and assignments, which are the same for \(c\) and \(\pi\cdot c\). An
assignment at a stage with cover \(J'\) has to belong to
\(H(\omega_2)^{N_{J'}}\), and this model is the same on both sides.
\item \emph{Clauses with \(\Pos\).} By \eqref{s:D1}, the remaining
requirements for \(c\in\B_k(a;H)\) are instances of \eqref{s:index} and of
Definition~\ref{def:complete} at stages of the forcing, with smaller
indices. They are invariant by the induction hypothesis at these stages.
\item \emph{\(\Pos_k\) and completeness.} The map \(c\mapsto\pi\cdot c\) is
a bijection between the codes on the two sides, with inverse
\(c\mapsto\pi^{-1}\cdot c\). By (a) and (b) it preserves membership in
\(\B_k(a;H)\), and by the induction hypothesis at the final stage it
preserves the statement which is forced in \eqref{s:poss}; the forcing
relation is that of the same forcing over the same model. This gives the
second line of \eqref{s:inv}. Completeness is a conjunction of instances of
\(\Pos_k\) and of their negations.
\end{enumerate}
The versions given by names are invariant by the same argument below
conditions. If \(a\) is the data at the stage with cover \(J\), then
\(\pi\cdot a=a\), and \eqref{s:inv} gives
\phantomsection
\begin{equation*}\tag{\(\ast\ast\)}\label{s:inv0}
 \Pos_k(a;H,A)^{\bar G}\ \Longleftrightarrow\ \Pos_k(a;H,A)^{\pi^*[\bar G]}.
\end{equation*}

By the ccc, the coding areas used by a witness are contained in
\(J\cup Z\) for a fixed countable set \(Z\in L\) disjoint from \(J\), and
some condition of \(\Cc_J\) forces that there is a witness which uses
coding areas from \(Z\). Given any condition \(q\in\Cc_J\), move \(Z\) and
the support of this condition, by a permutation \(\pi\) as above, to a new
countable set above the prescribed set and above both supports. The moved
condition is compatible with \(q\), and it forces that the translated
witness exists and uses coding areas from \(\pi[Z]\). So the conditions which force that there is a witness whose coding
areas avoid the prescribed set are dense, which gives the first
conclusion. For the second one, let \(c\) force \(\Pos_k(a;H,A)\) and let
\(q\) be any condition. Choose \(\pi\) such that \(\pi^*(c)\) is compatible
with \(q\). By \eqref{s:inv0}, \(\pi^*(c)\) forces \(\Pos_k(a;H,A)\) as
well, so \(q\) does not force its negation. The countable sets of ordinals which occur in this
argument belong to \(N_J\), so it suffices to bound each set which has to
be avoided by an ordinal. The size of the check name of a complete
assignment does not add coordinates to a condition.
\end{proof}

\subsection{Products and complete assignments}

Two \(m\)-allowable forcings \(\Qq^1\) and \(\Qq^2\) over the same stage,
with suitable pair \((M,J)\), are \emph{disjoint} if they use coding
areas from sets \(Z_1,Z_2\in L\) with \(Z_1\cap Z_2=\emptyset\). Let
\(G_1\times G_2\) be \(\Qq^1\times\Qq^2\)-generic over \(M\), let
\(G_i\res\xi:=G_i\cap\Qq^i_\xi\), and consider the intermediate models
\[
 M_{\beta,\gamma}:=M[G_1\res\beta][G_2\res\gamma].
\]
By Lemma~\ref{lem:interface}\,\ref{T:prod}, the product can be presented as an
allowable forcing over \((M,J)\) in either order, also for initial
segments, with the same reals added at the stages of \(\Qq^1\) and of
\(\Qq^2\). For pairs of stages \(v=(\beta,\gamma)\) and
\(w=(\beta',\gamma')\),
\begin{equation}\label{s:intersection}
                 \R^{M_v}\cap\R^{M_w}=\R^{M_{v\wedge w}},
\end{equation}
where \(v\wedge w\) is the coordinatewise minimum. Indeed, if
\(v\not\le w\) and \(w\not\le v\), say \(\beta<\beta'\) and \(\gamma'<\gamma\),
then \(M_v\) and \(M_w\) are mutually generic extensions of \(M_{v\wedge w}\)
by \(G_2\res[\gamma',\gamma)\) and \(G_1\res[\beta,\beta')\), and two mutually
generic extensions have only the reals of the ground model in common
\cite{solovay}. In
particular, a real
in \(M_{\beta,\gamma}\) which belongs to \(M[G_1]\) already belongs to
\(M[G_1\res\beta]\). We now show that both presentations also satisfy the
restrictions on the sets \(B_k\) and on the entries of the assignments.

\begin{prop}\label{s:product}
For every \(k\ge3\) the following hold.
\begin{enumerate}[leftmargin=*,topsep=.4em,itemsep=.2em,label=\textup{(\roman*)}]
\item If \(\Qq^1,\Qq^2\in\B_k(a;H)\) are disjoint, then their product,
presented as a two-step iteration in either order, is \(k\)-allowable, and
both presentations have the same final stage, including all assignments
below level \(k\). One family of assignments works for all pairs of
stages whose coordinates belong to the relevant sets \(C_\ell\).
\item Let \(\Qq\in\B_k(a;H)\) have length \(\delta\), let
\(a':=\da_k(\delta)\), and let \(H'\) be a level-\(k\) assignment at the
final stage of \(\Qq\) with \(H\preceq_M H'\). Then, for every
assertion \(A\) in normal form with parameters from \(M\),
\begin{equation}
             \Pos_k(a;H,A)\ \longleftrightarrow\
             \Pos_k(a';H',A). \label{s:change}
\end{equation}
In particular, \(\Pos_k\) is invariant under \(k\)-allowable forcing if
the assignment is kept fixed.
\end{enumerate}
\end{prop}

\begin{figure}[tp]
\centering
\begin{tikzpicture}[
  font=\footnotesize, >={Stealth[length=2.2mm]},
  bx/.style={draw, rounded corners=2pt, align=flush center, inner sep=3pt},
  dep/.style={->, semithick},
]
\node[bx, text width=5.6cm] (top) at (0,0)
  {\(\B_k(a;H)\) and \(\Pos_k(a;H,A)\)\\ at the index \((k,\eta)\), where \(\eta=\htR(R)\)};
\node[bx, text width=3.9cm] (same) at (-4.7,-2.3)
  {\eqref{s:index} for an entry \(r\in R\):\\ the index \((k,\rank(r))\),
   \(\rank(r)<\eta\)};
\node[bx, text width=3.9cm] (lower) at (0,-2.3)
  {completions and \eqref{s:index} at levels \(\ell<k\), and
   \(\cert_{k-1}\) at the final stage: indices \((\ell,\zeta)\)};
\node[bx, text width=3.9cm] (base) at (4.7,-2.3)
  {allowability over suitable pairs (Lemma~\ref{lem:interface}\,\ref{T:codes});
   \(\cert_2\) is truth in \((\omega,\R^M)\)};
\draw[dep] (same) -- (top);
\draw[dep] (lower) -- (top);
\draw[dep] (base) -- (top);
\node[anchor=west, font=\small] at (-6.8,0.55) {(a)};
\end{tikzpicture}

\bigskip
\begin{tikzpicture}[
  font=\footnotesize, >={Stealth[length=2.2mm]},
  bx/.style={draw, rounded corners=2pt, align=flush center, inner sep=3pt,
    text width=2.6cm, minimum height=9mm},
  fbx/.style={bx, dashed},
  dep/.style={->, semithick},
  num/.style={font=\scriptsize, fill=white, inner sep=1pt},
]
\node[fbx] (Ap) at (0,0)     {part A at level \(k-1\)};
\node[fbx] (Bp) at (4.5,0)   {part B at level \(k-1\)};
\node[fbx] (Cp) at (9,0)     {part C at level \(k-1\)};
\node[bx]  (A)  at (0,-2.9)  {part A at level \(k\)\\ (A1)--(A5)};
\node[bx]  (B)  at (4.5,-2.9){part B at level \(k\)\\ (B1)--(B4)};
\node[bx]  (C)  at (9,-2.9)  {part C at level \(k\)\\ (C1)--(C6)};
\node[fbx] (An) at (0,-5.8)  {part A at level \(k+1\)};
\node[fbx] (Bn) at (4.5,-5.8){part B at level \(k+1\)};
\node[fbx] (Cn) at (9,-5.8)  {part C at level \(k+1\)};
\draw[dep] (Ap) -- (Bp);
\draw[dep] (Bp) -- (Cp);
\draw[dep] (Ap.north east) to[bend left=12] (Cp.north west);
\draw[dep] (An) -- (Bn);
\draw[dep] (Bn) -- (Cn);
\draw[dep] (An.north east) to[bend left=12] (Cn.north west);
\draw[dep] (A) -- node[num]{(3)} (B);
\draw[dep] (B) -- node[num]{(5)} (C);
\draw[dep] (A.north east) to[bend left=12] node[num]{(4)} (C.north west);
\draw[dep] ([xshift=-8mm]Ap.south) -- node[num]{(1)} ([xshift=-8mm]A.north);
\draw[dep] (Cp.south) |- ++(0,-0.85) -| node[num, pos=0.25]{(2)}
  ([xshift=8mm]A.north);
\draw[dep] ([xshift=-8mm]A.south) -- node[num]{(1)} ([xshift=-8mm]An.north);
\draw[dep] (C.south) |- ++(0,-0.85) -| node[num, pos=0.25]{(2)}
  ([xshift=8mm]An.north);
\node[anchor=west, font=\small] at (-1.6,0.85) {(b)};
\end{tikzpicture}

\bigskip
\begin{tikzpicture}[
  font=\footnotesize, >={Stealth[length=2.2mm]},
  sb/.style={draw, align=flush center, inner sep=3pt, minimum height=8mm},
  dep/.style={->, semithick},
]
\node[sb] (A12) at (0,0)   {(A1), (A2)};
\node[sb] (A3)  at (3.4,0) {(A3) at \(\eta\)};
\node[sb] (A4)  at (6.4,0) {(A4) at \(\eta\)};
\node[sb] (A5)  at (9.2,0) {(A5)};
\node[sb, text width=4.6cm] (Hr) at (3.4,-1.7)
  {(A3) and (A4) for the stored \(H_r\), \(r\in R\),
   at the heights \(\htR(R_r)<\eta\)};
\draw[dep] (A12) -- (A3);
\draw[dep] (A3) -- (A4);
\draw[dep] (A4) -- (A5);
\draw[dep] (Hr) -- (A3);
\node[anchor=west, font=\small] at (-1.6,0.75) {(c)};
\end{tikzpicture}
\caption{Dependencies in the recursion (D1)--(D4) and in the proof of
Proposition~\ref{s:product}. An arrow from \(X\) to \(Y\) means that the
definition or the proof of \(Y\) uses \(X\).
(a)~By (D1), \(\B_k(a;H)\) and \(\Pos_k(a;H,A)\) at the index \((k,\eta)\)
only use smaller indices in the lexicographic order and the base. Along
this order, (D3) defines the tables, and (D4)\,(a)--(c) and
\eqref{s:inv} in the proof of Lemma~\ref{s:symmetry} are proved by
induction.
(b)~The induction on \(k\) in the proof of Proposition~\ref{s:product}.
(1)~Assertion~(ii) at level \(k-1\) gives the persistence of
certificates \eqref{s:A4a}.
(2)~Part~C at level \(k-1\), applied with \(m=k\), provides the families
valid through level \(k-1\) which part~A uses for products of forcings
in \(\B_k(a;H)\).
(3)~(B3) uses assertion~(ii) at level \(k\).
(4)~(C1) uses \eqref{s:A3} and assertion~(ii) at level \(k\).
(5)~(C4) applies part~B at the new pairs of stages.
For \(k=3\) there is no lower level: the family is empty, and
\eqref{s:A4a} is Shoenfield absoluteness.
(c)~Inside part~A, \eqref{s:A3} and \eqref{s:A4} are proved by induction
on \(\eta=\htR(R)\). The proofs use Lemma~\ref{s:restriction} in (A3),
(A4), (A5), (B3) and (C1), and Lemma~\ref{s:symmetry} in (A4) and (B2).
Outside the proposition, Lemma~\ref{s:extend} uses part~B at every level,
and Lemma~\ref{s:truth} uses assertion~(ii) at every
level.}\label{fig:deps}
\end{figure}
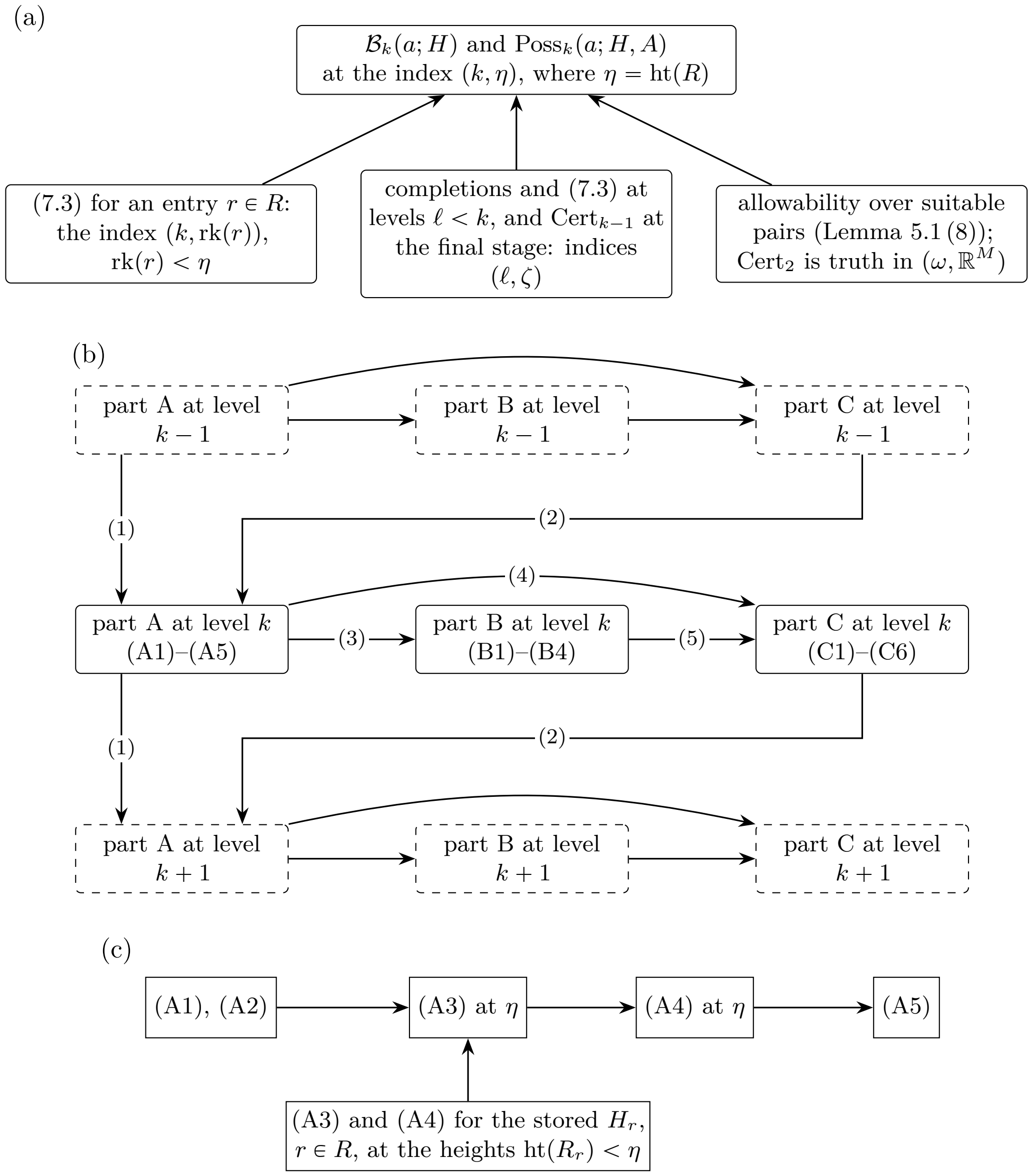

\begin{proof}
We argue by induction on \(k\). For each \(k\) we first compare \(\Pos_k\)
on a product, assuming that assignments at the lower levels on the
product have been chosen (part~A). Then we construct the level-\(k\)
assignments which are needed for the next step of the induction (parts~B
and~C). Figure~\ref{fig:deps}\,(b) and~(c) show how the parts depend on
each other.

\paragraph{Products of two forcings.}
Let \(\Qq^1,\Qq^2\in\B_k(a;H)\) be disjoint, of lengths \(\delta_1,\delta_2\),
with \(H=(K,R,F)\), and let \(G_1\times G_2\) be
\(\Qq^1\times\Qq^2\)-generic over \(M\). The data of \(\Qq^i\) alone carry
the superscript \(i\): \(C^i_\ell\), \(J^i_\xi\), \(X^i_\xi\),
\(B^{i,\xi}_\ell\), \(U^{i,\xi}_\ell(p)\), \(H^{i,\xi}_\ell\),
\(\da^i_k(\xi)\). Pairs of stages
are ordered componentwise. For \(v=(\beta,\gamma)\le(\delta_1,\delta_2)\)
let \(M_v=M[G_1\res\beta][G_2\res\gamma]\) as above, and
\begin{equation}\label{s:pdata}
 J_v:=J^1_\beta\cup J^2_\gamma,\quad
 X_v:=X^1_\beta\cup X^2_\gamma,\quad
 B^v_\ell:=B^{1,\beta}_\ell\cup B^{2,\gamma}_\ell,\quad
 U^v_\ell(p):=U^{1,\beta}_\ell(p)\cup U^{2,\gamma}_\ell(p).
\end{equation}
For \(i\in\{1,2\}\) and \(\zeta\le\delta_{3-i}\) let
\(v_1(\xi,\zeta):=(\xi,\zeta)\) and \(v_2(\xi,\zeta):=(\zeta,\xi)\), and let
\(\Qq^i[\zeta]\) be \(\Qq^i\), carried out over \(M_{v_i(0,\zeta)}\); its
stage \(\xi\) is the pair \(v_i(\xi,\zeta)\). The presentation
\(\Qq^1_\beta*\Qq^2_\gamma\) of \(\Qq^1_\beta\times\Qq^2_\gamma\) runs
through \(\Qq^1[0]\) up to \(\beta\) and then through \(\Qq^2[\beta]\) up
to \(\gamma\), and symmetrically for \(\Qq^2_\gamma*\Qq^1_\beta\). The
data of a step of \(\Qq^i\) are given by \(\Qq^i\)-names over \(M\), so
every step of \(\Qq^i\) involves only reals and nodes of \(M[G_i]\).

A \emph{family of lower assignments} is a family of level-\(\ell\)
assignments at the pairs \(v\),
\[
 \mathcal H=(H^v_\ell)_{3\le\ell<k,\ v\in C^1_\ell\times C^2_\ell},\qquad
 H^{(\beta,0)}_\ell=H^{1,\beta}_\ell,\quad H^{(0,\gamma)}_\ell=H^{2,\gamma}_\ell,
\]
where the equalities on the axes are equalities of full assignments,
including the stored \(H_r\). Along \(\Qq^i[\zeta]\),
\(\zeta\in C^{3-i}_\ell\), the level-\(\ell\) assignment in effect at a pair
\(v\) is formed from \(\mathcal H\) and \(C^i_\ell\) as in~\ref{N:assign},
with \(B^v_{\ell+1}\) from \eqref{s:pdata}. An assignment
\(H'=(K',R',F')\) at \(v\) \emph{agrees on old pairs} with an assignment
\(H''=(K'',R'',F'')\) at \(v''\le v\) if
\begin{equation}\label{s:oldpairs}
 K'\res\dom K''=K''\qquad\text{and}\qquad
 \{r\in R':\parm(r)\in M_{v''}\}=R''.
\end{equation}
The family is \emph{valid through level \(k-1\)} if the following hold for
all \(3\le\ell<k\), \(i\in\{1,2\}\), \(\zeta\in C^{3-i}_\ell\) and \(p\),
where \(v(\xi):=v_i(\xi,\zeta)\):
\begin{enumerate}[leftmargin=*,topsep=.4em,itemsep=.2em,label=\textup{(V\arabic*)}]
\item\label{V:complete} for \(\xi\in C^i_\ell\), \(H^{v(\xi)}_\ell\) is a
complete level-\(\ell\) assignment at \(v(\xi)\) which agrees on old pairs
with \(H^{v(\xi')}_\ell\) for all \(\xi'\in C^i_\ell\cap\xi\);
\item\label{V:fresh} for \(\xi\in C^i_\ell\setminus\{0\}\), if
\((p,d,H_r)\in R^{v(\xi)}_\ell\) and \(p\) has no entry in
\(R^{v(\xi')}_\ell\) for any \(\xi'\in C^i_\ell\cap\xi\), then
\(d\notin U^{v(\xi)}_\ell(p)\);
\item\label{V:steps} every step of \(\Qq^i\) of level \(\ell\) at a stage
\(\xi\) respects the level-\(\ell\) assignment in effect at \(v(\xi)\);
\item\label{V:finite} \(U^{v(\delta_i)}_\ell(p)\) is finite.
\end{enumerate}
No requirement of level \(k\) occurs in \ref{V:complete}--\ref{V:finite}.
For \(v\in C^1_{k-1}\times C^2_{k-1}\) we have
\(v\in C^1_\ell\times C^2_\ell\) for all \(\ell<k\), as the sets
\(C^i_\ell\) are nested, so
\[
 \da_k(v):=\bigl(M_v,\,J_v,\,X_v,\,(B^v_\ell)_{4\le\ell\le k},\,
   (H^v_\ell)_{3\le\ell<k}\bigr),\qquad
 \da_k(\beta,0)=\da^1_k(\beta),\quad \da_k(0,\gamma)=\da^2_k(\gamma).
\]
Let \((\beta,\gamma)\in C^1_{k-1}\times C^2_{k-1}\). The restriction of
\(\mathcal H\) to \(\Qq^1_\beta\times\Qq^2_\gamma\) is again valid. In the
presentation \(\Qq^1_\beta*\Qq^2_\gamma\), the set \(C_\ell\),
\(2\le\ell<k\), consists of the pairs \((\xi,0)\) with
\(\xi\in C^1_\ell\cap[0,\beta]\) and the pairs \((\beta,\zeta)\) with
\(\zeta\in C^2_\ell\cap[0,\gamma]\), and the level-\(\ell\) assignment at
such a pair is the one given by \(\mathcal H\); symmetrically for
\(\Qq^2_\gamma*\Qq^1_\beta\). As \(0,\beta\in C^1_\ell\) and
\(0,\gamma\in C^2_\ell\) for \(\ell<k\), these sets are closed and nested
and contain the initial and the final pair. So
\ref{V:complete}--\ref{V:finite} for the paths \(\Qq^1[0]\),
\(\Qq^2[\beta]\), \(\Qq^2[0]\) and \(\Qq^1[\gamma]\) give the requirements
of the changing assignments for both presentations of
\(\Qq^1_\beta\times\Qq^2_\gamma\). For the agreement on old pairs of the
assignments at \((\xi,0)\) and \((\beta,\zeta)\) we pass through
\((\beta,0)\), as agreement on old pairs is transitive: if
\(v'''\le v''\le v\), then \(M_{v'''}\subseteq M_{v''}\). For the fresh
indices, an entry which appears for the first time at a pair of the
presentation also appears there for the first time along the path which
contains this pair, so \ref{V:fresh} applies.

\paragraph{A. The fixed level \(k\) for a given lower family.}
Assume that the proposition holds below \(k\). Let \(\Qq^1,\Qq^2\) be as
above, and let \(\mathcal H\) be valid through level \(k-1\); for \(k=3\),
\(\mathcal H=\varnothing\). For \(k>3\), part~C at level \(k-1\), applied
with \(m=k\), gives for any two disjoint \(\Qq^1,\Qq^2\in\B_k(a;H)\)
\[
 \mathcal H\text{ valid through }k-2\ \Longrightarrow\
 \exists\,\mathcal H'\supseteq\mathcal H\text{ valid through }k-1;
\]
assertion~(i) at level \(k-1\) alone yields only levels \(\le k-2\).

\emph{(A1) New nodes.} Let \(\Qq^i\) put a node \(v\) into \(B_\ell\) at
its stage \(\xi\), and let \(\zeta\le\delta_{3-i}\). Then
\begin{equation}\label{s:A1}
 \lc(v)\in\R^{M_{v_i(\xi+1,0)}}\setminus\R^{M_{v_i(\xi,\zeta)}}
 \qquad\text{and}\qquad v\notin M[G_{3-i}].
\end{equation}
Indeed, \(\lc(v)\) is the generic real of the stage \(\xi\) of \(\Qq^i\),
so \(\lc(v)\notin M_{v_i(\xi,0)}\). Now
\(v_i(\xi+1,0)\wedge v_i(\xi,\zeta)=v_i(\xi,0)\) and
\(v_i(\xi+1,0)\wedge v_i(0,\delta_{3-i})=(0,0)\), so \eqref{s:intersection}
gives both assertions. With \(v\le_T w(v)\) for leaves
(Section~\ref{sec:separators}) we obtain, for all \(\eta\le\delta_i\) and
\(\ell\),
\begin{equation}\label{s:A1c}
 \bigl(B^{i,\eta}_\ell\setminus B^{(0,0)}_\ell\bigr)\cap M[G_{3-i}]=\varnothing,
 \qquad
 \bigl\{w(v):v\in B^{i,\eta}_4\setminus B^{(0,0)}_4\bigr\}\cap M[G_{3-i}]=\varnothing.
\end{equation}
For a step of level \(\ell\) at the pair \(u\), write
\[
 F^u_\ell=
 \begin{cases}
     B^u_{\ell+1},&3\le\ell<k,\\
     F,&\ell=k .
 \end{cases}
\]
These are the forbidden sets in Definition~\ref{def:mallowable}(ii) for
the fixed assignment \(H=(K,R,F)\). At \(u=v_{3-i}(\xi',\eta)\), a step
of \(\Qq^{3-i}\) involves only nodes and reals of \(M[G_{3-i}]\).
At each lower level, the active factor respects its own current forbidden
set, including the initial nodes. By \eqref{s:A1c}, the nodes added by the other
factor and the reals of the leaves it adds all lie outside \(M[G_{3-i}]\). At level \(k\) both factors use the same \(F\).
Together with \eqref{s:A1}, this gives
\begin{equation}\label{s:A1b}
 \begin{aligned}
 &\lc(v')\notin M_u,\qquad
     \operatorname{parent}(v')\notin F^u_{\ell'}\quad
                       \text{if \(v'\) has a parent},\\
 &w(v'')\text{ coded at }u\ \Longrightarrow\ v''\notin F^u_3 ,
 \end{aligned}
\end{equation}
where \(v'\) is put into \(B_{\ell'}\) by the step and \(v''\) is a leaf.
In particular the condition on leaves uses \(F\) when \(k=3\) and
\(B^u_4\) when \(k>3\). Each factor also respects the fixed \(K\).

\emph{(A2) Used indices.} For all \(p\) and \(\ell\le k\),
\begin{equation}\label{s:A2}
 U^{(\delta_1,\delta_2)}_\ell(p)=U^{1,\delta_1}_\ell(p)\cup
 U^{2,\delta_2}_\ell(p)\ \text{is finite},\qquad
 \parm(p)\notin M[G_i]\ \Longrightarrow\ U^{i,\delta_i}_\ell(p)=U^{(0,0)}_\ell(p).
\end{equation}

\emph{(A3) and (A4).} We prove them by induction on \(\eta:=\htR(R)\),
simultaneously for all data, with (A3) before (A4) at each \(\eta\). Since
\[
 \htR(R_r)=\rank(r)<\eta\qquad(r\in R),
\]
(A3) at \(\eta\) may use (A3) and (A4) for every \(H_r\), also for limit
\(\eta\), and (A4) at \(\eta\) uses (A3) at \(\eta\).
\begin{enumerate}[leftmargin=*,topsep=.4em,itemsep=.2em,label=\textup{(A\arabic*)},start=3]
\item For every valid family and all
\((\beta,\gamma)\in C^1_{k-1}\times C^2_{k-1}\),
\begin{equation}\label{s:A3}
 \Qq^1_\beta*\Qq^2_\gamma\in\B_k(a;H)\quad\text{and}\quad
 \Qq^2_\gamma*\Qq^1_\beta\in\B_k(a;H),
\end{equation}
with the data \eqref{s:pdata} and \(\mathcal H\).
\item For every \(\Tt\in\B_k(a;H)\) of length \(\delta\) and every assertion
\(A\) in normal form with parameters from \(M\),
\begin{equation}\label{s:A4}
 \Pos_k(a;H,A)\ \longleftrightarrow\ \Pos_k(\da_k(\delta);H,A).
\end{equation}
\end{enumerate}

\emph{Proof of \eqref{s:A3}.} By \ref{V:complete}--\ref{V:finite}, \eqref{s:A1b} and
\eqref{s:A2}, only \eqref{s:index} remains to be checked for the steps
governed by entries of \(R\); if \(R=\varnothing\), there is nothing to
check. Let \(r=(p,d,H_r)\in R\) govern a step of \(\Qq^1\) at its stage
\(\xi\) which concerns \((x,i)\) and is carried out at \((\xi,\zeta)\),
where \(\zeta\in C^2_{k-1}\), namely \(\zeta=0\) or \(\zeta=\gamma\). The
case of \(\Qq^2\) is symmetric. Then \(\xi\in C^1_{k-1}\),
\(x\in M_{(\xi,0)}\), and in \(\Qq^1\)
\[
 i=1\ \Longleftrightarrow\ \Pos_k(\da_k(\xi,0);H_r,\,x\in A^k_{p,1}).
\]
As \(H_r\preceq H\), Lemma~\ref{s:restriction} gives
\(\Qq^1,\Qq^2\in\B_k(a;H_r)\), with the same valid family, and
\(\htR(R_r)<\htR(R)\). By \eqref{s:A3} for \(H_r\),
\(\Qq^1_\xi*\Qq^2_\zeta\in\B_k(a;H_r)\). By \eqref{s:remaining}, its
remaining part \(\Qq^2[\xi]\) up to \(\zeta\) belongs to
\(\B_k(\da_k(\xi,0);H_r)\), with final data \(\da_k(\xi,\zeta)\). So
\eqref{s:A4} for \(H_r\) gives
\[
 \Pos_k(\da_k(\xi,0);H_r,x\in A^k_{p,1})\ \longleftrightarrow\
 \Pos_k(\da_k(\xi,\zeta);H_r,x\in A^k_{p,1}),
\]
and the step satisfies \eqref{s:index} at \((\xi,\zeta)\).

\emph{Proof of \eqref{s:A4}.} Write \(A=\exists b\,P(b)\) and
\(a':=\da_k(\delta)\).
(\(\rightarrow\)) Let \(\Qq\) witness \(\Pos_k(a;H,A)\). By
Lemma~\ref{s:symmetry} we may assume that \(\Tt\) and \(\Qq\) are
disjoint. We use the setting above with \(\Qq^1:=\Tt\) and \(\Qq^2:=\Qq\),
of lengths \(\delta\) and \(\delta_2\), write \(G_\Tt:=G_1\) and
\(G_\Qq:=G_2\), and put
\[
 s_0:=(0,\delta_2),\qquad s_1:=(\delta,\delta_2).
\]
By \eqref{s:A3} and \eqref{s:remaining}, \(\Qq^2[\delta]\) belongs to
\(\B_k(a';H)\) and has the final data \(\da_k(s_1)\). By agreement on the
axes, \(\da_k(s_0)\) are the final data of \(\Qq\) alone, so there is
\(b\in M_{s_0}\) with \(\cert_{k-1}(\da_k(s_0),P(b))\); fix it. It suffices to
show
\begin{equation}\label{s:A4a}
 \cert_{k-1}\bigl(\da_k(s_0),P(b)\bigr)\ \longrightarrow\
 \cert_{k-1}\bigl(\da_k(s_1),P(b)\bigr),
\end{equation}
since then \(\Qq^2[\delta]\) witnesses \(\Pos_k(a';H,A)\). For \(k=3\),
\(P(b)\) is a \(\Pi^1_2\) matrix and \(M_{s_0}\subseteq M_{s_1}\), so
\eqref{s:A4a} is Shoenfield absoluteness. Let \(k>3\). By \eqref{s:A3},
\eqref{s:remaining} and Lemma~\ref{s:restriction},
\[
 \Qq^1[\delta_2]\resk{k-1}H^{s_0}_{k-1}\in
 \B_{k-1}\bigl(\da_{k-1}(s_0);H^{s_0}_{k-1}\bigr),
\]
with final data \(\da_{k-1}(s_1)\). Moreover
\(H^{s_0}_{k-1}\preceq_{M_{s_0}}H^{s_1}_{k-1}\): by
Definition~\ref{def:complete}, \(\dom K^{s_0}_{k-1}=\omega\times\R^{M_{s_0}}\);
by~\ref{V:complete} along \(\Qq^1[\delta_2]\), \(H^{s_1}_{k-1}\) agrees with
\(H^{s_0}_{k-1}\) on old pairs; and \(\lc(v)\notin M_{s_0}\) for
\(v\in B^{s_1}_k\setminus B^{s_0}_k\) by \eqref{s:A1}.
So \eqref{s:change} at level \(k-1\) gives, for every
conjunct \(\neg A'\) of \(P(b)\),
\[
 \Pos_{k-1}\bigl(\da_{k-1}(s_0);H^{s_0}_{k-1},A'\bigr)\ \longleftrightarrow\
 \Pos_{k-1}\bigl(\da_{k-1}(s_1);H^{s_1}_{k-1},A'\bigr),
\]
that is, \(\cert_{k-1}(s_0,\neg A')\leftrightarrow\cert_{k-1}(s_1,\neg A')\)
by \eqref{s:cert}. This proves \eqref{s:A4a}.
The direction \eqref{s:A4a} can also be checked directly. Write
\(H_i:=H^{s_i}_{k-1}\). If \(\mathbb W\in\B_{k-1}(\da_{k-1}(s_1);H_1)\)
witnessed \(\Pos_{k-1}(\da_{k-1}(s_1);H_1,A')\), then
\(\mathbb W\in\B_{k-1}(\da_{k-1}(s_1);H_0)\) by Lemma~\ref{s:restriction},
as \(H_0\preceq H_1\). The restricted remaining part
\(\Qq^1[\delta_2]\resk{k-1}H_0\), followed by a name for \(\mathbb W\),
would then witness \(\Pos_{k-1}(\da_{k-1}(s_0);H_0,A')\), contradicting
\(\cert_{k-1}(s_0,\neg A')\).

(\(\leftarrow\)) Let \(G_\Tt\) be \(\Tt\)-generic over \(M\) with
\(\Pos_k(a';H,A)\) in \(M[G_\Tt]\). Choose \(t\in G_\Tt\) and a
\(\Tt\)-name \(\dot\Qq'\) such that \(t\) forces that \(\dot\Qq'\in\B_k(\dot
a';H)\) witnesses \(\Pos_k(\dot a';H,A)\). Then
\((\Tt\res t)*\dot\Qq'\in\B_k(a;H)\) witnesses \(\Pos_k(a;H,A)\).

\emph{(A5) Change of assignment.} Let \(\Tt\), \(a'\) and \(A\) be as in
\eqref{s:A4}, and let \(H'=(K',R',F')\) be a level-\(k\) assignment at the
final stage of \(\Tt\) with \(H\preceq_M H'\).
(\(\rightarrow\)) Let \(\Qq\) and \(\Qq^2[\delta]\in\B_k(a';H)\) be as in
the proof of \eqref{s:A4}. As \(H'\in N_{J_\delta}\subseteq M[G_\Tt]\) and
\(H\preceq_M H'\),
\[
 \parm(p),\ \parm(r),\ u\ \in M[G_\Tt]\setminus M\qquad
 \bigl(p\in\dom K'\setminus\dom K,\ r\in R'\setminus R,\ u\in F'\setminus F\bigr).
\]
Since \((\delta,0)\wedge(0,\delta_2)=(0,0)\), \eqref{s:intersection} gives
that none of these reals and nodes belongs to \(M[G_\Qq]\). Nor does
the real \(w(u)\) of a leaf \(u\in F'\setminus F\), since \(u\le_T w(u)\);
this holds even if \(\lc(u)\) is old. The steps of
\(\Qq^2[\delta]\) involve only reals and nodes of \(M[G_\Qq]\). So
\(H'\) imposes no restriction on them beyond those of \(H\). The lower
data do not depend on the fixed level-\(k\) assignment, by
\hyperref[obs:completion-locality]{Locality of completions}, and neither does \(\cert_{k-1}\) at the final stage.
So \(\Qq^2[\delta]\in\B_k(a';H')\) witnesses \(\Pos_k(a';H',A)\).
(\(\leftarrow\)) If \(\Qq'\in\B_k(a';H')\) witnesses \(\Pos_k(a';H',A)\),
then \(\Qq'\in\B_k(a';H)\) by Lemma~\ref{s:restriction}, and
\eqref{s:A4}~(\(\leftarrow\)) gives \(\Pos_k(a;H,A)\).

Thus (A3) proves assertion~(i) for every family valid through level
\(k-1\), and (A4), (A5) prove assertion~(ii) at level \(k\).

\paragraph{B. Completing level \(k\) at one stage.}
Assume that \eqref{s:change} holds at level \(k\), i.e.\ part~A has been
carried out at level \(k\). In this part \(s_0,s_1\) denote stages of an
\(m\)-allowable forcing with \(m>k\), of the main iteration, or pairs of
stages of a product as above, with the notation
\ref{N:models}--\ref{N:newpar}. Let \(s_1\) be complete through level
\(k-1\), with lower assignments in \(N_{J_{s_1}}\), and let
\[
 N:=N_{J_{s_1}},\quad a:=\da_k(s_1),\quad P:=\omega\times\R^{M_{s_1}},\quad
 B:=B^{s_1}_{k+1},\quad U:=U^{s_1}_k,
\]
where \(U(p)\) is finite for every \(p\in P\). The \emph{inherited data} are
sets \(\bar K,\bar R\) such that \((\bar K,\bar R,B)\) is a level-\(k\)
assignment over \((M_{s_1},J_{s_1})\) and
\begin{equation}\label{s:B0}
 \bar P:=\dom\bar K\subseteq P,\qquad U(p)\subseteq\bar K(p)\quad(p\in\bar P).
\end{equation}

\emph{(B1) The completion.} Let \(P^1:=P\setminus\bar P\), fix a
well-order \(\lhd\in N\) of \(P^1\), and let
\[
 d_p:=\min(\omega\setminus U(p)),\qquad
 K:=\bar K\cup\bigl\{(p,\,U(p)\cup\{d_p\}):p\in P^1\bigr\}.
\]
By recursion on \(\nu\) let \(H_\nu:=(K,R_\nu,B)\), where \(R_0:=\bar R\),
\(R_\lambda:=\bigcup_{\nu<\lambda}R_\nu\) for limit \(\lambda\), and, as long
as \(E_\nu\neq\varnothing\),
\begin{equation}\label{s:B1}
 \begin{gathered}
 E_\nu:=\bigl\{p\in P^1:\ p\text{ has no entry in }R_\nu,\ \
   \neg\Pos_k(a;H_\nu,I^k_p)\bigr\},\\
 R_{\nu+1}:=R_\nu\cup\bigl\{(p_\nu,d_{p_\nu},H_\nu)\bigr\},\qquad
 p_\nu:=\min\nolimits_\lhd E_\nu .
 \end{gathered}
\end{equation}
Let \(\nu^*\) be least with \(E_{\nu^*}=\varnothing\), and put
\(R:=R_{\nu^*}\) and \(H_k:=(K,R,B)\).

The pairs \(p_\nu\) are pairwise distinct elements of \(P^1\). The model
\(M_{s_1}\) is an extension of \(W\) by an allowable forcing, so by
Lemma~\ref{lem:interface}\,\ref{T:size} it has the cardinals of \(L\) and
satisfies \(\CH\). As
\(L\subseteq N\subseteq M_{s_1}\) and \(\R^N=\R^{M_{s_1}}\) by
Lemma~\ref{lem:interface}\,\ref{T:loc},
\[
 \omega_1^N=\omega_1^{M_{s_1}},\qquad \omega_2^N=\omega_2^{M_{s_1}},\qquad
 N\models\CH,\qquad \nu^*<\omega_2 ;
\]
an injection in \(N\) of \(\omega_2^N\) into \(\R\) would collapse
\(\omega_2^{M_{s_1}}\) in \(M_{s_1}\).

\emph{(B2) New pairs.} \(H_k\in N\) is a level-\(k\) assignment over
\((M_{s_1},J_{s_1})\) with
\begin{equation}\label{s:B2}
 \begin{gathered}
 \dom K=P,\qquad K\res\bar P=\bar K,\qquad U(p)\subseteq K(p)\ \ (p\in P),\\
 R\setminus\bar R\ \text{consists of entries for pairs in }P^1,\\
 (p,d,H_r)\in R,\ p\in P^1\ \Longrightarrow\
   d=d_p\notin U(p)\ \wedge\ \neg\Pos_k(a;H_r,I^k_p),\\
 p\in P^1\text{ without entry in }R\ \Longrightarrow\ \Pos_k(a;H_k,I^k_p).
 \end{gathered}
\end{equation}
So conditions (i) and (ii) of Definition~\ref{def:complete} hold for the
pairs in \(P^1\), and the index of every new entry is not in \(U(p)\).

Indeed, the displayed properties hold by \eqref{s:B1} and since
\(E_{\nu^*}=\varnothing\). For \(H_k\in N\): by
Lemma~\ref{s:symmetry}, for every \(\nu\) with \(H_\nu\in N\) the trivial
condition of \(\Cc_{J_{s_1}}\) decides over \(N\) every instance
\(\Pos_k(a;H_\nu,I^k_p)\), and
\(M_{s_1}=N[G^1\res(\omega_1\setminus J_{s_1})]\) by
Lemma~\ref{lem:interface}\,\ref{T:loc}. So the recursion \eqref{s:B1}, carried out
in \(N\) with \(\neg\Pos_k\) replaced by ``the trivial condition of
\(\Cc_{J_{s_1}}\) forces \(\neg\Pos_k\)'', gives the same sequence
\((H_\nu)_{\nu\le\nu^*}\), by induction on \(\nu\). For
Definition~\ref{def:assignment}: each \(p_\nu\) receives exactly one entry,
pairs in \(\bar P\) have entries only in \(\bar R\), \(d_{p_\nu}\in
K(p_\nu)\), and \(H_\nu\preceq H_k\) as \(R_\nu\subseteq R\). If
\(r\in R_{\nu+1}\setminus R_\nu\), then \(R_r=R_\nu\), so by induction on
\(\nu\) every \(R_\nu\) is closed under predecessors and
\(\htR(R_\nu)\le\htR(\bar R)+\nu<\omega_2\). Finally \(|K|,|R|\le\aleph_1\)
in \(N\), so \(H_k\in H(\omega_2)^N\).

\emph{(B3) Old pairs.} Let \(s_0\) be a stage which is complete through
level \(k\), let \(\Tt\in\B_k(\da_k(s_0);H^{s_0}_k)\) have the final stage
\(s_1\), i.e.\ the final data \(a\), and assume
\begin{equation}\label{s:B3}
 K^{s_0}_k\subseteq\bar K,\quad
 \{r\in\bar R:\parm(r)\in M_{s_0}\}=R^{s_0}_k,\quad
 B^{s_0}_{k+1}\subseteq B,\quad
 \lc(v)\notin M_{s_0}\ \ (v\in B\setminus B^{s_0}_{k+1}).
\end{equation}
Then \(U(p)\subseteq K^{s_0}_k(p)=\bar K(p)\) for
\(p\in\omega\times\R^{M_{s_0}}\). If moreover \eqref{s:B0} holds, then
\(H_k\) agrees on old pairs with \(H^{s_0}_k\) in the sense of
\eqref{s:oldpairs}, and satisfies (i) and (ii) of
Definition~\ref{def:complete} for all \(p\in\omega\times\R^{M_{s_0}}\).

Indeed, the steps of \(\Tt\) for \((p,d)\) at level \(k\) respect
\(K^{s_0}_k(p)\) by Definition~\ref{def:mallowable}(ii),
\(U^{s_0}_k(p)\subseteq K^{s_0}_k(p)\) by completeness at \(s_0\), and
\(U\) is read from \(a\) by~\ref{N:coding}; this gives the first assertion.
As \(\dom K^{s_0}_k=\omega\times\R^{M_{s_0}}\) and
\(P^1\cap\dom K^{s_0}_k=\varnothing\), \eqref{s:B3} and \eqref{s:B2} give
\eqref{s:oldpairs} for \(H_k\) and \(H^{s_0}_k\), and
\begin{equation}\label{s:B3a}
 H^{s_0}_k\preceq_{M_{s_0}}H_k ,
\end{equation}
since new pairs and new entries have parameters outside \(M_{s_0}\), and
\(v\notin M_{s_0}\) for \(v\in B\setminus B^{s_0}_{k+1}\) by the last
part of \eqref{s:B3}. Let \(p\in\omega\times\R^{M_{s_0}}\).
\begin{enumerate}[leftmargin=*,topsep=.4em,itemsep=.2em,label=\textup{(\roman*)}]
\item If \(r=(p,d,H_r)\in R^{s_0}_k\), then \(H_r\preceq H^{s_0}_k\), so
\(\Tt\in\B_k(\da_k(s_0);H_r)\) by Lemma~\ref{s:restriction}, and
\eqref{s:change} with \(H=H'=H_r\) gives
\[
 \neg\Pos_k\bigl(\da_k(s_0);H_r,I^k_p\bigr)\ \longrightarrow\
 \neg\Pos_k\bigl(a;H_r,I^k_p\bigr).
\]
\item If \(p\) has no entry in \(R^{s_0}_k\), hence none in \(R\), then
\eqref{s:change} with \(H=H^{s_0}_k\) and \(H'=H_k\), by \eqref{s:B3a}, gives
\[
 \Pos_k\bigl(\da_k(s_0);H^{s_0}_k,I^k_p\bigr)\ \longrightarrow\
 \Pos_k\bigl(a;H_k,I^k_p\bigr).
\]
\end{enumerate}
In both cases the left side holds by completeness at \(s_0\).

\emph{(B4) Summary.} Let \(S\) be a set of stages \(s_0\) as in (B3) with
\(\bar P=\bigcup_{s_0\in S}\omega\times\R^{M_{s_0}}\). Then \eqref{s:B0}
holds by (B3), and by (B2) and (B3)
\begin{equation}\label{s:B4}
 H_k\ \text{is complete at }s_1,\qquad
 H_k\ \text{agrees on old pairs with }H^{s_0}_k\ \ (s_0\in S),
\end{equation}
and the index of every entry first added at \(s_1\) is not in \(U(p)\). For
\(S=\varnothing\), i.e.\ \(\bar K=\bar R=\varnothing\), this is (B1) alone;
this applies at the initial stage. Along a single forcing
\(\Qq\in\B_m(\da_m(0);H)\), \(m>k\), for \(s_0<s_1\) in \(C_k\) and
\((\bar K,\bar R)=(K^{s_0}_k,R^{s_0}_k)\), we have
\(\Tt:=\Qq_{s_0,s_1}\resk{k}H^{s_0}_k\in\B_k(\da_k(s_0);H^{s_0}_k)\) by
Lemma~\ref{s:restriction}, the first three parts of \eqref{s:B3} are
trivial, and the last one holds, as
the last coordinate of a node which is put into \(B_{k+1}\) at a stage
\(\xi\ge s_0\) is the generic real of the stage \(\xi\), which does not
belong to \(M_\xi\supseteq M_{s_0}\) by Lemma~\ref{lem:interface}\,\ref{T:fresh}.

\paragraph{C. From level \(k-1\) to level \(k\) on products.}
Let \(m\ge k+1\) and let \(\Qq^1,\Qq^2\in\B_m(a;H)\) be disjoint. We use
the notation of the paragraph on products of two forcings with \(m\) in
place of \(k\). For \(j<m\), a family of lower assignments is \emph{valid
through level \(j\)} if \ref{V:complete}--\ref{V:finite} hold for
\(3\le\ell\le j\). The proofs of (A1) and (A2) only use
\eqref{s:intersection} and Definition~\ref{def:mallowable}, so
\eqref{s:A1}--\eqref{s:A2} hold for \(\Qq^1,\Qq^2\) at all levels
\(3\le\ell\le m\), with \(m\) as the fixed level in \eqref{s:A1b}.
In particular \(F^u_k=B^u_{k+1}\) since \(k<m\).
Let \(H_*=(K_*,R_*,B^{(0,0)}_{k+1})=H^{1,0}_k=H^{2,0}_k\) be
the common initial level-\(k\) assignment. Assuming parts~A and~B at level
\(k\), we show
\begin{equation}\label{s:C}
 \mathcal H\text{ valid through }k-1\ \Longrightarrow\
 \exists\,\mathcal H'\supseteq\mathcal H\text{ valid through }k,
\end{equation}
where \(\mathcal H'\) adds level-\(k\) assignments \(H^v_k\),
\(v\in C^1_k\times C^2_k\), with \(H^{(\beta,0)}_k=H^{1,\beta}_k\) and
\(H^{(0,\gamma)}_k=H^{2,\gamma}_k\). Fix \(\mathcal H\) valid through
\(k-1\). As \((\delta_1,0)\wedge(\beta,\gamma)=(\beta,0)\) and
\((\delta_1,0)\wedge(0,\delta_2)=(0,0)\), \eqref{s:intersection} gives
\begin{equation}\label{s:C0}
 \R^{M[G_1]}\cap\R^{M_{(\beta,\gamma)}}=\R^{M_{(\beta,0)}},\qquad
 \R^{M[G_2]}\cap\R^{M_{(\beta,\gamma)}}=\R^{M_{(0,\gamma)}},\qquad
 \R^{M[G_1]}\cap\R^{M[G_2]}=\R^{M}.
\end{equation}
Every step of \(\Qq^i\) involves only reals and nodes of \(M[G_i]\). We
treat the steps of \(\Qq^1\); the steps of \(\Qq^2\) are treated
symmetrically.

\emph{(C1) Entries in the product.} Let \(r=(p,d,H_r)\in R^{1,\xi}_k\)
govern a step of \(\Qq^1\) at its stage \(\xi\) which concerns the input
\(x\), and let \(\zeta\in C^2_{k-1}\). Then
\begin{equation}\label{s:C1}
 \Pos_k\bigl(\da_k(\xi,0);H_r,\,x\in A^k_{p,1}\bigr)\ \longleftrightarrow\
 \Pos_k\bigl(\da_k(\xi,\zeta);H_r,\,x\in A^k_{p,1}\bigr),
\end{equation}
so the step satisfies \eqref{s:index} at \((\xi,\zeta)\).

Indeed, \(\xi\in C^1_{k-1}\) by Definition~\ref{def:mallowable}(ii). By
Lemma~\ref{s:restriction}, \(\Qq^i_{0,\beta}\resk{k}H_*\in\B_k(\da^i_k(0);H_*)\)
for \(\beta\in C^i_{k-1}\). These restrictions have the same lower data by
\eqref{s:restr1}, so \(\mathcal H\) is valid for them, and \eqref{s:A3} at
level \(k\) and \eqref{s:remaining} show that the part \(\Tt\) of
\(\Qq^2[\xi]\) up to \(\zeta\), restricted to \(H_*\), satisfies
\[
 \Tt\in\B_k\bigl(\da_k(\xi,0);H_*\bigr),\qquad\text{with final data }
 \da_k(\xi,\zeta).
\]
We have \(H_r\preceq H^{1,\xi}_k\), and \(H^{1,\xi}_k\) agrees on old pairs
with \(H^{1,0}_k=H_*\). The assignment \(H_*\) is complete, as it is
part of \(a\), so \(\dom K_*=\omega\times\R^M\). So the old part of \(H_r\),
\[
 H_r^\circ:=\bigl(K_r\res(\omega\times\R^M),\ \{r'\in R_r:\parm(r')\in M\},\
 F_r\cap B^{(0,0)}_{k+1}\bigr),
\]
is a level-\(k\) assignment with \(H_r^\circ\preceq H_*\), as the stored
assignments of the entries in \(R_*\) are \(\preceq H_*\). The other pairs and entries of \(H_r\)
have parameters in \(M[G_1]\setminus M\), hence outside \(M[G_2]\) by
\eqref{s:C0}, and
\(F_r\setminus B^{(0,0)}_{k+1}\subseteq B^{1,\xi}_{k+1}\setminus
B^{(0,0)}_{k+1}\) is disjoint from \(M[G_2]\) by \eqref{s:A1c}, also for
the reals of leaves if \(k=3\). By Lemma~\ref{s:restriction},
\(\Tt\in\B_k(\da_k(\xi,0);H_r^\circ)\), and as \(H_r\) imposes on the steps
of \(\Qq^2\) exactly the restrictions of \(H_r^\circ\),
\(\Tt\in\B_k(\da_k(\xi,0);H_r)\). Now
\eqref{s:change} with \(H=H'=H_r\) gives \eqref{s:C1}.

\emph{(C2) Steps respect the assignments at pairs.} Let
\(v=(\beta',\gamma')\in C^1_k\times C^2_k\) be a pair at which \(H^v_k\)
is defined and agrees on old pairs with \(H^{(\beta',0)}_k\). Let a step of
\(\Qq^1\) of level \(k\) at a stage \(\xi\ge\beta'\) be carried out at a pair
\((\xi,\gamma'')\) with \(\gamma''\ge\gamma'\), \(\gamma''\in C^2_{k-1}\).
Then the step respects
\begin{equation}\label{s:C2}
 K^v_k,\qquad R^v_k\ \text{together with \eqref{s:index}},\qquad
 B^{(\xi,\gamma'')}_{k+1}.
\end{equation}
Indeed, let the step be for \((p,d)\), and let
\(\xi^*:=\max(C^1_k\cap[0,\xi])\ge\beta'\). If \(\parm(p)\notin M_v\), then
\(p\notin\dom K^v_k\), and \(H^v_k\) imposes no restriction. If
\(\parm(p)\in M_v\), then \(\parm(p)\in M_{(\beta',0)}\) by
\eqref{s:C0}, and the agreement on old pairs of \(H^v_k\) with
\(H^{1,\beta'}_k\), and of \(H^{1,\xi^*}_k\) with \(H^{1,\beta'}_k\), give
\[
 K^v_k(p)=K^{1,\beta'}_k(p)=K^{1,\xi^*}_k(p),
\]
and \(p\) has the same entry, or no entry, in \(R^v_k\), \(R^{1,\beta'}_k\)
and \(R^{1,\xi^*}_k\).
In \(\Qq^1\) alone the step respects \(K^{1,\xi^*}_k(p)\) and the entry
for \(p\), and \eqref{s:index} holds at \((\xi,\gamma'')\) by (C1). The
condition on parents and, for \(k=3\), the condition on leaves hold by
\eqref{s:A1b}.

For \(v=(\beta',\gamma')\le w=(\beta_w,\gamma_w)\) in
\(C^1_k\times C^2_k\), let \(\Tt_{v,w}\) be the part of
the product from \(v\) to \(w\) along either of the two coordinate paths,
restricted to level \(k\), with the lower assignments of \(\mathcal H\) and
the fixed assignment \(H^v_k\). Its steps of \(\Qq^1\) are carried out at
pairs \((\xi,\gamma'')\) with \(\gamma''\in\{\gamma',\gamma_w\}\subseteq
C^2_k\), and symmetrically for \(\Qq^2\). So (C2), the validity of
\(\mathcal H\) through level \(k-1\), \eqref{s:A1b}, \eqref{s:A2},
Lemma~\ref{lem:interface}\,\ref{T:prod} and~\ref{T:seg}
give
\begin{equation}\label{s:C2a}
 \text{\ref{V:steps} at level \(k\)},\qquad
 \Tt_{v,w}\in\B_k\bigl(\da_k(v);H^v_k\bigr)\ \text{with final data }
 \da_k(w).
\end{equation}
For \ref{V:steps}, take \(v:=(\xi^*,\zeta)\), whose assignment is in effect at
\((\xi,\zeta)\) by~\ref{N:assign}.

\emph{(C3) Inherited data at a new pair.} We choose \(H^w_k\) by recursion
along the componentwise order. Let \(w=(\beta,\gamma)\in C^1_k\times C^2_k\)
with \(\beta,\gamma>0\), let
\(S_w:=\{v\in C^1_k\times C^2_k:v<w\}\), and assume that \(H^v_k\) has been
chosen for \(v\in S_w\) such that \(H^v_k\) agrees on old pairs with
\(H^{v'}_k\) for \(v'\le v\) in \(S_w\). Let
\[
 \bar K_w:=\bigcup_{v\in S_w}K^v_k,\qquad
 \bar R_w:=\bigcup_{v\in S_w}R^v_k .
\]
Then \((\bar K_w,\bar R_w,B^w_{k+1})\) is a level-\(k\) assignment over
\((M_w,J_w)\), and for \(v\in S_w\)
\begin{equation}\label{s:C3}
 K^v_k\subseteq\bar K_w,\qquad
 \{r\in\bar R_w:\parm(r)\in M_v\}=R^v_k .
\end{equation}
Indeed, for \(v,v'\in S_w\) also \(v\wedge v'\in S_w\), and
\(\R^{M_v}\cap\R^{M_{v'}}=\R^{M_{v\wedge v'}}\) by \eqref{s:intersection}.
If \(p\in\dom K^v_k\cap\dom K^{v'}_k\), then \(\parm(p)\in M_{v\wedge v'}\),
so \(K^v_k(p)=K^{v\wedge v'}_k(p)=K^{v'}_k(p)\). If \(r\in R^{v'}_k\) and
\(\parm(r)\in M_v\), then \(r\in R^{v\wedge v'}_k\subseteq R^v_k\). Each
\(R^v_k\) is closed under predecessors, so a descending chain which starts
in \(R^v_k\) stays in \(R^v_k\), and as \(S_w\) is countable,
\(\htR(\bar R_w)<\omega_2\). Finally \(J_v\subseteq J_w\), so
\(H^v_k\in N_{J_v}\subseteq N_{J_w}\) for \(v\in S_w\). Fix an injection
\(e\in N_{J_w}\) of \(H(\omega_2)^{N_{J_w}}\) into the ordinals. Then
\((e(H^v_k))_{v\in S_w}\) is a countable sequence of ordinals in \(M_w\),
hence an element of \(N_{J_w}\) by Lemma~\ref{lem:interface}\,\ref{T:loc}, and so
are \((H^v_k)_{v\in S_w}\), \(\bar K_w\) and \(\bar R_w\).

\emph{(C4) Completion at the new pair.} We apply part~B with \(s_1:=w\),
\((\bar K,\bar R):=(\bar K_w,\bar R_w)\) and \(S:=S_w\). The pair \(w\) is
complete through level \(k-1\) by \ref{V:complete} for \(\mathcal H\), and
\(U^w_k(p)\) is finite by \eqref{s:A2}. For \(s_0:=v\in S_w\), the first
two parts of \eqref{s:B3} are \eqref{s:C3}, the third one holds by
\eqref{s:pdata}, and the last one by \eqref{s:A1}: a node
\(u\in B^w_{k+1}\setminus B^v_{k+1}\) is put into \(B_{k+1}\) by \(\Qq^i\)
at a stage \(\xi\) which is not below the \(i\)-th coordinate of \(v\), so
\(\lc(u)\notin M_{v_i(\xi,\zeta)}\supseteq M_v\), where \(\zeta\) is the
other coordinate of \(v\). For \(\Tt\) take \(\Tt_{v,w}\) from
\eqref{s:C2a}. By (B4), the assignment \(H^w_k:=H_k\) is complete at \(w\)
and agrees on old pairs with \(H^v_k\) for all \(v\in S_w\). This
continues the recursion and gives \ref{V:complete} at level \(k\).

\emph{(C5) New indices.} Let \(\xi\in C^1_k\setminus\{0\}\),
\(\zeta\in C^2_k\), let \(r=(p,d,H_r)\in R^{(\xi,\zeta)}_k\), and assume
that \(p\) has no entry in \(R^{(\xi',\zeta)}_k\) for
\(\xi'\in C^1_k\cap\xi\). Let \(y:=\parm(p)\). As \(0\in C^1_k\cap\xi\) and
\(H^{(\xi,\zeta)}_k\) agrees on old pairs with \(H^{(0,\zeta)}_k\), we have
\(y\notin M_{(0,\zeta)}\), so
\(y\notin M\cup M[G_2]\) by \eqref{s:C0}. As \(y\notin M\), no step for
\(p\) occurs before the initial stage, so \(U^{(0,0)}_k(p)=\varnothing\),
and \eqref{s:A2} and \eqref{s:pdata} give
\[
 U^{2,\zeta}_k(p)\subseteq U^{2,\delta_2}_k(p)=U^{(0,0)}_k(p)=\varnothing,
 \qquad U^{(\xi,\zeta)}_k(p)=U^{1,\xi}_k(p).
\]
If \(y\notin M[G_1]\), then \(U^{1,\xi}_k(p)=\varnothing\) in the same
way. If \(y\in M[G_1]\), then \(y\in M_{(\xi,0)}\) by \eqref{s:C0}, so
\(r\in R^{1,\xi}_k\) by the agreement with \(H^{(\xi,0)}_k=H^{1,\xi}_k\),
and \(p\) has no entry in \(R^{1,\xi'}_k\) for \(\xi'\in C^1_k\cap\xi\),
since such an entry would belong to \(R^{(\xi',\zeta)}_k\). So \(r\) is
new at the stage \(\xi\) of \(\Qq^1\), and \(d\notin U^{1,\xi}_k(p)\) by
the requirement on new entries after Definition~\ref{def:complete},
applied in \(\Qq^1\). In both cases
\(d\notin U^{(\xi,\zeta)}_k(p)\), which is \ref{V:fresh} at level \(k\).
\ref{V:finite} at level \(k\) is \eqref{s:A2}.

\emph{(C6) Names.} For \(v=(\beta,\gamma)\), the forcing
\(\Qq^1_\beta\times\Qq^2_\gamma\) decides whether
\(v\in C^1_k\times C^2_k\), and the data which are used in (C3) and (B1) at
\(v\) belong to \(M_v\). Enumerate the countable set of pairs below
\((\delta_1,\delta_2)\) in an order which extends the componentwise order.
Along this enumeration choose, by the maximum principle, a
\(\Qq^1_\beta\times\Qq^2_\gamma\)-name for \(H^v_k\) as in (C3) and (C4),
given the earlier names; on the axes take the names of \(\Qq^1\) and
\(\Qq^2\) alone. A \(\Qq^1_\beta\times\Qq^2_\gamma\)-name is also a
\(\Qq^1_{\beta'}\times\Qq^2_{\gamma'}\)-name for \(\beta\le\beta'\),
\(\gamma\le\gamma'\). So the restriction of \(\mathcal H'\) to
\(\Qq^1_\beta\times\Qq^2_\gamma\) is the family constructed for this
product, in either presentation. Together with (C4), (C5) and
\eqref{s:C2a}, this proves \eqref{s:C}.

\emph{Conclusion.} For \(k=3\), part~A with \(\mathcal H=\varnothing\)
proves (i) and~(ii) at level three, and parts~B and~C at level three give
the families valid through level three. If the families valid through
level \(k-1\) are given by part~C at level \(k-1\), then part~A proves (i)
and (ii) at level \(k\), and parts~B and~C at level \(k\) give the families
valid through level \(k\). This proves Proposition~\ref{s:product} by
induction on \(k\). The construction of part~C applies in the same way to
initial segments of forcings of higher level, using
Lemma~\ref{s:restriction} at each level \(k\), and each higher level
extends the family already chosen.

\end{proof}

\subsection{The iteration}\label{sec:iteration}

The iteration \(\Pp_{\omega_1}\) is built in \(\omega_1\) steps, each of
which adds a segment of countable length; as at the end of
Section~\ref{sec:long}, its initial segments are allowable over
\((W,\varnothing)\) with covers \(J_\beta\). The stages at which the
segments begin, and their limits, are the \emph{boundary stages}. We use
the notation \ref{N:models}--\ref{N:newpar} for the stages of
\(\Pp_{\omega_1}\). Along the iteration we keep closed sets of stages
\[
 \omega_1=C_2\supseteq C_3\supseteq C_4\supseteq\cdots,\qquad
 \sigma\in C_\ell\quad(\sigma\text{ a boundary stage},\ \ell\ge2),
\]
with complete level-\(\ell\) assignments at the stages in \(C_\ell\). So
all finite levels are complete at a boundary stage. At the initial
boundary stage \(0\), i.e.\ in \(W\) with \(X_0=\varnothing\),
\(B^0_k=\varnothing\) and \(J_0=\varnothing\), we apply (B1) and (B4) in
the proof of Proposition~\ref{s:product}, with \(S=\varnothing\), for
\(k=3,4,\ldots\) in this order; \(\Pos_k\) is evaluated with
\(\da_k(0)\), which only uses the levels below \(k\). As
\(U^0_k(p)=\varnothing\), all initial values are \(K^0_k(p)=\{0\}\). As in
(C3) in the proof of Proposition~\ref{s:product}, the sequence
\((H^0_k)_{k\ge3}\) belongs to \(N_{J_0}\).

\begin{lem}\label{s:extend}
Let \(\sigma\) be a boundary stage and \(m\ge3\).
\begin{enumerate}[leftmargin=*,topsep=.4em,itemsep=.2em,label=\textup{(\roman*)}]
\item Every \(\Qq\in\B_m(\da_m(\sigma);H^\sigma_m)\) can be used as the
segment which starts at \(\sigma\), with its sets \(C_\ell\), \(\ell<m\),
and its assignments, and its final stage can be made a boundary stage.
\item A countable limit of boundary stages is a boundary stage.
\item For all stages \(\alpha\le\beta\) in \(C_{m-1}\) and every
level-\(m\) assignment \(H'\preceq H^\alpha_m\) over
\((M_\alpha,J_\alpha)\),
\begin{equation}\label{s:E}
 \Pp_{\alpha,\beta}\resk{m}H'\in\B_m\bigl(\da_m(\alpha);H'\bigr).
\end{equation}
In particular this holds for boundary stages \(\alpha\le\beta\) and
\(H'=H^\alpha_m\).
\end{enumerate}
\end{lem}

\begin{proof}
The restriction \(\Pp_{\alpha,\beta}\resk{m}H'\) is defined as in
Definition~\ref{def:restriction}, with the stages of \(\Pp_{\omega_1}\). We
prove (i)--(iii) by induction on the boundary stages: (i) at \(\sigma\)
uses (B4) and induction on the finite levels at the new endpoint, (iii) for
\(\beta\) inside or at the end of the segment starting at \(\sigma\) uses
(i) at \(\sigma\) and (iii) up to \(\sigma\), and (ii) at a limit
\(\lambda\) uses (iii) for all \(\beta<\lambda\).

\emph{(E1) Proof of (i).} Let \(\Qq\) have length \(\delta\) and put
\(\sigma':=\sigma+\delta\). Inside \([\sigma,\sigma']\) let
\begin{equation}\label{s:E1}
 C_\ell\cap[\sigma,\sigma']=\sigma+C^{\Qq}_\ell\ \ (\ell<m),\qquad
 C_\ell\cap[\sigma,\sigma']=\{\sigma,\sigma'\}\ \ (\ell\ge m),
\end{equation}
with the assignments of \(\Qq\) at the stages of \(\sigma+C^{\Qq}_\ell\),
where \(\sigma+C:=\{\sigma+\gamma:\gamma\in C\}\). These sets are closed and
nested. As \(\Qq\) puts nodes only into sets \(B_k\) with \(k\le m\), the
set \(B_{m+1}\) and the level-\(m\) assignment \(H^\sigma_m\) do not change
before \(\sigma'\), so the restrictions of \(\Qq\) at level \(m\) are those
of \(H^\sigma_m\). The lower assignments of \(\Qq\), and with them the
instances of \(\cert_{m-1}\) at \(\sigma'\), are kept, by~\ref{N:data}
and \hyperref[obs:completion-locality]{Locality of completions}. At
\(\sigma'\) we complete the levels \(j=m,m+1,\ldots\) in this order. For
\(j\ge m\), the segment with its data through level \(j\) satisfies
\begin{equation}\label{s:E1a}
 \Tt_j\in\B_j\bigl(\da_j(\sigma);H^\sigma_j\bigr),\qquad\text{with final
 data }\da_j(\sigma').
\end{equation}
For \(j=m\) this is \(\Qq\). For \(j>m\), the levels \(m,\ldots,j-1\) have
been completed at \(\sigma'\), so the sets \eqref{s:E1} are as required in
Definition~\ref{def:mallowable}(i); \(\Qq\) makes no step of level \(j\),
so \(H^\sigma_j\) imposes no restriction; and
\(U^{\sigma'}_j=U^\sigma_j\). Now (B4), with \(S=\{\sigma\}\),
\((\bar K,\bar R)=(K^\sigma_j,R^\sigma_j)\) and \(\Tt:=\Tt_j\), gives a
complete level-\(j\) assignment \(H^{\sigma'}_j\) which agrees on old pairs
with \(H^\sigma_j\); the last part of \eqref{s:B3} is void, as
\(B^{\sigma'}_{j+1}=B^\sigma_{j+1}\) for \(j\ge m\). As in (C3), the
sequence \((H^{\sigma'}_j)_{j\ge m}\) belongs to \(N_{J_{\sigma'}}\). So
\(\sigma'\) is a boundary stage. A segment of length zero is represented
by a trivial iterand.

\emph{(E2) Proof of (iii).} Let \(\alpha\le\beta\) in \(C_{m-1}\) and
\(H'\preceq H^\alpha_m\). By~\ref{N:assign} and as the data of old pairs
are kept along \(C_m\),
\begin{equation}\label{s:E2}
 H'\preceq H^\alpha_m\preceq H^\gamma_m\qquad(\alpha\le\gamma\le\beta),
\end{equation}
as in \eqref{s:restr2}. Let \([\sigma,\sigma']\) be a segment of level
\(m_0\), i.e.\ \(\Qq\in\B_{m_0}(\da_{m_0}(\sigma);H^\sigma_{m_0})\), and
let \([\alpha',\beta']:=[\sigma,\sigma']\cap[\alpha,\beta]\neq\varnothing\).
Then \(\alpha',\beta'\in C_{m-1}\), as boundary stages belong to all
\(C_\ell\). If \(\alpha'=\beta'\), the part \(\Pp_{\alpha',\beta'}\resk{m}H'\)
is the trivial forcing, which belongs to \(\B_m(\da_m(\alpha');H')\). Let
\(\alpha'<\beta'\), so \(\sigma\le\alpha'<\sigma'\). Then
\(\Pp_{\alpha',\beta'}\resk{m}H'\) satisfies \eqref{s:E}, with \(\alpha'\)
in place of \(\alpha\):
\begin{enumerate}[leftmargin=*,topsep=.4em,itemsep=.2em,label=\textup{(\alph*)}]
\item if \(m_0>m\), by Lemma~\ref{s:restriction}, applied to \(\Qq\), as
\(H'\preceq H^{\alpha'}_m\) by \eqref{s:E2};
\item if \(m_0=m\), by \eqref{s:remaining} and the second part of
Lemma~\ref{s:restriction}, as \(H'\preceq H^{\alpha'}_m=H^\sigma_m\) by
\eqref{s:E2} and since \(C_m\cap[\sigma,\sigma')=\{\sigma\}\);
\item if \(m_0<m\), then \([\alpha',\beta']=[\sigma,\sigma']\), as
\(C_{m-1}\cap[\sigma,\sigma']=\{\sigma,\sigma'\}\) by \eqref{s:E1} with
\(m_0\le m-1\) in place of \(m\), and
\eqref{s:E1a} with \(m_0\) in place of \(m\) and \(j:=m\) gives
\(\Tt_m\in\B_m(\da_m(\sigma);H^\sigma_m)\), hence
\(\Tt_m\in\B_m(\da_m(\sigma);H')\) by Lemma~\ref{s:restriction}.
\end{enumerate}
These parts have the same fixed assignment \(H'\), so their concatenation
is \(m\)-allowable, see the remarks after Definition~\ref{def:mallowable};
at the limit boundary stages in \((\alpha,\beta]\), completeness and the
finiteness of the sets \(U_\ell(p)\), \(\ell\le m\), hold by (ii). This
proves \eqref{s:E}.

\emph{(E3) Proof of (ii).} Let \(\lambda\) be a countable limit of
boundary stages. We put \(\lambda\) into every \(C_\ell\); these sets stay
closed, as they contain all boundary stages below \(\lambda\). Let
\(X_\lambda:=\bigcup_{\gamma<\lambda}X_\gamma\),
\(B^\lambda_k:=\bigcup_{\gamma<\lambda}B^\gamma_k\), and, by
Lemma~\ref{lem:interface}\,\ref{T:closure}, let the cover \(J_\lambda\) be a countable ordinal
with \(J_\gamma\subseteq J_\lambda\) for \(\gamma<\lambda\). For
\(j\ge3\) and \(p=(e,y)\),
\begin{equation}\label{s:E3}
 \begin{gathered}
 U^\lambda_j(p)\subseteq K^\sigma_j(p)\ \ \text{for the least boundary stage
 }\sigma<\lambda\text{ with }y\in M_\sigma,\\
 U^\lambda_j(p)=\varnothing\ \ \text{if there is no such }\sigma.
 \end{gathered}
\end{equation}
Indeed, if a step for \(p\) occurs at a stage \(\xi<\lambda\), then
\(y\in M_\xi\subseteq M_{\sigma''}\) for the next boundary stage
\(\sigma''<\lambda\), so \(\sigma\) exists. An index used before \(\sigma\),
which is possible if \(y\) first appears inside a segment, belongs to
\(U^\sigma_j(p)\subseteq K^\sigma_j(p)\), by completeness at \(\sigma\). For
a step of level \(j\) for \(p\) at a stage \(\xi\) with
\(\sigma\le\xi<\lambda\), let \(\beta:=\xi+1\). Then \(\beta\in C_{j-1}\) by
Definition~\ref{def:mallowable}(ii) and \eqref{s:E1}, and \(\beta<\lambda\),
so (iii) for \([\sigma,\beta]\) and \(H^\sigma_j\) shows that the step
respects \(K^\sigma_j(p)\). This holds for all \(j\), also if the levels of
the segments below \(\lambda\) are unbounded in \(\omega\).

We complete the levels \(j=3,4,\ldots\) at \(\lambda\) in this order. Let
\(S:=C_j\cap\lambda\) and
\[
 \bar K:=\bigcup_{\alpha\in S}K^\alpha_j,\qquad
 \bar R:=\bigcup_{\alpha\in S}R^\alpha_j ,
\]
which is coherent, as the data of old pairs are kept along \(C_j\). For
\(\alpha\in S\), let \(\Tt:=\Pp_{\alpha,\lambda}\resk{j}H^\alpha_j\). By
(iii) for \(\beta<\lambda\), the remarks after
Definition~\ref{def:mallowable}, the completeness at \(\lambda\) of the
levels below \(j\) and \eqref{s:E3},
\(\Tt\in\B_j(\da_j(\alpha);H^\alpha_j)\), with final data
\(\da_j(\lambda)\). The first two parts of \eqref{s:B3} hold since the
data of old pairs are kept along \(C_j\), the third one is trivial, and
the last one holds along a single forcing, see (B4). By (B4), we obtain a
complete level-\(j\) assignment \(H^\lambda_j\) which agrees on old pairs
with \(H^\alpha_j\) for all \(\alpha\in S\); a pair whose parameter first
appears at \(\lambda\) is treated by (B1), with \(U^\lambda_j(p)=\varnothing\).
As in (C3), \(\bar K\), \(\bar R\) and the sequence
\((H^\lambda_j)_{j\ge3}\) belong to \(N_{J_\lambda}\). So \(\lambda\) is a
boundary stage, and (iii) extends to \(\beta=\lambda\) by (E2).
\end{proof}

The iteration has \(\omega_1\) steps. We schedule the following tasks for
all finite levels, codes of formulas and tuples of reals in intermediate
models.
\begin{enumerate}[leftmargin=*,topsep=.4em,itemsep=.2em,label=\textup{(\arabic*)}]
\item At a boundary stage with current \(H_m\), for an existential
assertion \(A\) of level \(m\ge3\) in normal form, decide whether
\(\Pos_m(a;H_m,A)\) holds. If it does, use a witness as the next segment,
by Lemma~\ref{s:extend}(i). The assertions \(I^m_p\) are included.
\item For an existing leaf \(v\notin B_4\) whose real \(w(v)\) has not been
coded yet, code \(w(v)\). For an existing node \(v\) of level \(k\ge5\)
with \(v\notin B_k\) and without a child in \(B_{k-1}\), put a child
\(v\cat u\) into \(B_{k-1}\), where \(u\) is the new real added by a coding
of \(w_{\rm aux}\). These steps concern reals of leaves and nodes which are
not roots, so the assignments do not restrict them. The only possible
restriction on parents has been checked.
\item For an entry \(r=(p,d,H_r)\) at level \(m\) and an input \(x\), carry
out the step for \((p,d)\) and \(x\) with the index \(i\) given by
\eqref{s:index}: for \(m=3\), code \(w_3(p,d,x,i)\), and for \(m\ge4\), put
the root \((m,p,d,x,i,t)\) into \(B_m\), with a new real \(t\). Such a root
has no parent. The reals \(w_3(p,d,x,i)\) are different from all reals of
leaves. The index \(d\) belongs to \(K(p)\).
\end{enumerate}

A single nontrivial step of level \(k\) in (2) or (3), carried out at a
boundary stage \(\sigma\), is a forcing \(\mathbb S\) of length one with
\(C^{\mathbb S}_\ell:=\{0,1\}\) for \(2\le\ell<k\). It respects the
conditions on parents and leaves by the choice in (2), and
\(K^\sigma_k(p)\) and \eqref{s:index} by the choice of \(d\) and \(i\) in
(3). At its final stage we provide the lower data by (B4), with
\(S=\{\sigma\}\), for \(\ell=3,\ldots,k-1\) in this order: restricted to
level \(\ell<k\), \(\mathbb S\) has the same iterand and counts as a coding
of \(w_{\rm aux}\), which is unrestricted, since every step which puts a
node into some \(B_k\) codes \(w_{\rm aux}\); and no set \(U_\ell(p)\),
\(\ell<k\), grows. So \(\mathbb S\in\B_k(\da_k(\sigma);H^\sigma_k)\), and
Lemma~\ref{s:extend}(i) applies to \(\mathbb S\).

At each boundary stage we choose a name for an enumeration of its reals in
order type \(\omega_1\). Each earlier boundary stage, tuple of indices,
kind of task and code of a formula is scheduled cofinally often after that
stage. Such a bookkeeping exists, as \(\omega_1\) can be partitioned into
\(\aleph_1\) unbounded sets. In (3), the pair and the input are specified
by \(m,e,y,x\), and \(d\) and \(H_r\) are given by the current assignment.
Names, and the conditions below which we restrict, can be chosen by the
forcing theorem and the axiom of choice. We bound the length of each
segment which is given by a name by a countable ordinal at its initial
stage, extend it by trivial iterands at the end if necessary, and use at
least one stage for each step of the iteration. Trivial iterands at the
end of a boundary stage repeat its data. The finite support iteration
\(\Pp_{\omega_1}\) has length \(\omega_1\). Every initial segment of
\(\Pp_{\omega_1}\) below \(\omega_1\) has countable length, and every real
of the final model belongs to the model at some boundary stage, and the
iteration preserves all cardinals and \(\CH\), by
Lemma~\ref{lem:interface}\,\ref{T:size}. Presented as a name over \(\Pp^0\times\Pp^1\), it yields a set
forcing over \(L\).

\subsection{The separator in the final extension}

Let \(V_{\mathrm{sep}}:=W[G_{\omega_1}]\), where \(G_{\omega_1}\) is
\(\Pp_{\omega_1}\)-generic over \(W\), and write \(V:=V_{\mathrm{sep}}\).

\begin{lem}\label{s:truth}
For every \(m\ge3\), every boundary stage \(s\), and every \(\Sigma^1_m\)
assertion \(A\) in normal form with parameters from the model at \(s\),
\begin{equation}
               V\models A\quad\longleftrightarrow\quad
               \Pos_m(a;H_m,A), \label{s:trutheq}
\end{equation}
where \(a\) and \(H_m\) are the data at \(s\). Consequently,
\(V\models\theta\) if and only if \(\cert_m(s,\theta)\), for every
\(\Pi^1_m\) assertion \(\theta\) in normal form with parameters from the
model at \(s\).
\end{lem}
\begin{proof}
Between two boundary stages, we use Lemma~\ref{s:extend}(iii) to regard the
part of the iteration in between as an \(m\)-allowable forcing with
respect to the assignment \(H_m\) at the first of them. The assignment
\(H_m\) at the second stage adds values of \(K\) and entries of \(R\) only
for new parameters, and it adds only forbidden parents which contain new
reals. So, by Proposition~\ref{s:product}, \(\Pos_m\) is constant for old
parameters, even over intervals which contain steps of unbounded finite
level.

Write \(A=\exists b\,P(b)\). If \(A\) holds in \(V\), choose a later
boundary stage whose model contains a witness \(b\). For \(m=3\), its
\(\Pi^1_2\) matrix is already true there, by Shoenfield absoluteness. For
\(m>3\), each conjunct of the matrix satisfies \(\cert_{m-1}\) there, by
the induction hypothesis. So the trivial forcing witnesses \(\Pos_m\) at
this stage, and by constancy \(\Pos_m\) also holds at \(s\). Conversely,
if \(\Pos_m\) holds at \(s\), then it still holds when task (1) for this
assertion and this tuple is carried out. The chosen \(m\)-allowable
forcing gives a particular \(b\) whose matrix satisfies \(\cert_{m-1}\),
and its final stage is a boundary stage with the same data through level
\(m-1\). Shoenfield absoluteness for \(m=3\), and the induction hypothesis
for \(m>3\), give \(V\models P(b)\). The statement about negations follows
from \eqref{s:cert}. This is an ordinary induction on the finite levels of
formulas, carried out in the metatheory; no truth predicate for all levels
is needed. The inductions on \(k\) in the proofs of
Proposition~\ref{s:product} and Lemma~\ref{s:extend} are inductions inside
\(\ZFC\), see the end of (D4).
\end{proof}

\begin{lem}\label{s:nodes}
Let \(B^\infty_k\) be the union of the sets \(B_k\) along the iteration.
For every \(k\ge4\) and every node \(v\) of level \(k\) in \(V\),
\[
             V\models T_k(v)\quad\longleftrightarrow\quad
                    v\in B^\infty_k.
\]
\end{lem}
\begin{proof}
A node which exists at some stage and does not belong to \(B_k\) there is
never put into \(B_k\) later, since its last coordinate would then be the
new generic real of that later stage. The real of a leaf in \(B^\infty_4\)
had not been coded when the leaf was put into \(B_4\), and it is never
coded afterwards. For a leaf not in \(B^\infty_4\), task (2) codes its
real if necessary. So Lemma~\ref{lem:interface}\,\ref{T:exact} proves the assertion for
\(k=4\). A node in \(B^\infty_k\) has no child in \(B^\infty_{k-1}\),
neither before nor after it was put into \(B_k\). Every node not in
\(B^\infty_k\) eventually has a child in \(B^\infty_{k-1}\), by task (2)
after the node appears. Using \eqref{s:decoder}, this gives the induction
step.
\end{proof}

In task (1), the iteration carries out witnesses for the current complete
assignment. The assignment stored in an entry may allow further witnesses.
For a single membership statement with respect to such an assignment we
therefore use the implication from membership in the final model to
\(\Pos\), together with the fact that \(\Pos\) cannot hold for both indices
at the same input.

\begin{proof}[Proof of Theorem~\ref{main:sep}]
Fix \(j\), \(e\) and \(y\) in \(V\), and suppose that the two sets given by
\(p=(e,y)\) are disjoint. At a boundary stage whose model contains \(y\),
\(\Pos_j(a;H_j,I^j_p)\) fails, by Lemma~\ref{s:truth}. By completeness,
there is an entry \(r=(p,d,H_r)\). Let \(\sigma\) be the first stage in
\(C_j\) at which this entry appears; it may lie inside a segment which
was carried out earlier. Then \(d\notin U_j(p)\) at \(\sigma\), and
\(\Pos_j(a;H_r,I^j_p)\) fails at \(\sigma\). The entry and \(K(p)\) are
kept at every later stage in \(C_j\). By Lemma~\ref{s:extend}(iii), every
later part of the iteration which begins at \(\sigma\), or which lies
between later stages in \(C_{j-1}\), is \(j\)-allowable with respect to
the stored \(H_r\). So \(\Pos_j(a;H_r,I^j_p)\) fails at every such later
stage, and for an old \(x\), the truth value of
\(\Pos_j(a;H_r,x\in A^j_{p,i})\) is constant between these stages.

At any such stage whose model contains \(x\), the statement
\(\Pos_j(a;H_r,x\in A^j_{p,i})\) cannot hold for both \(i<2\). Otherwise,
move the two witnesses to disjoint sets of coding areas and use the
product of Proposition~\ref{s:product}. Both instances of \(\cert_{j-1}\)
for the matrices hold at the common final stage: for \(j=3\) by
Shoenfield absoluteness, and for \(j>3\) by \eqref{s:change} at level
\(j-1\), along the remaining forcing in each order, as in part A of the
proof of Proposition~\ref{s:product}. So \(x\) and the two witnesses show
that \(\Pos_j(a;H_r,I^j_p)\) holds at this stage. Since \(x\) is part of
the witness of \(I^j_p\), the only real parameter from the pair is \(y\).
So the contradiction transfers back to \(\sigma\), even if \(x\) was added
later.

If \(x\in A^j_{p,i}\) in \(V\), take a boundary stage after the
appearance of the entry whose model contains \(x\) and a witness to
\(x\in A^j_{p,i}\). By Lemma~\ref{s:truth}, and by Shoenfield absoluteness
for \(j=3\), the matrix satisfies \(\cert_{j-1}\) for this witness there.
The trivial forcing shows that \(\Pos_j(a;H_r,x\in A^j_{p,i})\) holds,
since \(\cert_{j-1}\) at its final stage depends only on the data of the
lower levels at that stage. Consequently, \eqref{s:index} selects the
correct index \(i\) for every element of \(A^j_{p,i}\), and, by the
invariance of the stored instance of \(\Pos_j\), it gives a single value
for each input, also outside both sets. This argument also covers an
earlier step for \((p,d)\) inside a segment, after the entry appeared: the
stage immediately before the step belongs to \(C_{j-1}\), so the part of
the iteration from there to the later boundary stage is \(j\)-allowable
with respect to \(H_r\).

At \(\sigma\) no step for \((p,d)\) had occurred. All later steps for
\((p,d)\) use the prescribed index, and task (3) carries out a step for
every input. Lemma~\ref{lem:interface}\,\ref{T:exact} at level three and
Lemma~\ref{s:nodes} at the higher levels now give
\[
       A^j_{p,i}\subseteq S_i^j(p,d,\cdot)\quad(i<2),\qquad
       S^j_0(p,d,\cdot)=\R\setminus S^j_1(p,d,\cdot).
\]
Both sets have the \(\Sigma^1_j(y)\) definitions \eqref{s:separators}. So each
of them is \(\Delta^1_j(y)\), defined from \(y\) and the natural numbers
\(e\) and \(d\). Although \(d\) depends on the pair, it occurs in each
definition as a numeral, so the only real parameter is \(y\). All these
definitions are evaluated in \(V_{\mathrm{sep}}\).
\end{proof}

\paragraph{The examples of Section~\ref{f:examples} in \(V_{\mathrm{sep}}\).}
We return to the pair of the toy problem in Section~\ref{f:toy}; recall
that the range of the
tuple coding \(\langle i,z,t,j\rangle\) used there avoids \(w_*\) and the
ranges of \(c\) and of the injections \(e_j\). Let \(w\) be a real of \(V\)
which is neither the real of a leaf nor of the form \(w_3(p,d,x,i)\). At
a boundary stage \(s\) whose model contains \(w\), the forcing
\(\Code(w,\eta)\) of length one, with a coding area \(\eta\) outside the
cover, belongs to \(\B_3(\da_3(s);H^s_3)\), since the coding of such a real
is unrestricted. It forces \(D_3(w)\) by Lemma~\ref{lem:interface}\,\ref{T:exact}. So
\(\Pos_3(\da_3(s);H^s_3,D_3(w))\) holds, and \(V\models D_3(w)\) by
Lemma~\ref{s:truth}. In particular every real \(\langle i,z,t,0\rangle\)
is coded in \(V\), and
\[
 V\models\neg P_i(z)\quad\text{for all }i,z,\qquad
 A_0=A_1=\varnothing .
\]
No attempt of the kind described in the toy problem survives in \(V\): a
real which is not coded in \(V\) is the real of a leaf in \(B^\infty_4\) or
of the form \(w_3(p,d,x,i)\).

A pair can, however, read the candidate separators themselves. Fix
\(d_0\in\omega\). By the recursion theorem \cite[Chapter~7]{mosch} there is a lightface
pair
\(p=(e,0^\omega)\) with
\[
 A^4_{p,0}=\{x:S^4_1(p,d_0,x)\},\qquad A^4_{p,1}=\{x:S^4_0(p,d_0,x)\},
\]
so that \(x\) belongs to \(A^4_{p,0}\) as soon as \(x\) is put on the side
\(1\) of the candidate with index \(d_0\). At the initial stage
\(K^0_4(p)=\{0\}\), and this value is kept.
\begin{enumerate}[leftmargin=*,topsep=.4em,itemsep=.2em,label=\textup{(\roman*)}]
\item If \(d_0\ne0\), no step for \((p,d_0)\) at level four is ever
allowed. No root \((4,p,d_0,x,i,t)\) is put into \(B_4\), and both sets
are empty in \(V\) by Lemma~\ref{s:nodes}.
\item Let \(d_0=0\), and let \(H=(K,R,F)\) be a level-four assignment at a
stage which is complete through level three, with \(0\in K(p)\) and
without an entry for \(p\). Let \(\Qq\) consist of two steps for \((p,0)\)
which put roots \(v_i=(4,p,0,x,i,t_i)\), \(i<2\), into \(B_4\) for the same
\(x\), with the lower data provided as for the single steps in
Section~\ref{sec:iteration}. Then \(\Qq\in\B_4(a;H)\): the index \(0\) is
allowed, no entry prescribes the side, and roots have no parents. At the
final stage \(v_0,v_1\in B_4\), so no forcing which is \(3\)-allowable with
respect to the complete level-three assignment there codes \(w(v_0)\) or
\(w(v_1)\). By Lemma~\ref{lem:interface}\,\ref{T:exact}, \(\cert_3\) holds there for
\(\neg D_3(w(v_0))\) and for \(\neg D_3(w(v_1))\), i.e.\ for the matrices of
\(x\in A^4_{p,1}\) and \(x\in A^4_{p,0}\). Hence
\[
 \Pos_4(a;H,I^4_p)\quad\text{for every such }H .
\]
In the recursion \eqref{s:B1} the pair \(p\) therefore never belongs to a
set \(E_\nu\) and receives no entry. By
Definition~\ref{def:complete}(ii) and Lemma~\ref{s:truth}, the two sets
have a common element in \(V\).
\end{enumerate}
This is the test announced in Section~\ref{f:toy}. A pair which can read the side
we are about to choose is dealt with before any side is chosen for it:
either the index it reads is never used, or a common point is possible
and is forced.

\section{Uniformization at a fixed projective level}\label{sec:up}

\paragraph{Two candidates at level four.}
Fix \(x\) and a \(\Pi^1_4\) relation \(A\) on \(\R^2\), and write
\(A_x=\{v:(x,v)\in A\}\). Suppose we provisionally choose \(y\in A_x\).
We want to leave \(h(x,y)\) uncoded at level four and to code \(h(x,z)\)
for a competing value \(z\), where \(h\) is a fixed recursive pairing
with recursive inverses. By \eqref{eq:common-decoder}, coding the
competitor requires a witness \(t\) such that \(e_3(h(x,z),t)\) stays
uncoded at level three. What happens if we allow ourselves to give up
this witness whenever the current choice drops out of \(A_x\)?

For a concrete example take \(y=0^\omega\), \(z=1^\omega\), and put
\(v_{2n}=y\), \(v_{2n+1}=z\). Fix a uniformly recursive sequence of
pairwise distinct reals \(t_{n,j}\), \(n,j<\omega\). Let
\[
 b_{n,j}(u)=e_3(h(u,v_n),t_{n,j}),\qquad
 B_n(u)\ \longleftrightarrow\ \forall j<\omega\;D_3(b_{n,j}(u)).
\]
Call block \(n\) complete at \(u\) when \(B_n(u)\) holds. Put
\(V_0=\{y,z\}\) and \(V_n=\{v_n\}\) for \(n>0\). Define \(A\) by
\[
 (u,v)\in A\quad\longleftrightarrow\quad
 \exists n<\omega\,
 \bigl[(\forall k<n\;B_k(u))\land\neg B_n(u)\land v\in V_n\bigr].
\]
Thus \(A_u=V_n\) for the first incomplete block \(n\).
When every block is complete, \(A_u\) is empty. Each \(B_n\) is
\(\Sigma^1_3\), since one real codes the sequence of witnesses for its
\(D_3\) assertions. Nonmembership in \(A\) says that every block is
complete or that \(v\notin V_n\) for the first incomplete block \(n\).
This description and the displayed membership test are both \(\Sigma^1_4\).
Thus this \(\Pi^1_4\) example already lies in \(\Delta^1_4\).

Start over \(W\) with \(x\in W\). Initially \(A_x=\{y,z\}\). To keep \(y\),
code the reals \(b_{0,j}(x)\) at level three, one at a time on fresh
coding areas, and reserve \(b_{1,0}(x)\) as a witness to the exclusion
of \(z\). After \(\omega\) steps block zero is complete and \(A_x=\{z\}\).
Switch to \(z\) and begin block one. The first real of this block,
\(b_{1,0}(x)\), is the very real we had reserved to exclude \(z\). Reserve \(b_{2,0}(x)\)
to exclude \(y\). Completing block one puts \(y\) back into \(A_x\) and
removes \(z\).
Continue in this way. At \(\omega^2\) every block is complete and
\(A_x=\varnothing\). Every provisional witness for excluding either
candidate has now been cancelled.

To evaluate the two \(D_4\) predicates, continue this toy iteration to
\(\omega_1\). Schedule the coding of \(e_3(h(x,v),t)\) for \(v=y,z\)
and every real \(t\) as it appears. All reservations above were
provisional, so each such real can be coded on a fresh area. In the
final model \eqref{eq:exact} and \eqref{eq:common-decoder} give
\[
 A_x=\varnothing,\qquad \neg D_4(h(x,y)),\qquad \neg D_4(h(x,z)).
\]
Reading the graph from the uncoded pairs alone would give both \(y\)
and \(z\) at \(x\). Intersecting with \(A\), as in \eqref{u:graph},
gives the empty section here. For a nonempty section we must also
ensure that exactly one of its members remains uncoded.

The switches in the example are cofinal in a countable limit. For a
general \(\Pi^1_4\) relation, write \((x,v)\notin A\) as
\(\exists a\,P_3(x,v,a)\) with \(P_3\) a \(\Pi^1_3\) matrix. If this
holds in the final model of our \(\omega_1\)-iteration, a witnessing
\(a\) appears at a bounded stage. Its \(\Pi^1_3\) matrix holds at every
subsequent stage, since a witness to its \(\Sigma^1_3\) negation would
persist by Shoenfield absoluteness.
During the construction we certify exclusion once we have a particular
\(a\) for which no allowable level-three continuation can force
\(\neg P_3(x,v,a)\). This certificate persists through later coding
by \eqref{u:Eprops}.

The first level of the thinning at which no allowable continuation can
produce such a certificate gives \((x,v)\) its rank. At the fixed point,
if the rank is undefined, we carry out a witnessing continuation and
thereby put \(v\) permanently outside \(A_x\). For a nonempty final
section we choose a value of least attainable rank. If no new real
can attain that minimum, the coherent choice function selects one of
the old values and the rank rule keeps its pair uncoded. Otherwise we
force a new real attaining the minimum, forbid coding its pair as soon
as it appears, and restart the thinning with that prohibition in effect.
Tasks~\ref{u:T1}--\ref{u:T3} thus choose once for each nonempty final
section and eventually code all competing values.

For \(n=3\), Theorem~\ref{main:up} follows by applying
\cite[Theorem~3.1 and Corollary~3.2]{up} to the lightface universal
relation, with the parameter as one of the coordinates of the input. Fix
\(n\ge4\). We work over \(W=L[G^0][G^1]\). The final iteration
\(\Pp_{\omega_1}\) has initial segments which are allowable over
\((W,\varnothing)\), as at the end of Section~\ref{sec:long}, and we use
the notation \ref{N:models}--\ref{N:par} for its stages: \(s\) denotes a
stage below \(\omega_1\), with model \(M_s=W[G_s]\), cover \(J_s\) and
\(N_s:=N_{J_s}\), and \(t\) denotes the new generic real which is added by
the coding forcing at a given stage. Write \(u=\pa{m,p,x}\) for an input of
\(A_n\). For each \(u\) with nonempty final section the forcing will code
\(h(u,y)\) at level \(n\) for every \(y\) except the selected value.
The pairs in \(A_n\) left uncoded form the \(\Pi^1_n\) graph
\eqref{u:graph}.

\subsection{Coding iterations over \texorpdfstring{\(W\)}{W}}\label{sec:uptransfer}

The rank and product arguments of \cite[\S\S2.6--3]{up} are formulated for
the coding iterations of \cite[\S2.4]{up}. There, every stage \(i\) first
forces with a block \(\Cc_{E_i}=\operatorname{Fn}_{<\omega_1}(E_i\times
\omega_1,2)^L\) of the Cohen product computed in \(L\), and then with an
almost disjoint coding forcing, and \cite[(1) in \S2.4]{up} presents an
iteration of length \(\lambda\le\omega_1\) as
\(\Cc_{E_{<\lambda}}*\dot{\mathbb R}_\lambda\), with all Cohen factors in
front. For \(\lambda=\omega_1\), the Cohen part is forcing equivalent to
\(\Pp^1\). So these iterations are, up to the presentation, iterations of
almost disjoint coding forcings over \(W\); the stage models of
\cite{up} contain only the Cohen coordinates added so far, while the models
\(M_s\) contain all of them, with the same reals by
Lemma~\ref{lem:interface}\,\ref{T:loc}. By \cite[Remark after Theorem~3.1]{up}, the arguments of
\cite[Lemmas~2.12--2.16 and Theorem~3.1]{up} use only the
properties of coding iterations in \cite[Lemmas~2.9--2.11]{up} and two
absoluteness facts. Lemma~\ref{lem:interface} provides these properties
for allowable forcings over suitable pairs: its items \ref{T:seg},
\ref{T:lift} and \ref{T:prod} correspond to
\cite[Lemma~2.9\,(1)--(3)]{up}, item \ref{T:size} to
\cite[Lemma~2.10]{up}, and item \ref{T:exact} to \cite[Lemma~2.11]{up}.
Moreover, \eqref{u:Eprops} below replaces the absoluteness facts. The model which is called \(W\) in
\cite{up} is our \(L[G^0]\). In Section~\ref{sec:ranks} the rank and
product lemmas are proved in full in the present notation; the proofs
follow \cite[Lemmas~2.12--2.16]{up}.

\emph{Choice functions.} At each stage \(s\) we use a choice function
\[
 \begin{gathered}
 \chi_s:\bigl(\mathcal P(\R)\setminus\{\varnothing\}\bigr)^{N_s}
   \longrightarrow\R^{M_s},\qquad \chi_s\in N_s,\qquad \chi_s(B)\in B,\\
 s\le s'\ \Longrightarrow\ \chi_{s'}\res\dom\chi_s=\chi_s ,
 \end{gathered}
\]
and likewise at the intermediate models of the forcings of countable
length and of their products. Coherent choice functions on a product grid, extending prescribed
ones, are provided by \ref{CF:ext} and \ref{CF:ops} below; the proof
follows \cite[Lemma~2.12]{up}. As \(\chi_s\in N_s\), the choice
functions are invariant under the translations of
Lemma~\ref{lem:interface}\,\ref{T:trans}, and so are all the definitions below.

On a product they are obtained as follows. Let \(\Qq\) and \(\mathbb R\)
be as in Lemma~\ref{lem:interface}, with \(K\cap K'=\varnothing\), and let
\(G_\Qq\times G_{\mathbb R}\) be generic over \(M\). For pairs
\(a=(i,j)\le(\delta,\delta')\), ordered coordinatewise, put
\[
 M_a:=M[G_\Qq\res i][G_{\mathbb R}\res j],\quad
 J_a:=J^{\Qq}_i\cup J^{\mathbb R}_j,\quad N_a:=N_{J_a},\quad
 D_a:=\bigl(\mathcal P(\R)\setminus\{\varnothing\}\bigr)^{N_a}.
\]
\begin{enumerate}[leftmargin=*,topsep=.4em,itemsep=.2em,label=\textup{(CF\arabic*)}]
\item\label{CF:meet} \(D_a\cap D_b=D_{a\wedge b}\) for all pairs \(a,b\).
\item\label{CF:ext} Let \(\mathcal I\) be a set of pairs which is closed
downwards, and for \(c=(i,j)\in\mathcal I\) let \(\chi^0_c\in N_c\) be a
choice function on \(D_c\) which is given by a
\(\Qq_i\times\mathbb R_j\)-name, with \(\chi^0_c\res D_{c'}=\chi^0_{c'}\)
for \(c'\le c\). Then there are choice functions
\[
 \chi_a\colon D_a\to\R^{M_a},\quad \chi_a\in N_a,\quad \chi_a(B)\in B,
 \qquad a'\le a\ \Longrightarrow\ \chi_a\res D_{a'}=\chi_{a'},
\]
which are given by names in the same way, with \(\chi_c=\chi^0_c\) for
\(c\in\mathcal I\).
\item\label{CF:ops} The functions of \ref{CF:ext} can be chosen for the
restrictions of \(\Qq\times\mathbb R\) below conditions and for its
quotients by \(\Qq_{i_0}\times\mathbb R_{j_0}\), and they can be extended,
keeping all functions already chosen, when a continuation with its own
coherent choice functions is appended to \(\Qq_{i_0}\) as in
Lemma~\ref{lem:interface}\ref{T:prod}.
\end{enumerate}
The axes \(\{(i,0):i\le\delta\}\cup\{(0,j):j\le\delta'\}\) are closed
downwards, so coherent choice functions along \(\Qq\) and along
\(\mathbb R\) which agree at \((0,0)\) extend to the product.

For \ref{CF:meet}, if \(a'\le a\) then \(N_{a'}\subseteq N_a\) and
\(D_{a'}\subseteq D_a\), and for incomparable pairs the assertion is
Lemma~\ref{lem:interface}\ref{T:stage}. For \ref{CF:ext} let, for
\(B\in\bigcup_aD_a\),
\[
 \ell(B):=\bigl(\min\{i:\exists j\ B\in D_{(i,j)}\},\
               \min\{j:\exists i\ B\in D_{(i,j)}\}\bigr).
\]
Taking the meet of two pairs which witness the two minima,
\ref{CF:meet} gives
\phantomsection
\begin{equation*}\tag{\(\ast\)}\label{u:ell}
 B\in D_{\ell(B)},\qquad B\in D_a\ \Longleftrightarrow\ \ell(B)\le a .
\end{equation*}
For every pair \(a=(i,j)\) fix a \(\Qq_i\times\mathbb R_j\)-name for a
well-order \(\lhd_a\in N_a\) of \(\R^{N_a}\), and for \(B\in D_a\) let
\[
 \chi_a(B):=
 \begin{cases}
   \chi^0_{\ell(B)}(B), & \ell(B)\in\mathcal I,\\[2pt]
   \min\nolimits_{\lhd_{\ell(B)}}B, & \text{otherwise.}
 \end{cases}
\]
By \eqref{u:ell} the value does not depend on \(a\), so
\(\chi_a\res D_{a'}=\chi_{a'}\) for \(a'\le a\). If \(c\in\mathcal I\) and
\(B\in D_c\), then \(\ell(B)\le c\), hence \(\ell(B)\in\mathcal I\) and
\(\chi_c(B)=\chi^0_{\ell(B)}(B)=\chi^0_c(B)\). The sequences
\((D_c)_{c\le a}\), \((\lhd_c)_{c\le a}\) and
\((\chi^0_c)_{c\le a,\,c\in\mathcal I}\) are countable sequences of elements
of \(N_a\), and they belong to \(M_a\), as they are given by names over
initial segments of \(\Qq_i\times\mathbb R_j\). As in (C3) in the proof of
Proposition~\ref{s:product}, with an injection in \(N_a\) of a set which
contains these elements into the ordinals, the three sequences belong to
\(N_a\), and so does \(\chi_a\), which is defined from them. The
definition of \(\chi_a\) is uniform in the generic filter, which gives the
names.

For \ref{CF:ops}, a restriction below \((q,r)\) is the product of
\(\Qq\res q\) and \(\mathbb R\res r\), with the same models. For
\(s=(i_0,j_0)\), the functions \(\chi_a\), \(a\ge s\), are coherent choice
functions on the grid of the quotient. If a continuation \(\dot\Tt\) is
appended to \(\Qq_{i_0}\), with coherent choice functions along \(\dot\Tt\)
which extend \(\chi_{(i_0,0)}\), apply \ref{CF:ext} to
\((\Qq_{i_0}*\dot\Tt)\times\mathbb R\), where \(\mathcal I\) is the set of
all pairs \((i,j)\) with \(i\le i_0\) or \(j=0\), with the functions
already chosen and those of the continuation. This set is closed
downwards, and the prescribed functions are coherent.

Along a single forcing, and along the final iteration, the same
definition with one coordinate gives coherent choice functions: a set
\(B\) is given the value chosen at the first stage \(s\) with
\(B\in(\mathcal P(\R))^{N_s}\).

\subsection{Codings at the lower levels}\label{sec:levelcodings}

At each stage \(s\) we keep countable sets \(Z_j(s)\), \(3\leq j<n\), of
reals which may not be coded at level \(j\); initially they are empty. In
Section~\ref{sec:sep}, the corresponding information was given by the set
\(B_{j+1}\). The elementary steps are the following.

\begin{enumerate}[leftmargin=*,topsep=.4em,itemsep=.2em,label=\textup{(\roman*)}]
\item Coding a real \(d\notin Z_3\) at level three is the coding of \(d\)
by \(\Code(d,\eta)\), with a coding area \(\eta\) outside the current
cover \(J_s\).
\item Let \(4\leq j\leq n\), and let the real \(d\) be chosen before the
stage. If \(j<n\), we require \(d\notin Z_j\). If \(j=n\), we require that
\(d=h(u,y)\) for some \(u\) and \(y\), and the restrictions for
uniformization below apply. We code \(w_*\) with a coding area outside
the current cover, let
\(t\) be the generic real added by this coding forcing, say that \(d\) is
coded at level \(j\), and put \(e_{j-1}(d,t)\) into \(Z_{j-1}\).
\end{enumerate}
Every coding of the auxiliary real \(w_*\) counts as a coding of \(w_*\)
at level three. A stage at which a real is coded at a level above three
codes exactly this one real at this level, in addition to the coding of
\(w_*\). At limit stages we take unions. A real which is newly put into
some \(Z_i\) has not been coded at level \(i\) before, since it computes
the new real \(t\), whereas the reals coded earlier already existed. The
real coded at this stage is \(w_*\), which does not belong to the range of
\(e_3\). The new prohibition applies to all later stages. This shows that,
at each level, the set of reals coded at this level and the set of reals
forbidden at this level are disjoint.

The predicates \(D_j\) of \eqref{eq:common-decoder} read off these
codings.

\begin{lem}\label{u:decode}
Suppose that for each \(3\leq j<n\) and every real \(d\) of the final
model, coding \(d\) at level \(j\) is requested cofinally often and is
carried out whenever it is allowed. Then, for \(3\leq j\leq n\),
\(D_j(d)\) holds in the final model if and only if \(d\) has been coded at
level \(j\) at some stage.
\end{lem}
\begin{proof}
As every real appears at a stage below \(\omega_1\), the requests ensure
that the reals coded at level \(j<n\) are exactly the reals of the final
model outside \(Z_j^\infty:=\bigcup_s Z_j(s)\). For \(j=3\), the assertion
is \eqref{eq:exact}, by Lemma~\ref{lem:interface}\ref{T:exact}. If it holds for some \(j<n\), then
\[
 D_{j+1}(d)\ \longleftrightarrow\
 \exists t\ \bigl(e_j(d,t)\in Z_j^\infty\bigr)
 \ \longleftrightarrow\ d\text{ has been coded at level }j+1.
\]
The last equivalence uses the injectivity of \(e_j\) and the fact that
reals are put into \(Z_j\) only by the rule above. This step also applies
for \(j+1=n\), where the requests are made only at the preceding level.
\end{proof}

\subsection{Certificates for projective truth}

We use recursive normal forms down to level two. Starting with
\(P_n=A_n\), we fix matrices such that, for \(3\leq k\leq n\),
\begin{equation}\label{u:normal}
       \neg P_k(\bar v)\ \longleftrightarrow\
                 \exists z\,P_{k-1}(\bar v,z),\qquad
                 P_{k-1}\in\Pi^1_{k-1}.
\end{equation}
The tuple \(\bar v\) includes the witnesses for the outer quantifiers. We
need the certificates below only up to \(k=n-1\). The same definitions
apply to any fixed normal form with its real parameters.

A forcing of countable length starting at a stage \(s\), which is
allowable over \((M_s,J_s)\) and consists of the steps of
Section~\ref{sec:levelcodings}, is called \emph{unrestricted} if it is
considered without the restrictions for uniformization which are
introduced below. For \(3\leq k<n\), let \(\B_k(s)\) be the class of
unrestricted forcings of countable length, starting at \(s\), which code
reals only at levels at most \(k\). So \(Z_3,\ldots,Z_{k-1}\) may grow,
while \(Z_k\) stays fixed. The sets \(Z_i\) with \(i>k\) do not affect
which steps are allowed in \(\B_k(s)\). Every forcing in \(\B_k(s)\) may
be carried out in every phase of the final iteration, since all
restrictions for uniformization introduced below concern only codings at
level \(n\).

Let \(M_s\) be the model at the stage \(s\). We write
\(\Pos_{\mathcal D}(s,\theta)\) if some condition of some forcing in
\(\mathcal D(s)\) forces \(\theta\) at its final stage \(s'\); the trivial
forcing is included. In a forcing statement, data or parameters at a final
stage denote their canonical names. Let
\begin{align}
  S_2(s,P_2(\bar v))&\ \longleftrightarrow\ M_s\models P_2(\bar v),
                                                        \label{u:sbase}\\
  S_k(s,P_k(\bar v))&\ \longleftrightarrow\
    \neg\Pos_{\B_k}\bigl(s,\exists z\,
         S_{k-1}(s',P_{k-1}(\bar v,z))\bigr)
           \quad(3\leq k<n).                           \label{u:cert}
\end{align}
These predicates are defined in set theory, using the forcing relation. If
the presentation of the model at \(s\) is fixed, then \(S_k\) depends only
on \(Z_3,\ldots,Z_k\). This follows by recursion from the definition of
\(\B_k\), starting with \eqref{u:sbase}. The choice functions and the
phases introduced below do not enter into these predicates.

\begin{lem}\label{u:forget}
Let \(3\leq k<n\), and let \(\Qq\) be an unrestricted forcing of countable
length starting at \(s\). If every coding at a level above \(k\) in \(\Qq\)
is replaced by the mere coding of \(w_*\) at the same stage, the result is
a forcing in \(\B_k(s)\) with the same iterands, the same intermediate
models and the same sets \(Z_3,\ldots,Z_{k-1}\) at every stage. Every
\(S_k\), \(2\leq k<n\), persists through unrestricted forcings of
countable length.
\end{lem}
\begin{proof}
We keep the original names for all iterands, including every coding of
\(w_*\), and the original names for the codings at levels at most \(k\).
Along the resulting forcing, exactly the same codings change the sets
\(Z_i\), \(i<k\), with the same new reals. The set \(Z_k\) is now
constant, equal to \(Z_k(s)\), which is a subset of its old value at every
later stage. So every remaining coding at level \(k\) is still allowed.
This proves the first assertion at successor stages, and taking unions
proves it at limit stages.

For the persistence of \(S_2\) we use Shoenfield absoluteness. Suppose
that \(S_k(s,P_k)\) holds, but that some condition of an unrestricted
forcing of countable length forces that it fails at the final stage. Below
this condition, append a name for a forcing in \(\B_k\) which forces a
witness to the failure, so that its final stage satisfies
\(\exists z\,S_{k-1}(s',P_{k-1}(\bar v,z))\). Replace the codings above
level \(k\) in the whole composition as above. The final model and every
set \(Z_i\) used by \(S_{k-1}\) stay the same. So the same condition and
the same name for the witness contradict \eqref{u:cert} at \(s\). This uses
equality only up to \(Z_{k-1}\). In particular, the argument applies to
forcings in either presentation of a product and to forcings carried out
in a later phase.
\end{proof}

\begin{lem}\label{u:truth}
Suppose that, for each \(3\leq k<n\), every tuple is considered at level
\(k\) cofinally often, and that whenever \(S_k\) fails for it, a
witnessing forcing in \(\B_k\) is carried out below a condition, which
produces a particular \(z\) with \(S_{k-1}(P_{k-1}(\bar v,z))\). Let
\(M_\infty\) be the final model. Along such an iteration, for
\(2\leq k<n\),
\[
 M_\infty\models P_k(\bar v)
 \quad\longleftrightarrow\quad
 S_k(s,P_k(\bar v))\text{ holds at some stage \(s\) whose model contains }\bar v .
\]
A true instance satisfies \(S_k\) at some stage after any given stage, and
\(S_k\) persists.
\end{lem}
\begin{proof}
For \(k=2\) we use that the parameters belong to the model at some stage
below \(\omega_1\), and Shoenfield absoluteness between the intermediate
models and the final model. Suppose that the assertion holds for \(k-1\).
If \(S_k\) holds and \(P_k\) fails in the final model, choose \(z\) in the
final model as in \eqref{u:normal}. At a later stage, the induction
hypothesis gives \(S_{k-1}\) for this \(z\), and the trivial forcing in
\(\B_k\) at this stage contradicts the persistence of \(S_k\). Conversely,
if \(S_k\) fails when the tuple is considered, we obtain \(z\) with
\(S_{k-1}\), hence \(\neg P_k\) in the final model, by the induction
hypothesis. So an instance which is true in the final model satisfies
\(S_k\) each time it is considered.
\end{proof}

At the final stage of a product we only need persistence, by
Lemma~\ref{u:forget}. Agreement with projective truth is needed only along
the complete iteration with its bookkeeping, in Lemma~\ref{u:truth}. In
particular, the proof does not require the intermediate models to compute
higher projective truth correctly.

\subsection{Ranks and phases}\label{sec:ranks}

We follow \cite[\S2.6]{up} and the proof of \cite[Theorem~3.1]{up}, with
the following dictionary; there, the triples \((m,x,y)\) are coded by
\(\Phi\), here the pairs \((u,y)\) are coded by \(h(u,y)\) at level \(n\):
\[
\begin{array}{@{}l@{\qquad}l@{}}
 \cite{up} & \text{here}\\ \hline
 \neg A(b) & E(s,u,y)\\
 \text{hierarchy }\xi,\ \text{level }(\xi,\beta) & \text{phase }\xi,\ \text{level }d=(\xi,\alpha)\\
 \mathcal K^\xi_\beta(M),\ K_\xi(M) & \A_d(s),\ \K_\xi(s)\\
 b\notin P^\xi_\beta(M) & T_d(s,u,y)\\
 \rho_{\le(\xi,\beta)},\ \Lambda^\xi_\beta,\ B,\ j^\xi_\beta,\ J_M
   & \rank_d,\ \Lambda_d,\ H_d,\ j_d,\ \chi_s\\
 F_\xi & F_\xi
\end{array}
\]

The proofs below are given in full; the dictionary is meant for a
comparison with \cite{up}.

\emph{(R1) The exclusion predicate.} A pair \((u,y)\) is excluded at the
final stage \(s\) of a forcing of countable length if
\begin{equation}\label{u:exclude}
 E(s,u,y)\ \longleftrightarrow\ \exists a\ S_{n-1}\bigl(s,P_{n-1}(u,y,a)\bigr).
\end{equation}
Let \(V_n\) be the final model. By Lemmas~\ref{u:forget} and~\ref{u:truth},
the second along the final iteration with its bookkeeping,
\begin{equation}\label{u:Eprops}
 \begin{gathered}
 E(s,u,y)\ \Longrightarrow\ E(s',u,y)\quad(s'\text{ the final stage of an
   unrestricted forcing starting at }s),\\
 V_n\models\neg A_n(u,y)\ \Longleftrightarrow\ \exists s\ E(s,u,y)
   \ \Longleftrightarrow\ E(s,u,y)\text{ for cofinally many }s .
 \end{gathered}
\end{equation}
The first line gives the persistence of a witness in a product, see
(P3) in the proof of Lemma~\ref{u:preserve}, and the second line is used in
both directions in the proof of Theorem~\ref{main:up}. In
\cite[\S2.6]{up} the corresponding statements for \(\neg A(b)\), and the
passage from an intermediate model to the final model, follow from
Shoenfield absoluteness.

\emph{(R2) The classes.} Let \(F_\xi\) be the cumulative set of pairs
\((u,y)\) for which the coding of \(h(u,y)\) at level \(n\) is forbidden
in phase \(\xi\), with \(F_\zeta=\bigcup_{\xi<\zeta}F_\xi\) for limit
\(\zeta\). Phase \(\xi\) begins with
\[
 \begin{aligned}
 \A_{(\xi,0)}(s):=\bigl\{\Qq\text{ unrestricted over }s:\ &
 \Qq\in\K_\eta(s)\ (\eta<\xi),\\
 &\Qq\text{ codes no }h(u,y)\text{ at level }n\text{ with }(u,y)\in F_\xi\bigr\},
 \end{aligned}
\]
see \cite[Definition~9.9]{mixed}. Earlier levels at a later stage are
evaluated with the classes and the forbidden pairs of their own phase and
with the current coding data; evaluating a predicate never begins a new
phase. For a level \(d=(\xi,\alpha)\) let
\[
 T_d(s,u,y)\ :\longleftrightarrow\ \Pos_{\A_d}\bigl(s,E(s',u,y)\bigr),
\]
and let \(\rank_d(s,u,y)\) be the lexicographically least \(e\le d\) with
\(\neg T_e(s,u,y)\), where earlier phases contribute their levels up to
their fixed levels; \(\rank_d(s,u,y)\) is undefined if there is no such
\(e\). For an input \(u\in M_s\) let
\[
 \Lambda_d(s,u):=\bigl\{\lambda:\ \exists\Qq\in\A_d(s)\ \exists q\in\Qq\
 \exists\dot y\ \ q\Vdash_\Qq\rank_d(s',u,\dot y)=\lambda\bigr\}.
\]
If \(\lambda:=\min\Lambda_d(s,u)\) exists, it is \emph{stable} if no such
\(q\) forces \(\dot y\notin\R^{M_s}\). In this case
\[
 \begin{gathered}
 H_d(s,u):=\bigl\{y\in\R^{M_s}:\ \exists\Qq\in\A_d(s)\ \exists q\in\Qq\
 \ q\Vdash_\Qq\rank_d(s',u,\check y)=\lambda\bigr\}\neq\varnothing,\\
 j_d(s,u):=\chi_s\bigl(H_d(s,u)\bigr),
 \end{gathered}
\]
otherwise \(j_d(s,u)\) is undefined. Here \(H_d(s,u)\neq\varnothing\): a
condition which witnesses the stable minimum does not force
\(\dot y\notin\R^{M_s}\), so it can be strengthened to decide
\(\dot y=\check y\) for some \(y\in\R^{M_s}=\R^{N_s}\). Moreover
\(H_d(s,u)\in N_s\): by (R3) below the attainment formula is invariant, so
by Lemma~\ref{lem:interface}\ref{T:stage} the trivial condition of
\(\Cc_{J_s}\) decides it for each \(y\in\R^{N_s}\), and separation in
\(N_s\) gives \(H_d(s,u)\). The next class is
\begin{equation}\label{u:requirement}
 \begin{aligned}
 \Qq\in\A_{(\xi,\alpha+1)}(s)\ :\longleftrightarrow\ &
 \Qq\in\A_{(\xi,\alpha)}(s)\ \wedge\\
 &\bigl[\,h(u,y)\text{ coded at level }n\text{ at }t\ \wedge\
   j_{(\xi,\alpha)}(t,u)\downarrow\ \Longrightarrow\
   y\neq j_{(\xi,\alpha)}(t,u)\,\bigr],
 \end{aligned}
\end{equation}
for all stages \(t\) of \(\Qq\), below every condition, and
\(\A_{(\xi,\lambda)}:=\bigcap_{\alpha<\lambda}\A_{(\xi,\alpha)}\) for limit
\(\lambda\). Codings at levels below \(n\) are subject only to the rules of
Section~\ref{sec:levelcodings}.

\emph{(R3) Invariance.} All definitions in (R2) are invariant under the
translations of Lemma~\ref{lem:interface}\,\ref{T:trans} at the stage \(s\), i.e.\ under a
permutation \(\pi\in L\) of \(\omega_1\) with \(\pi\res J_s=\mathrm{id}\),
which translates coding areas and covers as in the proof of
Lemma~\ref{s:symmetry}. We prove this by induction on the level \(d\),
simultaneously for all named stages, in the order of the recursion in
(R2).
\begin{enumerate}[leftmargin=*,topsep=.4em,itemsep=.2em,label=\textup{(\alph*)}]
\item \emph{Base.} The definition of an unrestricted forcing reads
\(G^1\) through the coding forcings and the covers; by
Lemma~\ref{lem:interface}\,\ref{T:trans} it is invariant, with the covers
translated by \(\pi\). The finite recursion
\eqref{u:sbase}--\eqref{u:cert} then gives the invariance of \(S_k\),
\(2\le k<n\), and of \(E\). The sets \(Z_j\) and \(F_\xi\) are countable
sets of reals, so they are fixed.
\item \emph{Successor.} Assume that \(\A_d\) is invariant at every named
stage. A translation and its inverse carry the forcings in \(\A_d(s)\),
their conditions and their names for reals to forcings in \(\A_d(s)\),
conditions and names, in both directions. Hence they preserve \(T_d\),
the first level at which \(T\) fails, i.e.\ \(\rank_d\), the attainable
ranks \(\Lambda_d(s,u)\), and whether a new real can attain the minimum.
Now Lemma~\ref{lem:interface}\ref{T:stage} applies to \(H_d(s,u)\), which
belongs to \(N_s\) and is therefore fixed. Transporting the names for the
choice functions with the presentation gives the same \(j_d\), hence the
same restriction \eqref{u:requirement}, and \(\A_{(\xi,\alpha+1)}\) is
invariant.
\item \emph{Limits and new phases.} At limit levels we take
intersections. A new phase retains the earlier fixed classes and the
cumulative set \(F_\xi\), which are invariant.
\end{enumerate}

\begin{lem}\label{u:preserve}
\begin{enumerate}[leftmargin=*,topsep=.4em,itemsep=.2em,label=\textup{(\roman*)}]
\item \emph{Fixed points} \cite[Lemma~2.13]{up}. For each phase \(\xi\)
there is an ordinal \(\theta_\xi\), independent of \(s\), with
\(\A_{(\xi,\alpha)}(s)=\A_{(\xi,\theta_\xi)}(s)=:\K_\xi(s)\) for
\(\alpha\ge\theta_\xi\). The classes are closed under restrictions below
conditions, initial segments, quotients and concatenation, see the
beginning of the proof of \cite[Lemma~2.14]{up}.
\item \emph{One phase} \cite[Lemma~2.14]{up}. Let \(\Qq,\mathbb R\in\A_d(s)\)
use disjoint sets of coding areas, and let \(s_{\mathbb R}\) be the final
stage of \(\mathbb R\). Then \((\Qq\times\mathbb R)/G_\Qq\in\A_d\) over the
final stage of \(\Qq\), and symmetrically. For old pairs and inputs,
\[
 \rank_d(s_{\mathbb R},u,y)=\rank_d(s,u,y),\qquad
 \Lambda_d(s_{\mathbb R},u)=\Lambda_d(s,u),
\]
including undefinedness. A common minimum is stable over \(s\) if and only
if it is stable over \(s_{\mathbb R}\), and in the stable case
\(H_d(s_{\mathbb R},u)=H_d(s,u)\) and \(j_d(s_{\mathbb R},u)=j_d(s,u)\).
\item \emph{Across phases} \cite[Lemma~2.15]{up}. Let \(s\) be a stage in
phase \(\eta\), let \(\Pp\in\K_\eta(s)\) be a part of the final iteration
of countable length starting at \(s\), and assume
\begin{equation}\label{u:newforbidden}
 (u,y)\in F_\xi\setminus F_\eta\ \text{forbidden during }\Pp\
 \Longrightarrow\ y\notin M_s ,
\end{equation}
where \(\xi\) is the phase at the end of \(\Pp\). Let \(\Qq\in\K_\eta(s)\)
use coding areas disjoint from those of \(\Pp\). Then
\((\Qq\times\Pp)/G_\Pp\) belongs to all intermediate fixed classes and to
\(\K_\xi\). For an initial segment \(\Qq_i\), the models
\(S:=M_s[G\res\Qq_i]\) and \(N:=S[G_\Pp]\), and the ranks through
\((\eta,\theta_\eta)\) over \(S\) and through \((\xi,\theta_\xi)\) over
\(N\), the ranks agree on the pairs of \(S\). For \((u,\cdot)\) in \(S\),
if the minimum over \(S\) in \(\K_\eta\) exists, then it is the minimum
over \(N\) in \(\K_\xi\), with the same stability, and in the stable case
with the same set \(H\) and the same choice.
\end{enumerate}
\end{lem}
\begin{proof}

We give the proofs of \cite[Lemmas~2.13--2.15]{up} in the present
notation, with \(E\) in place of the statement \(\neg A(b)\) used there.
Four points have no counterpart there and are treated first: the
definability of the recursion over the stages, the coding areas, the sets
\(Z_j\) of Section~\ref{sec:levelcodings}, and the explicit prohibitions.
For unrestricted forcings \(\Qq\) and \(\mathbb R\) over a stage \(s\) with
disjoint sets of coding areas we write \(s_\Qq\) and \(s_{\mathbb R}\) for
their final stages, and \(s_{\Qq\times\mathbb R}\) for the final stage of
their product.

\emph{(P1) Definability.} The finite recursion
\eqref{u:sbase}--\eqref{u:cert} is fixed before the derivation, and the
sets \(Z_3,\ldots,Z_{n-1}\) and \(F_\xi\) are countable sets of reals.
The forcings of countable length are given by their codes in \(H(\omega_1)\),
see Lemma~\ref{lem:interface}\,\ref{T:codes} and (D2) before
Lemma~\ref{s:symmetry}; together with the
names for reals, sets of reals, choice functions, conditions and
embeddings of initial segments they form a set. So the recursion of
(R2) can be carried out with forcing relations, as in the proof of (i)
below. By (R3),
these definitions are invariant under the translations of
Lemma~\ref{lem:interface}\,\ref{T:trans}.

\emph{(P2) Coding areas.} We apply Lemma~\ref{lem:interface}\ref{T:disj}
with \(\varphi(\Qq,q)\) asserting both the relevant forcing statement and
the membership of \(\Qq\) in the relevant class \(\A_d\); this formula is
invariant by (R3), so the relocated witness belongs to the same class.
Whenever a forcing which is given by an existential statement has to
be multiplied with a given forcing, we first \emph{relocate} it in this
way, with \(Z\) the set of coding areas of the given forcing, and
Lemma~\ref{lem:interface}\ref{T:prod} gives the products and
continuations of \cite[Lemma~2.9(3)]{up}. In the proof across phases we
apply this over \(S=M_s[G\res\Qq_i]\) to the witnessing continuation,
avoiding the coding areas of \(\Pp\). This keeps the prescribed initial
segment \(\Qq_i\) and its choices, and
Lemma~\ref{lem:interface}\ref{T:lift}--\ref{T:prod} lets us append the
relocated tail; this is the use of a prescribed initial segment in
\cite[Lemmas~2.14 and~2.15]{up}.

\emph{(P3) The sets \(Z_j\) in products.} In \(\Qq\times\mathbb R\) put
\(Z_j:=Z^{\Qq}_j\cup Z^{\mathbb R}_j\). A prohibition
\(e_{j-1}(d,t)\in Z_{j-1}\) which is new in \(\Qq\) at its stage \(\gamma\)
contains the generic real \(t\) of this stage, and
\[
 t\notin M_s[G_{\mathbb R}]\,[G\res\Qq_\gamma],
\]
by \eqref{s:intersection} for the product of \(\Qq_{\gamma+1}\) and
\(\mathbb R\), factored over \(\Qq_\gamma\). So no coding of \(\mathbb R\)
violates it, in either order, and the two orders have the same sets
\(Z_j\) at the common final stage. Hence \(\mathbb R\) is an
unrestricted forcing over \(s_\Qq\) with final stage
\(s_{\Qq\times\mathbb R}\), and likewise \(\Qq\) over \(s_{\mathbb R}\).
By the first line of \eqref{u:Eprops}, for \((u,y)\) in the respective
model,
\phantomsection
\begin{equation*}\tag{\(\dagger\)}\label{u:Eprod}
 E(s_\Qq,u,y)\ \Longrightarrow\ E(s_{\Qq\times\mathbb R},u,y),\qquad
 E(s_{\mathbb R},u,y)\ \Longrightarrow\ E(s_{\Qq\times\mathbb R},u,y).
\end{equation*}

\emph{(P4) Explicit prohibitions across phases.} A coding in \(\Qq\) cannot
mention the value \(y\) of a pair which is newly forbidden during \(\Pp\):
by \eqref{u:newforbidden} and \eqref{s:intersection},
\[
 y\in M_s[G\res\Qq_i]\cap M_s[G_\Pp]=M_s\quad\text{would contradict}\quad
 y\notin M_s .
\]

\emph{Proof of (i).} Let \(s_0\) be the stage at which phase \(\xi\) begins.
In \(M_{s_0}\) let \(X_\xi\) be the set of all triples
\((\mathbb H,\dot\Tt,p)\), where \(\mathbb H\) is one of the presentations
of (P1) of an unrestricted forcing over \(s_0\), with its names for the
data of its stages, \(\dot\Tt\) is an \(\mathbb H\)-name for an unrestricted
forcing of a fixed countable length over the final stage of \(\mathbb H\),
and \(p\in\mathbb H\). Products and quotients are among these
presentations. By recursion on \(\alpha\) we define
\[
 D_\alpha:=\bigl\{(\mathbb H,\dot\Tt,p)\in X_\xi:\
   p\Vdash_{\mathbb H}\dot\Tt\in\A_{(\xi,\alpha)}\bigr\}.
\]
For \(\alpha=0\) the triple has to force that \(\dot\Tt\) codes no
\(h(u,y)\) at level \(n\) with \((u,y)\in F_\xi\) and belongs to
\(\K_\eta\) for \(\eta<\xi\); the latter is given by the sets
\(D_{\theta_\eta}\) of the earlier phases, along the embeddings from the
stages at which these phases began. Given \(D_\alpha\), the forcing
theorem expresses, by conditions and names over the forcings
\(\mathbb H*\dot\Tt_t\),
\[
 T_{(\xi,\alpha)}(t,\dot u,\dot y),\qquad
 \rank_{(\xi,\alpha)}(t,\dot u,\dot y)=\lambda,\qquad
 \lambda\in\Lambda_{(\xi,\alpha)}(t,\dot u),\qquad
 \dot y=j_{(\xi,\alpha)}(t,\dot u)
\]
at the stages \(t\) of \(\dot\Tt\): each of them quantifies over triples
in the sets \(D_\beta\), \(\beta\le\alpha\), and in the corresponding sets
of the earlier phases, and uses the forcing relation for \(E\), which is
definable by (P1). So \eqref{u:requirement} defines \(D_{\alpha+1}\), and at limits
we take intersections. The sets \(D_\alpha\) decrease. Every triple leaves
them at most once, so there is \(\theta_\xi\) with
\(D_\alpha=D_{\theta_\xi}\) for \(\alpha\ge\theta_\xi\). If \(G\subseteq
\mathbb H\) is generic and \(s\) is the final stage of \(\mathbb H\), then
\[
 \dot\Tt^G\in\A_{(\xi,\alpha)}(s)\ \Longleftrightarrow\
 \exists p\in G\ \ (\mathbb H,\dot\Tt,p)\in D_\alpha .
\]
Every stage reached from \(s_0\) and every unrestricted forcing over it are
presented in this way, so \(\A_{(\xi,\alpha)}(s)=\A_{(\xi,\theta_\xi)}(s)\)
for \(\alpha\ge\theta_\xi\) at all these stages.

The requirements which define the classes concern the single stages \(t\)
of a forcing, below every condition, and \(j_{(\xi,\alpha)}(t,u)\) only
depends on the stage \(t\). So, by induction on the level \(d\), for
\(\Tt\in\A_d(s)\) of length \(\delta\), \(t\in\Tt\) and \(\gamma\le\delta\),
with \(s_\gamma\) the stage \(\gamma\) of \(\Tt\),
\begin{gather*}
 \Tt\res t\in\A_d(s),\qquad \Tt_\gamma\in\A_d(s),\qquad
 \Tt_{\gamma,\delta}\in\A_d(s_\gamma),\\
 \Tt\in\A_d(s)\ \wedge\ \Tt\Vdash\dot{\mathbb U}\in\A_d(s_\Tt)
 \ \Longrightarrow\ \Tt*\dot{\mathbb U}\in\A_d(s).
\end{gather*}
For the quotient we use Lemma~\ref{lem:interface}\ref{T:seg}: its conditions
are the restrictions \(p\res[\gamma,\delta)\) of the conditions \(p\) of
\(\Tt\) with \(p\res\gamma\in G\res\gamma\). A name which is given below a
condition is first lifted by Lemma~\ref{lem:interface}\ref{T:lift}; trivial
iterands satisfy all requirements.

\emph{Proof of (ii).} We argue by induction on \(d\) in the lexicographic
order of the levels, where earlier phases contribute their levels up to
their fixed levels. The stages of \(\Qq\times\mathbb R\) carry the choice
functions of \ref{CF:ext}, which extend those of \(\Qq\) and of
\(\mathbb R\).

\emph{Products.} We show \((\Qq\times\mathbb R)/G_\Qq\in\A_d(s_\Qq)\); the
other order is symmetric. By (P3), \(\mathbb R\) is an unrestricted
forcing over \(s_\Qq\), and it codes the same reals at level \(n\) as over
\(s\). So it codes no \(h(u,y)\) with \((u,y)\in F_\xi\), and it belongs to
\(\K_\eta(s_\Qq)\) for \(\eta<\xi\) by the induction hypothesis at the
levels \((\eta,\theta_\eta)\). Let \(e=(\xi,\beta)<d\), and let \(h(u,y)\)
be coded at level \(n\) at the stage \(j\) of \(\mathbb R\). Let \(t\) be
this stage of \(\mathbb R\) over \(s\), and \(t'\) this stage of
\(\mathbb R\) over \(s_\Qq\). Then \(t'\) is the final stage of
\((\Qq\times\mathbb R_j)/G_{\mathbb R_j}\), which belongs to \(\A_e(t)\) by
the induction hypothesis at \(e\) for \(\Qq\) and \(\mathbb R_j\). The
induction hypothesis at \(e\), applied over \(t\) to this forcing, gives
\[
 j_e(t',u)=j_e(t,u),
\]
including definedness. As \(\mathbb R\in\A_d(s)\), we have
\(y\ne j_e(t,u)\) if \(j_e(t,u)\) is defined. So \eqref{u:requirement}
holds at \(t'\).

\emph{Ranks.} Let \((u,y)\in M_s\) and \(e\le d\); note
\(\mathbb R\in\A_e(s)\). We show
\[
 T_e(s_{\mathbb R},u,y)\ \Longleftrightarrow\ T_e(s,u,y),
\]
which gives \(\rank_d(s_{\mathbb R},u,y)=\rank_d(s,u,y)\), including
undefinedness. Let \(\Tt\in\A_e(s_{\mathbb R})\) and \(t\in\Tt\) with
\(t\Vdash E(s',u,y)\). Choose names \(\dot\Tt,\dot t\), decide the length
of \(\dot\Tt\) below some \(r\in G_{\mathbb R}\), and lift \(\dot\Tt\). Then
\(\mathbb R\res r*\dot\Tt^\uparrow\in\A_e(s)\) by (i), and \((r,\dot t)\)
forces \(E\) at its final stage, which is the final stage of \(\Tt\). So
\(T_e(s,u,y)\) holds. Conversely, let \(\Qq'\in\A_e(s)\) and \(q\in\Qq'\)
with \(q\Vdash E(s_{\Qq'},u,y)\). By (P2) we may assume that the coding
areas of \(\Qq'\) are disjoint from those of \(\mathbb R\). Then
\((\Qq'\times\mathbb R)/G_{\mathbb R}\in\A_e(s_{\mathbb R})\) by the
assertion on products, and by \eqref{u:Eprod} \(q\) forces \(E\) at its
final stage \(s_{\Qq'\times\mathbb R}\). So \(T_e(s_{\mathbb R},u,y)\)
holds.

\emph{Attainable ranks.} If \(\Tt\in\A_d(s_{\mathbb R})\), \(t\) and
\(\dot y\) witness \(\lambda\in\Lambda_d(s_{\mathbb R},u)\), then, as
before, \(\mathbb R\res r*\dot\Tt^\uparrow\), \((r,\dot t)\) and \(\dot y\)
witness \(\lambda\in\Lambda_d(s,u)\). Conversely, let \(\Qq'\in\A_d(s)\),
\(q\) and \(\dot y\) witness \(\lambda\in\Lambda_d(s,u)\), with \(\Qq'\)
relocated by (P2). By the part on ranks, applied over \(s_{\Qq'}\) to
\((\Qq'\times\mathbb R)/G_{\Qq'}\in\A_d(s_{\Qq'})\),
\[
 q\Vdash\rank_d(s_{\Qq'\times\mathbb R},u,\dot y)
       =\rank_d(s_{\Qq'},u,\dot y)=\lambda ,
\]
and \(s_{\Qq'\times\mathbb R}\) is the final stage of
\((\Qq'\times\mathbb R)/G_{\mathbb R}\in\A_d(s_{\mathbb R})\). So
\(\lambda\in\Lambda_d(s_{\mathbb R},u)\), and
\(\Lambda_d(s_{\mathbb R},u)=\Lambda_d(s,u)\).

\emph{Stability.} Let \(\lambda\) be the common minimum. If \(\Qq',q,
\dot y\) are as before and \(q\Vdash\dot y\notin\R^{M_s}\), then
\(\dot y\in M_s[G_{\Qq'}]\setminus M_s\) and
\(\dot y\notin M_s[G_{\mathbb R}]=M_{s_{\mathbb R}}\) by
\eqref{s:intersection}; so \(\lambda\) is not stable over
\(s_{\mathbb R}\). If \(\Tt,t,\dot y\) are as before and
\(t\Vdash\dot y\notin\R^{M_{s_{\mathbb R}}}\), then the composed witness
forces \(\dot y\notin\R^{M_s}\); so \(\lambda\) is not stable over \(s\).

\emph{Minimizers and choices.} Let \(\lambda\) be stable. By the part on
ranks and the trivial forcing,
\[
 H_d(s,u)=\{y\in\R^{M_s}:\rank_d(s,u,y)=\lambda\},\qquad
 H_d(s_{\mathbb R},u)=\{y\in\R^{M_{s_{\mathbb R}}}:
    \rank_d(s_{\mathbb R},u,y)=\lambda\},
\]
so \(H_d(s,u)=H_d(s_{\mathbb R},u)\cap M_s\). Let
\(y\in H_d(s_{\mathbb R},u)\), and let \(\dot y\) be an
\(\mathbb R\)-name for \(y\). Every sufficiently strong
\(r\in G_{\mathbb R}\) forces \(\rank_d(s_{\mathbb R},u,\dot y)=\lambda\), so
\(\mathbb R\res r\), \(r\) and \(\dot y\) witness \(\lambda\) over \(s\). By
the stability over \(s\), no such \(r\) forces \(\dot y\notin\R^{M_s}\).
Hence \(y\in M_s\). So both sets are equal to a set \(B\), which belongs
to \(N_s\) by (R2), and
\[
 j_d(s_{\mathbb R},u)=\chi_{s_{\mathbb R}}(B)=\chi_s(B)=j_d(s,u)
\]
by the coherence of the choice functions.

\emph{Proof of (iii).} Write \(d_\eta:=(\eta,\theta_\eta)\), and let \(S\)
and \(N\) also denote the stages with these models, i.e.\ the stage \(i\)
of \(\Qq\) over \(s\) and over the final stage \(s_\Pp\) of \(\Pp\). If
\(\xi=\eta\), the assertions of (iii) are those of (ii), applied to
\(\Qq\) and \(\Pp\) over \(s\) and to \(\Pp\) over \(S\). So let
\(\eta<\xi\). A \emph{later level} is a level \(e=(\xi',\beta)\) with
\(\eta<\xi'\le\xi\) and \(\beta\le\theta_{\xi'}\). By induction on the
later levels \(e\) we prove, simultaneously for all \(\Qq\in\K_\eta(s)\) whose coding areas are
disjoint from those of \(\Pp\) and for all their initial segments
\(\Qq_i\):
\begin{enumerate}[leftmargin=*,topsep=.4em,itemsep=.2em,label=\textup{(\alph*\(_e\))}]
\item \((\Qq\times\Pp)/G_\Pp\in\A_e(s_\Pp)\);
\item \(\rank_e(N,u,y)=\rank_{d_\eta}(S,u,y)\) for \((u,y)\in M_S\),
including undefinedness;
\item if \(\lambda=\min\Lambda_{d_\eta}(S,u)\) exists, then
\(\lambda=\min\Lambda_e(N,u)\), with the same stability, and in the stable
case \(H_e(N,u)=H_{d_\eta}(S,u)\) and \(j_e(N,u)=j_{d_\eta}(S,u)\).
\end{enumerate}
For \(e=(\xi,\theta_\xi)\) these are the assertions of (iii). On the grid
of \(\Qq\times\Pp\), and on the grids of the continuations below, we use
any choice functions as in \ref{CF:ext} and \ref{CF:ops}; the proof only
uses their coherence. Let
\(\mathbb V:=(\Qq_i\times\Pp)/G_{\Qq_i}\), which is \(\Pp\) carried out over
\(S\). By (ii), \(\mathbb V\in\K_\eta(S)\) and
\phantomsection
\begin{equation*}\tag{\(\ddagger\)}\label{u:iii0}
 \rank_{d_\eta}(N,u,y)=\rank_{d_\eta}(S,u,y)\qquad((u,y)\in M_S).
\end{equation*}

\emph{Level zero.} Let \(e=(\xi',0)\). The forcing
\((\Qq\times\Pp)/G_\Pp\) belongs to \(\K_\eta(s_\Pp)\) by (ii), and to
\(\K_{\xi''}(s_\Pp)\) for \(\eta<\xi''<\xi'\) by (a) at the fixed levels of
these phases. It codes no \(h(u,y)\) with \((u,y)\in F_\eta\), as
\(\Qq\in\K_\eta(s)\), and none with \((u,y)\in F_{\xi'}\setminus F_\eta\),
by (P4). This is (a\(_e\)).

\emph{From \textup{(a\(_e\))} to \textup{(b\(_e\))}.} Let \((u,y)\in M_S\).
If \(\rank_{d_\eta}(S,u,y)\) is defined, it is equal to
\(\rank_{d_\eta}(N,u,y)\) by \eqref{u:iii0}, and hence to
\(\rank_e(N,u,y)\). Suppose that it is undefined. Then there are
\(\Tt\in\K_\eta(S)\) and \(t\in\Tt\) with \(t\Vdash E(s',u,y)\). By (P2)
we relocate \(\Tt\) over \(S\) so that its coding areas are disjoint from
those of \(\Pp\). We restrict \(\Tt\) below \(t\), choose its name and
decide its length below a condition of \(G_{\Qq_i}\), and lift it. By (i),
\[
 \Qq':=\Qq_i*\dot\Tt^\uparrow\in\K_\eta(s),
\]
and its coding areas are disjoint from those of \(\Pp\). By (a\(_e\)) for
\(\Qq'\), by (i) and by Lemma~\ref{lem:interface}\ref{T:prod},
\[
 \mathbb U:=(\Qq'\times\Pp)/(G_{\Qq_i}\times G_\Pp)\in\A_e(N),
\]
and \(\mathbb U\) is \(\Tt\) carried out over \(N\). By \eqref{u:Eprod} for
\(\Qq'\) and \(\Pp\), \(\mathbb U\) forces \(E\) at its final stage. So
\(T_{e'}(N,u,y)\) holds for all later levels \(e'\le e\), as
\(\A_e(N)\subseteq\A_{e'}(N)\), and \(\rank_{d_\eta}(N,u,y)\) is undefined
by \eqref{u:iii0}. Hence \(\rank_e(N,u,y)\) is undefined.

\emph{From \textup{(a\(_e\))} and \textup{(b\(_e\))} to \textup{(c\(_e\))}.}
Let \(\lambda=\min\Lambda_{d_\eta}(S,u)\), and let \(\Tt\in\K_\eta(S)\),
\(t\) and \(\dot y\) witness it. Define \(\Qq'\) and \(\mathbb U\) as
before. By (b\(_e\)) for \(\Qq'\) and its full length,
\[
 \rank_e(N',u,\dot y)=\rank_{d_\eta}(S',u,\dot y)=\lambda ,
\]
where \(S'\) and \(N'\) are the final stages of \(\Qq'\) over \(s\) and
over \(s_\Pp\); and \(N'\) is the final stage of \(\mathbb U\) over \(N\).
So \(\lambda\in\Lambda_e(N,u)\). Now let \(\mu\le\lambda\) belong to
\(\Lambda_e(N,u)\), witnessed by \(\mathbb U'\in\A_e(N)\), a condition and a
name \(\dot y'\). Lift a name for \(\mathbb U'\). As
\(\A_e(N)\subseteq\K_\eta(N)\), (i) gives
\[
 \mathbb V*\dot{\mathbb U}'^{\uparrow}\in\K_\eta(S).
\]
As \(\mu\le\lambda\le d_\eta\), the composed witness forces
\(\rank_{d_\eta}(s',u,\dot y')=\mu\). So \(\mu\in\Lambda_{d_\eta}(S,u)\),
hence \(\mu=\lambda\) and \(\lambda=\min\Lambda_e(N,u)\). Moreover, if the
witness over \(N\) forces \(\dot y'\notin\R^{M_N}\), then the composed
witness forces \(\dot y'\notin\R^{M_S}\). Conversely, if
\(t\Vdash\dot y\notin\R^{M_S}\), then \(\dot y\in M_S[G_\Tt]\setminus M_S\)
and \(\dot y\notin M_N\) by \eqref{s:intersection}, so \(\mathbb U\)
forces \(\dot y\notin\R^{M_N}\). So \(\lambda\) has the same stability
over \(S\) and over \(N\).

Let \(\lambda\) be stable. By (ii) and the trivial forcing,
\(H_{d_\eta}(S,u)\) and \(H_e(N,u)\) are the sets of the reals of \(M_S\),
respectively of \(M_N\), with rank \(\lambda\), and the first set is
contained in the second by (b\(_e\)). Let \(y\in H_e(N,u)\), and let
\(\dot y\) be a \(\mathbb V\)-name for \(y\). Then
\(\rank_{d_\eta}(N,u,y)=\lambda\), so for every sufficiently strong
\(v\in G_{\mathbb V}\), the forcing \(\mathbb V\res v\), \(v\) and
\(\dot y\) witness \(\lambda\) over \(S\). By the stability over \(S\), no
such \(v\) forces \(\dot y\notin\R^{M_S}\). Hence \(y\in M_S\), and
\(y\in H_{d_\eta}(S,u)\) by \eqref{u:iii0}. So the two sets are equal, and
the choices agree by coherence.

\emph{From \textup{(a\(_e\))}, \textup{(b\(_e\))}, \textup{(c\(_e\))} to
the next level.} Let \(e=(\xi',\beta)\), and let \(h(u,y)\) be coded at
level \(n\) at the stage \(i\) of \(\Qq\). As \(\Qq\in\K_\eta(s)\), we have
\(y\ne j_{d_\eta}(S,u)\) if this value is defined. If
\(\Lambda_{d_\eta}(S,u)\ne\varnothing\), then
\(j_e(N,u)=j_{d_\eta}(S,u)\), including definedness, by (c\(_e\)). If
\(\Lambda_{d_\eta}(S,u)=\varnothing\), then \(\rank_{d_\eta}(S,u,y)\) is
undefined, by the trivial forcing, and \(\rank_e(N,u,y)\) is undefined by
(b\(_e\)), whereas \(j_e(N,u)\) has a rank over \(N\) if it is defined. In
both cases \eqref{u:requirement} holds at \(N\) for the level \(e\).
Together with (a\(_e\)) this is (a) at the level \((\xi',\beta+1)\). At
limit levels we take intersections, and a new phase starts with the step
of level zero.\qedhere

\end{proof}

\emph{(R4) Earlier requirements.} As after \cite[Lemma~2.14]{up}: if a
requirement \eqref{u:requirement} of an earlier level \(e\le d\), or of an
earlier phase, forbids a coding at \(u\) over \(M_s\), then
\[
 j_e(s,u)\downarrow\ \Longrightarrow\ j_d(s,u)=j_e(s,u),
\]
with the same stable minimum and the same set \(H\), for the current fixed
level \(d\). Indeed, the old values already have the minimal rank at \(s\),
so the trivial forcing keeps them available, and a smaller rank, or a
further value of the same rank, in the smaller class would already have
been available in the larger one. If the current minimum can be attained by
a new real, then no earlier requirement of this kind is defined at \(u\)
over \(M_s\). The explicit prohibitions in \(F_\xi\) still apply.

\begin{lem}\label{u:newmin}
Let \(s\) be a stage of phase \(\xi\), let \(d=(\xi,\theta_\xi)\) and
\(u\in M_s\), and suppose that \(\lambda=\min\Lambda_d(s,u)\) exists and
is not stable. Then there are \(\Qq\in\K_\xi(s)\) of length \(\delta>0\)
and a \(\Qq\)-name \(\dot y\) with
\[
 \Qq\Vdash\ \dot y\notin\bigcup_{\gamma<\delta}\R^{M_s[\dot G\res\gamma]}
 \ \wedge\ \rank_d(s',u,\dot y)=\lambda .
\]
\end{lem}
\begin{proof}
This is \cite[Lemma~2.16]{up}. Choose \(\Qq^0\in\K_\xi(s)\) of length
\(\delta_0\), \(q\in\Qq^0\) and \(\dot y\) with
\[
 q\Vdash\ \dot y\notin\R^{M_s}\ \wedge\ \rank_d(s',u,\dot y)=\lambda .
\]
Strengthen \(q\) to decide the ordinal
\(\delta:=\min\{\gamma\le\delta_0:\dot y\in M_s[\dot G\res\gamma]\}\) and
to force \(\dot y=\dot z\) for a \(\Qq^0_\delta\)-name \(\dot z\). Then
\(\delta>0\), and \(q\) forces
\(\dot z\notin\bigcup_{\gamma<\delta}\R^{M_s[\dot G\res\gamma]}\). The
remaining forcing \(\Qq^0_{\delta,\delta_0}\) belongs to \(\K_\xi\) over
the stage \(\delta\) of \(\Qq^0\), by Lemma~\ref{u:preserve}(i). So, by
Lemma~\ref{u:preserve}(ii), \(q\) forces that \((u,\dot z)\) has the rank
\(\lambda\) at this stage. Both assertions concern the extension by
\(\Qq^0_\delta\). If some \(r\le q\res\delta\) in \(\Qq^0_\delta\) forced
the negation of one of them, then
\(r\cup q\res[\delta,\delta_0)\le q\) by
Lemma~\ref{lem:interface}\ref{T:seg}, a contradiction. So
\(\Qq:=\Qq^0_\delta\res(q\res\delta)\), which belongs to \(\K_\xi(s)\) by
Lemma~\ref{u:preserve}(i), and \(\dot z\) are as required.
\end{proof}

\subsection{The iteration and the uniformizing function}

At a stage \(s\) of the final iteration, let \(\xi_s\) be the current phase
and abbreviate
\[
 \begin{aligned}
 \K_s&:=\K_{\xi_s}(s),&
 \rank_s&:=\rank_{(\xi_s,\theta_{\xi_s})}(s,\cdot),\\
 \Lambda_s(u)&:=\Lambda_{(\xi_s,\theta_{\xi_s})}(s,u),&
 j_s(u)&:=j_{(\xi_s,\theta_{\xi_s})}(s,u).
 \end{aligned}
\]
In addition to the codings of Lemma~\ref{u:decode} and the tasks of
Lemma~\ref{u:truth}, we carry out the following tasks cofinally often, for
all relevant arguments; they are the instructions (1), (3) and (2) of the
proof of \cite[Theorem~3.1]{up}. The value \(y_*(u)\) is defined at most
once.
\begin{enumerate}[leftmargin=*,topsep=.4em,itemsep=.2em,label=\textup{(T\arabic*)}]
\item\label{u:T1} If \(\rank_s(u,y)\) is undefined, append \(\Qq\res q\)
with \(\Qq\in\K_s\), \(q\Vdash_\Qq E(s',u,y)\).
\item\label{u:T2} If \(y_*(u)\) is undefined and
\(\Lambda_s(u)\neq\varnothing\), let \(\lambda:=\min\Lambda_s(u)\). If
\(\lambda\) is stable, put \(y_*(u):=j_s(u)\). Otherwise, by Lemma~\ref{u:newmin},
there are \(\Qq\in\K_s\) of length \(\delta>0\) and \(\dot y\) with
\[
 \Qq\Vdash\ \dot y\notin\bigcup_{\gamma<\delta}\R^{M_s[\dot G\res\gamma]}
 \ \wedge\ \rank_{(\xi_s,\theta_{\xi_s})}(s',u,\dot y)=\lambda .
\]
Append \(\Qq\), put \(y_*(u):=\dot y^{G_\Qq}\), and begin the phase
\(\xi_s+1\) with
\(F_{\xi_s+1}:=F_{\xi_s}\cup\{(u,y_*(u))\}\), before the next iterand.
\item\label{u:T3} If \(\K_s\) permits it, code \(h(u,y)\) at level \(n\).
\end{enumerate}

The bookkeeping is the one of the proof of \cite[Theorem~3.1]{up}. At each
stage we choose a name for an enumeration of its reals in order type
\(\omega_1\), and each kind of task, with indices into earlier
enumerations, is scheduled cofinally often; the tasks for the lower
predicates include every \(3\le j<n\) and every tuple of
\eqref{u:normal}. A task of countable length is given by a name at its
first stage, chosen below the elements of a maximal antichain which decides
its length. At a limit \(\zeta\) of phases we put
\(F_\zeta=\bigcup_{\xi<\zeta}F_\xi\) and
\(\A_{(\zeta,0)}=\bigcap_{\xi<\zeta}\K_\xi\), and between tasks we extend
the choice functions coherently. A task without effect gets one stage with
a trivial iterand. So the starting stages \(\gamma_\nu\) of the tasks
satisfy
\[
 \gamma_{\nu+1}=\gamma_\nu+\ell_\nu\ \ (0<\ell_\nu<\omega_1),\qquad
 \gamma_\delta=\sup_{\nu<\delta}\gamma_\nu\ \ (\delta\text{ limit}),\qquad
 \sup_{\nu<\omega_1}\gamma_\nu=\omega_1 .
\]
The initial segments of \(\Pp_{\omega_1}\) are allowable over
\((W,\varnothing)\), as at the end of Section~\ref{sec:long}, and for
\(s<s'<\omega_1\)
\begin{equation}\label{u:segments}
 \Pp_{s,s'}\in\K_s,\qquad
 (u,y)\in F_{\xi_{s'}}\setminus F_{\xi_s}\ \Longrightarrow\ y\notin M_s ,
\end{equation}
the second by the stopping point in \ref{u:T2}. By Lemma~\ref{lem:interface}
\ref{T:size}, \(\Pp^0\times\Pp^1\) followed by \(\Pp_{\omega_1}\) is a set
forcing in \(L\) which preserves all cardinals and \(\CH\).

\begin{proof}[Proof of Theorem~\ref{main:up} for \(n\geq4\)]

Let \(V_n\) be the final model. By \eqref{u:segments} and
Lemma~\ref{u:preserve}(ii)--(iii),
\[
 \rank_{s'}\res(\omega\times\R^2)^{M_s}=\rank_s\qquad(s\le s').
\]
For \((u,y)\in M_s\), \eqref{u:Eprops} and \ref{u:T1} give
\[
 \begin{aligned}
 \rank_s(u,y)\downarrow\ &\Longrightarrow\ \neg E(s',u,y)\ \text{for all
   }s'\ge s\ \Longrightarrow\ V_n\models A_n(u,y),\\
 \rank_s(u,y)\uparrow\ &\overset{\ref{u:T1}}{\Longrightarrow}\
   E(s',u,y)\ \text{for some }s'\ge s\ \Longrightarrow\
   V_n\models\neg A_n(u,y);
 \end{aligned}
\]
for the first line, if \(E(s',u,y)\) held, the part \(\Pp_{s,s'}\in\K_s\),
restricted below a suitable condition, would witness the instance of \(T\)
which fails at \(s\). Hence \(V_n\models A_n(u,y)\) if and only if
\(\rank_s(u,y)\) is defined. With the trivial forcing,
\[
 A_{n,u}^{V_n}\neq\varnothing\ \Longrightarrow\ \exists s,y\
 \rank_s(u,y)\downarrow\ \Longrightarrow\ \exists s\
 \Lambda_s(u)\neq\varnothing\ \overset{\ref{u:T2}}{\Longrightarrow}\
 y_*(u)\downarrow\ \wedge\ V_n\models A_n(u,y_*(u)).
\]

Let \(A_{n,u}^{V_n}\neq\varnothing\), let \(y_*:=y_*(u)\) be chosen at the
stage \(\gamma\), and let \(\mathcal C_n\) be the set of reals coded at
level \(n\). If \(y_*\) was chosen at a stable minimum, then for every
stage \(s\) with \(u,y_*\in M_s\), by Lemma~\ref{u:preserve} and (R4),
\[
 \rank_s(u,y_*)\downarrow\ \Longrightarrow\ \Lambda_s(u)\neq\varnothing
 \ \Longrightarrow\ j_s(u)=j_\gamma(u)=y_*,
\]
so \eqref{u:requirement} gives \(h(u,y_*)\notin\mathcal C_n\). Otherwise,
by \ref{u:T2},
\[
 s<\gamma\ \Longrightarrow\ y_*\notin M_s,\qquad
 s\ge\gamma\ \Longrightarrow\ (u,y_*)\in F_{\xi_s},\qquad\text{so}\quad
 h(u,y_*)\notin\mathcal C_n .
\]
Now let \ref{u:T3} consider \((u,y)\) at a stage \(s\ge\gamma\) with
\(y\neq y_*\). In the stable case, (R4) gives \(j_e(s,u)=y_*\) for every
level \(e\) whose requirement \eqref{u:requirement} is defined at \(u\).
In the unstable case, let \(\gamma_\nu\) be the start of the task which
chose \(y_*\); then \(y_*\notin M_{\gamma_\nu}\), and
Lemma~\ref{u:preserve}(iii) over \(M_{\gamma_\nu}\) shows that the minimum
at \(u\) can still be attained by a new real at \(s\), so \(j_e(s,u)\) is
undefined for every such level \(e\), by (R4). As
\(\{z:(u,z)\in F_{\xi_s}\}\subseteq\{y_*\}\), task \ref{u:T3} codes
\(h(u,y)\) at level \(n\). By Lemma~\ref{u:decode}, in \(V_n\),
\[
 D_n(h(u,y))\ \Longleftrightarrow\ h(u,y)\in\mathcal C_n\
 \Longleftrightarrow\ y\neq y_*(u)\qquad(A_{n,u}\neq\varnothing).
\]
So the lightface relation
\begin{equation}\label{u:graph}
 U_n(u,y)\ \longleftrightarrow\ A_n(u,y)\ \wedge\ \neg D_n(h(u,y))
\end{equation}
satisfies \(\{y:U_n(u,y)\}=\{y_*(u)\}\) if \(A_{n,u}\neq\varnothing\), and
\(\{y:U_n(u,y)\}=\varnothing\) otherwise. As \(D_n\in\Sigma^1_n\),
\(U_n\in\Pi^1_n\). Substituting \(u=\pa{m,p,x}\) gives a uniformizing
function for the relation with parameter \(p\), defined from the same
parameter; in particular, every lightface relation has a lightface
uniformizing function. Cardinals and \(\CH\) are preserved by
Lemma~\ref{lem:interface}\ref{T:size}.

\end{proof}

\paragraph{The two candidates revisited.}
Let \(n=4\), let \(A\) be the relation of the example at the beginning of
this section, and let \(u=\pa{m,p,x}\) be an input of \(A_4\) with
\(A_4(u,\cdot)=A_x\). A real is put into \(Z_3\) only in the form
\(e_3(d,t)\), where \(t\) is the generic real of a stage, and by
Lemma~\ref{lem:interface}\,\ref{T:fresh} such a \(t\) does not belong to the model of
its stage. The reals \(t_{k,j}\) are recursive. So no real
\(b_{k,j}(x)=e_3(h(x,v_k),t_{k,j})\) belongs to \(Z^\infty_3\), and by the
proof of Lemma~\ref{u:decode} all of them are coded at level three:
\[
 V_4\models\forall k,j<\omega\ D_3(b_{k,j}(x)),\qquad A_x^{V_4}=\varnothing .
\]
By the proof of Theorem~\ref{main:up}, \(\rank_s(u,v)\) is undefined for
every stage \(s\) and every \(v\in M_s\), the value \(y_*(u)\) is never
defined, and \(\Lambda_s(u)=\varnothing\) at every stage. By (R4), no
requirement \eqref{u:requirement} is then defined at \(u\), and no pair
\((u,v)\) belongs to a set \(F_\xi\). So task~\ref{u:T3} codes \(h(u,v)\)
at level four for every \(v\), and in particular
\[
 V_4\models D_4(h(u,y))\wedge D_4(h(u,z)).
\]
In the toy iteration both candidates ended uncoded; here both are coded,
and in both cases \eqref{u:graph} has an empty section at \(u\). The
reservation which was cancelled in the example cannot be made in the
construction itself. The witness for \(D_4(h(u,v))\) is a real
\(e_3(h(u,v),t)\) with a generic real \(t\). It is forbidden at level
three from the stage at which it appears, and no relation can name it in
advance.

\section{Reduction and the separation pattern}\label{sec:consequences}

The reduction property \(\Red(\Gamma)\) was defined in the introduction.
The following facts are standard, see
\cite[Exercises~4B.11--4B.12]{mosch} and \cite[\S22]{kechris} for
reduction and separation.

\begin{prop}\label{prop:duality}
For every projective pointclass \(\Gamma\), uniformly relative to a real
parameter,
\[
  \UP(\Gamma)\ \Longrightarrow\ \Red(\Gamma)\
       \Longrightarrow\ \Sep(\check\Gamma).
\]
For \(k\ge1\), \(\UP(\Pi^1_k(a))\) also implies
\(\UP(\Sigma^1_{k+1}(a))\). If one of the lightface pointclasses
\(\Gamma=\Sigma^1_k,\Pi^1_k\) has the reduction property, then there are
two disjoint lightface \(\Gamma\) sets which cannot be separated by a
\(\RRrel{\Delta}{k}\) set.
\end{prop}

Proposition~\ref{prop:duality}, applied to Theorem~\ref{main:sep} and
to Theorem~\ref{main:up} in the form \(\forall a\in\R\ \UP(\Pi^1_n(a))\)
(see the remark after it), gives the following.

\begin{cor}\label{cor:models}
In \(V_n\), for every \(a\in\R\), both \(\Red(\Pi^1_n(a))\) and
\(\Sep(\Sigma^1_n(a))\) hold, and \(V_n\) satisfies
\(\UP(\Sigma^1_{n+1}(a))\). Each of the pointclasses \(\Pi^1_n\) and
\(\Sigma^1_{n+1}\) contains a disjoint pair of lightface sets which cannot
be separated by a boldface set of the corresponding \(\Delta\)-level.
In \(V_{\mathrm{sep}}\), both \(\Sigma^1_k\)- and
\(\Pi^1_k\)-uniformization fail for every \(k\ge3\), in their lightface
as well as in their boldface forms.
\end{cor}

In particular, the fixed-level models realize the separation pattern
attributed to Harrington in \cite[Theorem~1.21]{sami}. The same
proposition explains the endpoint \(n=2\). Effective Kond\^o
uniformization \cite{mosch} gives \(\UP(\Pi^1_1)\), hence
\(\UP(\Sigma^1_2)\) and \(\Red(\Sigma^1_2)\), and therefore the failure of
\(\Sep(\Sigma^1_2)\). If \(\UP(\Pi^1_2)\) held, it would imply
\(\Sep(\Sigma^1_2)\). More generally, two adjacent \(\Pi\)-levels cannot
both have the uniformization property: \(\UP(\Pi^1_k)\) gives
\(\Red(\Sigma^1_{k+1})\), whereas \(\UP(\Pi^1_{k+1})\) gives
\(\Sep(\Sigma^1_{k+1})\).

\section{Open questions}\label{sec:questions}

We close with some questions which are suggested by the results of this
paper.

\paragraph{Reduction.}
Theorem~\ref{main:up} gives \(\Pi^1_n\)-reduction at a single level
\(n\). Two adjacent \(\Pi\)-levels cannot both have the uniformization
property (Section~\ref{sec:consequences}), but we know of no such
obstruction for reduction.

\begin{question}\label{q:adjacent}
Is \(\Pi^1_3\)-reduction consistent with \(\Pi^1_4\)-reduction? More
generally, can \(\Pi^1_n\)- and \(\Pi^1_{n+1}\)-reduction hold
simultaneously for some \(n\geq3\), for instance in a set forcing
extension of \(L\)?
\end{question}

In such a model \(\Pi^1_n\)-uniformization fails: it would give
\(\Sigma^1_{n+1}\)-uniformization and hence \(\Sigma^1_{n+1}\)-reduction,
which is incompatible with \(\Pi^1_{n+1}\)-reduction by
Proposition~\ref{prop:duality}. So for \(n=3\) such a model has to
combine \(\Pi^1_3\)-reduction with the failure of
\(\Pi^1_3\)-uniformization, as in \cite{Ho3}.

\begin{question}\label{q:globalred}
Is there a model of \(\ZFC\) in which \(\Pi^1_n\)-reduction holds for
every \(n\geq3\)? Can such a model be obtained by forcing over \(L\)?
\end{question}

The universal-set argument after Theorem~\ref{main:up} applies to
reduction as well, so lightface \(\Pi^1_n\)-reduction implies
\(\Pi^1_n(a)\)-reduction for every real \(a\). Hence, by
Proposition~\ref{prop:duality}, a model of Question~\ref{q:globalred}
satisfies the conclusion of Theorem~\ref{main:sep}, and both
\(\Sigma^1_n\)- and \(\Pi^1_n\)-uniformization fail in it for every
\(n\geq3\). As \(\Sigma^1_3\)-separation does not imply
\(\Pi^1_3\)-reduction \cite{SepRed}, a positive answer would strengthen
Theorem~\ref{main:sep}.

\paragraph{Uniformization.}
By Section~\ref{sec:consequences}, a global version of
Theorem~\ref{main:up} can concern at most every other level, and
\(\Pi^1_k\)-uniformization implies \(\Sigma^1_{k+1}\)-uniformization.
This leaves two global patterns: the pattern of projective determinacy,
and its shift by one level.

\begin{question}\label{q:pd}
Is there a set forcing extension of \(L\) in which
\(\Pi^1_{2k+1}\)-uniformization holds for every \(k\geq1\), so that the
uniformization pattern of projective determinacy holds? Is there one in
which \(\Pi^1_{2k}\)-uniformization holds for every \(k\geq2\)?
\end{question}

In such a model the reals have no good projective wellorder, since a good
\(\Delta^1_m\)-wellorder yields \(\Sigma^1_k\)-uniformization for every
\(k\geq m\) \cite{addison}, and \(\Sigma^1_k\)- and
\(\Pi^1_k\)-uniformization exclude each other by
Proposition~\ref{prop:duality}. The first pattern gives a projective
uniformization of every projective relation, so together with Lebesgue
measurability and the Baire property for all projective sets it requires
large cardinals \cite{woodin82,hauserschindler}. The second main theorem
of \cite{mixed} realizes over \(L\) the uniformization pattern of the
canonical inner model with two Woodin cardinals: \(\Pi^1_3\)-uniformization
together with \(\Sigma^1_m\)-uniformization for all \(m\geq4\). We do not
know whether \(\Pi^1_3\)- and \(\Pi^1_5\)-uniformization can hold
simultaneously in an extension of \(L\), nor whether the weaker separation
pattern of projective determinacy, \(\Sigma^1_n\)-separation for odd
\(n\geq3\) and \(\Pi^1_n\)-separation for even \(n\geq4\), can be forced
over \(L\).

\paragraph{Combining the two theorems.}
The conclusions of Theorems~\ref{main:sep} and~\ref{main:up} cannot hold
in one model: \(\Pi^1_n\)-uniformization implies
\(\Sigma^1_{n+1}\)-uniformization and hence \(\Sigma^1_{n+1}\)-reduction,
so \(\Sigma^1_{n+1}\)-separation fails by Proposition~\ref{prop:duality}.
This is the only obstruction we know of.

\begin{question}\label{q:combine}
Is there, for some \(n\geq3\), a model of \(\Pi^1_n\)-uniformization in
which \(\Sigma^1_m\)-separation holds for every \(m\geq3\) with
\(m\neq n+1\)?
\end{question}

\paragraph{Separation on both sides.}
Under the axiom of determinacy, exactly one of a nonselfdual pointclass
which is closed under continuous preimages and its dual has the
separation property \cite{vanwesep,steel81}. For projective subsets of
\(\omega\), lightface separation can hold on both sides of a fixed level,
see the discussion of \cite{kl} in the introduction.

\begin{question}\label{q:bothsides}
Is there a model of \(\ZFC\) in which, for some \(n\geq3\), both
\(\Sigma^1_n\)- and \(\Pi^1_n\)-separation hold?
\end{question}

By Proposition~\ref{prop:duality}, neither \(\Sigma^1_n\)- nor
\(\Pi^1_n\)-reduction could hold in such a model.

\paragraph{Regularity properties.}
All models of this paper have the same \(\omega_1\) as \(L\), so by
\cite{shelah84} some \(\RRrel{\Sigma}{3}\) set is not Lebesgue measurable
in them. It is natural to ask how the patterns above interact with the
regularity properties of projective sets.

\begin{question}\label{q:sepreg}
Is the conclusion of Theorem~\ref{main:sep} consistent with Lebesgue
measurability and the Baire property for all projective sets? If so,
what is the consistency strength of this combination? Is an inaccessible
cardinal sufficient?
\end{question}

Projective determinacy fails in such a model, as it implies
\(\Sigma^1_{2k}\)-reduction and hence the failure of
\(\Sigma^1_{2k}\)-separation for every \(k\geq2\).

\begin{question}\label{q:upreg}
Let \(n\geq3\). What is the consistency strength of
\(\Pi^1_n\)-uniformization together with Lebesgue measurability and the
Baire property for all \(\RRrel{\Sigma}{n-1}\) sets? Is it consistent
relative to \(\ZFC\) for \(n=3\), and relative to an inaccessible cardinal
for \(n\geq4\)?
\end{question}

By Theorem~\ref{main:up}, \(\Pi^1_n\)-uniformization alone needs no large
cardinals, and for \(n\geq4\) the regularity part needs an inaccessible
cardinal \cite{shelah84}. For odd \(n\), the combination follows from
projective determinacy \cite{mosch71,mosch}. For even \(n\), projective
determinacy refutes \(\Pi^1_n\)-uniformization, and we do not know
whether the combination is consistent at all.

\end{document}